\documentclass[11pt]{article}
\usepackage{bbm}
 \usepackage{amssymb}
\usepackage{amssymb, amsthm, amsmath, amscd}
\usepackage[usenames,dvipsnames]{color}

\newtheorem{thm}{Theorem}[section]
\newtheorem{lem}[thm]{Lemma}
\newtheorem{prop}[thm]{Proposition}
\newtheorem{cor}[thm]{Corollary}
\theoremstyle{definition}
\newtheorem{NN}[thm]{}
\theoremstyle{definition}
\newtheorem{df}[thm]{Definition}
\theoremstyle{definition}
\newtheorem{rem}[thm]{Remark}

\theoremstyle{definition}

\newcommand{\red}{\textcolor{red}}
\newcommand{\blue}{\color{blue}}

\renewcommand{\phi}{\varphi}

\newcommand{\A}{\mathbb{A}}
\newcommand{\aff}{\rm aff}

\newcommand{\N}{\mathbb{N}}
\newcommand{\Z}{\mathbb{Z}}
\newcommand{\Q}{\mathbb{Q}}
\newcommand{\R}{\mathbb{R}}
\newcommand{\C}{\mathbb{C}}
\newcommand{\T}{\mathbb{T}}

\numberwithin{equation}{section}

\newcommand{\Aff}{\operatorname{Aff}}

\newcommand{\id}{\operatorname{id}}

\newcommand{\cpc}{c.p.c.~map}
\newcommand{\morp}{contractive completely positive linear map}

\newcommand{\hm}{homomorphism}
\newcommand{\dt}{\delta}
\newcommand{\ep}{\varepsilon}

\newcommand{\td}{\tilde}

\newcommand{\lr}{\longrightarrow}
\newcommand{\ld}{\lambda}
\newcommand{\cd}{\cdots}

\newcommand{\Green}{\color{Green}}

\newcommand{\p}{\mathfrak{p}}
\newcommand{\q}{\mathfrak{q}}

\newcommand{\LD}{\Lambda}

\newcommand{\DT}{\Delta}

\newcommand{\la}{\langle}
\newcommand{\ra}{\rangle}
\newcommand{\andeqn}{\,\,\,{\rm and}\,\,\,}
\newcommand{\rforal}{\,\,\,{\rm for\,\,\,all}\,\,\,}
\newcommand{\CA}{$C^*$-algebra}
\newcommand{\SCA}{$C^*$-subalgebra}

\newcommand{\af}{{\alpha}}
\newcommand{\bt}{{\beta}}
\newcommand{\dist}{{\rm dist}}

\newcommand{\diag}{{\rm diag}}
\newcommand{\wilog}{without loss of generality}
\newcommand{\Wlog}{Without loss of generality}

\newcommand{\D}{\mathbb D}
\newcommand{\beq}{\begin{eqnarray}}
\newcommand{\eneq}{\end{eqnarray}}
\newcommand{\tforal}{\,\,\,\text{for\,\,\,all}\,\,\,}
\newcommand{\tand}{\,\,\,\text{and}\,\,\,}
\newcommand{\zo}{{\cal Z}_0}

\newcommand{\LAff}{{\rm LAff}}

\usepackage{amsfonts}
\usepackage{mathrsfs}
\usepackage{textcomp}
\usepackage[all]{xy}

\title{On classification  of  non-unital 
amenable  simple  C*-algebras, III, 
Stably projectionless \CA s.}
\author{Guihua Gong and  Huaxin Lin
 }
\date{
}

\begin{document}

\maketitle

\begin{abstract}
We show, based on previous results,   that 
two separable simple stably projectionless amenable ${\cal Z}$-stable
\CA s 
 which satisfy the UCT are isomorphic if and only if
they have the same Elliott invariant.
\end{abstract}

\section{Introduction}

We will present a unified classification of separable finite simple \CA s
of  finite nuclear dimension
which satisfy the Universal Coefficient Theorem.

{{In 1989, Elliott initiated a program aimed at the complete classification of simple separable amenable $C^*$-algebras by $K$-theoretical invariant {{(the Elliott invariant).}}
{{The success of the program has  deep impacts in the study of operator algebras as well as
its applications in dynamical systems and non-commutative geometry
(\cite{GPS}, \cite{EL}, \cite{GMPS}, \cite{SW}, \cite{FR}, \cite{BEGJ},  \cite{BK}, \cite{ER04}, \cite{Yu2},
to name a few from the beginning of  an incomplete list).
 The program inspires  a great deal of
research  over last three decades (see  \cite{Rbk} and \cite{ETbull} for some earlier assessments).}}
As suggested by Elliott (\cite{Ellicm}),
{{separable simple
\CA s are}}
naturally divided into three  cases
according to their
 $K_0$-groups.\\
Case 1:
$K_0(A)=K_0(A)_+$, and $V(A)\not=\{0\}$, where $V(A)$ is the set of Murry-von Neumann  equivalence classes of projections in {{the stabilization $A\otimes {\cal K}.$}}\\
Case 2: $K_0(A)\not=\{0\}$ is an ordered group.\\
Case 3.
$K_0(A)_+=V(A)=\{0\}$.

The Jiang-Su algebra ${\cal Z}$ (\cite{JS}) constructed  in 1998 during the development of the Elliott program is an
infinite dimensional unital nuclear \CA\, in the UCT class with the feature
that its Elliott invariant
 is exactly the same as that of the complex field.
It turns out that $A$ and $A\otimes {\cal Z}$ have the same tracial structure and
$K_i(A)=K_i(A\otimes {\cal Z})$ ($i=0,1$). Moreover $A\otimes {\cal Z}$ and $A$ are in the same  one of the three cases
mentioned above.
 Since $A\otimes {\cal Z}$ is ${\cal Z}$-stable, i.e.,
$A\otimes {\cal Z}\cong (A\otimes {\cal Z})\otimes {\cal Z},$
naturally one studies  simple ${\cal Z}$-stable \CA s.  Fortunately,
a separable nuclear  simple \CA\, is ${\cal Z}$-stable if and only if it has finite nuclear dimension (see
\cite{Winter-Z-stable-02}, \cite{MS},
\cite{MS2},
 \cite{SWW}, \cite{CETWW}, \cite{T-0-Z},  and \cite{CE}).

For {{ Case 1 and  2, the Elliott  program is now completed for  separable
simple nuclear ${\cal Z}$-stable $C^*$-algebras in the UCT class.}}
(Note that the  {{non-untial}} $C^*$-algebras in {{Case}} 2 can be reduced to the unital case, by considering ${{p(A\otimes {\cal K})p}}$ for a {{non-zero}} projection $p\in {{A\otimes {\cal K}}}$.)
These are the results of decades of work by many  mathematicians (see \cite{KP}, \cite{Pclass}, \cite{GLN},
\cite{TWW} and \cite{EGLN} for the historical discussion there, also \cite{Rcl}, \cite{EG} and \cite{Lininv}).
These {{progress}}  could be summarized  briefly as the following:
Two 
{{unital}} separable simple  nuclear ${\cal Z}$-stable \CA s $A$ and $B$
with some non-trivial projections in {{$A\otimes {\cal K}$} }(or {{$B\otimes {\cal K}$}})}} which satisfy
the UCT are isomorphic if and only if  their Elliott invariant ${\rm Ell}(A)$ and ${\rm Ell}(B)$ are
isomorphic. {{It is long stand conjecture that all separable nuclear $C^*$-algebras satisfy UCT---the universal coefficient theorem for their KK groups {{(for some recent progress on this subject, see \cite{Willett-Yu1}, \cite{Willett-Yu2}
and  {{and \cite{BL}}}).}}). Moreover,
{{all weakly unperforated Elliott invariant can be  achieved by a separable simple nuclear ${\cal Z}$-stable  \CA s in {{the}} UCT class.}}

In this {{paper,}} we will study the {{Case 3,}} the case that $K_0(A)_+=V(A)=\{0\},$ i.e.,
the case that \CA s are stably projectionless.
{{The study of stably projectionless simple \CA s has a long history (see \cite{K80}, for example).}}
{{Stably projectionless simple \CA s can naturally occur in the study of
flow actions (see \cite{K80}, \cite{KK} and
\cite{Rl}).   Following Razak work in \cite{RzW}, Tsang (\cite{Tsang}) showed that any metrizable Choquet simplex
can be the tracial state space for  some  stably projectionless simple  nuclear \CA s.
In fact, in \cite{GLII}, we show that there is a unique
separable stably projectionless simple \CA\, $\zo$ with finite nuclear dimension in the UCT class
and with $K_0(\zo)=\Z,$ $K_1(\zo)=\{0\}$ and with a unique tracial state (such \CA s was known to exist).
It turns out, for any separable simple  nuclear \CA\, $A,$ $A\otimes \zo$ is stably projectionless.}} In fact
we show  (in \cite{GLII}) that, for any abelian group $G_0,$ any compact metrizable Choquet simplex $\Delta$ and
any \hm\, $\rho: G_0\to \Aff(\Delta),$ the space of continuous real affine functions on $\Delta$ such that $\rho(G_0)\cap \Aff_+(\Delta)=\{0\},$  there exists a  separable  stably projectionless simple nuclear  ${\cal Z}$-stable \CA\, $A$
such that $(K_0(A), \rho_A)=(G_0, \rho)$ (and with arbitrarily given $K_1(A)$).}}

 {{Let us {{point}} out that  {{purely}} infinite $C^*$-algebras in  Case 1 {{are}} of real rank zero, which have {{a rich structure of}} projections, {{a $C^*$-algebra}} $A$ in
{{Case}} 2  {{admits}}
at least one {{non-zero}} projection in
${{A\otimes {\cal K}}}$   so that
$A\otimes U$ have plenty of projections {{for any infinite dimensional UHF algebra $U.$}}
It is worth  mentioning  that one successful study of non-simple \CA s is to consider
\CA s with ideal property, i.e., all {{ideals}} of the $C^*$-algebras are generated by the projections inside the ideal (see {{
\cite{Pa},
\cite{Pa-Ph1}, 
\cite{Pa-R1}, \cite{Pa-R2}}}).
The total absence  of non-zero projections in $A\otimes {\cal K}$ causes a great deal of different  challenges.
Nevertheless,  significant results in this direction
have  been made in \cite{eglnkk0} and \cite{GLII} (early work may also be found in  \cite{Rl} and \cite{RzW}). The former classifies the subclass of separable simple \CA s
of finite nuclear dimension
which are $KK$-contractible.
The latter classifies the class
of separable simple \CA s with finite nuclear dimension in UCT class whose $K_0$-groups vanish on traces.
%
%
%

This paper  studies  the case $K_0(A)_+=\{0\}$ but $\rho_A(K_0(A))$ may not be zero.
  {{It should be noted that almost all the  technical results in the unital case
cannot be applied to the stably projectionless case, as one could not even find a single
non-zero projection in $A\otimes {\cal K}.$}}
{{However, while}} it is  largely independent from the unital case,
much of the work of this paper depends on \cite{eglnp}, \cite{eglnkk0} and \cite{GLII}.
 We will  also present a unified form
of classification for {{finite}} separable simple amenable ${\cal Z}$-stable  \CA s
in the UCT class.
{{In this paper, we will}} show that any two  separable finite simple amenable ${\cal Z}$-stable  \CA s 
in the UCT class are isomorphic if and only if they have the same Elliott invariant (see Theorem \ref{TisomorphismC}).
 Combining with the classification theorem of Kirchberg and Phillips, {{this gives}}  the complete classification of all simple separable
 nuclear ${\cal Z}$-stable $C^*$-algebras
 which satisfy UCT.
Since
given a stably finite (or purely infinite) separable simple amenable \CA\, $A$ {{with weakly unperforated K-theory,}}
$A\otimes {\cal Z}$  and $A$ have exactly the same Elliott invariant, as mentioned above,
we limit ourselves
to the class of separable simple amenable ${\cal Z}$-stable \CA s, or, equivalently,
the class of separable simple \CA s with finite nuclear dimension. The  nuclear dimension of {{certain $C^*$-algebras}} associated with a discrete metric space is related to  asymptotical dimension  of the underline space, and the concept of asymptotical dimension  has fundamental applications to geometry and topology (see \cite{Yu1} and \cite{Yu2}).

The paper is organized as follows.
Section 2 serves as preliminaries, {{some notations are directly taken  from  \cite{GLrange}.}} Section 3 discusses the strict comparison
for mapping tori as \SCA s of $C([0,1], A)$ for some simple \CA\, $A.$ 
Section {{4}} 
 contains some facts about $U(\td B)/CU(\td B).$
Section {{5  and 6}} 
 contain some versions of existence theorems.
Section {{7}} 
provides an isomorphism theorem for the class of separable simple
\CA s with the form $A\otimes U,$ where $A$ is a separable simple {{\CA}}\, with generalized tracial rank at most one
in the UCT class, and $U$ is an infinite dimensional  UHF-algebra.  Section {{8}} 
studies \CA s which are generated by $C_0((0,1))$ and $C,$ where
$C$ is an 1-dimensional NCCW complexes. Section {{9}} 
contains some homotopy lemmas for
 \CA s classified in section {{7}}. 
 In Section {{10}}, 
 we study maps from certain \CA s to the mapping
tori associated with two \hm s from one separable simple  \CA\,
to another   such \CA.
 Section {{11}} 
 discusses some stable results related to the homotopy lemmas.
Section {{12 and 13}} 
study asymptotic unitary equivalence of \hm s from one separable simple \CA\, to another
and  related invariant including the rotation maps. In Section {{14}}, 
we state the non-unital version
of an important result {{in \cite{Wlocal}.}} 
Applying this 
theorem  together with results in Section {{12 and 13}}, 
among others, we  present the main
theorem of the paper: the isomorphic theorem for separable {{finite}} simple \CA s 
{{with}} finite nuclear dimension in the UCT class.
{{It is}}  a uniform version for both unital and non-unital
cases.

{{{\bf{Acknowledgement}}:    This research
began when both authors stayed
in the Research Center for Operator Algebras in East China Normal University
in the summer of 2016 and 2017.
    Both authors acknowledge the support by the Center
 which is in part supported  by NNSF of China (11531003)  and Shanghai Science and Technology
 Commission (13dz2260400)
 and  Shanghai Key Laboratory of PMMP.
The second named author was also supported by NSF grants (DMS 1665183 and DMS 1954600).}}
The main part of this research was first reported in East Coast Operator Algebra Symposium, October, 2017.
Both authors would like to express their sincere gratitude to George A. Elliott for inspiration
as well as for  his moral support throughout  the three decade of research.
Much of the work of this paper is based on \cite{eglnp} and \cite{eglnkk0}. We would also like to take the opportunity
to thank  G.A. Elliott and Z. Niu for the collaboration.
\section{Preliminaries}

{{For readers' convenience, we will repeat some notations   {{in}} \cite{GLrange}.}}

{{First, throughout of the paper, a c.p.c. map is a completely positive contractive linear map.}}

\begin{df}\label{Dher}
Let $A$ be a \CA.  Denote by $A^{\bf 1}$  {{the}} unit ball of $A.$
Let  $a\in A_+.$  Denote by  ${\rm Her}(a)$ the hereditary \SCA\, $\overline{aAa}.$
If $a, b\in A_+,$ we write $a\lesssim b$ ($a$ is Cuntz smaller than $b$),
if there exists a sequence of $x_n\in A$ such that $a=\lim_{n\to\infty} x_n^*x_n$
and $x_nx_n^*\in {\rm Her}(b).$  {{If both $a\lesssim b$ and $b\lesssim a$, then we say $a$ is Cuntz equivalent to $b$.}} The Cuntz equivalence class represented by $a$ will
be denoted by $\la a\ra.$
A projection $p\in M_n(A)$ defines {{an element}} $[p]\in K_0(A)_+$.  {{We}} will also {{write  $[p]$ for}} the Cuntz equivalence class represented by $p.$ 

\end{df}

\begin{df}\label{DTtilde}
Let $A$ be a \CA.  Denote by $T(A)$ the tracial state {{space}} of $A$ {{(which could be an empty set).}}
{{Denote by $T_f(A)$ the set of faithful traces of $A$
(if $a\in A_+\setminus \{0\}$ and
$\tau\in T_f(A),$ then $\tau({{a}})>0$).}}
Let $\Aff(T(A))$ be the space of all {{real}} affine continuous functions on $T(A)$.
Let ${\tilde{T}}(A)$ be the cone of densely defined,
positive lower semi-continuous traces on $A$ equipped with the topology
of point-wise convergence on elements of the Pedersen ideal  ${\rm Ped}(A)$ of $A.$
Let $B$ be another \CA\, with $\td T(B)\not=\{0\}$
and let $\phi: A\to B$ be a \hm.

{{\it{In what follows we will also {{write}}  $\phi$ for $\phi\otimes \id_{M_k}: M_k(A)\to M_k(B)$
whenever it is convenient.}}}
We will write   $\phi_T: \td T(B)\to \td T(A)$ for the induced continuous affine map.
Denote by $\td T^{b}(A)$ the subset of $\td T(A)$ which are bounded on $A.$
Of course $T(A)\subset \td T^b(A).$ Set $T_0(A):=\{t\in \td T(A): \|{{t}}\|\le 1\}.$ It is a compact convex subset of $\td T(A).$

Let $r\ge 1$ be an integer and $\tau\in {\tilde T}(A).$
We will continue to {{write}} $\tau$  on $A\otimes M_r$ for $\tau\otimes {\rm Tr},$ where ${\rm Tr}$ is the standard
trace on $M_r.$
Let  $S$ be a convex subset (of a convex topological cone with Choquet simplex as a base).
{{We assume  that the convex cone contains 0. Denote by
$\Aff(S)$ the set of affine continuous real functions on $S$ with the property that, if  $0\in S,$ then
$f(0)=0$ for all $f\in \Aff(S).$}}
Define (see \cite{Rl})
\beq
\Aff(S)_+&=&\{f: C(S, \R)_+: f \,\, 
{{\rm affine}}, f(\tau)\ge 0\},\\
\Aff_+(S)&=&\{f: C(S, \R)_+: f \,\, 
{{\rm affine}}, f(\tau)>0\,\,{\rm for}\,\,\tau\not=0\}\cup \{0\},\\
{\rm LAff}(S)_+&=&\{f:S\to [0,\infty]: \exists \{f_n\}, f_n\nearrow f,\,\,
 f_n\in \Aff(S)_+\},\\
{\rm LAff}_+(S)&=&\{f:S\to [0,\infty]: \exists \{f_n\}, f_n\nearrow f,\,\,
 f_n\in \Aff_+(S)\}\andeqn\\
{{{\rm LAff}_{{+}}^{\sim}}}(S)
 &=&\{f_1-f_2: f_1\in {\rm LAff}_+(S)\andeqn f_2\in
 {\rm Aff}_+(S)\}.
 \eneq
 For the most part of this paper, $S={\tilde T}(A)$ {{or}} $S=T(A),$ or $S=T_0(A)$
in the above definition will be used.
Recall {{that}} $0\in {\tilde T}(A)$ and if ${{g}}\in \LAff(\td T(A)),$ then ${{g}}(0)=0.$


\end{df}

%
%

\begin{NN}\label{DElliott+}

{{Recall from Definition 2.3 of \cite{GLrange}, a (simple)}} ordered group paring is a triple $(G, T,  \rho),$ where
$G$ is a countable abelian group, $T$ is  {{a convex  topological cone}} with a Choquet simplex as its base,
and  $\rho: G\to \Aff(T)$  {{is a \hm.}}
Define $G_+=\{g\in G: \rho(g)>0\}\cup \{0\}.$ If $G_+\not=\{0\},$
then $(G, G_+)$ is an ordered group {{in the sense that $G_+\cap (-G_+)=\{0\}$ and $G_+-G_+=G$}}. It has the property that if $ng>0$ for some integer $n>0,$
then $g>0.$  In other words, $(G, G_+)$ is weakly unperforated.

The ordered group pairing {{above}} satisfies the following condition: either
 $G_+=\{0\},$ or
$(G, G_+)$ is a simple ordered group, i.e.,  every element $g\in G_+\setminus \{0\}$
is  an order unit.

A scaled simple ordered group paring is a quintuple $(G, \Sigma(G), T, s, \rho)$ such that
$(G, T, \rho)$ is a simple ordered group paring,
 $s\in \LAff_+(T)\setminus \{0\}$  {{and}}
{{\beq
\Sigma(G):=\{g\in G_+: \rho(g)<s\}\,\,{\rm or}\,\,\,
\Sigma(G):=\{g\in G_+: \rho(g)<s\}\cup\{u\},
\eneq
where {{$\rho(u)=s.$}}}}
We allow $\Sigma(G)=\{0\}.$   Note {{that}} $s(\tau)$ could be infinite for some $\tau\in T.$
It is called unital scaled simple ordered group paring, if $\Sigma(G)=\{g\in G_+: \rho(g)<s\}\cup \{u\}$
with $\rho(u)=s,$ in which case, $u$ is called the unit for $G.$  Note {{also}} that, in this {{case,}}
$u$ is the maximum element of {{$\Sigma(G),$
one may}} write
$(G, u, T, \rho)$ for $(G, \Sigma(G), T, s,\rho).$  {{On the other hand,
$\Sigma(G)$ is determined by $s$ if it has no unit.}}
One may {{then}} write $(G, T, s, \rho)$ for $(G, \Sigma(G), T, s, \rho)$
(see Theorem {{5.2 of \cite{GLrange}.}}) 

Let $(G_i, \Sigma(G_i), T_i, s_i, \rho_i),$  $i=1,2,$ be {{scaled}} simple ordered group parings.
A map
$$\Gamma_0: (G_1, \Sigma(G_1), T_1, s_1, \rho_1)\to (G_2, \Sigma(G_2), T_2, s_2, \rho_2)$$
is said to be a \hm, if there is a group \hm\,
$\kappa_0: G_1\to G_2$ and a continuous {{affine cone}} map $\kappa_T: T_2\to T_1$
such that
\beq\label{july12-2021}
{{\rho_2}}(\kappa_0(g))(t)={{\rho_1(g)}}(\kappa_T(t))\rforal g\in G_1\andeqn t\in T_2,\andeqn\\ \label{july12-2021-1}
\kappa_0(\Sigma(G_1))\subset \Sigma(G_2)), \andeqn   s_1(\kappa_T(t))\le s_2(t) \rforal t\in T_2.
\eneq
We say a \hm\, $\Gamma_0$ is an isomorphism
if $\kappa_0$ is an isomorphism, $\kappa_0(\Sigma(G_1))={{\Sigma(G_2)}},$ $\kappa_T$ is a cone
homeomorphism, and $s_1(\kappa_T(t))=s_2(t)$ for all $t\in T_2.$
\end{NN}

\begin{df}\label{Dfep}
For any $\ep>0,$ define $f_\ep\in {{C([0,\infty))_+}}$ by
$f_\ep(t)=0$ if $t\in [0, \ep/2],$ $f_\ep(t)=1$ if $t\in [\ep, \infty)$ and
$f_\ep(t)$ is linear in $(\ep/2, \ep).$

Let $A$ be a \CA\, and $\tau$ {{a}} quasitrace.
For each $a\in A_+,$ define $d_\tau(a)=\lim_{\ep\to 0} {{\tau(f_\ep(a)).}}$
Note that $f_\ep(a)\in {\rm Ped}(A)$ for all $a\in A_+.$

{{Let $S$ be a  convex subset of $\tilde{T}(A)$ and $a\in M_n(A)_+.$
The function $\hat{a}(s)=s(a)$ (for $s\in S$) is an  affine function.
Define $\widehat{\la a\ra }(s)=d_s(a)=\lim_{\ep\to 0}s(f_\ep(a))$ (for $s\in S$),
which is a lower semicontinuous function.
If $a\in {\rm Ped}(A),$ then $\hat{a}$ is in $\Aff_+(S)$ and
$\widehat{\la a\ra}\in \LAff_+(S)$  {{in general}} (see \ref{DTtilde}). {{Note
that $\hat{a}$ is different from $\widehat{\la a\ra}.$}}
 In most cases, $S$ is ${\tilde T}(A),$ $T_0(A),$ or $T(A).$}}
Note {{also}} that, there is a nature map from ${\rm Cu}(A)$ to {{${\rm LAff}_+(\tilde{T}(A))$}} by sending $\la a \ra$ to $\widehat{\la a \ra}$.

\end{df}

\begin{NN}\label{range4.1}
If $A$ is a unital \CA\, and $T(A)\not=\emptyset,$ then
there is a canonical  
\hm\, $\rho_A: K_0(A)\to \Aff(T(A)).$

{{Now consider the case that $A$ is not unital.
Let $\pi_\C^A: \td A\to \C$ be the quotient map.
Suppose that $T(A)\not=\emptyset.$
Let $\tau_\C:=\tau_\C^A: \td A\to \C$ be
the tracial state which factors through $\pi_\C^A.$
Then
\beq
T(\td A) = \big\{t{{\tau_{_\C}^A}}+(1-t)\tau:~ t\in [0,1],~ \tau\in {{T(A)}}\big\}.
\eneq
The map ${{T(A)}}\hookrightarrow {{T(\td A)}}$ induces a map {{${{\Aff(T(\tilde{A}))}}\to \Aff(T(A))$}}. Then
 the  map {{$\rho_{\td A}: K_0({{\td A}})\to \Aff(T(\td A))$}} induces a \hm\,  $\rho':~ K_0(A)\to \Aff (T(A))$ by
\beq
{{\rho':~K_0(A)\to  K_0({{\td A}})\stackrel{\rho_{\td A}}{\longrightarrow}
\Aff(T(\td A))   \to \Aff(T(A)).}}
\eneq
However, in the case that $A\not={\rm Ped}(A),$ we will not use $\rho'$  in general, as it is possible
that $T(A)=\emptyset$ but ${\td T}(A)$ is rich (consider the case $A\cong A\otimes {\cal K}$).}}
\end{NN}

\begin{df}{{[Definition 2.6 of \cite{GLrange}]}}\label{Dparing}
{{Let $A$ be a \CA\, with $\td T(A)\not=\{0\}.$
If $\tau\in \td T(A)$ is bounded on $A,$ then $\tau$ can be extended naturally
to {{a trace}} on $\td A.$
Recall that $\td T^b(A)$ is the set of bounded traces on $A.$
Denote by $\rho_A^b: K_0(A)\to \Aff(\td T^b(A))$ the \hm\, defined by
$\rho_A^b([p]-[q])=\tau(p)-\tau(q)$ for all $\tau\in \td T^b(A)$
and for projections $p, q\in M_n(\td A)$ (for some integer $n\ge 1$)
and $\pi_\C^A(p)=\pi_\C^A(q).$ {{Note $p-q\in  M_n(A).$
Therefore $\rho_A^b([p]-[q])$ is continuous on $\td T^b(A).$}}
In the case that $\td T^b(A)=\td T(A),$ for example, $A={\rm Ped}(A),$
we write $\rho_A:=\rho_A^b.$

Let $A$ be a $\sigma$-unital \CA\, with a strictly positive element $0\le e\le 1.$
Put  $e_n:=f_{1/2^n}(e).$ Then $\{e_n\}$ forms an approximate identity
for $A.$ Note $e_n\in {\rm Ped}(A)$ for all $n.$
Set $A_n={\rm Her}(e_n):=\overline{e_nAe_n}.$
Denote by $\iota_n:A_n\to A_{n+1}$   and $j_n: A_n\to A$ the embeddings. It extends to
$\iota_n^\sim: \td A_n\to \td A_{n+1}$ and $j_n^\sim :\td A_n\to \td A$ {{unitally.}}
Note that $e_n\in {\rm Ped}(A_{n+1}).$ Thus  $\iota_n$  and $j_n$ induce
continuous cone maps ${\iota_n}_T^b: \td T^b(A_{n+1})\to
\td T^b(A_n)$ and ${j_n}_T: \td T(A)\to \td T^b(A_n)$
(defined by ${\iota_n}_T^b(\tau)(a)=\tau(\iota_n(a))$
for $\tau\in {{\td T^b}}(A_{n+1}),$ and ${j_n}_T(\tau)(a)=\tau(j_n(a))$
for all $\tau\in \td T(A)$ and all $a\in A_n$), respectively.
Denote by $\iota_n^\sharp: \Aff(\td T^b(A_n))
\to {{\Aff(\td T^b(A_{n+1}))}}$  and $j_n^\sharp: \Aff(\td T^b(A_n))\to \Aff(\td T(A))$ the induced continuous linear maps.
Recall that $\cup_{n=1}^\infty A_n$ is dense in ${\rm Ped}(A).$
A direct computation shows that one may obtain  the following inverse direct limit of {{convex  topological cones}} (with continuous cone maps):
\beq\label{Drho-trace}
\td T^b(A_1)\stackrel{{\iota_1}_T^b}
{\longleftarrow} \td T^b(A_2)\stackrel{{\iota_2}_{T}^b}{\longleftarrow} \td T^b(A_3)
\cd\longleftarrow\cd\longleftarrow \td T(A).
\eneq
}}
which induces the following commutative diagram.
\beq
\Aff(\td T^b(A_1))\stackrel{\iota_1^\sharp}
{\longrightarrow} \Aff (\td T^b(A_2))\stackrel{{\iota_2}^\sharp}{\longrightarrow} \Aff (\td T^b(A_3))
\cd\longrightarrow\cd\longrightarrow \Aff (\td T(A)).
\eneq
Hence one also has the following commutative diagram:
\begin{displaymath}
    \xymatrix{
        K_0(A_1) \ar
        [d]_{\rho_{A_1}}\ar[r]^{\iota_{1*o}} & K_0(A_2) \ar[r]^{\iota_{2*o}} \ar
        [d]_{\rho_{A_2}}& K_0(A_3) \ar[r]
        \ar
        [d]_{\rho_{A_3}}& \cd K_0({{A}} )\\
        \Aff(\td T^b({{A}}_1))
        \ar[r]^{\iota^{\sharp}_{1,2}}
        &
         \Aff(\td T^b({{A}}_2)) \ar[r]^{\iota_2^{\sharp}}
         &
         \Aff(\td T^b({{A}}_3))
         \ar[r]
         & \cd
         \Aff(\td T({{A}})).}
\end{displaymath}
Thus one obtains a \hm\, $\rho: K_0(A)\to \Aff(\td T(A)).$
It should be noted that, when $A$ is simple, $\td T^b(A_n)=\td T(A_n)$ for all $n$
{{(see Definition 2.6 of \cite{GLrange} for more details). Moreover, the map $\rho$ does not depend
on the choice of $\{e_n\}.$}} {{We will write $\rho_A:=\rho.$ In  the case that $T(A)$ generates ${\tilde T}(A)$
such as the case that $A={\rm Ped}(A),$
we may also write $\rho_A: K_0(A)\to \Aff(T(A))$  by   restricting
$\rho_A(x)$  on $T(A)\subset {\tilde T}(A)$ and for all $x\in K_0(A).$}}
\end{df}


\begin{df}\label{DElliott}
{{Let $A$ be a \CA\, with $\td T(A)\not=\{0\}.$  In 2.6 of \cite{GLrange} (see 2.6 above),
we define a paring $\rho_A: K_0(A)\to \Aff(\td T(A))$ which is a \hm.
If every trace in $\td T(A)$ is bounded, then it coincides with the usual paring
(see the detail in 2.6 of \cite{GLrange}).}}
{{Recall from 2.6  of \cite{GLrange}, we {\it will write} $\pi_{\aff}^{\rho, A}: \Aff(\td T(A))\to \Aff(\td T(A))/\overline{\rho_A(K_0(A))}$
for the quotient map.   This may be simplified to $\pi_{\aff}^\rho$ if $A$ is clear. When $T(A)\not=\emptyset,$ we will use the same notation for
the quotient map $\Aff(T(A))\to \Aff(T(A))/\overline{\rho_A(K_0(A))}.$
In this case, we also {{write}}  $\rho_A^\sim:  K_0(\td A)\to \Aff(T(A))$ for the map
defined by $\rho_A^\sim([p])(\tau)=\tau(p)$ for projections $p\in M_l({{\td A}})$ (for all integer $l$).}}

{{Recall, from 2.6 and 2.7 of \cite{GLrange},  {{that,}} with the  pairing $\rho_A: K_0(A) \to {{\Aff (\td T(A))}}$ defined in  2.6  of \cite{GLrange}, }} the Elliott invariant  for separable simple \CA s {{(see {{\cite{Ellicm} and}} \cite{point-line}),}}
%
%
%
{{for}} the case  ${{{\tilde T}(A)}}\not=\{0\},$
 is {{described by}} the six-tuple:
$$
{\rm{Ell}}(A):=((K_0(A), \Sigma(K_0(A)), {\tilde T}(A), \Sigma_A, \rho_A), K_1(A)),
$$
where $\Sigma(K_0(A))=\{x\in K_0(A): x=[p]\,\, {\rm for\,\, some\,\, projection},\,\, p\in A\},$
and $\Sigma_A$ is a function in $\LAff_+(\td T(A))$ defined by
\beq
\Sigma_A(\tau)=\sup\{\tau(a): a\in {\rm Ped}(A)_+,\,\,\|a\|\le 1\}.
\eneq
 Let $e_A\in A$ be a strictly positive element.
Then $\Sigma_A(\tau)=\lim_{\ep\to 0} \tau(f_\ep(e_A))$ for all $\tau\in {{\td T(A),}}$ which
is independent of the choice of $e_A.$

Let $B$ be another separable  \CA.  A {{\hm\,}} $\Gamma: {\rm Ell}(A)\to {\rm Ell}(B)$ consists of a {{\hm\,}}   $\Gamma_0:(K_0(A), \Sigma(K_0(A)), {\tilde T}(A), \Sigma_A, \rho_A)$ to $(K_0(B), \Sigma(K_0(B)), {\tilde T}(B), \Sigma_B, \rho_B)$ (as
in \ref{DElliott+}, see (\ref{july12-2021}) and (\ref{july12-2021-1}) also) and a homomorphism  $\kappa_1:K_1(A)\to K_1(B)$. We say that $\Gamma$ is an isomorphism if  both $\Gamma_0$ and $\kappa_1$ are isomorphisms.

In the case that $\rho_A(K_0(A))\cap {\rm LAff}_+(\td T(A))=\{0\},$ we often consider
the reduced case that $T(A)$ is compact which gives a base for $\td T(A).$ In that  case, we
may write ${\rm Ell}(A)=(K_0(A), T(A), \rho_A, K_1(A)).$ Note {{that}}, in the said situation,
$\Sigma(K_0(A))=\{0\}$ {{and}}
$\td T(A)$ is determined by $T(A)$ and $\Sigma_A(\tau)=1$ for all $\tau\in T(A).$

\end{df}

\begin{df}\label{Dastablerk1} Let $A$ be a \CA. We say $A$ has almost stable rank one,
if $A$
has the following property:
the set of  invertible elements of $\td B$ of
every hereditary \SCA\, $B$ of {{$A$}} is dense in $B.$
$A$ is said {{stably to have}} 
almost stable rank one,
if $M_n(A)$ has almost stable rank one for all integer $n\ge 1.$
\end{df}

\begin{df}\label{Dlambdas}
Let $A$ be a \CA\, with $T(A)\not=\emptyset.$
Suppose that $A$ has a strictly positive element $e_A\in {\rm Ped}(A)_+$ with $\|e_A\|=1.$
Then $0\not\in \overline{T(A)}^w,$ the closure of $T(A)$ in ${\tilde T}(A)$
(see Theorem 4.7  of \cite{eglnp}).
Define
$$
\lambda_s(A)=\inf\{d_\tau(e_A): \tau\in {{\overline{T(A)}^w}}\}{{=\lim_{n\to\infty}\inf\{\tau(f_{1/n}(e_A)):\tau\in T(A)\}.}}
$$
 Let $A$ be a \CA\, with $T(A)\not=\{0\}.$
There is an affine  map
$r_{\aff}: A_{s.a.}\to \Aff(T_0(A))$
defined by
$$
r_{\aff}(a)(\tau)=\hat{a}(\tau)=\tau(a)\tforal \tau\in T_0(A)
$$
and for all $a\in A_{s.a.}.$ Denote by $A_{s.a.}^q$ the space  $r_{\aff}(A_{s.a.}),$
$A_+^q=r_{\aff}(A_+)$ {{and $A_+^{q, {\bf 1}}=r_{\aff}(A_+^{\bf 1}).$}}

\end{df}

%
%

%

\begin{df}\label{DkappaJ}
Let $A$ be a unital \CA.    Denote by $U(A)$ the unitary group of $A$ and
$U_0(A)$ the connected component of $U(A)$ containing $1_A.$
Denote by $CU(A)$ the closure of the commutator subgroup of $U(A).$
Let $\phi: A\to B$ be a unital \hm\,(assuming $B$ is also unital). Then $\phi$ induces a continuous \hm\,
$\phi^\dag: U(A)/CU(A)\to U(B)/CU(B).$
%
%
%
%
%
%
%

Suppose that $A$ is a unital \CA\, with stable rank $k.$ {{Let $\Pi_{cu}^{A}: U(M_k(A))/CU(M_k(A))\to K_1(A)$ be the canonical map.}}
By 3.2 of \cite{Thomsen} (see 2.16 of \cite{GLN}),
there is a  split short exact sequence
\beq\label{Dcu-4}
\hspace{-0.2in}0\to \Aff(T(A))/\overline{\rho_A(K_0(A))}\to U(M_k(A))/CU(M_k(A))
{{\overset{\Pi^A_{cu}}{\underset{J^A_{cu}}{\rightleftharpoons}}}}K_1(A)\to 0.
\eneq
In  what follows, for  each unital \CA\, $A$ of stable rank $k,$  we will {\it fix} one splitting map  $J_{cu}^A,$
and, {{we will identify $\Aff(T(A))/\overline{\rho_A(K_0(A))}$ with a
subgroup of  $U(M_k(A))/CU(M_k(A)).$}}

For the most part, $A$ will be stable rank one. So $k=1$ in the above diagram.

For each continuous and piecewise smooth path $\{u(t): t\in [0,1]\}\subset U(M_k(A)),$
or $\{u(t):t\in [0,1]\}\subset U(M_k(\td A)),$ if $A$ is not unital, but $T(A)\not=\emptyset,$
 define 
$$
D_A(\{u(t)\})(\tau)={1\over{2\pi i}}\int_0^1 \tau({du(t)\over{dt}}u^*(t))dt,\quad \tau\in T(A)\,\, ({\text{or}} \,\,T(\td A)).
$$
Note that {{here}} we view $T(A)$ as a convex subset of $T(\td A),$ if $A$ is not unital
(by extending tracial states of $A$).
Let us consider the non-unital case.
For each $\{u(t)\},$ the map $D_A(\{u(t)\})$ is a real continuous affine function on $T(\td A).$\index{$D_A,$ $\overline{D}_A$}
Let us recall de la Harpe and
Skandalis determinant
$$\overline{D}_A: U_0(M_k(\td A))/CU(M_k(\td A))\to \Aff(T(\td A))/\overline{\rho_{\td A}(K_0(\td A))}$$
{{is given by,}}
 for any ${\bar{u}}\in U_0(M_k(\td A))/CU(M_k(\td A))$ (represented by $u$),
$$
\overline{D}_A(\bar u)=D_A(\{{{u(t)}}\})+\overline{\rho_{\td A} (K_0(\td A))},
$$
where $\{u(t): t\in [0,1]\}\subset M_k(\td A)$ is a continuous and piecewise
smooth path of unitaries with $u(0)=1$ and  $u(1)=u.$)  It is known that the de la {{Harpe}} and Skandalis determinant is independent
of the choice of representative for $\bar u$ and the choice of  path $\{u(t)\}.$

\end{df}


\begin{df}\label{DLddag}
Let $A$ be a unital \CA.
 Suppose  $x\in A$ such
 {{that}}
$\|xx^*-1\|<1$ and $\|x^*x-1\|<1.$ Then $x|x|^{-1}$ is a unitary.
Let us write $\lceil x\rceil$ for the unitary $x|x|^{-1}.$

For any finite subset ${\cal U}\subset U(A),$ there exists $\dt>0$ and  a finite subset ${\cal G}\subset  A$
satisfying the following:
If $B$ is another unital \CA\, {{and}} $L: A\to B$ is {{a ${\cal G}$-$\dt$}}-
multiplicative \cpc,
then
$\overline{\lceil L(u)\rceil}$ is a {{well-defined}} element in $U(B)/CU(B)$  for all $u\in {\cal U}.$
{{We may assume that $[L]|_{{\cal S}}$ is well defined, where ${\cal S}$ is the image of ${\cal U}$
in $K_1(A)$ (see, for example, 2.12 of \cite{GLN}).}}
Let $G({\cal U})$ be the subgroup generated by ${\cal U}.$
{{Suppose that $1/2>\ep>0$ is given.  By  Appendix in \cite{LinLAH},  we may assume that
there is a \hm\, $L^\dag: G({\cal U})\to U(B)/CU(B)$ such that
\beq\label{Ddag-n1}
{\rm dist}(L^\dag(\bar u),\overline{\lceil L(u)\rceil})<\ep\rforal u\in {\cal U}.
\eneq}}
{{Moreover, as in  Definition 2.17 of  \cite{GLN},  we may also assume that
\beq\label{Ddag-n2}
L^{\dag}((G({\cal U})\cap U_0(A))/CU(A))\subset U_0(B)/CU(B).
\eneq}}
{{It follows that
${{\Pi_{cu}^{B}}}\circ {{L^{\dag}}}(\overline{u})=[L]\circ {{\Pi_{cu}^{A}}}([u])\rforal u\in G({\cal U}),$
where ${{\Pi_{cu}^{A}}}$ and ${{\Pi_{cu}^{B}}}$ are defined as in \ref{DkappaJ}
(see Definition 2.17 of \cite{GLN}).}}
In what follows,  {{when $1/2>\ep>0$ is given, whenever we write
 $L^{\dag},$ we mean that $\dt$ is small enough and ${\cal G}$ is large enough so that $L^{\dag}$ is
 defined,  {{and hence}} \eqref{Ddag-n1} and \eqref{Ddag-n2} hold (see 2.17 of \cite{GLN}).}}
 Moreover, for an integer $k\ge 1,$ we will also {{write}}  $L^{\dag}$
 {{for the map on  some given subgroup of}} $U(M_k(A))/CU(M_k(A))$ induced by $L\otimes {\rm id}_{M_k}.$  In particular, when $L$ is a unital \hm, the map $L^{\dag}$ is well defined
 on $U(M_k(A))/CU(M_k(A)).$

{{If $A$ is not unital, $L^\dag$ is defined to be ${\tilde L}^\dag,$ where ${\tilde L}: {\tilde A}\to {\tilde B}$
 is the unital extension of $L.$}}

\end{df}

\begin{df}\label{Dstrongaue}
Let $A$ and $B$ be two  \CA s.  A sequence of
linear maps $L_n: A\to B$ is said be approximately multiplicative
if
$$
\lim_{n\to\infty}\|L_n(a)L_n(b)-L_n(ab)\|=0\rforal a, b\in  A.
$$
Let
$\phi, \psi: A\to B$ be \hm s. We say $\phi$ and $\psi$ are asymptotically unitarily
equivalent if there is a continuous path of unitaries $\{u(t): t\in [1, \infty)\}$ in
$B$ (if $B$ is not unital, $u(t)\in \td B$) such
that
$$
\lim_{t\to\infty} u^*(t)\phi(a) u(t)=\psi(a)\rforal a\in A.
$$
We say $\phi$ and $\psi$ are
strongly  asymptotically unitarily equivalent if $u(1)\in  U_0(B)$ (or in $U_0(\td B)$).
\end{df}

\begin{df}\label{THfull}
Let $A$ and $B$ be \CA s, and let $T: A_+\setminus\{0\}\to \N\times \R_+\setminus \{0\}$
defined by $a\mapsto (N(a), M(a)),$ where $N(a)\in \N$ and $M(a)\in \R_+\setminus\{0\}.$
Let ${\cal H}\subset A_+\setminus \{0\}.$ A map $L: A\to B$ is said to be
$T$-${{\cal H}}$-full
if, for  any $a\in {\cal H}$ and any $b\in  B_+$ with $\|b\|\le 1,$  any $\ep>0,$ there are $x_1, x_2,...,x_N\in B$
with $N\le N(a)$ and $\|x_i\|\le M(a)$  such that
\beq
\|\sum_{j=1}^N x_i^*L(a)x_i-b\|\le \ep.
\eneq
$L$ is said to {{be}}  exactly $T$-${{\cal H}}$-full, if $\ep=0$ in the above formula
{{(see also 5.5 of \cite{eglnp}).}}

\end{df}

 \begin{df}
 \label{DfC1}
{\rm
Let $A$ and $B$ be \CA s and $\phi_0, \phi_1: A\to B$ be \hm s.
By  mapping torus  $M_{\phi_0, \phi_1},$ we mean the following
\CA:
\beq\label{dmapping}
M_{\phi_0, \phi_1}=\{(f,a)\in C([0,1], B)\oplus A: f(0)=\phi_0(a)\andeqn f(1)=\phi_1(a)\}.
\eneq

{{One has the
short exact sequence
\begin{equation*}\label{Mtoruses}
0\to SB\stackrel{\imath}{\to}M_{{{\phi, \psi}}} \stackrel{\pi_e}{\to} A\to 0,
\end{equation*}
where $\imath: SB\to M_{{{\phi, \psi}}}$ is the embedding and $\pi_e$ is the
quotient map from $M_{{{\phi, \psi}}}$ to $A.$
Denote by  $\pi_t: M_{{{\phi, \psi}}}\to B$  the point evaluation at $t\in [0,1].$}}

Let $F_1$ and $F_2$ be two finite dimensional \CA s.
Suppose that there are  (not {{necessarily}} unital)  \hm s
$\phi_0, \phi_1: F_1\to F_2.$
Denote the mapping torus $M_{\phi_1, \phi_2}$ by
$$
A=A(F_1, F_2,\phi_0, \phi_1)
=\{(f,g)\in  C([0,1], F_2)\oplus F_1: f(0)=\phi_0(g)\andeqn f(1)=\phi_1(g)\}.
$$

Denote by ${\cal C}$ the class of all  \CA s of the form $A=A(F_1, F_2, \phi_0, \phi_1).$
These \CA s {{are called  Elliott-Thomsen building blocks as well as
one  dimensional non-commutative CW complexes  (see  \cite{ET-PL} and  \cite{point-line}).}}

{{Recall that ${{\cal C}_0}$ is the class of all $A\in {\cal C}$
with $K_0(A)_+=\{0\}$ such that $K_1(A)=0$   and
$\lambda_s(A)>0,$ and ${\cal C}_0^{(0)}$ the class of all $A\in {\cal C}_0$ such that $K_0(A)=0.$  Denote by {{${\cal C}'$,}} ${\cal C}_0'$  and ${\cal C}_0^{0'}$ the class of all full hereditary \SCA s of \CA s in {{$\cal C$,}} ${\cal C}_0$ and
{{${\cal C}_0^{{(0)}},$}} respectively.}}
}
\end{df}


 \begin{df}\label{DD0}{{(cf. 8.1 and 8.2 of \cite{eglnp})}}
{{Recall}} the definition of class ${\cal D}$ and ${\cal D}_0.$

 Let $A$ be a non-unital simple \CA\,  with a strictly positive element $a\in A$ such that 
 $\|a\|=1.$   Suppose that there exists
$1> \mathfrak{f}_a>0,$ for any $\ep>0,$  any
finite subset ${\cal F}\subset A$ and any $b\in A_+\setminus \{0\},$  there are ${\cal F}$-$\ep$-multiplicative \cpc s $\phi: A\to A$ and  $\psi: A\to D$  for some
\SCA\, $D\subset A$ with $D\in {\cal C}_0'$ (or ${\cal C}_0^{0'}$) {{such that}}  $D\perp \phi(A)$ and
\beq\label{DNtr1div-1++}
&&\|x-(\phi(x)+\psi(x))\|<\ep\rforal x\in {\cal F}\cup \{a\},\\\label{DNtrdiv-2}
&&c\lesssim b,\\\label{DNtrdiv-4}
&&t(f_{1/4}(\psi(a)))\ge \mathfrak{f}_a\rforal t\in T(D),
\eneq
where $c$ is a strictly positive element of $\overline{\phi(A)A\phi(A)}.$
  Then
 we say $A\in {\cal D}$ (or ${\cal D}_0$).

 {{ Note, by  Remark 8.11 of \cite{eglnp},
 $D$ can {\it always} be chosen to be in ${\cal C}_0$ (or ${\cal C}_0^{{(0)}}$).}}

{{When $A\in {\cal D}$ and $A$ is separable, {{then}}  $A=\mathrm{Ped}(A)$ (see 11.3  of \cite{eglnp}).
Let $a\in A_+$ with $\|a\|=1$ be a strict positive element.
Put
\beq\label{1225dd}
d=\inf\{\tau(f_{1/4}(a)): \tau\in T(A)\}.
\eneq
Then, for any $0<\eta<d,$ $\mathfrak{f}_a$ can be chosen to  be
$d-\eta$ (see  Remark 9.2 of \cite{eglnp}).
One may also assume that $f_{1/4}(\psi(a))$ is
{{full}} in $D.$
Furthermore,  there exists a map:  $T: A_+\setminus \{0\}\to \N\times \R_+\setminus \{0\}$
which is independent of
${\cal F}$ and $\ep$ such that, for any  finite subset ${\cal H}\subset A_+\setminus \{0\},$ we can further  require that $\psi$ is exactly $T$-${\cal H}$-full (see 8.3
and 9.2
of \cite{eglnp}).
}}
For any $n\ge 1,$ one can choose a strictly positive element $b\in A$ with $\|b\|=1$ such that
$f_{1/4}(b)\ge f_{1/n}(a).$ Therefore, if $A$ has continuous scale, $0<d<1$ can be chosen to be arbitrarily close to $1,$
if the strictly positive element is chosen accordingly.

In \cite{eglnp}, it is proved that  if $A$ a separable simple \CA\, in ${\cal D},$
then $A$ has stable rank one, {{is stably projectionless with}}
 ${\rm Cu}(A)=\LAff_+(\td T(A)),$  and every 2-quasitrace
on $A$ is a trace {{(see 9.3 and 11.11 of \cite{eglnp}).}}

Let $A$ be a separable stably projectionless simple \CA.  Recall that  $A$ has generalized tracial rank at most
one and {{is written}} $gTR(A)\le 1,$ if there exists $e\in {\rm Ped}(A)_+$ with $\|e\|=1$ such that
$\overline{eAe}\in {\cal D}$ (see  11.6 of \cite{eglnp}).  It should be noted that, in  the definition of $D$ above,
if we assume that $A$ is unital, and replace ${\cal C}_0$ by ${\cal C},$ then $gTR(A)\le 1$ (see 9.1,
9.2 and 9.3 of \cite{GLN}). But {{the}} condition \eqref{DNtrdiv-4} and constant $\mathfrak{f}_a$ are not needed.
In  the case $K_0(A)_+\not=\{0\}$ but $A$ is not unital,  we may define $gTR(A)\le 1,$ if
for some nonzero projection $e\in M_k(A)$ (for some $k\ge 1$) $gTR({{eM_k(A)e}})\le 1$ (see \cite{GLN}).

 \end{df}

 \begin{df}\label{DD1}
 Let $A\in {\cal D}$ {{be}} as defined in \ref{DD0}.
 If, in  addition,
 for any integer $n,$ {{we can choose $D$ and $\psi$ to satisfy the following condition:}}
 $D=M_n(D_1)$ for some $D_1\in {\cal C}_0$ such that
\vspace{-0.14in} \beq\label{DD1-1}
 \psi(x)=\diag(\overbrace{\psi_1(x), \psi_1(x),...,\psi_1(x)}^n)\rforal x\in {\cal F},
 \eneq
where $\psi_1: A\to D_1$ is an ${\cal F}$-$\ep$-multiplicative \cpc, then we say $A\in {\cal D}^d.$

Note that here, as in  8.3 and 9.2 of \cite{eglnp}, {{there exists a map $T$ as mentioned in \ref {DD0}
{{(independent of $\ep$ and ${\cal F}$)}} such that $\psi$ {{can be chosen}} exactly $T$-${\cal H}$-full, for a pre-given set ${\cal H}\subset A_+\setminus \{0\},$}}
and
$\mathfrak{f}_a$ can be also
chosen as $d-\eta$ for any $\eta>0$ with $d$ as in \eqref{1225dd} for a certain strictly positive element $a.$

\end{df}

\begin{rem}\label{RDd}
  It follows  from
  {{10.4 and 10.7 of \cite{eglnp}}} that, if $A\in {\cal D}_0,$  then $A\in {\cal D}^d.$ Moreover,
  $D_1$ can be chosen in ${\cal C}_0^{(0)}$.
  If $A$ is a separable simple \CA\, in ${\cal D}$ and $A$ is tracially approximate divisible  {{(in the sense of
  10.1 of \cite{eglnp}),}}
  then $A\in {\cal D}^d.$
\end{rem}

\begin{df}\label{DWtrace}
Throughout the paper,
${\cal W}$ is the separable simple \CA\, with a unique tracial state  {{and}}  is an inductive limit
of \CA s in ${\cal C}_0^{{(0)}}$
(see \cite{RzW}).
{{In fact }}${\cal W}$ is the unique separable simple \CA\, with finite
nuclear dimension which is $KK$-contractible and with a unique tracial state {{(see \cite{eglnkk0}).}}
Denote by $\tau_W$ the unique tracial state of ${\cal W}.$

Let $A$ be a \CA\, and let $\tau$ be a nonzero trace of $A.$
We say $\tau$ is a ${\cal W}$-trace, if there exists a sequence of approximately multiplicative \cpc s
$\phi_n: A\to {\cal W}$ such that
\beq
\lim_{n\to\infty}\tau_W\circ \phi_n(a)=\tau(a)\rforal a\in A.
\eneq
\end{df}

\begin{df}\label{Dcompatible}
Let $A$ and $B$ be two
separable \CA s.
%
For convenience, let us assume that $U(M_k(\td A))/CU(M_k(\td A))=U(\td A)/CU(\td A)$ here
{{for all $k\ge 1.$}}  Otherwise, we will
replace $U(\td A)/CU(\td A)$ by
$U(M_\infty(\td A))/CU(M_\infty(\td A)).$
Let $\kappa_i: K_i(A)\to K_i(B)$ be a \hm\, ($i=0,1$), {{$\kappa_T:
{\tilde T}(B)\to {\tilde T}(A)$ be a continuous affine map such that
 $\kappa_T(T(B))\subset T(A)$}} 
   and
$\kappa_{cu}: U(\td A)/CU(\td A)\to U(\td B)/CU(\td B)$
(in case $A$ or $B$ is unital, we replace $\td A$ by $A,$
or $\td B$ by $B$) be a continuous \hm.
We say $\kappa_0$ and $\kappa_T$ are compatible,
if $\rho_B(\kappa_0(x))(\tau)=(\rho_A(x))(\kappa_T(\tau))$
for all {{$\tau\in  {\tilde T}(B)$}} and $x\in K_0(A).$
We say $\kappa_1$ and
$\kappa_{cu}$ are compatible, if {{$\Pi_{cu}^B(\kappa_{cu}(x))=\kappa_{1}(\Pi_{cu}^A(x))$
for all $x\in U(\td A)/CU(\td A).$}}

In the case that $A$ and $B$ are {{non-unital}} separable \CA s
we have $K_0({\td A})=\{(x,n): x\in K_0(A), n\in \Z\}$ and $T({\td A})=\{\af\tau+(1-\af)\tau_\C^A: \tau \in T(A), 0\leq \af\leq 1\}$, where $\tau_\C^A$ is the trace factoring through $\pi_\C^A: \td A\to \C$ (see \ref{range4.1}). Hence $\kappa_0$ induces a map from $K_0({\td A})$ to $K_0({\td B})$ (still denoted by $\kappa_0$), and ${\kappa_T}|_{T(B)}$ induces an affine map from $T({\td B})$ to $ T({\td A})$ (still denoted by $\kappa_T$). Furthermore $\kappa_T$ induces a map $\kappa^\sharp:
\Aff( T({\td A}))\to \Aff( T({\td B}))$.
We say $\kappa_{cu}$ and $\kappa_T$ are compatible,
if {{$\kappa_{cu}(\pi_{\aff}^{\rho,{\td A}}(f))=\pi_{\aff}^{\rho, {\td B}}(\kappa^{\sharp}(f))$}}
for all {{$f\in \Aff( T({\td A}))$}} (see  \ref{DElliott} above), 
{{where $\pi_{\aff}^{\rho,{\td A}}(f)$ is in $ \Aff(T(\td A))/\overline{\rho_{\td A}(K_0(\td A))},$ which is identified with a {{subgroup}} of $U(\td A)/CU(\td A)$ (see Definition \ref{DkappaJ}). }}
 We say $\kappa_1, \kappa_{cu}$ and {{$\kappa_T$}} are compatible
if $\kappa_1$ and $\kappa_{cu}$ are compatible and $\kappa_{cu}$ and {{$\kappa_T$}}
are compatible. We say $\kappa_0, \kappa_1, {{\kappa_T}}, \kappa_{cu}$ are
compatible if $\kappa_0$ and {{$\kappa_T$}} are compatible and
$\kappa_1, \kappa_{cu}$ and {{$\kappa_T$}} are {{compatible.}}
If $\kappa\in KL(A,B)$ (or $KK(A,B)$ 
) which induces $\kappa_i: K_i(A)\to K_i(B)$ ($i=0,1$),
we say $\kappa, {{\kappa_T}},\kappa_{cu}$ are compatible, if $\kappa_0, \kappa_1, {{\kappa_T}}, \kappa_{cu}$ are
compatible.  {{In this paper, this definition often applies to the case
that 
{{$T(B)$ is}}  non-empty  and  {{generates}} 
{{${\tilde T}(B)$.}}
In these cases, $\kappa_T$ may be regarded as an affine map from
$T(B)$ to $T(A).$}}

\end{df}

\begin{df}[9.3 of \cite{GLII}]\label{DpropertyW}
%
Let $A$ be a separable \CA.  We say $A$ has property (W), if there is a map $T: A_+^{\bf 1}\setminus \{0\}\to \N\times \R_+\setminus \{0\}$ and
a sequence of approximately multiplicative \cpc s $\phi_n: A\to {\cal W}$ such
that, for any finite subset ${\mathcal H}\subset A_+^{\bf 1}\setminus \{0\},$
there exists an integer $n_0\ge 1$ such that
$\phi_n$ is exactly $T$-${\mathcal H}$-full (see  {{\ref{THfull} above}} and 5.5   and 5.7  of \cite{eglnp}) for all $n\ge n_0.$
\end{df}

\section{Comparison}

\begin{lem}\label{Lnbhd}
Let $A$ be a  \CA.  Let $a, b\in C([0,1], A)_+\setminus \{0\}$. {{Suppose, for some $1>d>0$ and $t\in [0,1]$,}}
 $a(t)\lesssim (b(t)-d)_+$.
Then, for any $d/4>\ep>0$ (with $\ep<\|a\|, \|b\|$),  there exists $\dt>0$ such that
\beq
f_\ep(a)|_{[t-\dt, t+\dt]\cap [0,1]}\lesssim b|_{[t-\dt, t+\dt]\cap [0,1]}
\eneq
as elements in $B:=C([t-\dt, t+\dt]\cap [0,1], A).$
Moreover,
there is $x\in B$ such that
\beq
((a-\ep)_+)|_{[t-\dt, t+\dt]\cap [0,1]}=x^*x\andeqn xx^*\in {{{\rm Her}(b|_{[t-\dt, t+\dt]\cap [0,1]})}}.
\eneq
\end{lem}

\begin{proof}
Fix $a$ and $b$  as above.
Let $d/4>\ep>0.$   By the assumption, there is {{$y\in A$}} such
that $\|y^*y-a(t)\|<\ep/16$ and $yy^*\in {\rm Her}((b(t)-d)_+).$ Note $f_{d/4}(b(t))yy^*=yy^* f_{d/4}(b(t)){{=yy^*}}.$
There is $\dt>0$ such that, for $t'\in [t-\dt, t+\dt]\cap [0,1],$
\beq
\|a(t)-a(t')\|<\ep/16\andeqn {{\|b(t)-b(t')\|<\ep/16.}}
\eneq
By Lemma 2.2 of \cite{Rr2} there are {{$r_1, r_2\in C([t-\dt, t+\dt]\cap [0,1],A)$,  
such that
$(a(t')-\ep/4)_+=r_1(t')^*y^*yr_1(t')$ and
$(b(t)-d/8)_+=r_2 ((b(t')-d/16)_+)r_2^*$ for all $t'\in  [t-\dt_1, t+\dt_1]\cap [0,1].$}}
{{Write $c=b(t)$.}} View $y$  and $c$
as a constant function in $C([0,1], A).$  One has,  in $B=C([t-\dt, t+\dt]\cap [0,1], A),$
\beq
(a-\ep/4)_+\lesssim y^*y\andeqn yy^*\lesssim {{f_{d/4}(c)}}\lesssim f_{d/16}(b)\,\,\,{\rm on}\,\,  [t-\dt, t+\dt]\cap [0,1].
\eneq
It follows that there is $z\in B$ such that
\beq
\|(a-\ep/2)_+-z^*z\|<\ep/16\andeqn zz^*\in {\rm Her}(f_{d/16}(b))\,\,\, {\rm on}\,\, [t-\dt, t+\dt]\cap [0,1].
\eneq
By Lemma 2.2 of \cite{Rr2} again, there is $r\in B$
such that
$
(a-\ep)_+=r^*(z^*z)r.
$
Choose $x=zr\in B.$ Then  $(a-\ep)_+=x^*x$ and $xx^*=zrr^*z^*\in {\rm Her}(b)$ on $[t-\dt, t+\dt]\cap [0,1].$
%
%
%
%
%

%
%
\end{proof}

\begin{lem}\label{Lunitary}
Let $A$ be a separable simple \CA\, which has stable rank one.
Let $\ep>0.$ Suppose that $a\in A_+\setminus \{0\}$ and ${\cal F}\subset {{\rm Her}}(a)$
is a finite subset.
Suppose also that $0$ is a limit point of $sp(a)$ 
and $u\in {\widetilde{{{\rm Her}}(a)}}$ is a unitary.
Then there exists a unitary $v\in U_0({\widetilde{{{\rm Her}}(a)}})$ such that
\beq
\|vb-ub\|\le \ep,\,\|v^*bv-u^*bu\|\le \ep \tand  \|bv-bu\|\le \ep
\eneq
for all $b\in {\cal F}.$
Moreover, if there exists $d\in A_+\setminus \{0\}$ such that $da=ad=0,$
The above holds with $\ep=0.$
\end{lem}

\begin{proof}
Choose $0<d<\ep/2$ such
that
\beq
\|f_d(a)b-b\|<\ep/2\rforal b\in {\cal F}.
\eneq
Since $0$ is a limit point of $sp(a)$, there exists $c\in {\rm Her}(a)_+\setminus \{0\}$ 
such
that $cf_d(a)=f_d(a)c=0.$  Since $A$ is simple, ${{\rm Her}}(c)\otimes {\cal K}\cong A\otimes {\cal K},$
by \cite{Br}. Since $A$ has stable rank one, so does ${{\rm Her}}(c)$  (see Corollary 3.6  of \cite{BP}).
Therefore there exists $u_1\in {\widetilde{{{\rm Her}}(c)}}$ such that
$[u_1]=[u]$ in $K_1(A).$  We may write $u_1=1_{{\widetilde{{{\rm Her}}(c)}}}+x$ for some $x\in {{\rm Her}}(c).$
Put $u_2=1_{{\widetilde{{\rm Her}(a)}}}+x.$
Again, since $A$ has stable rank one, we conclude that $u_2^*u\in U_0({{{\widetilde{{\rm Her}(a)}}}}).$
Choose $v=u_2^*u.$ One checks $v$ satisfies the conclusion.
The last part of the statement also follows.
\end{proof}

\begin{thm}\label{Tcomparison}
Let $A$ be a  non-unital separable {{projectionless}}
simple \CA\, which has  stable rank one and let $a, \, b\in C([0,1], A)_+\setminus \{0\}.$
Suppose that, for some $0<d<1/2,$
\beq
a(t)\lesssim f_d(b(t))\tforal t\in [0,1].
\eneq
Then $a\lesssim b$ in $C([0,1], A).$
\end{thm}

\begin{proof}
\Wlog, we may assume that $\|a\|, \|b\|=1.$

Fix $d/8>\ep>0.$  For each $t\in [0,1],$  choose $\dt(t)>0$ so that the conclusion
of \ref{Lnbhd} holds.  By a standard compactness {{argument}}, there is a $\dt>0$ such
that the conclusion of \ref{Lnbhd} holds for all $t\in [0,1].$
There exists a
 partition
 $0=t_0<t_1<\cdots < t_m=1$ such that
 \beq
 \|a(t)-a(t_i)\|<\ep/64, \|b(t)-b(t_i)\|<\ep/64\rforal t\in (t_{i-1}, t_{1+i}),\,\,\, i=1,2,...,m.
 \eneq
We may assume that $0<t_i-t_{i-1}=2T<\dt$ for all $i.$
Let  $I_i=[t_{i-1}, t_i]$ and $B_i=C(I_i, A).$

By Lemma \ref{Lnbhd}, we may also assume that, there are  $x_i\in B_i$  with
$\|x_i\|\le 1$ such
that
\beq\label{17827-10}
(a(t)-\ep/64)_+|_{I_i}=x_i^*x_i \andeqn x_ix_i^*\in {{{\rm Her}(f_d(b)|_{I_i})}}.
\eneq
Put $C_i=Her(f_d(b(t_i))).$
Then
\beq
&&(x_i(t_i))^*x_i(t_i)=(x_{i+1}(t_i))^*x_{i+1}(t_i)\andeqn\\ \label{October 2,2020}
&&x_i(t_i)(x_i(t_i))^*{{\in C_i,~~\mbox{and} ~}}\,x_{i+1}(t_i)(x_{i+1}(t_i))^*\in  C_{i}.
\eneq
Thus, by Lemma 2.4 of \cite{CES}, for any $\ep/256>\eta>0,$
there is a unitary $u_i\in {\widetilde C_i}$
such
that
\beq\label{17827-40}
\|x_i(t_i)-u_ix_{i+1}(t_i)\|<\eta.
\eneq
{{Since $A$ is projectionless, $0$ is a limit point of $f_{d}(b(t_i))$, if $f_{d}(b(t_i))\not=0$.}} By  Lemma \ref{Lunitary},
we may assume that $u_i\in U_0(\widetilde{{{\rm Her}}(f_{d}(b(t_i)))}).$
Therefore there exists a continuous path of unitaries
$\{u_i(t): t\in [t_i+T/4, t_i+T/2]\}$ such that $u_i(t_i+T/4)=u_i$ and $u_i(t_i+T/2)=1.$
Now define $r(t)$ as follows.
On $[t_0, t_1],$ define {{$r(t)=x_1(t).$}} 
On $[t_i, t_{i+1}],$
define
\beq
r(t)=\begin{cases} ({{\frac{t_i+(T/4)-t}{T/4}}})x_i(t_i)+
({{\frac{t-t_i}{T/4}}})u_ix_{i+1}(t_i), & \,\,\, t\in [t_i, t_i+T/4];\\\nonumber
                                        u_i(t)x_{i+1}{{(t_i),}}  & \,\,\,t\in (t_i+T/4, t_i+T/2];\\
                                        x_{i+1}({{\frac{4}{3}}}(t-(t_i+T/2))+
                                        t_i),
    &\,\,\,t\in (t_i+T/2, t_{i+1}].
                                        \end{cases}
                                        \eneq
                                        Note that  $r(t_i)=x_i(t_i)$ and $r(t_{i+1})=x_{i+1}(t_{i+1}).$
                                        {{Moreover, if $x_i(t_i)=0,$ then $r(t)=0$ for $t\in [t_i, t_i+T/2].$}}
                                         It follows that
$r(t)\in C([0,1], A).$
On $[t_0, t_1],$ $r(t)^*r(t)={{x_1(t)^*x_1(t)}}=(a(t)-\ep/64)_+.$
On $[t_i, t_i+T/4],$  using \eqref{17827-40},
\beq\nonumber
r(t)^*r(t) &=& ({{\frac{t_i+(T/4)-t}{T/4}}})^2(x_i(t_i))^*x_i(t_i)+
({{\frac{t-t_i}{T/4}}})^2(x_{i+1}(t_i))^*u_i^*u_ix_{i+1}(t_i)\\
&&+ ({{\frac{t_i+(T/4)-t}{T/4}}})({{\frac{t-t_i}{T/4}}})(x_i(t_i))^*u_ix_{i+1}(t_i)\\
&&+ ({{\frac{t_i+(T/4)-t}{T/4}}})({{\frac{t-t_i}{T/4}}})(x_{i+1}(t_i))^*u_i^*x_i(t_i)\\
&\approx_{{4\eta}}&({{\frac{t_i+(T/4)-t}{T/4}}})^2(x_i(t_i))^*x_i(t_i)+
({{\frac{t-t_i}{T/4}}})^2(x_i(t_i))^*x_i(t_i)\\
&&+({{\frac{t_i+(T/4)-t}{T/4}}})({{\frac{t-t_i}{T/4}}})(x_i(t_i))^*x_i(t_i)\\
&&+({{\frac{t_i+(T/4)-t}{T/4}}})({{\frac{t-t_i}{T/4}}})(x_{i}(t_i))^*x_i(t_i)\\
&=&(x_i(t_i))^*x_i(t_i).
\eneq
On $[t_i+T/4, t_i+T/2],$
\beq
r(t)^*r(t)=(x_{i+1}(t_i))^*x_{i+1}(t_i).
\eneq
On $[t_i+T/2, t_{i+1}],$
\beq
r(t)^*r(t)= (x_{i+1}({{\frac{4}{3}}}(t-(t_i+T/2))+t_i))^* x_{i+1}({{\frac{4}{3}}}(t-(t_i+T/2))+t_i).
\eneq
It follows that
\beq\label{17827-1}
\|(a-\ep/64)_+-r(t)^*r(t)\|<\ep/64.
\eneq
It follows that there exists $r_0\in C([0,1], A)$ such that
\beq\label{17827-2}
r_0^*r_0=(a-\ep/16)_+\andeqn r_0r_0^*\in {{\rm Her}}(r^*r).
\eneq
Define $c(t)\in C([0,1], A)$ as follow{{s:}}
On $[t_0, t_1],$ {{define}} $c(t)=b(t),$
{{on}} $[t_i, t_{i+1}],$ define
\beq
c(t)=\begin{cases} b(t_i)&\,\,\, t\in [t_i, t_i+T/2];\\
                              b\big({{\frac{4}{3}}}(t-(t_i+T/2))+t_i\big) &\,\,\, t\in [{{t_i+T/2}}, t_{i+1}].
                              \end{cases}
\eneq
We have that
\beq
\|c-b\|<\ep/64.
\eneq
There is $r_1\in C([0,1], A)$ such that
\beq\label{17827-3}
r_1^*r_1=f_{d/4}({{c}})\andeqn r_1r_1^*\in {{{\rm Her}}}(b).
\eneq
Now consider $r(t)r(t)^*.$
On $[t_0, t_1], $
\beq\label{17827-30}
r(t)r(t)^*={{x_1(t)(x_1(t))^*.}}
\eneq
On $[t_i, t_i+T/4],$
\beq\nonumber
r(t)r(t)^*&=&({{\frac{t_i+(T/4)-t}{T/4}}})^2x_i(t_i)(x_i(t_i))^*+
({{\frac{t-t_i}{T/4}}})^2u_ix_{i+1}(t_i))(x_{i+1}(t_i))^*u_i^*\\
&&+({{\frac{t_i+(T/4)-t}{T/4}}})({{\frac{t-t_i}{T/4}}})(x_i(t_i))(x_{i+1}(t_i))^*u_i^*\\
&&+({{\frac{t_i+(T/4)-t}{T/4}}})({{\frac{t-t_i}{T/4}}})u_ix_{i+1}(t_i))(x_i(t_i))^*.
\eneq
By {{\eqref{October 2,2020}}} 
\beq
x_i(t_i)x_i(t_i)^*,\, x_{i+1}(t_i)(x_{i+1}(t_i))^*\in {{\rm Her}}(f_d(b(t_i))).
\eneq
Since $u_i\in {\widetilde{{{\rm Her}}(f_d(b(t_i)))}},$
we conclude that
\beq
r(t)r(t)^*(f_d(b(t_i)))^{1/n}\to r(t)r(t)^*
\eneq
uniformly on $[t_i, t_i+T/4]$ as $n\to\infty$ {{(note   that,
when $(f_d(b(t_i)))^{1/n}=0,$ $r(t)=0$ for $t\in [t_i, t_i+T/4]$).}}
In other words,
\beq\label{17827-11}
r(t)r(t)^*\in {{\rm Her}}(f_d(c(t))\rforal t\in [t_i, t_i+T/4].
\eneq
On $[t_i+T/4, t_i+T/2],$  since $c(t)=b(t_i),$
\beq
r(t)r(t)^*=u_i(t)x_{i+1}(t_i)(x_{i+1}(t_i))^*(u_i(t))^*\in {{\rm Her}}(f_d(c(t))).
\eneq
On $[t_i+T/2, t_{i+1}],$
\beq\label{17827-12}
&&\hspace{-0.5in}r(t)r(t)^*=
(x_{i+1}({{\frac{4}{3}}}(t-(t_i+T/2))+t_i)) (x_{i+1}({{\frac{4}{3}}}(t-(t_i+T/2))+t_i))^*
\eneq
{{which is in  ${{\rm Her}}(f_d(c(t))).$}}
Combining \eqref{17827-30}, \eqref{17827-11} and \eqref{17827-12},
we obtain
that
\beq\label{17827-13}
rr^*\in {{\rm Her}}(f_d(c)).
\eneq
By \eqref{17827-2}, \eqref{17827-13} and \eqref{17827-3}, we conclude
that
\beq
(a-\ep/16)_+\lesssim b.
\eneq
It follows that $a\lesssim b.$
\end{proof}

\begin{lem}\label{Tcompa2}
Let $A$ be a non-unital separable projectionless  simple \CA\, with stable rank one and with
strictly comparison for positive elements
and let $a, b\in C([0,1], A)_+\setminus \{0\}.$  Suppose that $\overline{T(A)}^w$ is a nonempty compact set with $0\not\in \overline{T(A)}^w.$
If $d_\tau(a(t))<d_\tau(b(t))$ for all $t\in [0,1]$ and for all $\tau\in \overline{T(A)}^w,$
then
$$
a\lesssim b.
$$

\end{lem}

\begin{proof}
We may assume that $\|a\|=\|b\|=1.$
Let $1/2>\ep>0.$
Then, for any $\tau\in \overline{T(A)}^w$ and any $t\in [0,1],$
\beq
d_\tau((a(t)-\ep)_+)=d_\tau(f_{2\ep}(a(t)))\le \tau(f_{\ep}(a(t)))<d_\tau(b(t)).
\eneq
Note that $[0,1]\times \overline{T(A)}^w$ is compact. Moreover,
$d_\tau(b(t))-\tau(f_{\ep}(a(t)))$ is lower semicontinuous on $[0,1]\times \overline{T(A)}^w.$
Consider $b_n(t)=f_{1/n}(b(t)))$ {{and define $g_n\in C([0,1]\times \overline{T(A)}^w)$ by $g_n(t, \tau)=\tau(b_n(t))$.}}
Then
$$g_n(t,\tau)-\tau(f_{\ep}(a(t)))\nearrow d_\tau(b(t))-\tau(f_{\ep}(a(t)))\,\,\,({\rm as}\,\,n\to\infty).$$ Therefore, for some $n\ge 1,$
$g_n(t, \tau)>0$ for all $(t, \tau)\in [0,1]\times \overline{T(A)}^w.$
Therefore  there is $\ep/2>d>0$ such that
\beq
d_\tau(f_{\ep}(a(t)))< d_\tau(f_d(b(t)))\rforal t\in [0,1]\andeqn \tau\in \overline{T(A)}^w.
\eneq
Since $A$ has strict comparison for positive elements, {{for any $t\in [0,1]$, $f_{\ep}(a(t))\lesssim f_d(b(t))$.  By \ref{Tcomparison}, we have $f_{\ep}(a)\lesssim b$. Consequently,}}
\beq
a\lesssim b.
\eneq

\end{proof}

\begin{thm}\label{CMphi0phi1}
Let $A$ and $B$ be two non-unital separable projectionless simple {{\CA s}} with stable rank one and with strict comparison
for positive elements.  Suppose that ${{Ped}}(A)=A$ and ${{Ped}}(B)=B.$
Let $\phi_0, \phi_1: A\to B$ be two \hm s such that
they {{map}} strictly positive elements to strictly positive elements.

Suppose that $(g_1, a), (g_2, b)\in (M_{\phi_0,\phi_1})_+\setminus \{0\}$ {{(recall \eqref{dmapping} in Definition \ref{DfC1})}} such that
\beq
&&d_\tau(g_1(t))<d_\tau(g_2(t))\tforal \tau\in  \overline{T(B)}^w\tand \tforal t\in [0,1],\tand\\
&&d_t(a)<d_t(b)\tforal t\in \overline{T(A)}^w.
\eneq

Then $(g_1,a)\lesssim (g_2, b)$ in $M_{\phi_0, \phi_1}.$

\end{thm}

\begin{proof}
Let $1/2>\ep>0.$
As in the proof of  {{Lemma}} \ref{Tcompa2}, there exists {{$\ep/2>d>0$}} 
such that
\beq
d_\tau(f_{\ep}(g_1(t)))<d_\tau(f_d(g_2(t)))\rforal t\in [0,1] \andeqn \tau\in \overline{T({{B}})}^w,\andeqn\\
d_t(f_{\ep}(a))<d_t(f_d(b))\rforal t\in \overline{T(A)}^w.
\eneq
Since $A$ {{has  stable}} rank one {{with strict comparison
for positive elements}}, there exists $r_a\in A$ such that
\beq
f_{\ep}(a)=r_a^*r_a\andeqn r_ar_a^*\in {{\rm Her}}(f_d(b)).
\eneq
Note that
\beq
\phi_0(f_{\ep}(a))=f_{\ep}(g_1(0)),\,\,\, \phi_1(f_{\ep}(a))=f_{\ep}(g_1(1)),\\
\phi_0(f_d(b))=f_d(g_2(0))\andeqn \phi_1(f_d(b))=f_d(g_2(1)).
\eneq
Let $\min\{d/4, \ep/4\}>\eta>0.$
Then, {{by \ref{Tcompa2},}} there exists $r_b(t)\in C([0,1], B)$ such that
\beq\label{17827-2-0}
\|f_\ep(g_1)-r_b^*r_b\|<\eta^2/4\andeqn  r_br_b^*\in {{\rm Her}}(f_d(g_2)).
\eneq
Replacing $r_b$ by $r_bf_{\ep}(g_1)^{1/k}$ for some large $k,$ we may assume
also that $r_b^*r_b\in {{\rm Her}}(f_{\ep}(g_1)).$
In particular,
\beq\label{17827-2-1}
\|\phi_0(r_a)^*\phi_0(r_a)-r_b^*(0)r_b(0)\|<\eta^2/4\andeqn
\|\phi_1(r_a)^*\phi_1(r_a)-r_b^*(1)r_b(1)\|<\eta^2/4.
\eneq
Also
\beq\label{17827-2-2}
&&r_b(0)^*r_b(0)\in {{\rm Her}}(f_\ep(g_1(0)))=\overline{\phi_0(r_a)^*B\phi_0(r_a)}\andeqn\\
&&r_b(1)^*r_b(1)\in {{\rm Her}}(f_\ep(g_1(1)))=\overline{\phi_1(r_a)^*B\phi_1(r_a)}.
\eneq
Moreover, {{by \eqref{17827-2-0},}}
\beq\label{17827-2-3}
&&r_b(0)r_b^*(0), \phi_0(r_a)\phi_0(r_a)^*\in {{\rm Her}}(f_d(g_2(0)))\andeqn\\
&& r_b(1)r_b^*(1), \phi_1(r_a)\phi_1(r_a)^*\in {{\rm Her}}(f_d(g_2(1))).
\eneq
By applying {{Lemma 2.4}} of \cite{CES},
we obtain a unitary $u_i\in {\widetilde{{{\rm Her}}(f_d(g_2(i)))}}$ ($i=0,1$) such
that
\beq\label{17827-2-3+}
\|r_b(i)-u_i\phi_i(r_a)\|<\eta/2,\,\,\, i=0,1.
\eneq
{{Since $B$ is projectionless, by}} applying \ref{Lunitary}, we may assume that $u_i\in U_0( {\widetilde{{{\rm Her}}(f_{d/2}(g_2(i)))}}.$

Choose $1/8>\dt_0>0$ such that, if $|t-t'|<3\dt_0,$
\beq
\|g_i(t)-g_i(t')\|<\eta/4.
\eneq
Define
\vspace{-0.15in}\beq
g_3(t)=\begin{cases} g_2(0) &\,\,\, t\in[0, \dt_0],\\
                                  g_2( (t-\dt_0)/(1-2\dt_0)) &\,\,\, t\in [\dt_0, 1-\dt_0];\\
                                  g_2(1) &\,\,\, t\in [1-\dt_0, 1].
                                  \end{cases}
                                  \eneq
Then
\vspace{-0.15in}\beq
\|g_3-g_2\|<\eta/4.
\eneq
Note {{that}} $(g_3,b)\in M_{\phi_0, \phi_1}.$  Thus
\beq\label{17828-2-20}
\|{{(g_3,b)-(g_2, b)}}\|<\eta/4.
\eneq
It follows that there is $x\in  M_{\phi_0, \phi_1}$ such that
\beq\label{17828-2-21}
x^*x=f_{d/4}((g_3, b))=(f_{d/4}(g_3), f_{d/4}(b))\andeqn xx^*\in {{\rm Her}}((g_2,b)).
\eneq
There are continuous paths of unitaries $\{u_0(t): t\in [0, \dt_0/2]\}$ in ${\widetilde{{{\rm Her}}(f_{d/2}(g_2(0)))}}$ and
$\{u_1(t): t\in [1-\dt_0/2, 1]\}$  in ${\widetilde{{{\rm Her}}(f_{d/2}(g_2(1)))}}$ such that
$u_0(0)=1,$ $u_0(\dt_0/2)=u_0,$ $u_1(1-\dt_0/2)=u_1$ and $u_1(1)=1.$
Define  $y=(y_b, y_a)\in M_{\phi_0, \phi_1}$ as follows
$y_a=r_a.$
\beq\label{17827-2-6}
y_b(t)=\begin{cases} u_0(t) \phi_0(r_a) & \,\,\, t\in [0, \dt_0/2];\\
                                  ({{\frac{\dt_0-t}{\dt_0/2}}})u_0\phi_0(r_a)
                                  +({{\frac{t-\dt_0/2}{\dt_0/2}}})r_b(0) &\,\,\, t\in (\dt_0/2, \dt_0];\\
                                  r_b((t-\dt_0)/(1-2\dt_0)) &\,\,\, t\in (\dt_0, 1-\dt_0];\\
                                  ({{\frac{{t-1+\dt_0}}{\dt_0/2}}})u_1\phi_1(r_a)
                                  +({{\frac{1-t-\dt_0/2}{\dt_0/2}}}) r_b(1) &\,\,\, t\in (1-\dt_0, 1-\dt_0/2];\\
                                  u_1(t)\phi_1(r_a)&\,\,\, t\in (1-\dt_0/2, 1].
                                  \end{cases}
\eneq
We check that $y\in M_{\phi_0, \phi_1}.$
On $[0, \dt_0/2],$
\beq\label{17828-2-7}
y_b(t)^*y_b(t) =\phi_0(r_a)^*\phi_0(r_a).
\eneq
On $[\dt_0/2, \dt_0],$ by \eqref{17827-2-0} and \eqref{17827-2-3+},
\beq
y_b^*(t)y_b(t) &=&  ({{\frac{\dt_0-t}{\dt_0/2}}})^2\phi_0(r_a)^*\phi_0(r_a)+
({{\frac{t-\dt_0/2}{\dt_0/2}}})^2r_b(0)^*r_b(0)\\
                      &&+({{\frac{\dt_0-t}{\dt_0/2}}})
                      ({{\frac{t-\dt_0/2}{\dt_0/2}}})\phi_0(r_a)^*u_0^*r_b(0)\\
                         &&+({{\frac{\dt_0-t}{\dt_0/2}}})
                         ({{\frac{t-\dt_0/2}{\dt_0/2}}}) r_b(0)^*u_0\phi_0(r_a)\\
                         &\approx_{\eta}& ({{\frac{\dt_0-t}{\dt_0/2}}})^2r_b(0)^*r_b(0)+
                         ({{\frac{t-\dt_0/2}{\dt_0/2}}})^2r_b(0)^*r_b(0)\\
                          &&+({{\frac{\dt_0-t}{\dt_0/2}}})
                          ({{\frac{t-\dt_0/2}{\dt_0/2}}})r_b(0)^*r_b(0)
                             +({{\frac{\dt_0-t}{\dt_0/2}}})
                             ({{\frac{t-\dt_0/2}{\dt_0/2}}}) r_b(0)^*r_b(0)\\\label{17828-2-8}
                             &=& r_b(0)^*r_b(0).
                         \eneq
On $[\dt_0, 1-\dt_0],$
\beq\label{17828-2-9}
y_b^*(t)y_b(t)=r_b^*((t-\dt_0)/(1-2\dt_0))r_b((t-\dt_0)/(1-2\dt_0)).
\eneq
On $[1-\dt_0, 1-\dt_0/2],$ as on $[\dt_0/2, \dt_0],$
\beq\label{17828-2-10}
y_b^*(t)y_b(t) &\approx_{\eta/2} r_b(1)^*r_b(1).
\eneq
On $[1-\dt_0/2, 1],$
\beq\label{17828-2-11}
y_b^*(t)y_b(t)=\phi_1(r_a)^*\phi_1(r_a).
\eneq
Combining \eqref{17828-2-7}, \eqref{17828-2-8}, \eqref{17828-2-9}, \eqref{17828-2-10}, \eqref{17828-2-11}
and \eqref{17827-2-0},  {{we obtain}}
\beq\label{17828-2-12}
\|y^*y-f_{\ep}(g_1,a)\|<2\eta.
\eneq
We also compute that
$y_b(t)y_b(t)^*\in {{\rm Her}}(f_{d/2}(g_3)).$
It follows that
\beq\label{17828-2-14}
yy^*\in {{\rm Her}}(f_{d/2}(g_3, b)).
\eneq
Let $z=xy.$
Then, by \eqref{17828-2-14} and \eqref{17828-2-21},
\beq
z^*z=y^*x^*xy=y^*f_{d/4}((b, g_3))y=y^*y\andeqn zz^*=xyy^*x^*\in {{\rm Her}}((b, g_2)).
\eneq
Thus
$\|z^*z-f_\ep((a, g_1))\|<2\eta.$
This implies that $f_\ep((g_1, a))\lesssim (g_2, b)$ for all $1/2>\ep>0.$
It follows that $(g_1, a)\lesssim (g_2, b).$
\end{proof}

%
%
%
%
%
%

\begin{lem}\label{Lasrunitary}
Let $A$ be a non-unital  but $\sigma$-unital \CA\, and let $u\in M(A).$
Suppose that $A$ has almost stable rank one.
Then, for any finite subset ${\cal F}\subset A$ and $\ep>0,$
there is a unitary $v\in {\tilde A}$ such that
\beq
vav^*\approx_{\ep} uau^* \rforal a\in {\cal F}.
\eneq

\end{lem}

\begin{proof}
\Wlog, we may assume that there
is
$e_1\in A$ with $0\le e_1\le 1$
such that
$e_1x=xe_1=x$ for all $x\in {\cal F}.$
We may also assume that ${\cal F}\subset A^{\bf 1}.$
Choose $\dt>0$ such that, for any pair of positive elements  $a, b\in C^{\bf 1}$
(in any \CA\, $C$), $\|a^{1/2}-b^{1/2}\|<\ep/4$
if $\|a-b\|<\dt.$ Put $\dt_1=\min\{\dt,\ep\}/4.$

Let $Z=ue_1^{1/2}\in A.$ Since $A$ has almost stable rank one, there is an invertible element $z\in {\tilde A}$
such that
$
\|z-Z\|<\dt_1/4.
$
Note that, for $x\in {\cal F},$
\beq\nonumber
zxz^*\approx_{\dt_1/2}ue_1^{1/2}xe_1^{1/2}u^*
=uxu^*.
\eneq
Also
$z^*z\approx_{\dt_1/2}e_1^{1/2}u^*ue_1^{1/2}=e_1.$
It follows that
\beq
\|(z^*z)^{1/2}-e_1^{1/2}\|<\ep/4.
\eneq
Write $z=v(z^*z)^{1/2}$ as a polar decomposition of $z$ in ${\tilde A}.$
One estimates
that, for all $x\in {\cal F},$
\beq\nonumber
vxv^*&=&ve_1^{1/2}xe_1^{1/2}v^*\approx_{\ep/2} v(z^*z)^{1/2}x(z^*z)^{1/2}v^*\\\nonumber
&=&zxz^*\approx_{\dt_1} uxu^*.
\eneq
\end{proof}
%
%
%
%
%
%
%
%
%
%
%
%
%
%
\vspace{-0.3in}
\section{The determinant map}

\begin{df}\label{Dtidelf}
Let $B$ be a \CA\, with $T(B)\not=\emptyset.$
For each $f\in \Aff(T(B)),$ define $\iota_B^\sharp: \Aff(T(B))\to \Aff(T(\td B))$ by
$\iota_B^\sharp(f)(\af t_\C+(1-\af) \tau)=(1-\af)f(\tau)$
for all $0\le \af\le 1,$ $\tau\in T(B)$ and
$t_\C\in T({\tilde B})$ such that $t_\C|_B=0.$
Put
$$
\Aff(T(B))^{\iota}=\{\iota_B^\sharp(f): f\in \Aff(T(B))\} \subset \Aff(T(\td B)).
$$
Therefore, if $B={\rm Ped}(B),$ then  $\rho_{\td B}(K_0(\td B))=\iota_B^\sharp(\rho_B(K_0(B)))+\Z.$
\end{df}

\begin{prop}\label{PK0divisible}
Let $A\in {\cal D}^d$ with continuous scale. Then, for any
integer $n\ge 1$  and $g\in \overline{\rho_A(K_0(A))},$
there is $g_1\in \overline{\rho_A(K_0(A))}$ such that $ng_1=g.$
Moreover, $\overline{\rho_{\td A}(K_0(\td A))}/\Z$ is divisible subgroup of $\Aff(T(\td A))/\Z$
and $\Aff(T(A))^\iota\cap \overline{\rho_{\td A}(K_0(\td A))}={\overline{\iota_A^\sharp(\rho_A(K_0(A)))}}=
\iota_A^\sharp(\overline{\rho_A(K_0(A))}).$
{{Furthermore,}} $\Aff(T(A))^\iota/\overline{\rho_{\td A}(K_0(\td A))}=\Aff(T(A))^\iota/\iota_A^\sharp(\overline{\rho_A(K_0(A))})$
is torsion free.

\end{prop}

\begin{proof}
To {{show}} the first part of the proposition,
{{it suffices to show that
it holds for $g\in \rho_A(K_0(A)).$}}

Let $g=[p]-[{\bar p}],$ where $p={\bar p}+x\in M_m(\tilde A)$ is  a projection,
${\bar p}\in M_m(\C \cdot 1_{\tilde A})\subset M_m({\tilde A})$
is a scalar matrix, and $x\in  {{M_m(A)_{s.a.}.}}$
Fix an integer $n\ge 1.$ Let $1/4m>\ep>0.$

Write $x=(x_{i,j})_{m\times m},$ where $x_{i,j}\in A.$
Put ${\cal F}_1=\{x_{i,j}: 1\le i, j\le m\}$ and ${\cal F}={\cal F}_1\cup{{\{xy: x, y\in {\cal F}_1\}}}.$
Note that $\|x_{i,j}\|\le 2.$

Since $A\in D^d,$ there are ${\cal F}$-$(\ep/64m)^2$-{{multiplicative}} \cpc s $\phi: A\to A$ and $\psi_1: A\to D_1$  with
$\phi(A)\perp D,$  $D:=M_n(D_1)$ is a \SCA\, of $A,$  and
\beq
y\approx_{\ep/2(64m)^2} \phi(y)+\psi(y)\andeqn
\psi(y)=\psi_1(y)\otimes 1_n
\eneq
for all $y\in {\cal F},$ and
\beq\label{Pk0div-1}
d_\tau(c)<\ep/16m\rforal \tau\in T(A),
\eneq
where $c$ is a strictly positive element of $\overline{\phi(A)A\phi(A)}.$
Let $e_d\in (D_1)^{\bf 1}_+$ be a  strictly positive element of $D_1.$
Let $\eta>0.$ Let $A_1=\overline{f_\eta(c)Af_\eta(c)}.$
Choosing a sufficiently small $\eta,$
and replacing $\phi(a)$ by $\phi'(a):=f_\eta(c)\phi(a)f_\eta(c)$ for all $a\in A,$
$\psi_1(a)$ by $\psi_1'(a):=f_\eta(e_d)\psi_1(a)f_\eta(e_d)$ for all $a\in {{A}}$
 if {{necessary,}}  we may
 write
 \beq
 y\approx_{\ep/(64m)^2} \phi'(y)+\psi'(y)\rforal y\in {\cal F}
 \eneq
 and
 assume
that
there are $f_0\in ({\overline{\phi(A)A\phi(A)}})^{\bf 1}_+$ and $f_1\in (D_1)_+^{\bf 1}$
such that
\beq
f_0\phi'(a)=\phi'(a)f_0\andeqn (f_1\otimes 1_n)\psi'(a)=\psi'(a)(f_1\otimes 1_n)\rforal a\in A.
\eneq

{{Set}} $x_{i,j,k}=\psi_1'(x_{i,j})\otimes e_k\in M_n(\overline{f_\eta(e_d)D_1f_{\eta}(e_d})),$
where $e_k:=e_{k,k}$ and $\{e_{k,k'}: 1\le k, k'\le {{n}}\}$ is a system of matrix unit for $M_n.$
(Note that $D_1$ is non-unital and $e_k\not\in M_n(D_1).$) Put $x_k=(x_{i,j,k})\in M_m(A)$  and
$x'=(\phi'(x_{i,j}))_{m\times m}\in M_m(A).$  Note {that,} since both
$\phi'$ and $\psi'$ are \cpc s, $x'$ and $x_k$ are self-adjoint.
Recall that, since $p$ is a projection,
$x{\bar p}+{\bar p}x +x^2=x.$
Let $y_0={\bar p}+x'.$
Then
\beq
&&(y_0+\sum_{k=1}^nx_k)\approx_{\ep/32^2}
(y_0+\sum_{k=1}^nx_k)^2\\
&&\approx_{\ep/64}{\bar p}^2+x'{\bar p}+{\bar p}x'+
x'^2+
\sum_{k=1}^n(x_k{\bar p}+{\bar p}x_k+x_k^2).
\eneq
It follows that
\beq
z:=x'{\bar p}+{\bar p}x'+
x'^2+
\sum_{k=1}^n(x_k{\bar p}+{\bar p}x_k+x_k^2)\approx_{\ep/32} x.
\eneq
Let $E_0=\diag(f_0, f_0,...,f_0),$
${\bar E}={{\diag(f_1\otimes 1_n, f_1\otimes 1_n, {{...,}} f_1\otimes 1_n)}},$
and $E_k=\diag(f_1\otimes e_k, f_1\otimes e_k,...,f_1\otimes e_k)$ {{($1\le k\le n$) be  $n+2$ many}}
 $m\times m$ diagonal elements in $M_m(A).$
Then
\beq
x'\approx_{(\ep/64m)^2}E_0x\approx_{\ep/32}E_0z&=&E_0x'{\bar p}+E_0{\bar p}x'+E_0x'^2+E_0(\sum_{k=1}^n {\bar p}x_k)\\
&=& x'{\bar p}+{\bar p}E_0x'+x'^2+\sum_{k=1}^n{\bar p}E_0x_k\\
&=& x'{\bar p}+{\bar p}x'+x'^2.
\eneq
Therefore
\beq
{\bar p}+x'\approx_{\ep/16} ({\bar p}+x')^2.
\eneq
A standard perturbation argument produces an element
$y'\in M_m(\overline{\phi(A)A\phi(A)})_{s.a.}$
such that $q_0:={\bar p}+y'$ is a projection and
\beq
\|q_0-({\bar p}+x')\|<\ep/4\andeqn \|x'-y'\|<\ep/4.
\eneq
The same argument shows that, for each $k,$
there is a self-adjoint element $y''\in M_m(D_1)\subset M_m(A)$
such that $q_k:={\bar p}+y''\otimes {{e_k}}$ is a projection {{($k=1,2,...,n$)}} and
\beq
\|q_k-({\bar p}+x_k)\|<\ep/4\andeqn \|y''\otimes e_k-x_k\|<\ep/4.
\eneq
We note
that
\beq
\|x-(y'+\sum_{k=1}^n y''\otimes e_k)\|<\ep/4.
\eneq
Combining with \eqref{Pk0div-1}, we have
\beq
\tau(x)\approx_\ep n\tau(y'')\rforal \tau\in T(A).
\eneq
Fix $k$ and choose $g_0=\rho_A([q_k]-[{\bar p}]).$
Then
\beq
|(\tau(p-{\bar p})-ng_0(\tau)|<\ep\rforal \tau\in T(A).
\eneq
Since $\ep$ is arbitrary, this shows that there is $g_1\in {\overline{\rho_A(K_0(A))}}$
such that $ng_1=g.$
The first part of the proposition then follows.

To see the second part, let $g\in \overline{\rho_{\td A}(K_0(\td A))}=\overline{\iota_A^\sharp(\rho_A(K_0(A)))+\Z}$
and $k\ge 2$ be an integer.
There are $f_n\in \iota_A^\sharp(\rho_A(K_0(A)))$ and $N(n)\in \Z$
such that $\lim_{n\to\infty}\|f_n-N(n)-g\|=0.$  From the first part of the proof, there are $g_n\in \iota_A^\sharp(\rho_A(K_0(A))$
such that $\lim_{n\to\infty}\|kg_n-f_n\|=0.$  Hence
$\lim_{n\to\infty}\|kg_n-N(n)-g\|=0.$  There are integers $D(n), r_n\in \Z$ with $0\le r_n<k$
such that {{$N(n)=kD(n)+r_n,$}} $n=1,{{2,....}}$ 
Hence  in $\Aff(T(\td A)),$
\beq
g/k=\lim_{n\to\infty} g_n-D(n)-r_n/k.
\eneq
By passing to a subsequence,  we may assume that, for some integer $0\le r<k,$
\beq
g/k+r/k=\lim_{n\to\infty} (g_n-D_n).
\eneq
Note that $\Z$ is a closed subgroup of $\Aff(T(\td A)).$
In the topological group $\Aff(T(\td A))/\Z,$ {$g/k+r/k$ gives an element $g_0:=[g/k+r/k]\in \overline{\iota_A^\sharp(\rho_A(K_0(A)))+\Z}/\Z$, where we {{temporarily}} use $[x]$ to denote the corresponding equivalent class of $x$ modulo $\Z$.} 
{{It}} follows that $ [g]=kg_0-k[r/k]=k{ g}_0,$ 
as desired.
(Warning: we do not know $g_0$ 
is in $\overline{\iota_A^\sharp(\rho_A(K_0(A)))}.$)
%

Now let $f\in \Aff(T(A))^\iota\cap \overline{\rho_{\td A}(K_0(\td A))}.$
Recall {{that}} $\overline{\rho_{\td A}(K_0(\td A))}=\overline{\iota_A^\sharp(\rho_A(K_0(A)))+\Z}.$
It follows that there are $g_n\in \iota_A^\sharp(\rho_A(K_0(A)))$ and $m(n)\in \Z$
such that $f(\tau)=\lim_{n\to\infty} g_n(\tau)+m(n)$ for all $\tau\in T(\td A).$
In particular, $0=f(\tau_\C)=\lim_{n\to\infty} g_n({{\tau_\C}})+m(n)=\lim_{n\to\infty} m(n).$
Thus, for some $N\ge 1,$ $m(n)=0$ for all $n\ge N.$  This implies that $g_n\to f$
uniformly on $T(\td A).$ In other words,
$f\in \overline{\iota_A^\sharp(\rho_A(K_0(A)))}.$  Since $T(A)\subset T(\td A),$
this implies that $\overline{\iota_A^\sharp(\rho_A(K_0(A)))}=\iota_A^\sharp(\overline{\rho_A(K_0(A))}).$

To see the  last part,  one notes that $\Aff(T(A))^\iota$ is torsion {{free}} and,  by the first part,
$\iota_A^\sharp(\overline{\rho_A(K_0(A))})$ is divisible.
If $x\in \Aff(T(A))^\iota$ and $y:=nx\in \iota_A^\sharp(\overline{\rho_A(K_0(A))})$ for some integer $n>1,$
then there exists $z\in \iota_A^\sharp(\overline{\rho_A(K_0(A))})$ such that $y=nz.$ Then
$n(x-z)=0.$ It follows that $x-z=0,$ or $x\in  \iota_A^\sharp(\overline{\rho_A(K_0(A))}).$
\end{proof}

\begin{NN}\label{DdertF}
Recall that ${\cal M}_1$ is the class of separable  stably projectionless simple
\CA s with continuous scale constructed in Theorem 4.31 or Remark 4.32 of \cite{GLrange} (see also  Definition 4.33 of \cite{GLrange}).
Let $B\in {\cal M}_1$ be as in (3) of Remark 4.32 of \cite{GLrange} written as $B=\lim_{{m\to \infty}}(B_m\oplus C_m\oplus D_m, \phi_{m,m+1})$,  with $B_m=B_1$.
Fix a finitely generated subgroup ${{F_1}}\subset K_1(B).$
We may write ${{F_1}}\subset K_1(B_m)$ and
$B_m=\lim_{k\to\infty} (M_{(k!)^2}A(W, \af_k)\oplus W_k, \Phi_{k,k+1})$
{{(see the statement of Theorem 4.31, (3) of Remark 4.32}} and description in 4.30 of \cite{GLrange}).
In particular, in Lemma \ref{L215}
below.
{{Let}} $E_n=M_{(n!)^2}(A(W, \af_n))$
{{be}}  as   $E_n$  in 11.4 of \cite{GLII}.   We will retain the notation in 11.3 and 11.4 of \cite{GLII}.
{{Put $\iota_n:=\phi_{1, \infty}\circ \Phi_{n, \infty}: E_n\to B.$}}
We may assume that ${{F_1\subset F}}\subset (\iota_n)_{*1}(K_1(A(W, \af_n)))$
for some $n\ge 1,$  {{where $F$ is a standard subgroup of $K_1(B)$ as mentioned in
11.3 of \cite{GLII}.}}
Write  $F={{\Z^{m_f}}}\oplus \Z/{k_1}\Z\oplus \Z/k_2\Z\oplus \cdots {{\Z/k_{m_t}\Z}}.$

Fix such $F.$ Let $J_{F,u}: F\to U({\tilde B})/CU({\tilde B})$ {{be}} given by the splitting map $J:=J_{cu}^{{{\td B}}}: K_1({\tilde B})\to
U({\tilde B})/CU({\tilde B})$ defined in \eqref{Dcu-4}.  {{We retain these notation for Lemma \ref{L215}.}}

%
%

\end{NN}


\begin{lem}\label{L215}
Let $C$ be a non-unital separable simple \CA\, in ${\cal D}^d$  with  continuous scale
and let
$B$ be as {{(3) in remark 4.32 (see Theorem 4.31 and Theorem 4.34 also)}} of \cite{GLrange}. 

Let $\ep>0, $ ${\cal F}\subset B$  {{and}}
${\cal P}\subset \underline{K}(B)$ be  {{finite subsets,}}  and let $1/2>\dt_0>0.$

For any finitely generated standard subgroup
$F$ (see \ref{DdertF}),  any finite subset $S\subset F,$
there exists an integer $n\ge 1$  with the following property:
for
any finite  subset  ${\cal U}\subset  U({\tilde B})$ such that
${\overline{\cal U}}\subset
J_{F,u}(F)\subset J_{F,u}({{{\iota_n}}}_{*1}(K_1(E_n)))$ and {{${{\Pi_{cu}^{\td B}}}(\overline{{\cal U}})=S$
(see \ref{DkappaJ})}}
{{for}}
any \hm\, \\ $\gamma: J_{F,u}({{(\iota_n)_{*1}}}(K_1(E_n)))\to \Aff(T({{\tilde C}}))/\overline{\rho_{\tilde C}(K_0({\tilde C})},$
such that $\gamma|_{{\rm Tor}(J_{F,u}((\iota_n)_{*1}(K_1(E_n)))}=0,$ and
any $c\in C_+$ with $\|c\|=1,$
there exists ${\cal F}$-$\ep$-multiplicative \cpc\,
$\Phi:  B
\to \overline{cCc}$
such that, in $U_0({\tilde C})/CU({\tilde C}){{\cong \Aff(T({\tilde C}))/\overline{\rho_{\tilde C}(K_0({\tilde C}))}}},$
\beq\label{L215-1}
[\Phi]|_{\cal P}=0\tand
{\rm dist}(\Phi^{\dag}({\bar z}),\gamma({\bar z}))<\dt_0\tforal z\in {\cal U}.
\eneq

({{Here}} we assume $\dist(\Phi^\dag({\bar z}), \overline{\lceil \Phi(z)\rceil})<\dt_0/4$
for all $z\in {\cal U}$--see \ref{DLddag} for the definition of $\Phi^{\dag}.$)

\end{lem}

\begin{proof}
{{As $B_n=B_1$ for all $n\geq 1$ in the assumption of (3) of Remark 4.32 of \cite{GLrange},}} we may {{further}} assume
{{$B=B_1$ since $K_1(C_n\oplus D_n)=\{0\}$}}. 
We will then reduce the lemma to 11.5 of \cite{GLII}.
Let  $C_1$ be a \CA\, {{which is an inductive limit of \CA s in ${\cal C}_0^{{(0)}}$}}
 with continuous scale such that $K_0(C_1)=\{0\}$ and
$T(C_1)=T(C).$  It follows from \cite{Rl} that there is a \hm\,
$j: C_1\to C$ which maps strictly positive elements to strictly positive elements and {{induces an
affine isomorphism}}
$j^{\sharp}: \Aff(T(C_1))\to \Aff(T(C)).$
Moreover, the map $j_T: T(C)\to T(C_1)$  is an affine homeomorphism.
%
%
%
Note that $\rho_{{\tilde C}_1}(K_0({\tilde C_1}))=\Z$ and (see \ref{Dtidelf})
\beq
\rho_{\tilde C}(K_0(\tilde C))=\iota_C^\sharp(\rho_C(K_0(C)))+\Z.
\eneq
Let $j_{\dag, o}: \Aff(T({\tilde C_1}))/\Z\to \Aff(T({\tilde C}))/\overline{\rho_{\tilde C}(K_0({\tilde C}))}$
be the map induced by $j^{\sharp}.$
{{Note that
\beq
{\rm ker}(j_{\dag,o})={j^{\sharp}}^{-1}(\overline{\{\iota_C^\sharp\circ \rho_C(x) \in \Aff(T({\tilde C})): x\in K_0(C)\}}),
\eneq
which is a divisible group by Proposition \ref{PK0divisible}.
Therefore
there {{exists  a \hm\,}}\\ $
{{\gamma_0:}} J_{F,u}(F)\subset J_{F,u}({{\iota_n}}_{*1}(K_1(E_n)))
\to \Aff(T(C_1))/\Z$ such that $j_{\dag,0}\circ 
{{\gamma_0=\gamma}}$.}}
{{Thus,}} it suffices to prove the lemma under the assumption
that $C=C_1.$ But that is exactly the same as 11.5 of \cite{GLII}.
\end{proof}

\begin{lem}[11.6 of \cite{GLII}]\label{Lderttorsion}
Let $C$ be a
non-unital separable \CA.
Suppose that $u\in U(M_s({\tilde C}))$ (for some integer $s\ge 1$) with $[u]\not=0$ in $K_1(C)$
but $u^k\in CU(M_s({\tilde C}))$
 for some $k\ge 1.$
 Suppose that $\pi_C(u)=e^{2{{\sqrt{-1}}} \pi \theta}$
 for some $\theta\in (M_s)_{s.a.},$  where
$\pi_C: {\tilde C}\to \C$ is the quotient map.
Then $k{\rm tr}(\theta)\in \Z,$ where ${\rm tr}$ is the tracial state
of $M_s.$

Let $B_1$ be a stably projectionless simple separable \CA\,
and with continuous scale and $B=B_1\otimes U$ for
some infinite dimensional UHF-algebra, or $B\in {\cal D}^d.$
For any $\ep>0,$ there exists $\dt>0$ and finite subset ${\cal G}\subset C$ satisfying the following:
If $L_1, L_2: C\to B$ are two ${\cal G}$-$\dt$-multiplicative \cpc s such that $[L_1](u)=[L_2](u)$
in $K_1(B),$
then
\beq\label{Lderttor-1}
{\rm dist}(\overline{\lceil L_1(u)\rceil}, \overline{\lceil L_2(u)\rceil} )<\ep.
\eneq


\end{lem}


\begin{proof}
Write  $u=e^{2{{\sqrt{-1}}}\pi \theta}+\zeta,$  where $\zeta\in {{M_s(C)}}$
and $\theta\in (M_s)_{s.a.}.$
Therefore, if $u^k\in CU({{M_s({\tilde C})}}),$ then $k{\rm tr}(\theta)\in \Z.$

Note {{that}} $L_i$ is originally defined on $C$ and the extension
$L_i: M_s({\tilde C})\to M_s({\tilde B})$ has the property  that
$L_i(u)=e^{2{{\sqrt{-1}}} \pi \theta}+L_i(\zeta),$  $i=1,2.$
To simplify notation,
\wilog, we may assume that $\lceil L_1(u)\rceil \cdot \lceil L_2(u^*)\rceil\in U_0(M_s({\tilde B})).$  Note that $$\pi_B(\lceil L_1(u)\rceil \cdot \lceil L_2(u^*)\rceil)=e^{2\sqrt{-1}  \pi \theta } e^{-2\sqrt{-1}  \pi \theta}=1$$
{{(where $\pi_B: M_s({\tilde B})\to M_s$  is the quotient map).}}
We may write (see Lemma 6.1 of \cite{GLII})
\beq\nonumber
\lceil L_1(u)\rceil \cdot \lceil L_2(u^*)\rceil &=&\prod_{j=1}^n \exp(2{{\sqrt{-1}}}  \pi h_j) \,\,\,\,\text{{{for}}\,\, some}\,\,h_1,h_2,...,h_n\in M_s({\tilde B})_{s.a}\,\,\text{with}\\
\label{1802}
\vspace{-0.1in}
\pi_B(h_j)&=&0\,\,\,~\mbox{and}~\pi_B( \exp(2\sqrt{-1} \pi h_j))=1 ~{\mbox{ for all}~j.}
\eneq

Recall $u^k\in CU(M_s({\tilde C})).$ {{It}} follows from 14.5  {{of \cite{LinLAH}}} {{(applied to $\lceil L_1(u^k)\rceil$ and $\lceil L_2((u^*)^k)\rceil $ separately)}}
that, by choosing small $\dt$ and large ${\cal G}$
(independent of $L_1$ and $L_2$)
there is $h_0\in {\tilde B}_{s.a.}$ such that
$\|h_0\|<{{\min\{1, \ep\}}}/2(k+1)$ and
\beq\label{e84-1}
((\exp(2{{\sqrt{-1}}}\pi h_0))(\prod_{j=1}^n \exp(2{{\sqrt{-1}}} \pi h_j)))^k\in CU(M_s({\tilde B})).
\eneq
{{By \eqref{1802}, $\pi_B(\exp(2i h_0))\in CU(M_s).$
Then $s\cdot t_\C^B(h_0)\in \Z,$
where $t_\C^B\in T({\tilde B})$ is defined by
$t_\C^B(b)={\rm tr}\circ \pi_B(b)$ for all $b\in {\tilde B}.$
However, since $\|h_0\|<1/4s(k+1),$
$t_\C^B(h_0)<1/4s(k+1).$
This implies that $t_\C^B(h_0)=0.$}}
Note also {{that}} $U_0(M_s({\tilde B}))/CU(M_s({\tilde B}))=\Aff(T({\tilde B}))/\overline{\rho_{\tilde B}(K_0({\tilde B}))}$ and
\beq
\overline{\rho_{\tilde B}(K_0({\tilde B}))}=\overline{\iota_B^\sharp(\rho_{B}(K_0(B))+\Z}
\eneq
(see  \ref{Dtidelf}   for $\iota_B^\sharp$).
Note, by  \eqref{e84-1}, $k(\sum_{j=1}^n h_j+h_0/k)^{\widehat{}} \in \Aff(T(B))^\iota\cap \overline{\rho_{\tilde B}(K_0({\tilde B}))}.$
It follows from \ref{PK0divisible}  that $g:=(\sum_{j=1}^n h_j+h_0/k)^{\widehat{}}\in \Aff(T(B))^\iota\cap \overline{\rho_{\tilde B}(K_0({\tilde B}))}.$
%
%
%
%
%
Therefore, since $g\in \overline{\rho_{\tilde B}(K_0(B))}$ and $\|h_0\|<\ep/2(k+1),$
\beq
{\rm dist}(\overline{\lceil L_1(u)\rceil}, \overline{\lceil L_2(u)\rceil} )<\ep.
\eneq

\end{proof}

We actually prove the following:

\begin{cor}\label{CCUtor}
Let $B_1$ be a stably projectionless simple separable \CA\,
 with continuous scale and $B=B_1\otimes U$ for
some infinite dimensional UHF-algebra, or $B\in {\cal  D}^d.$ Let $u=1_{M_s(\td B)}+x\in U_0(M_s(\td B)),$
where $x\in M_s( B)_{s.a.}$ such that $u^k\in CU(M_s(\td B)).$
Then $u\in CU(M_s(\td B)).$
\end{cor}

{{We end this section with the following lemma  for
the convenience in later sections.}}

\begin{lem}\label{Unitary}
Let $B$ be a non-unital separable simple \CA\, with stable rank one,
let
$\ep>0,$  {{${\cal F}\subset B$}} be a finite subset and $u\in {\tilde B}$ be a unitary.
Then there exists a unitary $w\in {\tilde B}$ such that
\beq
\|w^*aw-a\|<\ep\rforal a\in {\cal F}\andeqn wu\in CU({\tilde B}).
\eneq
In particular, $wu\in U_0({\tilde B})$ and $\overline{w^*}=\overline{u}$ in $U({\tilde B})/CU({\tilde B}).$

\end{lem}

\begin{proof}
We may assume that $0<\ep<1.$ \Wlog, we may assume that
$\|y\|\le 1$ for all ${{y}}\in {\cal F}.$
Consider an approximate identity $\{e_n\}$ of $B$ such that
$e_{n+1}e_n=e_n$ for all $n.$
 Write $u=\lambda\cdot 1_{\tilde B}+x$ for some $\lambda\in \T$ and $x\in B.$
Note $\|x\|\le 2.$ Replacing $u$ by ${\bar \lambda}u,$ if necessary, we may assume that $\lambda=1.$

Choose $n\ge 1$ such
that
\beq
\|(1_{\tilde B}-e_n){{y}}\|<\ep/{{128}}\andeqn \|{{y}}(1_{\tilde B}-e_n)\|<\ep/128\rforal {{y}}\in {\cal F}\cup \{{{x}}\}.
\eneq
Let $B_0=\overline{e_nAe_n}.$ Then $(1-e_{n+1})b=0$ for all $b\in B_0.$
 \Wlog, we may assume that $a\in B_0$ for all $a\in {\cal F}.$
Put $z=1_{\tilde B}+e_nxe_n.$  {{Then}}
 $\|z^*z-1_{\tilde B}\|<\ep/32$ and $\|zz^*-1_{\tilde B}\|<\ep/32.$
 Thus we obtain a unitary $v'=\lambda'\cdot 1_{\tilde B}+z'$ such that $\lambda'\in \T,$  $z'\in B_0$
 and
 \beq\label{Unitary-10}
 \|u-v'\|<\ep/16.
 \eneq
{{Note that, since we assume that $u=1_{\td B}+x,$
$|\lambda'-1|<\ep/16.$   Put $v=v'{\bar \lambda'}.$  Then $\|u-v\|<\ep/8$
and $v=1_{\td B}+z''$ for some $z''\in B_0.$}}
 Put $B_1=\overline{(e_{n+2}-e_{n+1})B(e_{n+2}-e_{n+1}}).$
 Since $B$ is separable and simple, $K_1(B_1)=K_1(B).$
 Since $B$ (and  so does $B_1$) has stable rank one, one obtains a unitary
 $u_1\in {\tilde B}_1$ such that $[u_1]=[u^*]$ in $K_1(B).$
 Write $u_1=\lambda_1\cdot 1_{\tilde B_1}+z_1$ with $\lambda_1\in \T$ and $z_1\in B_1.$
 Put $u_1'=\lambda_1\cdot 1_{\tilde B}+z_1.$ Then $u_1'$ is a unitary in ${\tilde B}$ and
 $[u_1']=[u_1]=[u^*]$ in $K_1(B).$ Replacing $u_1'$ by $u_1'{\bar \lambda_1},$
 we may assume that $u_1'=1_{\td B}+z_1.$

 Since we have assume that {{${\cal F}\subset  B_0,$}} $z_1^*{{b}}=0$ and ${{b}}z_1=0$ for all ${{b}}\in B_0,$  
 {{one has}}
 \beq
 u_1'a=au_1'\rforal a\in {\cal F}.
 \eneq
 Put $v_1=(u_1')^*v.$ Then
 $v_1\in U_0({\tilde B})$ and
 $v_1=1_{\tilde B}+z_2$ for ${{z_2}}\in \overline{e_{n+2}Be_{n+2}}.$

 Let us write $v_1=\exp({{\sqrt{-1}}}h_1)\exp({{\sqrt{-1}}} h_2)\cdots \exp({{\sqrt{-1}}} h_m).$
 Since $z_2\in \overline{e_{n+2}Be_{n+2}},$ we may assume that ${{h_j}}=\af_j\cdot 1_{\tilde B}+b_j,$
 where $\af_j\in \R$ and ${{b_j}}\in \overline{e_{n+2}Be_{n+2}},$ $j=1,2,...,m.$
 Since we assume that $v_1=1_{\tilde B}+z_2,$ $\sum_{j=1}^m \af_j=2k\pi$ for some
integer $k.$
 Therefore we may also write
 \beq
 v_1=\exp({{\sqrt{-1}}} b_1)\exp({{\sqrt{-1}}} b_2)\cdots \exp({{\sqrt{-1}}} b_m).
 \eneq
 Let $b_j=(b_j)_+-(b_j)_-,$ where $(b_j)_+, (b_j)_-\in \overline{{(e_{n+2}Be_{n+2})}}_+,$ $j=1,2,....$
 Put
 $$
 {{B_2= \overline{(e_{n+4}-e_{n+3})B(e_{n+4}-e_{n+3})} \andeqn {\rm{choose}}\,\,b_0\in (B_2)_+\setminus \{0\}.}}
 $$
Since $B$ is simple, {{as 3.4 of \cite{eglnp},}}
there exist $x_1, x_2,...,x_N\in B$ such that
\beq
\sum_{j=1}^N x_j^*b_0x_j=f_{1/2}(e_{2n+3}).
\eneq
 Note that $f_{1/2}(e_{2n+3})b=b$ for all $b\in \overline{e_{n+2}Be_{n+2}}.$
 Then
 \beq
 \sum_{i=1}^N (b_j)_+^{1/2}x_i^*b_0x_i(b_j)_+^{1/2}=(b_j)_+.
 \eneq
 Put $(b_j)_+'=\sum_{i=1}^Nb_0^{1/2}x_i(b_j)_+x_i^*b_0^{1/2}.$ Then $(b_j)_+'\in B_2$ and
 \beq
 \tau((b_j)_+')=\tau((b_j)_+)\rforal \tau\in T(B).
 \eneq
 This implies that there are $b_j'\in B_2$ such that $\tau(b_j')=\tau(b_j),$ $j=1,2,...,m.$
 Define
 $v_2=\exp({{\sqrt{-1}}}\sum_{j=1}^m b_j')$ in ${\tilde B}.$
 Put $v_3=v_2^*v_1.$
 Then, by Lemma 3.1 of \cite{Thomsen} and 3.11 of \cite{GLX},  $v_3\in CU({\tilde B}).$ Since $b_j'a=ab_j'$ for all $a\in {\cal F}$ (as we assume that $a\in B_0$),
 \beq
 v_2av_2^*=a\rforal a\in {\cal F}.
 \eneq
Put  $w_1=(u_1')v_2.$ Then $w_1^*v=v_2^*(u_1')^*v=v_2^*v_1=v_3\in CU({\tilde B})$ and $w_1^*u\in U_0({\tilde B}).$
Put $w_2=uv^*.$ Then $w_2\in U_0({\tilde B})$ and, by \eqref{Unitary-10},
\beq
\|w_2-1_{\tilde B}\|<{{\ep/8.}}
\eneq
 Put $w=w_2w_1.$ Then,
 $\|w^*aw-a\|=\|w_2^*aw_2-a\|<\ep/4\rforal a\in {\cal F}.$
 Moreover,
\beq\nonumber
w^*u=w_1^*w_2^*u=w_1^*vu^*u=w_1^*v\in CU({\tilde B}).
\eneq
\end{proof}


\section{Existence Theorems}

\begin{lem}\label{ExtAB}
Let $A$
be {{in ${\cal M}_1$}} as constructed in Theorem 4.31 {{or (3) of Remark 4.32}} of \cite{GLrange} 
{{(see also Theorem 4.34 there)}} with continuous scale and
let $B$   be  a
separabple  {{\CA}}  in ${\cal D}^d$ which
{{has}} continuous scale.
Suppose that there is $\kappa\in KL(A,B)$
and
an affine continuous map
$\kappa_T: T(B)\to T(A)$ such that $\kappa$ and $\kappa_T$ are compatible {{(see Definition \ref{Dcompatible})}}.
Then,
there exists a  sequence of approximate multiplicative  \cpc s $\phi_n: A\to M_2(B)$ such that
\beq
&&[\{\phi_n\}]=\kappa.
\eneq
\end{lem}

\begin{proof}
Let $\ep>0$ and
{{${\cal F}\subset A^{\bf 1}$}} be a finite {{subset.}}
Fix a finite subset ${\cal P}\subset \underline{K}(A).$

Choose $\dt>0$ and  {{a}} finite subset ${\cal G}\subset A$ so that $[L]|_{{\cal P}}$ is well defined
for any
${\cal G}$-$\dt$-multiplicative \cpc\, $L$ from $A.$
We may assume that $\dt<\ep$ and ${\cal F}
\subset {\cal G}.$
Since both $A$ and $B$ have continuous scales,
$T(A)$ and $T(B)$ are compact (5.3 of \cite{eglnp}).

By the assumption, {{there is a sequence of subalgebras $A_n=E_n\oplus C_n\oplus D_n$, as described  in Theorem 4.34 of \cite{GLrange},  such that $A=\overline{\cup_{n=1}^{\infty} A_n}$ and $\lim_{n\to\infty}{\rm dist} (x, A_n)=0~~\mbox{for any}~~x\in A$.}} 
\Wlog,
 we may assume that, for some large $n\ge 1,$
that
\beq
{\cal G}\subset  E_n\oplus C_n',
\eneq
where $C_n'=C_n\oplus D_n\in {\cal C}_0$ and $E_n$ {{is}} as in Theorem {{4.34 of \cite{GLrange}.}} 
Moreover, we may further assume ${\cal G}={\cal G}_0'\cup {\cal G}_1',$
where ${\cal G}_0'\subset E_n$ and ${\cal G}_1'\subset C_n'$ are finite subsets.
We may also assume, \wilog, that there are finite
subset ${\cal P}_0\subset \underline{K}(E_n)$ and
${\cal P}_1\subset\underline{K}(C_n')$ such that
${\cal P}\subset {{[\imath]}}({\cal P}_0\cup {\cal P}_1),$
where ${{\imath:=}}\imath_n: E_n\oplus C_n'\to A$ is the embedding.
Since $K_0(C_n')$ is finitely generated, we may assume that
${\cal P}_1\cap K_0(C_n')$ generates $K_0(C_n').$
Let  
{{$c_n\in C_n'$}} be a strictly positive element of $C_n'$ with $\|c_n\|=1.$

Note here we assume, as constructed {{in Theorem 4.34 of \cite{GLrange},}} 
{{that}} $K_0(E_n)$ is torsion
and $\iota_{*0}(K_0(C_n'))$ is free.
Denote by 
{{$\Psi_0: A\to E_n$}} and $\Psi_1: A\to C_n'$ (for sufficiently large $n$) two \cpc s
which are
${\cal G}$-$\dt/16$-multiplicative  and
\beq\label{ExtAB-10}
\|\Psi_0(b)-b\|<\dt/16\rforal b\in {\cal G}_0'\andeqn \|\Psi_1(c)-c\|<\dt/16\rforal c\in {\cal G}_1'.
\eneq
\Wlog, we may assume
that
\beq\label{2020-713-71}
([\Psi_0]+[\Psi_1])|_{{\cal P}}=[{\rm id}]|_{{\cal P}},
[\Psi_1]|_{{\cal P}\cap [\iota]({\cal P}_1)}=[{\rm \id}]|_{{\cal P}\cap [\iota]({\cal P}_1)}\andeqn
[\Psi_1]|_{{\cal P}\cap [\iota]({\cal P}_0)}=0.
\eneq

Let
$$
{\cal P}_2=[\imath]({\cal P}_0)\cup [\iota]{{({\cal P}_1).}}
$$
%
We may assume that, for some $m\ge 1,$
$$
{\cal P}_2\subset K_0(A)\bigoplus K_1(A)\bigoplus_{j=1}^m (K_0(A, \Z/j\Z)\oplus K_1(A, \Z/j\Z).
$$
Moreover,  we may also assume that $m!x=0$ for all 
$x\in {\rm Tor}(K_0(A))\cap {{\cal P}_2}.$
Let 
${{G_{0, {\cal P}_2}}}$be the subgroup generated by $K_0(A)\cap {\cal P}_2.$
We may write 
${{G_{0, {\cal P}_2}}}:=F_0\oplus G_0,$ {{where}} $F_{0}$ is free and
is generated by 
${{\iota_{*0}{{({\cal P}_1)}}}},$
and $G_0$ is generated by ${\cal P}_0\cap K_0(A)$ and
$G_0$ is a finite group.
In particular, $m!x=0$ for all $x\in G_0.$
Moreover, $F_0\subset (\iota)_{*0}(K_0(C_n)).$


Choose $0<\dt_1<\dt$  and {{a}} finite subset ${\cal G}_3\subset A$ such that
$[L']|_{{\cal P}_2}$ is well defined for
any ${\cal G}_3$-$\dt_1$-multiplicative \cpc\, from $A.$
We assume that {{${\cal G}\subset {\cal G}_3.$}}

Note that, by {{Theorem 3.3 of \cite{GLrange},}} 
$A$ satisfies the assumption of 9.8 of \cite{GLII} {{(see Definition 9.3 of \cite{GLII})}}.

It follows from Theorem 3.4 of \cite{GLrange} 
that
there exists a ${{{\cal G}_3}}$-$\dt_1/4$-multiplicative \cpc\,
$L: A\to B\otimes M_K$
 for some
integer $K$ such that
\beq\label{711Ext1-2}
{{[L]|_{{\cal P}_2}}}=
{{\kappa|_{{\cal P}_2}.}}
\eneq
\Wlog, we may assume that ${\cal G}_3\subset A^{\bf 1}.$

Let {{${\cal Q}\subset \underline{K}(B)$}} be a finite subset
which contains $[L]({\cal P}_2).$
We assume that
\beq\label{11Ext1-3+1}
{\cal Q}\subset K_0(B){{\oplus}} K_1(B){{\oplus}}\bigoplus_{i=0,1}\bigoplus_{j=1}^{m_1}K_i(B, \Z/j\Z)
\eneq
for some $m_1\ge 2.$ Moreover, we may assume
that $(m_1)!x=0$ for all $x\in {\rm Tor}(G_{0,b}),$ where $G_{0,b}$ is the subgroup
generated by ${\cal Q}\cap K_0(B).$ \Wlog, we may assume that
$m|m_1.$ Choose an integer $m_2$ such that
$m_1|m_2.$

Let $p_1, p_2,...,p_l\in M_r(\td C_n')$ be projections which generate
$K_0(\td C_n')_+$ {{(see Theorem 3.15 of \cite{GLN}),}}
 Let 
${{{\bar p}_i}}\in M_r(\C\cdot 1_{\td C_n'})$ be scalar projections  with rank
$R_i\ge 1$
such that $o([p_i])=[p_i]-[{\bar p}_i]\in K_0(C_n'),$ $i=1,2,...,k.$
Since $(\kappa, \kappa_T)$ is compatible, $R_i+\rho_{A}(\kappa(o([p_i])))(s)>0$
for all $s\in T(A).$  Set $R:=\max\{R_i: 1\le i\le l\}$ and
\beq\label{2020-713-n1}
\eta_1:=\min{{\{\inf \{|R_i+\rho_{A}(\kappa(o([p_i])))(s)|:
s\in T(A)\}: 1\le i\le l\}}}.
\eneq

Let $b_0\in B$ with $\|b_0\|=1$ such that
\beq\label{711Ext1-3}
d_\tau(b_0)<\min\{\dt_1,
\eta_1/R\}/16(K+1)m_2\rforal \tau\in T(B).
\eneq
Let $e_b\in   B\otimes M_K$
be a strictly positive element of $B\otimes M_K$
such that
\beq\label{11Ext1-3+}
\tau(e_b)>7/8\rforal \tau\in T(B\otimes M_K).
\eneq


Let ${\cal G}_b\subset B\otimes  M_K$ be a  finite subset and
$1/2>\dt_2>0$ be such that
$[\Phi]|_{\cal Q}$ is well defined for any ${\cal G}_b$-$\dt_2$-multiplicative \cpc\,  $\Phi$ from
$B\otimes M_K.$ Note that  $B\in {\cal D}^d.$

There are  ${\cal G}_b$-$\dt_2$-multiplicative \cpc s
$\phi_{0,b}: B\otimes M_K\to B_0:={{\overline{\phi_{0,b}(e_b)( B\otimes M_K)\phi_{0,b}(e_b)}}}$
(see 2.10 of \cite{eglnp})
and $\psi_{0,b}: B\otimes M_K\to D_b\subset B\otimes M_K$
with
$D_b\in {\cal C}_0$ such that
\beq\label{11Ext1-4}
\|b-\diag(\phi_{0,b}(b),\overbrace{\psi_{0,b}(b), \psi_{0,b}(b),..., \psi_{0,b}(b)}^{(m_2)!})\|<\min\{\dt_2, \ep/16, \eta/16\}\rforal b\in {\cal G}_b\\\label{11Ext1-4+}
\andeqn \phi_{0,b}(e_b)\lesssim b_0\andeqn t(\psi_{0,b}(e_b))>3/4\rforal t\in T(D_b).
\eneq

Note that $K_1(D_b)=\{0\}.$ Moreover, we may also assume
that
\beq\label{711Ext1-4+}
(m_2)![\psi_{0,b}]|_{{\rm Tor}(G_{0,b})}=0\andeqn
(m_2)![\psi_{0,b}]|_{{\cal Q}\cap K_i(B, \Z/j\Z)}=0,\,\,\, j=2,3,...,m_1.
\eneq
Therefore
\beq\label{711Ext1-5n}
&&[\phi_{0,b}]|_{{\rm Tor}(G_{0,b})}=[{\rm id}_B]|_{{\rm Tor}(G_{0,b})},
{[}\phi_{0,b}{]}|_{{\cal Q}\cap K_1(B)}=[{\rm id}_{B}]|_{{\cal Q}\cap K_1(B)}\andeqn\\\label{711Ext1-5n+2}
&&{[}\phi_{0,b}{]}|_{{\cal Q}\cap K_i(B, \Z/j\Z)}=[{\rm id}_B]|_{{\cal Q}\cap K_i(B, \Z/j\Z)},\,\,\, j=2,3,...,m_1.
\eneq
Let $G_{\cal P}$ be the subgroup
generated by ${\cal P}$ and let
$\kappa'=\kappa- \phi_{0,b}\circ [L]$
be defined on $G_{\cal P}.$
Then, by \eqref{711Ext1-2}, \eqref{711Ext1-5n} and \eqref{711Ext1-5n+2},
we compute that
\beq\label{711Ext1-101}
\kappa'|_{G_0}=0,
\kappa'|_{{\cal P}\cap K_1(A)}=0\andeqn
\kappa'|_{{\cal P}\cap K_i(A, \Z/j\Z)}=0,\,\,\, j=2,3,...,m.
\eneq
Let $\zeta:=\kappa'\circ \iota: K_0(C_n')\to K_0({{B}}).$
Then, by \eqref{711Ext1-3} and \eqref{2020-713-n1},  for all $s\in T(A),$
\beq
R_i+\rho_A{{(\zeta(o([p_i])))}}(s)
&=&R_i+\rho_A{{(\kappa(o([p_i])))}}(s)-
\rho_A{{([\phi_{0,b}]\circ [L](o([p_i])))}}(s)\\
&=&R_i+\rho_A{{(\kappa(o([p_i])))}}(s)-
\rho_A{{([\phi_{0,b}](\kappa(o([p_i]))))}}(s)\\
&\ge& R_i+\rho_A{{(\kappa(o([p_i])))}}(s)- (\eta_1/R16)\rho_A{{((\kappa(o([p_i]))))}}(s)\\
&>& R_i+\rho_A{{(\kappa(o([p_i])))}}(s)- \eta_1>0
\eneq
for $1\le i\le l.$ Since $\{p_1,p_2,..., p_l\}$ generates $K_0(\td C_n'),$
this implies that the unital extension $\zeta^\sim : K_0(\td C_n')\to K_0({{\td B}})$
is strictly positive.
By Theorem {{5.7 of \cite{GLrange},}} 
there is a \hm\, $h: C_n'\to {{B}}$ such that
$h_{*0}=\zeta.$
By \eqref{11Ext1-4+} and \eqref{711Ext1-3}, since $B$ {{stably}} has almost stable rank one (see Lemma 11.1
of \cite{eglnp}),
one obtains a unitary $U\in U(\widetilde{B\otimes M_K})$ such that
\beq\label{2020-910-1}
U^*B_0U\subset B.
\eneq
Define  $\phi: A\to M_2(B)$ by
$\phi(a)=({\rm Ad}\, U\circ \phi_{0,b}\circ L(a))\oplus h(\Psi_1(a))$
for $a\in A.$
Then, by choosing  sufficiently
large ${\cal G}_b,$
 $\phi$ is ${\cal G}$-$\ep$-multiplicative.
 Note that
 \beq
 \kappa'|_{{\cal G}_{\cal P}}=(\kappa-[\phi_{0,b}]\circ [L])|_{{\cal G}_{\cal P}}=(\kappa-[\phi_{0,b}]\circ \kappa)|_{{\cal G}_{\cal P}}=m_2![\psi_{0,b}]\circ \kappa|_{{\cal G}_{\cal P}}.
 \eneq
 One then checks (see \eqref{2020-713-71}, \eqref{711Ext1-4+}, \eqref{711Ext1-5n}, \eqref{711Ext1-5n+2},
 \eqref{711Ext1-101}) that
\beq
[\phi]|_{\cal P}=\kappa|_{\cal P}.
\eneq
The lemma follows.
\end{proof}



\begin{lem}\label{ExtTBA}
Let   $A$ and $B$ be
separable
simple amenable \CA s  in ${\cal D}$  with continuous scales.
Suppose  that   $B$ satisfies the UCT,
and that  there is $\af\in KL(B,A)$  and
an affine continuous map
$\af_T: T(A)\to T(B)$
such that {{$\af$}} and $\af_T$ are compatible {{(see Definition \ref{Dcompatible})}}.
Suppose also that there exists a sequence of \cpc s
$\Psi_n: B\to {{M_3(A)}}$ such that
\beq\label{2020-14-n1}
[\{\Psi_n\}]=\af\andeqn \lim_{n\to\infty}\|\Psi_n(ab)-\Psi_n(a){{\Psi_n(b)}}\|=0
\eneq
for all $a, b\in B.$
Then,  there exists
a  sequence of approximate multiplicative  \cpc s $\Phi_n: B\to A$  such that
\beq\label{ExtTBA-1}
&&
[\{\Phi_n\}]=\af,\\\label{ExtTBA-2}
&&\lim_{n\to\infty}\sup \{|\tau\circ \Phi_n(b)-\af_T(\tau)(b)|:\tau\in T(A)\}=0\tforal b\in B_{s.a.}.
\eneq
\end{lem}

\begin{proof}
Let $e_A$ and $e_B$ be  strictly positive elements for
$A$ and $B,$ respectively.

Fix $\ep>0,$ a finite subset ${\cal F}\subset B,$ $\eta>0,$ a finite subset
${\cal H}\subset B_{s.a.}$ and a finite subset ${\cal P}\subset \underline{K}(B).$
We assume that any ${\cal F}$-$\ep$-multiplicative \cpc\, $L$ defines $[L]|_{\cal P}$ {{well}}
and if $L_1$ and $L_2$ are both ${\cal F}$-$\ep$-multiplicative \cpc s and
$$
\|L_1(x)-L_2(x)\|<\ep\rforal x\in {\cal F},
$$
then $[L_1]|_{\cal P}=[L_2]|_{\cal P}.$

Put ${\cal F}_1={\cal F}\cup {\cal H}.$  We may assume that ${\cal F}_1\subset B^{\bf 1}.$
Let $\ep_0=\min\{\ep, \eta\}.$


Since $B\in {\cal D},$ there is a strictly  positive {{element}} $e_B\in B$ with $\|e_B\|=1$ and  there are  two sequences of mutually orthogonal \SCA s $B_{n,0}$ and $B_{n,1}$ of $B$
with $B_{n,1}\in {\cal C}_0,$ and two sequences of \cpc s $\phi_{n,i}: B\to B_{n,i}$ ($i=0,1$) such that
{{\beq\label{ExtTBA-10}
&&\lim_{n\to\infty}\|\phi_{n,i}(ab)-\phi_{n,i}(a)\phi_{n,i}(b)\|=0\rforal a, b\in B,\\\label{ExtTBA-10402}
&&\lim_{n\to\infty}\|a-\diag(\phi_{n,0}(a), \phi_{n,1}(a))\|=0\rforal a\in B,\\
&&t(f_{1/4}(\phi_{n,1}(e_B)))\ge 1/2\rforal t\in T(B_{n,1})\andeqn\\
&&\lim_{n\to\infty}\sup\{d_\tau(e_{n,0}): \tau\in T(B)\}=0,
\eneq}}
where $e_{n,0}$ is a strictly positive element of $B_{n,0}^{\bf 1}$ {{with
$\|e_{n,0}\|=1.$}}
We may assume, for all $n\ge 1,$ that
\beq\label{ExtTBA-10+1}
&&\|e_{n,0}\phi_{n,0}(x)e_{n,0}-\phi_{n,0}(x)\|<\ep/64\rforal x\in {\cal F}_1\cup \{ab: a, b\in {\cal F}_1\},\\\label{ExtTAB-10+2}
&&{[}\phi_{n,0}{]}|_{\cal P}+[\phi_{n,1}]|_{\cal P}=[{\rm id}_B]|_{\cal P},\\\label{ExtTAB-10+3}
&&\|x-\diag(\phi_{n,0}(x), \phi_{n,1}(x))\|<\ep_0/64\rforal x\in {\cal F}_1\andeqn\\\label{ExtTAB-10+4}
&&d_\tau(e_{n,0})+d_\tau(e_{n,1})\le 1\rforal \tau\in T(B),
\eneq
where $e_{n,1}$ is a strictly positive element in $B_{n,1}.$

Choose $n_0\ge 1$ such that $1/n_0<\ep_0/4.$  For some large $n_1,$
\beq\label{ExtTBA-11}
d_\tau(e_{n,0})<1/2n_0\rforal \tau\in T(B)
\eneq	
for all $n\ge n_1.$
For each $n\ge n_1,$ there are mutually {{orthogonal}} elements $a_{n,1}, a_{n,2},...,a_{n,2n_0}\in B$
and unitaries $u_i\in {\tilde B}$ such that
$u_i^*a_{n,1}u_i=a_{n,i},$ $i=1,2,...,2n_0,$ and $a_{n,1}=e_{n,0}.$
Since $B$ is stably {{projectionless}}, $sp(a_{n,1})=[0,1].$ 
 Therefore
$a_{n,1}$ and $a_{n,1}^{1/2}u_i,$ $i=1,2,...,2n_0,$ generate a \SCA\,
$C\cong C_0((0,1])\otimes M_{2n_0}{{\subset B}}$ which is semi-projective.
Therefore, there {{exist}} $n_2\ge n_1,$  $\dt>0$ and a finite subset ${\cal G}\subset  B$
satisfying the following:
for any $n\ge n_2,$ and any ${\cal G}$-$\dt$-multiplicative \cpc\, ${{L''}}: B\to D$ (for any \CA\, $D$),
there exists a \hm\, $h_n: C\to D$ such that
\beq\label{ExtTBA-12}
\|h_n(a)-{{L''}}(a)\|<\ep/16 \rforal a\in {\cal F}_1.
\eneq
Let $\{\Psi_n\}$ be as in the lemma.
Consider $L'=\Psi_m\circ \phi_{n,0}.$  By choosing $n\ge n_2$ and {{sufficiently}} large $m,$
we may assume that $L'$ is ${{{\cal F}_1}}$-$\ep/2$-multiplicative and
there exists a \hm\, $h: C\to {{M_3(A)}}$
such that
\beq\label{ExtTBA-14}
\|h(a_{n,i})-L'(a_{n,i})\|<\ep/16,\,\,\,i=1,2,...,2n_0.
\eneq
Therefore
\beq\label{ExtTBA-15}
f_{\ep/8}(L'(a_{n,1}))\lesssim h(a_{n,1})\,\,\,{\rm or}\,\,\, f_{\ep/8}(L'(e_{n,0}))\lesssim h(a_{n,1}).
\eneq
By choosing even larger $m,$ we may further assume
that  $L$ is ${\cal F}_1$-$\ep/2$-multiplicative, where
$L(x)=f_{\ep/8}(L'(e_{n,0}))L'(x)f_{\ep/8}(L'(e_{n,0}))$ for all 
{{$x\in B.$}}
By \eqref{ExtTBA-15},
\beq\label{ExtTBA-15+}
L(e_B)\lesssim f_{\ep/8}(L'(e_{n,0}))\lesssim h(a_{n,1}).
\eneq
It follows that
\beq\label{ExtTBA-15+2}
d_\tau(L(e_B))<1/2n_0\le \eta/2\rforal \tau\in T(A).
\eneq
Note
that (see \, \eqref{ExtTBA-10402} {{and \eqref{2020-14-n1}}})
\beq\label{ExtTBA-16}
[L]|_{\cal P}+\af\circ [\phi_{n,1}]|_{\cal P}={{\af|_{\cal P}.}}
\eneq
Let $\imath: B_{n,1}\to B$ be the embedding. Let $\af_T^{\sharp}: {\rm LAff}^\sim_{+}({\tilde T}(B))\to {\rm LAff}_+^\sim({\tilde T}(A))$
be the ordered semi-group \hm\, induced by $\af_T.$   In fact
(see
11.1 of \cite{GLII})
$\af_T^{\sharp}$
is a morphism in ${\bf Cu},$ since both $A$ and $B$ are stably projectionless (see 9.3 of \cite{eglnp}).
{{Recall that, by Theorem 6.2.3 of \cite{Rl} and Theorem 7.3 of \cite{eglnp}, the map defined in (6.6)
of \cite{Rl} (lines above Theorem 6.2.3 of \cite{Rl})  is an isomorphism in ${\bf Cu}$
for \CA s $A$, $B$ and {{$B_{n,1}.$}}}}
Let {{$\gamma: \LAff^\sim_+(\td T(B_{n,1})\to \LAff^\sim_+(\td A)$}} {{be defined}} by
$\gamma=\af_T^{\sharp}\circ {\rm Cu}^\sim(\imath){{|_{\LAff^\sim_+(\td T(B_{n,1}))}}}.$
Let $\gamma_0=\af\circ \imath_{*0}: K_0(B_{n,1})\to K_0(A).$  Since $(\af, \af_T)$ is compatible,
it is easy to check that $\gamma$ and $\gamma_0$ {{induce}}  a morphism
$\gamma^{\sim}: {\rm Cu}^{\sim}(B_{n,1})\to {\rm Cu}^{\sim}(A).$
It follows from 1.0.1 of \cite{Rl} that there is a \hm\, $H: B_{n,1}\to A$ such that
${\rm Cu}^{\sim}(H)=\gamma^{\sim}$ and
\beq\label{ExTBA-17}
d_\tau(H(e_{n,1}))=d_{\af_T(\tau)}(e_{n,1})\rforal \tau\in T(A).
\eneq
Define $\Phi: B\to {{M_4(A)}}$ by {{$\Phi(a)=L(a)\oplus H(\phi_{n,1}(a))$}} 
for all $a\in A.$ Then
$\Phi$ is ${\cal F}$-$\ep$-multiplicative.
{{Note that}}  $K_0(B_{n,1})$ is free and $K_1(B_{n,1})=\{0\}.$
By \eqref{ExtTAB-10+2}, \eqref{2020-14-n1}
and \eqref{ExtTBA-16},
one computes that
\beq
[\Phi]|_{\cal P}=\af|_{\cal P}.
\eneq
By \eqref{ExtTBA-15+2},  \eqref{ExtTAB-10+3},   \eqref{ExTBA-17}  and by \eqref{ExtTBA-15+2}, we have
\beq\label{2020-713-nn+10}
&&
|\tau(\Phi(x))-{{(\af_T(\tau))(x)}}|<\eta\rforal {{x\in}} {\cal H}\andeqn\\
&&|d_\tau(\Phi({{e_B}}))-1|<\eta\rforal \tau\in T(A).
\eneq

\noindent
This shows that we have a sequence of 
{{$\Phi_n: B\to {{M_4(A)}}$}}
which satisfies \eqref{ExtTBA-1} {{and}}
\eqref{ExtTBA-2}.

It remains to show that we can modify $\Phi_n$ so {{that}} it maps into $A$ instead of 
{{$M_3(A).$}}
Consider  positive elements $\Phi_n(e_B).$
Define 
$H_n: C_0:=C_0((0,1])\to {{M_3(A)}}$
by $H_n(f)=f(\Phi_n(e_B))$ for all $f\in C_0((0,1]).$
{{Let $\imath: C_0((0,1])\to B$ be defined by}}
$\imath(f)=f(e_B)$ for all $f\in C_0((0,1]).$ Then
 $\bt={{\af_T^{\sharp}}}\circ {{\rm Cu^\sim}}(\imath)$ gives a morphism from ${{\rm Cu}}^{\sim}(C_0((0,1]))$ to
 ${{\rm Cu}}^{\sim}(A).$  It follows from 1.0.1 of \cite{Rl} that there exists
a \hm\,  $H_0: C_0((0,1])\to A$  such that ${{\rm Cu^\sim}}(H_0)=\bt.$
Let
$$
\Delta(\hat{f})=(1/4)\inf\{\tau(H_0(f)): \tau\in T(A)\}\rforal f\in C_0((0,1])_+.
$$
Then,  as $\{\Phi_n\}$ satisfies \eqref{ExtTBA-2},
we have, for any fixed finite subset ${\cal H}_0\subset {C_0}_+^{\bf 1}\setminus \{0\}$ and
large $n,$
\beq
\tau(H_n)(f)\ge \Delta(f)\rforal f\in {\cal H}_0.
\eneq
Thus, by applying 
{{7.7 of \cite{eglnp},}}
there is a sequence of unitaries $u_n\in {{{\widetilde{M_{{4}}(A)}}}}$
such that
\beq
\lim_{n\to\infty}\|{\rm Ad}\, u_n\circ f(H_n(e_B))-f(H_0(e_B))\|=0\rforal f\in C_0((0,1]).
\eneq
Note that $f(H_0(e_B))\in A$ for all $f\in C_0((0,1]).$
Choose an approximate identity $\{a_n\}$ for $A.$
Then, there is a subsequence
$\{k(n)\}$ such that
\beq
\lim_{n\to\infty} \|a_{k(n)}u_n^*\Phi_n(b)u_na_{k(n)}-u_n^*\Phi_n(b)u_n\|=0\rforal b\in B.
\eneq
Define $\Phi_n': B\to A$ by $\Phi_n'(b)=a_{k(n)}u_n^*\Phi_n(b)u_na_{k(n)}$ for all $b\in B.$
We then replace $\Phi_n$ by $\Phi_n'.$ This lemma then follows.
\end{proof}

\begin{lem}\label{Lextcu}
Let $A$
{{in  ${\cal M}_1$ be as in}}
{{the form of}}  (3) of Remark 4.32
of \cite{GLrange} 
with continuous scale and
let $B$   be  a separable  simple \CA s in ${\cal D}^d$
with
continuous scale.
Suppose that there is $\kappa\in KL(A,B),$
an affine continuous map
$\kappa_T: T(B)\to T(A)$ such that $\kappa,$ $\kappa_T$ and  a continuous \hm\, $\kappa_{cu}: U({\tilde A})/CU({\tilde A})\to
U({\tilde B})/CU({\tilde B})$
are compatible {{(see Definition \ref{Dcompatible})}}.
Then,
there exists a  sequence of approximate multiplicative  \cpc s $\phi_n: A\to B$ such that
\beq
&&[\{\phi_n\}]=\kappa\tand\\
&&\lim_{n\to\infty}\sup \{|\tau\circ \phi_n(a)-\kappa_T(\tau)(a)|: \tau\in T(B)\}=0\tforal a\in A_{s.a.}\andeqn\\
&&\lim_{n\to\infty}{\rm dist}(\phi_n^{\dag}({\bar z}), \kappa_{cu}({\bar z}))=0{{\tforal}} {\bar z}\in U({\tilde A})/CU({\tilde A}).
\eneq
\end{lem}

\begin{proof}
Let $\ep>0,$ let $\eta>0$ and let $\sigma>0,$ {{let ${\cal P}\subset \underline{K}(A)$,
$S_u\subset U({\tilde A})/CU({\tilde A})$,}}
${\cal H}\subset {{A_{s.a.}}}$
and
 ${\cal F}\subset A$ be finite subsets.
\Wlog, we may assume that ${\cal F}\subset A^{\bf 1}$ and
$[L']|_{\cal P}$ and $(L')^{\dag}|_{S_u}$ are {{well defined}} for any
${\cal F}$-$\ep$-multiplicative \cpc\, from $A.$
Let $G_1\subset K_1(A)$ be the subgroup generated by ${\cal P}\cap K_1(A).$

 We now apply Theorem {{4.34 of \cite{GLrange}}} 
 and retain  notation used there such as
 $E_n$ and
$\{\imath_n\}.$

Fix $\dt>0$ and a finite subset ${\cal G}\subset A.$
We assume that $\dt<\min\{\ep/2, \eta/4, \sigma/16\}.$
To simplify notation, \wilog, we may assume
that $G_1\subset F\subset (\imath_{n_0})_{*1}(K_1(E_{n_0}))$ for some $n_0\ge 1,$
where $F$ is a finitely generated standard subgroup (see \ref{DdertF}).
We also choose $n_0$ larger than that required by {{Lemma}} \ref{L215} for $\dt$ (in place of $\ep$) ${\cal G}$
(in place of ${\cal F}$)
${\cal P}$ and $\sigma/16$ (in place of $\dt_0$).

\Wlog, we may write
\beq
S_u=S_{u,1}\sqcup S_{u,0},
\eneq
where $S_{u,1}\subset J_{F,u}(F)$ and $S_{u,0}\subset U_0({\tilde A})/CU({\tilde A})=
\Aff(T({\tilde A}))/\overline{\rho_{\tilde A}(K_0({\tilde A}))}$
and both $S_{u,1}$ and $S_{u,0}$ are finite subsets.
For $w\in S_{u,0},$
write
\vspace{-0.1in}\beq
w={{\prod_{j=0}^{l(w)}}}\exp(2i\pi h_{w,j}),
\eneq
where
$h_{w,0}\in \R,$  and
$h_{w,j}\in {{A_{s.a.}}},$ $j=1,2,...,l(w).$
Let
\vspace{-0.1in}\beq
{\cal H}_u=\{h_{w,j}: 1\le j\le l(w),\,\,\, w\in S_{u,0}\}\andeqn
M=\max\{\sum_{i=0}^{l(w)}\|h_{w,j}\|: w\in S_{u,0}\}.
\eneq
To simplify notation further, we may assume that $G_1=F.$

Write $G_1=\Z^{m_f}\oplus {\rm Tor}(G_1)$ and
$\Z^{m_f}$ is generated by cyclic and free generators $x_1,x_2,...,x_{m_f}.$
Let ${\rm Tor}(G)$ be generated by $x_{0,1}, x_{0,2},...,x_{0, m_t}.$
Let $u_1, u_2,..., u_{m_f}, u_{1,0},u_{2,0},...,u_{m_t,0}\in U({\tilde A})$
be unitaries such
that $[u_i]=x_i,$ $1\le i\le m_f,$
and $[u_{j,0}]=x_{0,j},$ $1\le j\le m_t.$
Recall that  ${{\Pi_{cu}^{\td A}}}: U({\tilde A})/CU({\tilde A})\to K_1(A)$ is the quotient map.
Let $G_u$ be the subgroup generated by $S_{u,1}.$
Since $(\kappa, \kappa_T, \kappa_u)$ is compatible, \wilog,
we may assume that ${{\Pi_{cu}^{\td A}}}(G_u)=\{x_1,x_2,...,x_{m_f}\}\cup \{x_{0,1}, x_{0,2},...,x_{0, m_t}\}$
and $S_{u,1}=\{{\bar u}_1, {\bar u}_2,...,{\bar u}_{m_f}, {\bar u}_{1,0},{\bar u}_{2,0},...,{\bar u}_{m_t,0}\}$
as described in \ref{DdertF}, in particular, $k_j{\bar u}_{j,0}=0$ in $U({\tilde A})/CU({\tilde A}),$
$j=1,2,...,m_t.$  Let $G_u':=J_{cu}^A(\iota_{*0}(K_1(E_{n_0}))),$ where
$J_{cu}^A: K_1(A)\to U({\td A})/CU(\td A)$ be a splitting map (for ${{\Pi_{cu}^{\td A}}}$). Note $G_u\subset G_u'.$

Let $\phi_n: A\to B$ be a sequence of approximately multiplicative \cpc s
given by \ref{ExtAB} and  \ref{ExtTBA}
such that
\beq\label{11Ext2-3}
[\{\phi_n\}]&=&\kappa\andeqn\\\label{11Ext2-3+}
\lim_{n\to\infty}\sup \{|\tau\circ \phi_n(a)-\kappa_T(\tau)(a)|\}&=&0\rforal a\in A_{s.a.}.
\eneq

By replacing 
{{$\phi_n(a)$ by $e_n\phi_n(a)e_n$}}
(for all 
{{$a\in A$}}) for some
{{$e_n\in B_+$}} with $\|e_n\|=1,$ we may assume that 
$\phi_n({{A}})^{\perp}\not=\{0\}.$
Choose 
$b_n\in (\phi_n({{A}})^\perp)_+$with $\|b_n\|=1$
and
\beq\label{711EXT2-11+}
d_{\tau}(b_n)<\min\{\eta, \sigma/16\}/2(M+1)\rforal 
{{\tau\in T(B),}} \,\, n=1,2,....
\eneq

Fix a sufficiently large $n.$
Define $\lambda=\kappa_{cu}|_{G_u'}-{\phi_n}^{\dag}|_{G_u'}: G_u\to U({\widetilde{B}})/CU({\widetilde{B}}).$
Write $G_u'=F(G_{u'})\oplus {\rm Tor}(G_u'),$   where $F(G_{u'})$ is the free part of $G_u'.$
Define $\lambda'|_{F(G_u')}=\lambda|_{F(G_u')}$ and $\lambda'|_{Tor(G_u')}=0.$
Since $(\kappa, \kappa_T, \kappa_u)$ is compatible,
${{\Pi_{cu}^{\td B}}}\circ \lambda({\bar u_i})=0$ and ${{\Pi_{cu}^{\td B}}}\circ \lambda({\bar u}_{0,j})=0,$
$i=1,2,...,{{m_f}}$ and $j=1,2,...,m_t.$

Let ${\cal F}_1={\cal F}\cup {\cal H}.$
It follows from {{Lemma}} \ref{L215} that there exists ${\cal F}_1$-$\min\{\ep/4, \eta/4\}$-multiplicative \cpc\,
$L: A\to \overline{b_nBb_n},$
such that
\beq\label{11Ext2-15}
[L]|_{\cal P}=0\andeqn {\rm dist}(L^{\dag}({\bar u}_j), \lambda'({\bar u}_j))<\sigma,\,\,\,j=1,2,...,m_f.
\eneq

Define $\Psi: A\to B$ by
\beq
\Psi(a)=\phi_n(a)\oplus L(a)\rforal a\in A.
\eneq
Then $\Psi$ is ${\cal F}$-$\ep$-multiplicative if $n$ is sufficiently large.

It follows from \eqref{11Ext2-3},
\eqref{11Ext2-15} and the definition of $\lambda'$ that
\beq
[\Psi]|_{\cal P}=\kappa|_{\cal P}\andeqn
{\rm dist}(\Psi^{\dag}({\bar u}_j), {{\kappa_{cu}}}({\bar u}_j))&<&\sigma,\,\,\, j=1,2,...,m_f.
\eneq
By \ref{Lderttorsion}, we may also have
\beq
{\rm dist}(\Psi^{\dag}({\bar u}_{j,0}), {{\kappa_{cu}}}({\bar u}_{j,0}))&<&\sigma,\,\,\, j=1,2,...,m_t.
\eneq
By the choice of $M$ and ${\cal H}_u,$  \eqref{711EXT2-11+}
and by the assumption
that $(\kappa, \kappa_T, {{\kappa_{cu}}})$ is compatible,
\beq
{\rm dist}(\Psi^{\dag}({\bar w}), {{\kappa_{cu}}}({\bar w}))<\sigma\rforal w\in S_{u,0}.
\eneq

Moreover, by \eqref{711EXT2-11+},  by \eqref{11Ext2-3+} and by choosing sufficiently large $n,$
\beq
\sup\{|\tau(\Psi(b))-\kappa_T(\tau)(b)|: \tau\in T(A)\}<\eta\rforal b\in {\cal H}.
\eneq
\end{proof}

\begin{thm}\label{Text1}
Let $A$ in ${\cal M}_1$
be
the form of  (3) of Remark 4.32
of \cite{GLrange}
with continuous scale,
let $B_1$   be  a non-unital separable  simple \CA s in ${\cal D}$ which
has continuous scale  and let $B=B_1\otimes U,$ where $U$ is an infinite dimensional UHF-algebra.
Suppose that there is $\kappa\in KL(A, B),$
an affine continuous map
$\kappa_T: T(B)\to T(A)$
and a continuous \hm\, ${{\kappa_{cu}}}: U({\tilde A})/CU({\tilde A})\to U({\widetilde{B}})/CU({\widetilde{B}})$ such that $(\kappa, \kappa_T, {{\kappa_{cu}}})$ is compatible {{(see Definition \ref{Dcompatible})}}.
Then there exists a  \hm\, $\phi: A\to B$
 such that
\beq\label{11ExtT1-1}
&&[\phi]=\kappa,\,\,\,
\tau\circ \phi(a)=\kappa_T(\tau)(a)\tforal a\in A_{s.a.}\tand
\phi^{\dag}={{\kappa_{cu}}}.
\eneq

\end{thm}

\begin{proof}
Let $e_a\in A$ be a strictly positive element of $A$ with $\|e_a\|=1.$
Since $A$ has continuous scale, \wilog, we may assume that
\beq\label{11ExtT1-5}
\min\{\inf\{\tau(e_a): \tau\in T(A)\},\,\inf\{\tau(f_{1/2}(e_a)):\tau\in T(A)\}\}>3/4.
\eneq
Let $T: A_+\setminus \{0\}\to \N\times \R_+\setminus \{0\}$ be given by
Theorem 5.7 of \cite{eglnp}.

By {{Lemma}} \ref{Lextcu}, there exists a sequence of approximately multiplicative  \cpc s $\phi_n: A\to B$ such that
\beq
&&[\{\phi_n\}]=\kappa\\
&&\lim_{n\to\infty}\sup \{|\tau\circ \phi_n(a)-\kappa_T(\tau)(a): \tau\in T(A)|\}=0\rforal a\in A_{s.a.}\andeqn\\
&&\lim_{n\to\infty}{\rm dist}({{\kappa_{cu}}}(z), \phi_n^{\dag}(z))=0\rforal z\in U({\tilde A})/CU({\tilde A}).
\eneq

Let $\ep>0$ and ${\cal F}\subset A$ be a finite subset.

We will apply   Theorem 5.3 of \cite{GLII} (see also 5.2 of \cite{GLII}).
Note that $B\in {\cal D}^d,$ $K_0({\tilde B})$
is weakly unperforated (see Proposition 5.5 of \cite{GLII}, also A.7 of \cite{eglnkk0} and
Theorem 16.10 of \cite{GLII}).
Let
 $\dt_{1,1}>0$ (in place of $\dt$), $\gamma_1>0$ (in place of $\gamma$),
$\eta_1>0$ (in place of $\eta$), let ${\cal G}_{1,1}\subset A$ (in place of ${\cal G}$) be a finite subset,
${\cal H}_{1,1}\subset A_+\setminus \{0\}$ (in place of ${\cal H}_1$) be a finite subset, ${\cal P}_1\subset \underline{K}(A)$ (in place of ${\cal P}$),
{{${\cal U}_1\subset U({\tilde A})$}} (in place of ${\cal U}$) with
$\overline{{\cal U}}={\cal P}\cap K_1(A)$  and let ${\cal H}_{1,2}\subset A_{s.a.}$ (in place of
${\cal H}_2$)  required by Theorem  5.3 of \cite{GLII}
for $T,$ $\ep$ and ${\cal F}$ (with $T(k,n)=n,$ see 5.2 of \cite{GLII}).

\Wlog, we may assume that ${\cal H}_{1,1}\subset A_+^{\bf 1}\setminus \{0\}$
and $\gamma_1<1/64.$

Let ${\cal G}_{1,2}\subset A$  (in place of ${\cal G}$)  be a finite subset and let
$\dt_{1,2}>0$ be required by Theorem 5.7 of \cite{eglnp}
for the above ${\cal H}_{1,1}$
(in place of ${\cal H}_1$).
Let $\dt_1=\min\{\dt_{1,1}, \dt_{1,2}\}$  and ${\cal G}_1={\cal G}_{1,1}\cup {\cal G}_{1,2}.$

Choose $n_0\ge 1$ such that
$\phi_n$ is ${\cal G}_1$-$\dt_1/2$-multiplicative,  for all $n\ge n_0,$
\beq
&&[\phi_n]|_{{\cal P}_1}=\kappa|_{{\cal P}_1},\\
&&\sup\{|\tau\circ \phi_n(a)-\kappa_T(\tau)(a)|:\tau\in T(A)\}<\gamma_1/2\rforal  a\in {\cal H}_{1,2},\\
&&\tau(f_{1/2}(\phi_n(e_a)))>3/8\rforal \tau\in T(A)\andeqn\\
&&{\rm dist}(\phi_n^{\dag}({\bar u}), {{\kappa_{cu}}}({\bar u}))<\eta/2\rforal u\in {\cal U}.
\eneq
{{By}} Theorem 5.7  of \cite{eglnp}, $\phi_n$ are all exactly $T$-${\cal H}_{1,1}$-full.
By applying Theorem 5.3 of \cite{GLII},
we obtain a unitary $u_n\in {\tilde B}$ (for each $n\ge n_0$)
such that
\beq\label{11ExtT1-15}
\|u_n^*\phi_n(a)u_n-\phi_{n_0}(a)\|<\ep\rforal a\in {\cal F}.
\eneq

Now let $\{\ep_n\}$ be an decreasing sequence of positive elements
such that $\sum_{n=1}^{\infty}\ep_n<\infty$ and let
$\{{\cal F}_n\}$ be an increasing sequence of finite subsets of $A$
such that $\cup_{n=1}^{\infty}{{{\cal F}_n}}$ is dense in $A.$

By what have been proved, we obtain a subsequence $\{n_k\}$ and
a sequence of unitaries $\{u_k\}\subset {\tilde A}$ such that
\beq\label{11ExtT1-16}
\|{\rm Ad}\, u_{k+1}\circ \phi_{n_{k+1}}(a)-{\rm Ad}\, u_{k}\circ \phi_{n_k}(a)\|<\ep_k\rforal a\in {\cal F}_k,
\eneq
$k=1,2,....$
Since $\cup_{n=1}^{\infty}{{{\cal F}_n}}$ is dense in $A,$  by \eqref{11ExtT1-16},
$\{{\rm Ad}\, u_k\circ \phi_{n_k}(a)\}$ is a Cauchy sequence.
Let
\beq
\phi(a)=\lim_{k\to\infty}{\rm Ad}\, u_k\circ \phi_{n_k}(a)\rforal a\in A.
\eneq
Then $\phi: A\to B$ is a \hm\, which satisfies \eqref{11ExtT1-1}.

\end{proof}

\section{Existence Theorem, continued}


In  \ref{C201GLN}, \ref{traces}, \ref{Kstate}, we will use a similar construction as in
section 20 of \cite{GLN}. These statements are taken from there and repeated with minimal
modification.

\begin{NN}\label{C201GLN}
Let $A\in D^d$ be a separable simple \CA\, with continuous scale
and let $e_A\in A$ be a strictly positive element with $\|e_A\|=1.$
Then
there are mutually orthogonal \SCA s $A_{n,0}$ and  $S_n\in {\cal C}_0$ of $A,$ two
\cpc s ${{\Psi_n}}: A\to A_{n,0}$ and $L_n: A\to S_n$ such that $A_{{n,0}}\perp S_n,$
{{\beq\label{C210GLN-n1}
&&\lim_{n\to\infty}\|a-\diag({{\Psi_n}}(a), L_n(a))\|=0\rforal a\in A,\\\label{C210GLN-n2}
&&\sup\{d_\tau(e_{n,0}):\tau\in T(A)\}<1/2^{n+4}\andeqn\\\label{1207-21-8-1}
&&t(f_{1/4}(L_n(e_A)))\ge 1/2\rforal  t\in T(S_n),
\eneq}}
where $A_{n,0}$ is a hereditary \SCA\, of $A$ generated by positive element $e_{n,0}$ {{for}} which
we also assume $\|e_{n,0}\|=1.$
{{Note that \eqref{1207-21-8-1} implies that $\lambda_s(S_n)\ge 1/2.$}}

{{By}} Theorem 3.15  of \cite{GLN},
the positive cone of the
$K_0({\tilde S_n})$ is finitely generated.  {{Denote by}} $\{s^n_1, s^n_2, ..., s^n_{r_n}\}\subseteq K_0({\tilde{S_n}})_+$ a set of generators of $K_0({\tilde{S_n}})_+.$  Let ${\cal P}'_n$ be a finite set of projections in the matrix algebras over ${\tilde{S_n}}$ which represent the elements $s^n_j$ for $j=1,2,\cdots, r_n$.

Let $\{x_1, x_2,....\}$ be a sequence of elements in the unit ball of $A$ such that it is
dense in $A^{\bf 1}.$
Fix a finite subset ${\cal P}\subset M_r({\tilde A})$ of projections
(for some integer $r\ge 1$), a finite subset ${\cal F}_0\subset A$ and $1/2>\dt_0>0.$
We assume that ${{{\cal F}_0}}$ 
is sufficiently large and $\dt_0$ is sufficiently small so that
any ${{{\cal F}_0}}$-${{\dt_0}}$ -multiplicative \cpc\, $L'$ from
$A$ defines $[L']|_{\cal P}$ {{well}} and gives a \hm\, on {{$G_0,$}} the subgroup
of $K_0({\tilde A})$ generated by ${\cal P}.$ {{Furthermore, we assume that if $L', L'': A\to B$ satisfy $\|L'(f)-L''(f)\|<2 \dt_0$ for any $f\in {\cal F}_0$, then $[L']|_{\cal P}=[L'']|_{\cal P}$. }}
Here we also {{write}}  $L'$ for the unitization $(L')^{\sim}: {\tilde A}\to {{\tilde B}}$ for $L': A\to B$
of $L'.$

Let ${\cal F}_n={\cal F}_0\cup\{x_1,x_2,...,x_n\}.$

Without loss of generality, passing  to  a subsequence, we will  assume that
\beq\label{Sept-10-2020-1}
\|a-\Psi_n(a)\oplus L_n(a)\|<\dt_0/2^{n+3}~~~\mbox{for all}~~a \in {\cal F}_n.
\eneq
We assume that $e_{n+1,0}\lesssim e_{n, 0},$ $n=1,2,....$
We also assume that $L_1$ is ${\cal F}_1$-$\dt_0/2$-multiplicative.
By passing to a subsequence, one may assume, since each $S_n$ is semi-projective,
that there is a unital \hm\, $h_n: {\tilde S_n}\to {\tilde S_{n+1}}$
(where we assume that $1_{\tilde S_n}=1_{\tilde A}$ for all $n$) {{with $h_n(S_n)\subset S_{n+1}$}}
such that
{{\beq
&&
{{\|h_n(g)-L_{n+1}(g)\|}}<\dt_0/2^{n+3}\rforal g\in L_n({\cal F}_n'){{\subset S_n,}}\\\label{62320-n11}
&&\sup\{|\tau((L_{n+1}(p))-L_n(p))|:\tau\in T(A)\}<\dt_0/2^{n+3}\andeqn\\\label{62320-nn11}
&&\|h_n(q_n(p))-L_{n+1}(q_n(p))\|<\dt_0/2^{n+3}\rforal p\in {\cal P},
\eneq}}
where $q_n(p)$ is a projection close to $L_n(p)$ and $[q_n(p)]=[L_n]([p]),$
where ${\cal F}_1'={\cal F}_1$ and ${\cal F}_n'={\cal F}_n\cup \cup_{j=1}^{n-1}(L_j({\cal F}_j').$ {{We further assume that
\beq\label{july-2021}
|\tau((L_{n+1}(p))-L_n(p))|<\dt_0/2^{n+3}\andeqn
&&\|h_n(p)-L_{n+1}(p)\|<\dt_0/2^{n+3},
\eneq
for all $\tau\in T(A)$ and $p\in {\cal P}_n$, where ${\cal P}_1={\cal P}'_1$ and ${\cal P}_n={\cal P}'_n\cup \bigcup_{i=1}^{n-1} h_{n-1}({\cal P}'_{n-1})$.}}

{{Let $\dt_1=\dt_0/4$. For $n>1$, there is $\dt_n<\min\{\dt_{n-1}/4,\dt_0/2^{n+1}\}$ such that if a completely positive map $\Theta: A \to E$ is ${\cal F}_n$-$\dt_n$ multiplicative ($E$ is any $C^*$-algebra), then $\Theta|_{{\cal P}_{n-1}}$ is well defined.}} We also assume that both ${{\Psi_n}}$ and $L_n$ are ${\cal F}_n$-$\dt_n$-multiplicative,
{{for $n=1,2,\cdots.$ Furthermore, we can strengthen (\ref{Sept-10-2020-1}) to the following
 \beq\label{July-31-2021}
\|a-\Psi_{n+k}(a)\oplus L_{n+k}(a)\|<\dt_{n+k}~(<\dt_n/2^{k+3})~~~\mbox{for all}~~a \in {\cal F}_n,
\eneq
for any positive integers $n, k$.}}

Define  $J_n: A\to A$ by $J_n(a)={{\Psi_n}}(a)\oplus L_n(a)$ for all $a\in A.$
Note that $J_n$ is ${\cal F}_n$-$\dt_n$-multiplicative.  We {{write}} $\imath_n: S_n\to A$ for the embedding.
Define $J_{m,n}=J_{n-1}\circ \cdots \circ J_m$ and $h_{m,n}=h_{n-1}\circ \cdots \circ  h_m: {{\tilde{S}_m\to \tilde{S}_n.}}$ 
Note also {{that}} $J_{m,n}$ is ${\cal F}_m$-$\dt_m$-multiplicative.

In what follows we will consider {{unitizations}}
${\tilde A}$ and  ${\tilde S}_n.$
We insist that $1_{\tilde A}=1_{\tilde S_n}.$
We will also use
$L_n, {{\Psi_n}},$ $J_n,$ $J_{m,n},$ $h_{m,n}$ for their  unitization and their extensions
on matrix algebras over $A$ (and over ${\tilde A}$). For example $L_n$ may also be viewed  a unital map
from ${{M_r}}({\tilde A})$ to $M_r({\tilde A}).$ Moreover,  it is important to keep in mind {{that}}{{$1_{\tilde A}=1_{\tilde S_m}=1_{\tilde S_n}$ and }}
$\imath_n(1_{\tilde S_n})=L_n(1_{\tilde A})={{\Psi_n(1_{\tilde A})=}}J_n(1_{{\tilde S_n}})=h_{m,n}(1_{\tilde S_m})
=J_{m,n}(1_{{\tilde A}}){{=1_{\tilde A}}}$ for all $n$ and $m.$ {{(Note that $\Psi_n(A)\perp L_n(A)$ does not contradicts with $L_n(1_{\tilde A})=\Psi_n(1_{\tilde A})=1_{\tilde A}$.)}} {{It follows from (\ref{Sept-10-2020-1}) that
\beq\label{Sept-10-2020-2}
[J_{m,n}]|_{\cal P}=[\id]|_{\cal P}, ~~\mbox{for any}~~m, n.
\eneq}}

In general, let $A$ and $B$ be non-unital \CA s and let $\phi, \psi: A\to B$ be \hm s
with orthogonal ranges. Suppose that ${\bar p}\in M_r(\C\cdot 1_{\tilde A})$ is a scalar projection and
$x\in M_r(A)_{s.a.}$ such that $p:={\bar p}+x$ is a projection.
Suppose that $\phi,\psi: {\tilde A}\to {\tilde B}$ are {{unitized so}} that $\phi(1_{\tilde A})=\psi(1_{\tilde A})=1_{\tilde B}.$
Consider also $H: {\tilde A}\to {\tilde B}$
defined by $H(1_{\tilde A})=1_{\tilde B}$ and $H(a)=\phi(a)+\psi(a)$ for all $a\in A.$
In $M_r({\tilde B}),$ if we identify scalar matrix in {{the}} obvious way,  then
$H(p)={\bar p}+H(x)={\bar p}+\phi(x)+\psi(x).$
On the other hand,
\beq
\phi(p-{\bar p})+\psi(p-{\bar p})&=&\phi({\bar p}+x)-{\bar p}+\psi({\bar p}+x)-{\bar p}\\
&=&{\bar p}+\phi(x)-{\bar p}+{\bar p}+\psi(x)-{\bar p}=\phi(x)+\psi(x).
\eneq
It follows that $[H]([p]-[{\bar p}]){{=[\phi]([p]-[{\bar p}])+[\psi]([p]-[{\bar p}])}}.$
\end{NN}

\begin{lem}\label{62320-1}
Let $A$ and $B$ be \CA s   such that  $T(A)$ and $T(B)$ are
compact
and $\phi_n: {\tilde A}\to {\tilde B}$ be a sequence {{of}} \cpc s
such that $\phi_n(1_{\tilde A})={{1_{\tilde B},}}$
{{$\phi(A)\subset B$}}  and
$\lim_{n\to\infty}\|\phi_n(a)\phi_n(b)-\phi_n(ab)\|=0$ for all $a, b\in A.$
Let $p:={\bar p}+x\in M_r({\tilde A})$ be a projection, where ${\bar p}\in M_r(\C\cdot 1_{\tilde A})$
is a scalar projection with rank $R$ and $x\in M_r(A)_{s.a.}.$
Suppose that
$
\lim_{n\to\infty}\sup\{|\tau(\phi_n(x))|: \tau\in T({{B}})\}=0.
$
Then
\beq
\lim_{n\to\infty} \sup\{|\tau([\phi_n]([p]))-R|: \tau\in T({{B}})\}=0.
\eneq
\end{lem}

\begin{proof}
Suppose that there is $\ep_0>0,$  $\tau_k\in T({{B}})$ and a subsequence $\{n_k\}$ such that
\beq
|\tau_k([\phi_{n_k}]([p]))-R|\ge \ep_0.
\eneq
Define $\Phi: {{\tilde{A}\to l^\infty(\tilde{B})}}$ be defined by
$\Phi(a)=\{\phi_{n_k}(a)\}$ for $a\in A.$
Let $\omega$ be a free ultra  filter and let $J_\omega=\{\{b_k\}; \lim_{\omega}\|b_k\|=0\}.$
Let $\Pi: l^\infty({{\tilde{B}}})\to l^\infty({{\tilde{B}}})/J_\omega$ be the quotient map.
Then $\Psi:=\Pi\circ \Phi$ is a unital \hm.
Therefore $\Psi(p)={\bar p}+\Psi(x),$ where we identify the scalar projection ${\bar p}$ with the one
in $M_r(\C\cdot 1_{\Pi(l^{\infty}({{\tilde{B}}})})).$
Note that $[\Pi]([\Phi]([p]))=(\Pi\circ \Phi)_{*0}([p]).$
Let $\tau_\omega(\Pi_\omega(\{b_k\}))=\lim_\omega \{\tau_k(b_k)\}.$
Then $\tau_\omega$ is a tracial state.
One computes that
\beq\nonumber
|\tau_\omega{{([\Pi]([\Phi]([p])))}}-R|&=&|\tau_\omega(\Psi(p))-R|
=|\tau_\omega(\Psi(x))|=
|\lim_\omega\{{{\tau_k}}(\phi_{n_k}(p))-R|\\\nonumber
&=&
|\lim_{\omega}\{\tau_k(\phi_{n_k}(x))\}|=0\andeqn\\\nonumber
|\tau_\omega{{([\Pi]([\Phi]([p])))}}-R|&=&\lim_\omega\{|\tau_k([\phi_{n_k}]([p]))-R|\}\ge \ep_0.
\eneq
A contradiction. This implies that
\beq
\lim_{k\to\infty}\sup\{|\tau(\phi_k(p))-R|:\tau\in T(A)\}=0.
\eneq

\end{proof}

\begin{lem}[Lemma 2.7 of \cite{Lnduke}]\label{traces}
Let $\mathcal P\subset \text{M}_r({\tilde A})$ (for some integer $r\ge 1$) be a finite set of
projections.
Assume that $\mathcal F_1$
is sufficiently large and $\delta_0$ is
sufficiently small such that $[L_{n}\circ J_{1,n}]|_{\mathcal P}$ and $[L_{n}\circ J_{1,n}]|_{G_0}$ are well defined, where $G_0$ is the subgroup generated by $\mathcal P$. Then
\beq\label{62320-10}
\lim_{n\to\infty}\sup_{\tau\in T(A)}|\tau([\iota_{n+1}\circ {{L_{n+1}}}\circ J_{1,n}]([p]))-\tau([p])|=0
\eneq
for any $p\in \mathcal P.$
%
%
{{Moreover,}}
for any $k\ge 1,$ we have
\beq\label{62320-10+}
|\tau(h_{k, k+n+1}\circ {{[L_{k}]}}([p]))-\tau{{(h_{k, k+n}\circ [L_{k}]([p]))}}|<(1/2)^{n+k}
\eneq
for all $\tau\in T(A)$ and
\beq\label{62320-10++}
\lim_{n\to\infty} \tau{{(h_{k,k+n}\circ [L_{k}]([p]))}}\ge (1-\sum_{i=1}^n 1/2^{i+k})\tau({{[L_{k}]}}([p]))>0
\eneq
for all  $p\in {\cal P}$
and $\tau\in T(A).$ {{Furthermore,
\beq\label{July-31-2021-1}
|\tau(h_{k, k+n+1}([p]))-\tau{{(h_{k, k+n}([p]))}}|<(1/2)^{n+k}~\mbox{and}
\\ \label{July-31-2021-2}
\lim_{n\to\infty} \tau{{(h_{k,k+n}([p]))}}\ge (1-\sum_{i=1}^n 1/2^{i+k})\tau([p])>0,
\eneq
for all $p\in {\cal P}_k$ and $\tau\in T(A)$.}} 

\end{lem}

\begin{proof}
Note {{that}}, by \eqref{Sept-10-2020-1}, we may assume that
$[J_{1,n}]([p])=[p]$ for $p\in {\cal P}.$ So, for \eqref{62320-10}, it suffices
to show $\lim_{n\to\infty}\sup_{\tau\in T(A)}|\tau([\iota_{n+1}\circ {L_{n+1}}]([p]))-\tau([p])|=0.$
Let us write  $p=\bar{p}+x$ with $\bar{p}\in M_r(\C\cdot 1_{\tilde{A}})$ and $x\in M_r(A).$
Note that all maps involved are unitized at $1_{\td A},$ we know $L_{n+1}(p)=\bar{p}+L_{n+1}(x)$ and $\Psi_{n+1}(p)=\bar{p}+\Psi_{n+1}(x)$.
By \eqref{Sept-10-2020-1} again, it suffices to show that, for all $p\in {\cal P},$
\beq\label{Ltrace-n1}
\lim_{n\to\infty}\sup_{\tau\in T(A)}|\tau{{([\Psi_{n+1}]([p]))}}-\tau({\bar p})|=0.
\eneq
Then,
by (\ref{C210GLN-n2}), $ |\tau(\Psi_{n+1}(x))|<r/2^{n+1}$. Hence \eqref{Ltrace-n1} follows from Lemma \ref{62320-1}.

For \eqref{62320-10+}, choose a projection {{$q$ with $[q]=[L_k]([p])$}} as in {{\eqref{62320-nn11}.}}
Write $q={\bar q}+x,$ where ${\bar q}\in {{M_r({\tilde S_{k}})}}$ is a scalar projection and
$x\in M_r({{S_{k}}})_{s.a.}$
Note {{that,}} again, all maps involved are unitized at $1_{\tilde A}.$

So, by  \eqref{62320-n11} and \eqref{62320-nn11},
\beq
\sup\{|\tau(h_{k, k+n+1}\circ {{[L_{k}]}}([p]))-\tau{{(h_{k, k+n}\circ [L_{k}]([p])|:\tau\in T(A)}})\}
<{{\frac{1}{2^{k+n}}}}.
\eneq
Thus \eqref{62320-10++} {{also}} holds. By (\ref{July-31-2021}), we may assume that
$[J_{k,n+k}]([p])=[p]$ for $p\in {\cal P}_k.$ {{The same argument shows that}}
 (\ref{July-31-2021-1}) and (\ref{July-31-2021-2}) {{also hold,}}  replacing \eqref{62320-n11} and \eqref{62320-nn11} by (\ref{july-2021}).
\end{proof}

\begin{rem}\label{Kstate}
Since $A$ is stably finite and assumed to be amenable, therefore exact, any positive state of $K_0(A)$ is the restriction of a tracial state of
$A$ ({{\cite{BR}  and \cite{Haagtrace}}}).  Thus, the lemma above still holds if one replaces the trace $\tau$ by any positive state $\tau_0$ on $K_0(A)$.  
\end{rem}

\begin{NN}\label{Constr-1508}
Fix ${\cal P}\subset K_0(A)$ with the form
${\cal P}=\{[p_i]-[{\bar p}_i]: i=1,2,...,j\},$
where $p_i\in  M_r({\tilde A})$ is a projection,
${\bar p}_i\in M_r(\C\cdot 1_{\tilde A})$ is a scalar projection such that $p_i-{\bar p}_i\in M_r(A)$
(for some $r\ge 1$)
and an integer $N\ge 1$ such that $[L_{N+i}]|_{\cal P},$
$[J_{N+i}]|_{\cal P}$ and $[\Psi_{N+i}]|_{\cal P}$ are all well defined.
Keep notation in \ref{C201GLN}.
Then, on ${\cal P},$ {{(note that we use ${\cal P}\subset K_0(A)$ and $\Psi_n(A)\perp L_n(A)$),}}
\beq\nonumber
[L_{N+1}\circ J_N]&=&[L_{N+1}\circ L_N]\oplus [L_{N+1}\circ \Psi_N]\\
&=& {{[h_{N}\circ L_N]}}\oplus [L_{N+1}\circ  \Psi_N], \andeqn
\eneq
\beq\nonumber
\hspace{-0.6in}[L_{N+2}\circ J_{N, N+2}]&=&{{[L_{N+2}\circ J_{N+1}\circ J_{N}]}}\\ \nonumber
&=&[L_{N+2}\circ L_{N+1}\circ J_N]\oplus [L_{N+2}\circ \Psi_{N+1}\circ J_N]\\ \nonumber
&=&[L_{N+2}\circ L_{N+1}\circ L_N]\oplus [L_{N+2}\circ L_{N+1}\circ \Psi_N]\\ \nonumber
&&\oplus [L_{N+2}\circ  \Psi_{N+1}\circ J_N]\\ \nonumber
\hspace{-0.1in}&=&{{[h_{N, N+2}]}}\circ [L_N]\oplus [L_{N+2}\circ L_{N+1}\circ \Psi_N]\oplus
{{[L_{N+2}\circ  \Psi_{N+1}\circ J_N].}}
\eneq
Moreover,  on ${\cal P},$
\beq\nonumber
[L_{N+n}\circ {{J_{N, N+n}}}] &=& {{[h_{N, N+n}]}}\circ [L_N]\oplus [L_{N+n}\circ \Psi_{N+n-1}\circ J_{N, N+n-1}]\\ \nonumber
&&\oplus [L_{N+n}\circ L_{N+n-1}\circ \Psi_{N+n-2}\circ J_{N, N+n-2}]\\ \nonumber
&&\oplus  [L_{N+n}\circ L_{N+n-1}\circ L_{N+n-2}\circ \Psi_{N+n-3}\circ J_{N, N+n-3}]\\ \nonumber
&&\oplus \cdots  \oplus{{[L_{N+n}\circ L_{N+n-1}\circ\cdots \circ L_{N+3}\circ \Psi_{N+2}\circ J_{N,N+2}]}}\\ \nonumber
&&{{\oplus}}[L_{N+n}\circ L_{N+n-1}\circ\cdots \circ L_{N+2}\circ \Psi_{N+1}\circ J_N]\\ \nonumber
&& \oplus [L_{N+n}\circ L_{N+n-1}\circ\cdots \circ L_{N+1}\circ \Psi_N].
\eneq

Set $\psi^N_{N}=L_N,$ $\psi^N_{N+1}=L_{N+1}\circ \Psi_N,$
$\psi^N_{N+2}=L_{N+2}\circ \Psi_{N+1},..., \psi^N_{N+n}=L_{N+n}\circ \Psi_{N+n-1},$
$n=1,2,....$
We may also assume $\psi_{k+j}^k(1_{\td A})=1_{\td S_{k+j}}.$

(Note that $\psi^N_{N+i}=\psi^{N+1}_{N+i}=\cdots=\psi^{N+i-1}_{N+i}$.  We insist {{on}} the notation $\psi^N_{N+i}$
{{in order}} to emphasize {{that}} our estimation {{begins}}  with a fixed index $N.$)
Moreover, for $j\geq 1$, 
\beq\label{1207-8-im-1}
{{d_\tau(\psi_{k+j}^k(e_A))<1/2^{k+j+1}\rforal \tau \in T(S_{k+j}).}}
\eneq

\end{NN}


\begin{NN}\label{Constr-Sn}
(a) For each $S_n$, since the abelian group $K_0({\tilde S_n})$ is finitely generated and torsion free, there is a set of free generators $\{e^{n}_1, e^{n}_2, ..., e^{n}_{l_n}\}\subseteq K_0({\tilde{S_n}})$.
We may assume that  $e^{n}_1=[1_{S^n}]$ and $\{e^n_2,e^n_3,\cdots,e^{n}_{l_n}\}\subset K_0(S_n)$ generate $K_0(S_n)$.
{{Recall {{that}} from \ref{C201GLN},}}
{{$K_0({\tilde S_n})_+$ is
generated}}
by $\{s^n_1, s^n_2, ..., s^n_{r_n}\}\subseteq K_0({\tilde{S_n}})_+$. Then there is an
$s_n\times l_n$ integer-valued matrix $R'_n$ such that
$${{\vec{s}_n}}=R'_n\vec{e}_n,$$ where $\vec{s}_n=(s^n_1, s^n_2, ..., s^n_{r_n})^{\mathrm T}$ and $\vec{e}_n=(e^{n}_1, e^{n}_2, ..., e^{n}_{l_n})^{\mathrm T}$. In particular, for any ordered group $H$, and any elements $h_1, h_2, ..., h_{l_n}\in H$, the map $e^n_i\mapsto h_i$, $i=1, ..., l_n,$ induces an abelian-group homomorphism $\phi: K_0({\tilde S_n})$ to $H$, and the map $\phi$ is positive (or strictly positive) if and only if
$$R'_n\vec{h}\in H^{r_n}_+ \quad \textrm{(or $R_n'\vec{h}\in {(H_+\setminus \{0\})^{r_n}}$ ),}$$
where $\vec{h}=(h_1, h_2, ..., h_{l_n})^{{T}} \in H^{l_n}$.
Moreover,  we may write $e^n_1=e^n_{1+},$
and, for $i>1,$  {{$e^n_i=e^n_{i,+}-e^n_{i,-}$}} for
{{$(e^n_i)_+, (e^n_i)_-\in K_0(\td S_n)_+$}}
and fix this decomposition.
Define  a $r_n\times (2l_n-1)$ matrix
\begin{displaymath}
R_n=
R'_n
\left(
\begin{array}{ccccccc}
1 & 0 & 0 & 0 &\cdots & 0 & 0 \\
0 & 1 & -1 & 0 & \cdots & 0 & 0 \\
\vdots & \vdots & \vdots & \vdots & \ddots & \vdots & \vdots \\
0 & 0 & 0 & 0 & \cdots & 1 & -1
\end{array}
\right).
\end{displaymath}
Then one has \beq\label{05-18-2021-1}\vec{s}_n=R_n\vec{e}_{n, \pm},\eneq where $\vec{e}_{n, \pm}=({{e^{n}_{1,+}, e^{n}_{2,+},
e^n_{2,-},..., e^{n}_{l_n,+}, e^{n}_{l_n,-}}})^{\mathrm T}$. Hence, for any ordered group $H$, and any elements $h_{1, +}, h_{2,+},h_{2, -}, ..., h_{l_n, +}, h_{l_n, -}\in H$, the map $e^n_1\mapsto h_{1+},$ $e^n_i\mapsto (h_{i, +}-h_{i, -})$, $i=2, ..., l_n$
(with  induces a positive (or strictly positive) homomorphism if and only if
$$R_n\vec{h}_{{\pm}}\in H^{r_n}_+\quad \textrm{ (or $R_n\vec{h}_{{\pm}}\in (H_+\setminus\{0\})^{r_n}$) },$$
where $\vec{h}_{\pm}=(h_{1, +}, h_{2, +},  h_{2,-},..., h_{l_n, +}, h_{l_n, -})^{{T}}   \in H^{2l_n-1}$.

{ {(b) Let $A\in D^d$ be a separable simple \CA\, with continuous scale  and let $B\in {\cal M}_1$
such that
$$((K_0(A), \{0\},   T(A),
\rho_A), K_1(A){{)}}\cong ((K_0(B), \{0\}, T(B), \rho_B), K_1(B)).$$
Let $\alpha\in KL(A, B)$ be
{{an}} element {{and $\kappa_T: T(B) \to T(A)$ be an affine homeomorphism}} which {{implements}} the isomorphism above. Then, for any $x\in K_0({A}), \tau \in T(B)$, we have $\rho_B(\af(x))(\tau)=\rho_A(x)(\kappa_T(\tau))$---that is $\tau(\af(x))=\kappa_T(\tau)(x).$

Let $e^n_1, e^n_i, e^n_{i,\pm} \in K_0(\td S_n)$ be as  above, $i=2,3,..., l_n.$
Let $s(0)=0, s(n)=\sum_{j=1}^n(2l_j-1)$.
Put
\beq\label{05-18-2021-2}
&&\af(\imath_n\circ  h_{j, n}(e^j_{1,+}))=g^{(n)}_{s(j-1)+1}(=1_{\tilde B}),\,\, {\rm and,}\,\, {\rm for}\,\, {{i\in \{2,3,\cdots l_j\}}},\\
&&\alpha(\imath_n\circ  h_{j, n}(e^j_{i,+}))=g^{(n)}_{s(j-1)+{{2i-2}}},
~~\alpha(\imath_n\circ h_{j, n}(e^j_{i,-}))=g^{(n)}_{s(j-1)+{{2i-1}}},\eneq
{and} {{$g^{(n)}_l=0~\mbox{if}~ l>s(n).$}}
Let $a^{(n)}_j=\rho_B(g^{(n)}_j) \in \Aff({{T(B))_+}}$ for $j=1,2,{ ...}.$ Then{{,}} by Lemma \ref{traces}, $\{a^n_j\}_n$ uniformly converge to {{some}} $a_j>0$ on $T(B)$. 

For $j\in \{1,2,{ ...}, n\}$,
let
$(s_1^j, s_2^j, { ...}, s_{r_j}^j) \in K_0(\td S_j)_{ +}\setminus \{0\}$
be
the generators of the positive cone {{$K_0(\td S_j)_+,$}}
and {{let $R_j$ be  the $r_j\times (2 l_j-1)$}}  {{matrix}} as {{in}} part (a). {{From (\ref{05-18-2021-1}) and {{(\ref{05-18-2021-2})}}, for $j\leq n$}}
\beq
&&R_j (g^{(n)}_{s(j-1)+1}, g^{(n)}_{s(j-1)+2}, { ...}, g^{(n)}_{s(j)})^T\\
&=&(\alpha(\imath_n\circ h_{j,n}(s_1^j)), \alpha(\imath_n\circ h_{j,n}(s_2^j)), { ...}, \alpha(\imath_n\circ h_{j,n}(s_{r_j}^j)))^T
\andeqn\\
&&R_j (a_{s(j-1)+1}, a_{s(j-1)+2}, { ...}, a_{s(j)})^T\\
&=&\lim_{n\to \infty} (\rho_B(\alpha(\imath_n\circ h_{j,n}(s_1^j))), \rho_B(\alpha(\imath_n\circ h_{j,n}(s_2^j))), { ...}, \rho_B(\alpha(\imath_n\circ h_{j,n}(s_{r_j}^j))))^T.
\eneq
By Lemma \ref{traces} {{and}} (\ref{July-31-2021-2}) (applied to $\kappa_T(\tau)$)}} 
for all $1\le j\le n,$ and  for all $\tau\in T(B),$ {{we obtain}}
$$
\tau (\alpha(\imath_n\circ h_{j,n}(s_i^j)))>(1-\sum_{k=1}^{\infty} 1/2^{j+k} )\tau{{(\af(\imath_j (s_i^j)))}}>0,{{\rforal}} i \in\{1,2,{ ...}, r_j\}.
$$
Hence each entry of $R_j (a_{s(j-1)+1}, a_{s(j-1)+2}, { ...}, a_{s(j)})^T$  is a strictly positive element of \\
$\Aff_+(T(B))\setminus\{0\}.$
Let ${\bar R}_n =\diag (R_1, R_2,..., R_n)$. Then
$$0\ll {\bar R}_n (a_1, a_2, ..., a_{s(n)})^T\in (\Aff({{T(B)}}))^{\sum_{i=1}^nr_k},$$
i.e.,  each coordinate  is strictly positive  on ${{T(B)}}$.

Furthermore, ${\bar R}_n(g^{(n)}_1, g^{(n)}_2, ..., g^{(n)}_{s(n)}) \in (K_0(B)_+\setminus \{0\})^{\sum_{i=1}^nr_k}.$ In particular, for each positive integer $N_0 <n$, we also have
$$\diag (R_{N_0+1},R_{N_0+1}, ..., R_n) (g^{(n)}_{s(N_0)+1}, g^{(n)}_{s(N_0)+2}, ..., g^{(n)}_{s(n)}) \in (K_0(B)_+\setminus \{0\})^{\sum_{i=N_0+1}^nr_k}.$$

(c)  {{Let $p\in M_r(\td S_n)$ be a projection and  ${\bar p}\in M_r(\C \cdot 1_{\td S_n})$  be a
(scalar) projection
such that $p-{\bar p}\in M_r(S_n).$}}
Put $z=[p]-[{\bar p}].$  {{Since $\{e^{n}_1, e^{n}_2, ..., e^{n}_{l_n}\}$ is a set of free generators of $K_0({\tilde S_n})$,}}
there is a unique  $l_n$-tuple of {\it integers} $m^n_1(z), ..., m^n_{l_n}(z)$ such that $z=m^n_1(z)e^n_1+m^n_2(z)e^n_2+\cdots+ m^n_{l_n}(z)e^n_{l_n}$.
 Hence for any homomorphism
$\tau: K_0({\tilde S_n})\to\mathbb R$, one has
$$
\tau(z)=\langle\vec{m}_n(z), \tau(\vec{e}_n)\rangle=\sum_{i=1}^{l_n} m_i^{{n}}(z)\tau(e^n_i)=m^n_1(z) e^n_1+
\sum_{i=2}^{l_n} m^n_i(z)\tau((e^n_{i, +})-m^n_i(z)\tau((e^n_{i,-}),$$
where $\vec{m}_n{{(z)}}={{(m^n_1(z), ..., m^n_{l_n}(z))}}$ and {{$\vec{e}_n=(e^{n}_1, e^{n}_2, ..., e^{n}_{l_n})^{T}.$}} 
Recall that $e^n_1=e^n_{1,+}=[1_{\td S_n}]$ and
$S_n$ is stably projectionless.
Note also that $m^n_1(z)=0,$ if $z\in K_0(S_n),$ and $e^n_i\in K_0(S_n)$ for $i\ge 2.$

Note that we assume $1_{\tilde{S}_n}=1_{\tilde A}$. There is a nature  {{rank}} map
$K_0(C) \to K_0(\C)=\Z$, for $C={\tilde S}_n$ or $C={\tilde A}$. Then $m^n_1(p)={\rm rank}(p)$. If $z=[p]-[{\bar p}]$, where ${\bar p}$, is a scalar matrix such that $p-{\bar p}\in M_r(S_n)$, 
then $m^n_1(z)=0$ and $m^n_i(z)=m^n_i(p)$ for $i\geq 2$.

For each $p\in M_m({\tilde A})$ {{(for some integer $m\ge1),$}} denote by $[\psi_{{k}, k+j}(p)]$ an element
in $K_0({{{\tilde S_{k+j}}}})$ associated with $\psi_{ k+j}^k(p).$  Let $\imath_n: S_n\to A$ be the imbedding.
Denote by
$${\overline{(\imath_n)_{*0}}}: {{\vec{e}_n}}\mapsto (((\imath_n)_{*0}(e_1^n), (\imath_n)_{*0}(e_2^n),...,(\imath_n)_{*0}(e_{l_n}^n)),$$
where $(\imath_n)_{*0}(e_1^n)=[1_{\tilde A}]$.

\end{NN}
Then, by Lemma \ref{traces} and Remark \ref{Kstate}, one has the following lemma.

\begin{lem}\label{rho-uniform}  
With the notion same as {{above,}} {{in particular ${\cal P}\subset K_0(A)$,}} for any $z\in\mathcal P,$  for each fixed $k$, one has that
$$
\tau(z)=
\lim_{n\to\infty}\sum_{{j=0}}^n({{\sum_{i=2}^{{l_{k+j}}} }} m_i^{{k+j}}([{{\psi^k_{k+j}}}(z)])\tau((\imath_{k+n}\circ h_{k+j,k+n})_{*0}({{e_{i,+}^{k+j}}}))
$$
$$~~~~~-m_i^{{k+j}}({{[\psi^k_{k+j}(z)]}})\tau((\imath_{k+n}\circ h_{k+j,k+n})_{*0}({{e_{i,-}^{{k+j}}}})))
$$
uniformly on  $T(A)$  ({{recall}} that $m_1^{{k+j}}([\psi^k_{k+j}(z)])=0$ since $z\in K_0(A)$).
Moreover,  $\rho_A\circ (\imath_n)_{*0}\circ h_{k+j,k+n}(e_{i,\pm}^{k+j})$
converge{{s}} to a strictly positive element in $\mathrm{Aff}(T(A))$
as $n\to\infty$ uniformly {{(recall $h_{k+j,k+n}(e_1^{k+j})=e_1^{k+n}=1_{\tilde A}$).}}
\end{lem}

\begin{proof}
For $z\in {\cal P},$ we write $z:=[p]-[{\bar p}]\in {\cal P},$ where
$p\in M_r({\tilde A})$ is a projection, ${\bar p}\in M_r(\C\cdot 1_{\tilde B})$ is scalar projection
and $p-{\bar p}\in M_r(A).$
We first compute that, if {{$j\geq 1$}} 
for $z\in {\cal P},$
\beq\label{rho-comput-1508}
&&\sum_{i=1}^{{l_{k+j}}} m_i^{{k+j}}([\psi^k_{k+j}](z))\tau((\imath_{k+n}\circ h_{k+j,k+n})_{*0}(e_i^{{k+j}}))\\
&=& \tau([L_{k+n}\circ\cdots \circ L_{k+j+1}\circ L_{k+j}\circ \Psi_{k+j-1}](z))
\eneq
(not the first term is zero, as $e^n_1=[1_{\td S_n}]$ and $e_i^{(k+1)}\in K_0(S_n)$ for $i\ge 2$.)
and, if {{$j=0$}} 
\beq\label{rho-compute-1508-1}
&&\sum_{i=1}^{{l_{k}}} {{m^{k}_i}}({{[\psi^k_{k}]}}(z))\tau((\imath_{k+n}\circ h_{{k,k+n}})_{*0}({{e^{k}_i}}))\\
&=& \tau([L_{k+n}\circ\cdots \circ L_{{k+1}}\circ L_k](z)).
\eneq
Thus  (see \ref{Constr-1508})
\beq\label{comput-1508-2}
&&{{\sum_{j=0}^n}} (\sum_{i=1}^{l_{k+j}} m^{k+j}_i([\psi^k_{k+j}(z)])\tau((\imath_n\circ h_{k+j,k+n})_{*0}(e^{k+j}_i)))\\
&=& \tau([L_{k+n}\circ\cdots \circ {{L_{k+1}}}\circ L_k](z))\\
&&+{{\sum_{j=1}^n}}\tau([L_{k+n}\circ\cdots \circ L_{k+j+1}\circ L_{k+j}\circ \Psi_{k+j-1}](z))\\
&=& \tau([L_{k+n}\circ\cdots \circ {{L_{k+1}}}\circ L_k](z))+{{\tau([L_{k+n}\circ\cdots \circ L_{k+1}\circ \Psi_k](z))}}\\
&&{{+\tau([L_{k+n}\circ\cdots \circ L_{k+2}\circ \Psi_{k+1}]\circ[J_k](z))}}\\
&&+{{\sum_{j=3}^n}}\tau([L_{k+n}\circ\cdots \circ L_{k+j+1}\circ L_{k+j}\circ \Psi_{k+j-1}\circ J_{k,k+j-1}](z))\\
&=&\tau([L_{k+n}\circ {{J_{k, n+k}}}](z)).
\eneq
Thus the first part of the lemma {{follows}} from \ref{traces}. The second part also follows.

\end{proof}

One then has the following
\begin{cor}\label{rho}
Let $\mathcal P$ be a finite subset of projections in a matrix algebra over ${\tilde A}$,  let $G_0$ be the subgroup of $K_0({\tilde A})$ generated by $\mathcal P$ {{ and let $k\ge 1$ be an integer}}.  Denote by $\tilde{\rho}: G_{{0}}\to\Pi\mathbb Z$ the map defined  (see the {{(c)}} of \ref{Constr-Sn}) by
\noindent
\begin{eqnarray}
&&\hspace{-0.7in}[p] \mapsto
({{{m}^{k}_1(q_0),}}  
~ {m}^{k}_2(q_0), -m^{k}_2(q_0),\cdots m^{k}_{l_k}(q_0), -m^k_{l_k}(q_0), \nonumber\\
&&{{m^{k+1}_1(q_1),}}  
~ m^{k+1}_2(q_1), -m^{k+1}_2(q_1), ...,m^{k+1}_{l_{k+1}}(q_1), -m^{k+1}_{l_{k+1}}(q_1),\nonumber\\
&&{{m^{k+2}_1(q_2),}}  
~ m^{k+2}_2(q_2), -m^{k+2}_2(q_2), ...,m^{k+2}_{l_{k+2}}(q_2), -m^{k+2}_{l_{k+2}}(q_2), \cdots), \label{rho-1}
\end{eqnarray}
where $q_i=[\psi^{{{k}}}_{k+i}(p)],$ $i=0,1,2,....$
If $\tilde{\rho}(g)=0$, then $\tau(g)=0$ for any trace over $A$ {{(recall that $m_1^{k+i}(q_i)={\rm rank}(p)$ for all $i$)}}
\end{cor}


%
%
%


\begin{NN}\label{D2010gln}
Let $S$ be a compact convex set, and let $\Aff(S)$ be the space of real affine continuous functions on $S$. Let $\mathbb{D}$ be an  ordered subgroup of $\Aff(S)$ with  {{the}} form $\Z\cdot 1+\mathbb{D}_0,$
where $\mathbb{D}_0$ is dense in $\R\D_0$ {{and $\mathbb{D}_0\cap \Aff_+(S)=\{0\}$.}}
Let $G$ be an abelian group  with the form
$G=\Z\cdot g_e+G_0$ for some $g_e\in G\setminus \{0\}$ and $G_0\subset G$ is a subgroup. 
Let $\rho: G\to \mathbb{D}$ be a {{surjective}} \hm\, such that $\rho(g_e)=1$ {{and $\rho(G_0)\subset \mathbb{D}_0$.}}
Define $G_+$ to be the set of those elements $g=m\cdot g_e+g_0,$ where $m\in \Z$ and $g_0\in G_0,$
such that $m>0$
and
$\rho(g)>0,$ and  {{the zero element $g=0.$}}
We further assume that $(G, G_+)$ is an ordered group.
In the next lemma, for $g\in G,$
we write $g=J(g)g_e+o(g),$ where
$J(g)\in \Z$ and $o(g)\in G_0,$ and for $d\in \D,$ we write $d=J(d)+o(d),$ where $J(d)\in \Z$ and $o(d)\in \D_0.$
{\it We assume that $g>0$ implies that $J(g)>0.$}
Let $r\in \N.$ Denote by $G^r=\{(g_1, g_2,...,g_r)^T: g_i\in G,1\le i\le r\}$ (as columns),
$\D^r=\{(d_1, d_2,...,d_r)^T: d_i\in \D, 1\le i\le r\},$ {{and}}  {{$\rho^r: G^r\to \D^r$}} 
the map defined by $\rho^r((g_1,g_2,...,g_r)^{{T}})=(\rho(g_1), \rho(g_2),...,\rho(g_r))^T$ for $(g_1, g_2,...,g_r)^T\in G^r.$ {{For convenience, for any $n> r$, we also {{write}}
$\rho^r$ for the map $\rho^r: G^n\to \D^r$ defined by $\rho^r((g_1,g_2,...,g_r,\cdots, g_n)^{T})=(\rho(g_1), \rho(g_2),...,\rho(g_r))^T$.}}
{{For any $n\geq r$, we also define $J^r: G^n\to \Z^r$ and $o^r: G^n\to G_0^r$  }} by $J^r({\tilde g})=(J(g_1), J(g_2), ...,J(g_r))^T$ and $o^r({\tilde g})=(o(g_1), o(g_2),...,o(g_r))^T,$
 if ${\tilde g}=(g_1,g_2,...,g_n)^T.$
Similarly, if ${\tilde d}\in {{\D^n}},$ $J^r({\tilde d})=(J(d_1), J(d_2),...,J(d_r))^T$
and $o^r({\tilde d})=(o(d_1), o(d_2), ...,o(d_r))^T,$ where ${\tilde d}=(d_1, d_2,..,d_{{n}})^T.$
In these cases, we may also identify $J^r(g)$ with $(J(g_1),...,{{J(g_r)}},0,...,0),$ for example,
if it is convenient. These notation will be used below.

\end{NN}

\begin{lem}[Lemma 3.4 of \cite{Lnduke}]\label{solveeq}
Let $S,$ $G,$ $G_0,$ $\D,$  $\D_0$  and $\rho$ be as above.
Let $\{{{x_{i,j}}}\}_{1\leq i\leq r, 1\leq j<\infty}$ be {{an}} $r\times\infty$ matrix having rank $r$ and with ${{x_{i,j}}}\in\mathbb Z$ for each $i,j$.  Let $g_j^{(n)}\in G$ be such that $\rho(g_j^{(n)})=a_j^{(n)}$, where $\{a_j^{(n)}\}$ is a sequence of positive elements in $\mathbb{D}$ such that $a_j^{(n)}\to a_j (>0)$ uniformly on $S$ as $n\to\infty$
and
$J(g_j^{(n)})\to J(a_j)(>0).$

{Further suppose that there is a sequence of integers $s(n)$ satisfying the following condition:}

Let $\widetilde{v_n}=(o(g_1^{(n)}), o(g_2^{(n)}),...,o(g_{s(n)}^{(n)}))$
be the part of $(o(g_1^{(n)}), o(g_2^{(n)}) , o(g_j^{(n)}),...)$
and let
$$\widetilde{y_n}=({{x_{i,j}}})_{r\times s(n)}\widetilde{v_n}.$$
Denote by $y_n=\rho^{(r)}(\widetilde{y_n}).$
Then there exists $z=(z_j)_{r\times1}$ such that  $y_n\to z$ on $S$ uniformly.


{With the above conditions,  there exist $1/4>\delta>0,$  and {{integers $K>0$  and $N>0$}}
satisfying the following:}

{For  any $n\ge N,$
if $M$ is a positive integer,  and if {{${{\tilde {z}'}}\in G^r$ satisfies}} $o(\tilde{z}')\in (K^3G_0)^r$ (i.e.,  there is $\tilde{z}''\in G_0^r$ such that $K^3({\tilde{z}''})=o({\tilde{z}'})$){{, $|| o(z)-M(o({{z'}})) ||<\delta$}}, where
$\td z'=(\td z'_1,\td z'_2,\cdots, \td z'_r)^T$ {{and $z'=(z'_1,z'_2,\cdots, z'_r)$}} with  $z'_j=\rho(\tilde{z}'_j),$
 then there is  ${\tilde{u}_0}=(c_j)_{s(n)\times1}\in G^{s(n)}$ 
  such that
  \beq\label{2020-7-10-n++}
  ({{x_{i,j}}})_{r\times s(n)}{o({\tilde{u}_{0})}}=\tilde{z}'.
  \eneq}

 {Moreover, if each $s(n)$ can be written as $s(n)=\sum_{k=1}^n l_k$, where $l_k$ are positive integers,  and for each $k$,   $R_k$ is {{an}} $r_k\times l_k$ matrix with entries in $\Z$ so that }
$$
{\td R}_n={\rm diag}(R_1, R_2,...,R_n)
$$
satisfies {{that}}
\beq\label{June14-1}
{\td R}_n({\td g}_n)\in (G_+\setminus \{0\})^{\sum_{k=1}^n r_k}\andeqn  {\td R}_n({\td a}_n)\in(\Aff_+(S)\setminus \{0\})^{\sum_{k=1}^nr_k},
\eneq
where
${\td g}_n=(g_1^n,g_2^n,...,g_{s(n)}^n)^T$ and ${\td a}_n=(a_1, a_2,...,a_{s(n)})^T,$
$n=1,2,{{...,}}$
then  there exist $\delta, K , {{N}}$ as described above but {{also}} with {{${\tilde u}_0$ above}} 
satisfying an extra condition
that
\beq\label{JUne14-2}
{\td R}_n({\tilde u})>0,
\eneq
{{where ${\tilde u}=J^{s(n)}({\tilde u}_0)+M\cdot o^{s(n)}({\tilde u}_0)=(J(a_1), J(a_2),...,J(a_r))+M\cdot o^{s(n)}({\td u}_0)
$.}} 
\end{lem}

\begin{proof}
The proof repeats  the argument of Lemma 3.4 of \cite{Lnduke}.
But we will also show that
$u=(\tilde{c}_j)_{s(n)\times1}$ can be chosen to make (\ref{JUne14-2}) hold (see also  20.10 of \cite{GLN}).

Keep in mind that $J(g)$  and $J(d)$ (for all $g\in G$ and $d\in \D$)  are integers and, in particular, $J(a_j)$ is
a positive integer. We will also identify these integers with the {{integer-valued}} affine functions in $\Aff(S).$

Without loss of generality, we may  assume  that $({{x_{i,j}}})_{r\times r}$ has rank $r$.  Choose integer $N_0$ such that $s(N_0)\geq r$. Write
\beq
{\td R}_{N_0} (a_1, a_2, {...,} a_{s(N_0)})^T=(b_1, b_2, {...,} b_{\bar r})^{{T}}\in (\Aff_+(S)\setminus \{0\})^{\bar r}\andeqn\\
{\td R}_{N_0}(J(a_1), J(a_2),..., J(a_{s(N_0)})^T=(J_1,J_2,...,J_{\bar r})^T\in {{\N^r\subset}}(\Aff_+(S)\setminus \{0\})^{\bar r}
\eneq
where ${\bar r}:=\sum_{i=1}^{N_0}r_i$.
 {{Note that $\min_{1\le i\le {\bar r}} \{J_i\}\geq 1$.}}
  {{Let $\ep_0=\min\{1/4, \min_{1\leq i\leq {\bar r}}\inf_{s\in S} \{b_i(s)\}\}>0.$}}

Choose $0<\dt_0<{\ep_0\over{\|{\td R}_{N_0}\|+1}}$
 {{such that if
$$\|(a'_1, a'_2, {...,} a'_{s(N_0)})^T-(a_1, a_2, {...,} a_{s(N_0)})^T\|<\dt_0, $$
then
$$\|{\td R}_{N_0}(a'_1, a'_2, {...,} a'_{s(N_0)})^T-(b_1, b_2, {...,} b_{\bar r})^{{T}}\|<\ep_0/4.$$}}
We further assume that $\dt_0< \frac14\min_{1\leq j\leq s(N_0)}\inf_{s\in S}\{a_j(s)\}.$
Consequently, if $(h_1, h_2, {...,}, h_{s(N_0)})\in G^{{{s(N_0)}}}$ { satisfies}
$$\|(\rho(h_1), \rho(h_2), {...,}, \rho(h_{s(N_0)}))^T-(a_1, a_2, {...,} a_{s(N_0)})^T\|<\dt_0,$$
then
\beq\label{2017-may-9-1}
&&(h_1, h_2, {...,} h_{s(N_0)})\in (G_+\setminus \{0\})^{s(N_0)}, \andeqn\\
&&{\td R}_{N_0}(h_1, h_2, {...,} h_{s(N_0)})\ge 3\ep_0/4\,\,\, {\rm in}\,\,
(\Aff_+(S)\setminus \{0\})^{\bar r}.
\eneq
Choose $N_1\ge N_0$ such that, for all $n\ge N_1,$  and for  $1\le j\le s(N_0),$
\beq\label{2020-7-10-n1}
\|a_j-a_j^{(n)}\|<\dt_0/8,\,\, \|J(a_j)-J(a_j^{(n)})\|<\dt_0/8\andeqn \|o(a_j)-o(a_j^{(n)})\|<\dt_0/8.
\eneq
Since $J(a_j)$ and $J(a_j^{(n)})$ are integers,  for all $n\ge N_1,$
\beq\label{2020-7-10-n2}
J(a_j)=J(a_j^{(n)}),\,\, j=1,2,...,s(N_0).
\eneq
In particular we have the following claim (which will be used at the end of the proof): For all $n\ge N_1,$
if for some $f_{j}\in G_+,$
\beq\label{2020-7-10-n3}
||J(f_{j})-(J(a_j^{(n)})||<\dt_0/8\andeqn \|o(f_j)-o(a_j^{(n)})\|<\dt_0/8,
\eneq
 $j=1,2,...,s(N_0),$
then
\beq\label{2020-7-10-n4}
{\td R}_{s(N_0)}((f_1, f_2,...,f_{s(N_0)})^T)\in (\Aff_+(S)\setminus \{0\})^{{\bar r}}.
\eneq

 {{Recall}} that we assume {{that}} $A:=(x_{i,j})_{r\times r}$ has rank $r$.  There is an invertible matrix $B\in M_{r}(\mathbb Q)$ with $BA=I_r$. There is an integer $K>0$ such that all entries of $KB$ and $K(B)^{-1}$ are integers.  Choose a positive number $\dt<\dt_0$ such that
$\|B\|\dt<\dt_0/8.$

{{Recall that $\widetilde{y_n}=({{x_{i,j}}})_{r\times s(n)}\widetilde{v_n}$, $\widetilde{v_n}=(o(g_j^{(n)}))_{s(n)\times1}$ , $\rho(g_j^{(n)})=a_j^{(n)}$ and
$y_n=\rho^{(r)}(\widetilde{y_n})$, we have}}
\beq\label{Sept-12,2020}
y_n=({{x_{i,j}}})_{r\times s(n)} (o(a^{(n)}_1), o(a^{(n)}_2), ..., o(a^{(n)}_{s(n)}))^T.
\eneq

Since {{$y_n \to z$}}
 (uniformly on $S$) as $n \to \infty,$
choose $N\ge N_1$ such that, if $n\ge N,$
$$
\|y_n-z\|<\dt/16.
$$
It follows that (for $n\ge N$)
\beq\label{2017-may-9-2}
\|B(y_n)-B(z)\|<\frac18\dt_0
\eneq

Let us show that $K$, $N$, and $\dt$ {{as defined above}} are as desired.

Put $A_n=({{x_{i,j}}})_{r\times s(n)}.$ Then  $BA_n=C_n,$ where $C_n=(I_r, D_n')$ for some $r\times (s(n)-r)$ matrix $D_n'$. Since all entries of {{$A_n$ and}} $KB$  are integers, $KD_n'$ is also a matrix with integer entries.
Put $D_n=(0_{r\times r}, D_n').$

{{Recall}}  that $\rho(g^{(n)}_j)=a^{(n)}_j$, {{and so}} from the first part of  (\ref{June14-1}), we have
$${\bar R}_n ((a^{(n)}_1, a^{(n)}_2, {...,} a^{(n)}_{s(n)})^T)\in (\Aff_+(S)\setminus \{0\})^{\sum_{k=1}^n r_k}.$$
{{For each $n\ge N$ (by the continuity of the linear maps),}} there is
{{$0<\dt_n<\dt/4$}}
such that if $(x_1, x_2, { ...}, x_{s(n)})\in \Aff(S)^{s(n)}$ satisfies
$$\|(x_1, x_2, {...,} x_{s(n)})-(a^{(n)}_1, a^{(n)}_2, {...,} a^{(n)}_{s(n)})\|<\dt_n, $$
then
\beq\label{gg17-02}
{{\|(0_{r\times r}, D_n')(x_1, x_2, {...,} x_{s(n)})^T-(0_{r\times r}, D_n')(a^{(n)}_1, a^{(n)}_2, {...,} a^{(n)}_{s(n)})^T\|<\frac{\dt_0}4,}}
\eneq
and ${\td R}_n ((x_1, x_2, {...,} x_{s(n)}))^T \in (\Aff_+(S)\setminus \{0\})^{\sum_{k=1}^n r_k}$. In particular, we have
\beq\label{2017-may-9}
\diag(R_{N_0+1}, R_{N_0+2}, {...,} R_n) (x_{s(N_0)+1}, x_{s(N_0)+2}, {...,} x_{s(n)})^T \in (\Aff_+(S))^{\sum_{k=N_0+1}^n r_k}.
\eneq

\noindent
Since $\mathbb D_0$  is dense in $\R\D_0$ in $\mathrm{Aff}(S)$ {{and $\rho: G\to \D$ is surjective}}, there are $\xi_n\in G^{s(n)}$ such that $\xi_n=(d_j^{(n)})_{s(n)\times 1}$ and,
for all $n\ge N,$
\beq\label{equ5001}
||K^3\rho(o({d}_j^{(n)}))-\frac{o(a_j^{(n)})}{M}||<\dt_n/16M\andeqn
J({d}_j^{(n)})=J(a_j^{(n)}),\,\,\, j=1,2,....,s(n).
\eneq
Let $\rho({d}_j^{(n)}):={\bar d}_j^{(n)},$ ${\bar \xi}_n:=({\bar d}_j^{(n)})_{s(n)\times 1}$ and
let $\tilde{w}_n:=J^{s(n)}(\xi_n)+K^3o({{ \xi_n}}).$ 
Hence
\beq\label{2020-7-10-n500-}
\tilde{w}_n=(J(d_1^{(n)})+K^3o(d_1^{(n)}), J(d_2^{(n)})+K^3o(d_2^{(n)}),...,
J(d_{s(n)}^{(n)})+K^3 o(d_{s(n)}^{(n)})).
\eneq

Let {{$\tilde {z}'\in G^r$ and}}  $\tilde{z}''\in {{G_0^r}}$ be such that
$K^3({\tilde{z}''})={{o({\tilde{z}'})}}$ and  $||{{o(z)-M o(z')}} ||<\delta,$
where
$\td z'=(\td z'_1,\td z'_2,\cdots, \td z'_r)^T$ and $z'_j=\rho(\tilde{z}'_j),$  as
 described in the lemma.

 Since both $KB$ and $KD_n$ are matrices over $\mathbb Z$,
 \beq\label{2020-7-10-nn+1}
 {{u':=}}J^r(a^{(n)})+KB(\tilde{z}'')-KD_no^{s(n)}(\xi_n)\in G^{r},
 \eneq
   where $J^r(a^{(n)})=(J(a_1^{(n)}),J(a_2^{(n)}),...,J(a_r^{(n)}))^T.$
%
 {{Then,}} for $n\ge N,$
\beq\label{equ4999}
\hspace{-0.3in}B\tilde{z}'-D_no(\tilde{w}_n)&=&
K^3B \tilde{z}''-K^3D_no^{s(n)}(\xi_n)
=K^2o^{r}(u').
\eneq
{{(See \ref{D2010gln} for notation.)}}

Write  $o(u')=(\tilde{c}_1, \tilde{c_2}, ..., \tilde{c}_{r})\in G^{r}$
and  $$\rho^{(r)}(u')=(c_1, c_2, ..., c_r)\in\mathbb D^r.$$
Set
\beq\nonumber
{\td u}_0:
=(J(a_1)+K^2 \tilde{c}_1,..., J(a_r)+K^2 \tilde{c}_r, J(d_{r+1}^{(n)})+K^3o(\tilde{d}_{r+1}^{({n})}),..., J(d_{s(n)}^{(n)})+K^3o(\tilde{d}_{s(n)}^{({n})}))^T.
\eneq
{{From \eqref{2020-7-10-nn+1}
and (\ref{2020-7-10-n500-}), we obtain}}
\beq\label{2020-7-10-n8}
{{(I_r,0) o({\td u}_0)}}=B\tilde{z}'-D_{{n}}o(\tilde{w}_n)\andeqn {{(0,D'_n) {\td u}_0=(0,D'_n) {\td w}_n=D_{{n}}\tilde{w}_n}}.
\eneq
%
%
By {{\eqref{2020-7-10-nn+1},}} \eqref{2020-7-10-n8} and \eqref{2020-7-10-n500-},
\beq
&&\hspace{-0.6in}A_no(\td u_0)=B^{-1}((I_r, D_n')o({\td u}_0))=B^{-1}({{(I_r,0)o({\tilde u}_0)}}+(0_{r\times r},D_n')o(\tilde{u}_0))\\
&&=B^{-1}(B(\tilde{z}')-D_no(\tilde{w}_n)+(0,D_n')o(\tilde{u}_0))=\tilde{z}'.
\eneq
Thus \eqref{2020-7-10-n++} holds.

Let ${\td u}:=J^{s(n)}(\td u_0)+Mo^{s(n)}(\td u_0):=(x_1', x_2', ...,x_{s(n)}').$
By \eqref{equ5001},
\beq
J^{s(n)}(\td u_0)=(J(a_1),...,J(a_r), J(a_{r+1}),...,J(a_{s(n)}).
\eneq

It remains to show ${\td R}_{s(n)}(\td u)>0$ when $n\ge N.$
By
\eqref{equ5001}, {{and}} the choice of $\xi_n$ and choice of $\dt_n$
(see also \eqref{2017-may-9}),
one has
\beq\label{2020-7-10-n20}
\diag(R_{N_0+1}, R_{N_0+2}, {...,} R_n) (x'_{s(N_0)+1}, x'_{s(N_0)+2}, {...,} x'_{s(n)})^T \in (\Aff_+(S))^{\sum_{k=N_0+1}^n r_k}.
\eneq
Put
$$g'=(J(a_1^{(n)}))+o(a_1^{(n)}), J(a_2^{(n)})+o(a_2^{(n)}),...,J(a_r^{(n)})+o(a_r^{(n)})).$$
Recall that the definition of ${\tilde w}_n$ and $u'$ (see (\ref{2020-7-10-nn+1})) give
\beq
u'=J^r(a^{(n)})+KBo^r(\tilde{z}'')-KD_n{{o^{s(n)}}}(\xi_n)\andeqn \tilde{w}_n:=J^{s(n)}(\xi_n)+K^3{{o^{s(n)}}}({{\xi_n}}).
\eneq
By  the choice of $\dt,$
\eqref{2017-may-9-2}, \eqref{equ5001} (twice),
 one has
\begin{eqnarray*}
|| MK^2o^r(\rho^r(u'))-o^r(\rho^r(g'))||&=&|| Bo^r(\rho^r(M\tilde{z}'))-D_no^{s(n)}{{(\rho^{s(n)}(M\tilde{w}_n))}}-
{{o^r(\rho^r(g'))}}||\\
&\leq&|| Bo^r(z)-D_no^{s(n)}{{(\rho^{s(n)}(M\tilde{w}_n))}}-o^r(\rho^r(g'))||+\dt_0/8\\
&\leq&||Bo^r(y_n)-D_no^{s(n)}{{(\rho^{s(n)}(M\tilde{w}_n))}}-
{{o^r(\rho^r(g'))}}||+\dt_0/4\\
&\le &\|M(B{{\rho^r}}(A_n{{o^{s(n)}}}(\tilde{w}_n))-D_no^{s(n)}(\rho^{s(n)}(\tilde{w})_n))
-{{o^r(\rho^r(g'))}}||+5\dt_0/16\\
&=&||M \rho^r((C_n-D_n)({{o^{s(n)}}}(\tilde{w}_n))-{{o^r(\rho^r(g'))}}) ||+5\dt_0/16\\
&=&|| MK^3o^r(\rho^r((\tilde{d_j^{(n)}})_{r\times 1}))-o^r(\rho^r(g') ) ||+6\dt_0/16\le
\dt_0/2.
\end{eqnarray*}
Therefore, combining \eqref{equ5001}, we have
$\|o(x_i')-o(a_i)\|<\dt_0,$ $i=1,2,...,s(n_0).$
It follows from the claim earlier that (see \eqref{2020-7-10-n4})
\beq\label{2020-7-10-n21}
{\td R}_{N_0}((x_1',x_2',...,x_{s(N_0)}')^T)\ge 3\ep_0/4\,\,({\rm in}\,\, (\Aff_+(S)\setminus \{0\})^{\bar r}.
\eneq
Combing \eqref{2020-7-10-n21} and \eqref{2020-7-10-n20}, one concludes
that
\beq
{\td R}(\td u)>0
\eneq
as desired.
\end{proof}

\begin{lem}\label{ExtBA}
Let $B$ be a non-unital separable simple amenable \CA\, in ${\cal D}^d$  with continuous scale which satisfies the UCT and let
$A=A_1\otimes U,$ where
$A_1\in {\cal M}_1$ and $U$ is an infinite dimensional UHF-algebra.
Suppose
that
\beq
((K_0(B),  T(B), \rho_B), K_1(B))\cong ((K_0(A), T(A), \rho_A), K_1(A))
\eneq
and suppose that $\af\in KL(B,A)$ is an element
which implements (part of) the isomorphism above.

Then,  there exists
a  sequence of approximate multiplicative  \cpc s $\phi_n: B\to A\otimes M_3$ such that
\beq
{{[\{\phi_n\}]=\af.}}
\eneq
\end{lem}

\begin{proof}
Let $\ep>0, \eta>0,$ ${\cal F}\subset B$ be a finite subset.
\Wlog, we may assume that ${\cal F}
\subset B^{\bf 1}.$
Fix a finite subset ${\cal P}_B\subset \underline{K}(B).$

Choose $\dt_1>0$ and {{a}} finite subset ${\cal G}\subset B$ so that ${{[L]|_{{\cal P}_B}}}$ is well defined
for any
${\cal G}$-$\dt_1$-multiplicative \cpc\, $L$ from $B.$
We may assume that $\dt_1<\ep$ and ${\cal F}\cup {\cal H}\subset {\cal G}.$
Since both $A$ and $B$ have continuous scale,
$T(A)$ and $T(B)$ are compact ({{by}} Theorem 5.3 of \cite{eglnp}).

Choose  $b_0\in B_+$ such that $\|b_0\|=1$ and
\beq\label{exBA-n2}
d_\tau(b_0)<\min\{\eta, \dt_1\}/4\rforal \tau\in T(B).
\eneq
Let  $e_0\in B$ be a strictly positive element of {{$B$}} with $\|e_0\|=1$
such that $\tau(e_0)>15/16$ for all $\tau \in T(B).$

Let $G({\cal P}_B)$ be the subgroup generated by ${\cal P}_B,$
$G_{oB}$ be the subgroup generated by ${\cal P}_B\cap K_0(B).$
Let ${\cal P}=\N \cdot [1_{\tilde B}]+{\cal P}_B\cap K_0(B)$ and
$G^{B0}=\Z\cdot [1_{\tilde B}]+G_{oB}.$  We may assume
that $G^{B0}$ is generated by projections $\{p_1,p_2,...,p_l\}$ in $M_{N_0}({\tilde B})$
for some integer $N_0\ge 1.$ We
may
{{also}} assume that $p_1=1_{\tilde B}.$

We may write $G_{oB}=G_{inf}\oplus G_{0,1},$ where $G_{inf}\subset {\rm ker}\rho_B$
and ${\rm ker}\rho_B\cap G_{0,1}=\{0\}.$ Let $k_0$ be an integer such that
$G({\mathcal P}_B)\cap\mbox{K}_i(B,\mathbb Z/k\mathbb Z)=\{0\}$ for any $k\geq k_0$, $i=0,1.$
Keep in mind that $G_{0,1}$ is free and $\rho_B(G_{oB})=\rho_B(G_{0,1}).$

By {{Theorem}} 3.3 and 3.4 of \cite{GLrange} 
there exists a subsequence $\{k(n)\}$  and
 a sequence of asymptotically multiplicative \cpc s
 ${{\Theta_n}}:B\to A\otimes M_{k(n)}$ such that
 \beq
 [\{{{\Theta_n}}\}]=\af.
 \eneq
In what follows we will identity $T(A)$ with $T(B)$ (both are assumed to be compact).

We now use the construction of \ref{C201GLN} and keep the notation used there.
Consider the map $\tilde\rho:G({\mathcal P})\cap\mbox{K}_0(B)\to l^\infty(\mathbb Z)$ defined in Corollary \ref{rho}.
 Let $\pi_B: {\tilde B}\to \C$ be the quotient map with ${\rm ker}\pi_B=B.$
 We may assume that $[p_i]=J([p_i])\cdot [1_{\tilde B}]+o([p_i]),$ where $J([p_i])\in \N$ and
 $\pi_B(p_i)=J([p_i])\cdot 1_{\tilde B},$ and $o([p_i])\in K_0(B).$
 We assume that $[p_1]=[1_{\tilde B}].$  Keep in mind that
 $o([p_1])=0.$
The
linear span of $\{\tilde\rho([p_1]), ..., \tilde\rho([p_l])\}$ over $\mathbb Q$ has finite rank, say {{$r+1$}}.
So,
we may assume that $\{ \tilde\rho([p_1]), \tilde\rho(o([p_2])), ..., \tilde\rho(o([p_{{r+1}}]))\}$ are linearly independent and the $\mathbb Q$-linear span of them  contains ${\tilde\rho}(G({\mathcal P})\cap\mbox{K}_0(B)).$
 Therefore, there is an integer $M$ such that for any $g\in\tilde\rho(G_{oB})$, the element $Mg$
is in the subgroup generated by $\{\tilde\rho(o([p_2]), ..., \tilde\rho(o([p_{{r+1}}])\}$.
{{Moreover, we may further assume that $J([p_j])=J([p_2])$ for
$2\le j\le {{r+1}}.$}}
Write ${{\td \rho(o([p_{i+1}])}}=(x_{i,1}, x_{i,2},...).$ In other words,
${{x_{i,j}={\tilde\rho(o([p_{i+1}]))_j},}}~{{i=1,2,\cdots, r.}}$ {{Let $s(0)=0, s(1)={{2l_k-1}}, s(2)={{(2l_k-1)+(2l_{k+1}-1)}},\cdots,$ $s(j)=\sum_{i=0}^{j-1}{{(2l_{k+i}-1)}}$.}} By Lemma \ref{rho},  {{we have $x_{i,s(j)+1}=m^{k+j}_1\psi^k_{k+j}(o([p_{i+1}]))$ and}}
\beq\label{05-20-2021}
&&\hspace{-0.3in}x_{i, {{s(j)+2l-2}}}= m^{k+j}_l(\psi^k_{k+j}(o([p_{i+1}]))~~~\mbox{and}~~~~x_{i, {{s(j)+2l-1}}}= -m^{k+j}_l(\psi^k_{k+j}(o([p_{i+1}])),
\eneq {{for $l=1,2, \cdots,l_{k+j}$.}} Let  $z_i=\rho_A(\af(o([p_{i+1}]))\in\mathbb{D}$, where $\mathbb{D}=\rho_A(K_0({\tilde A}))$ in
$\text{Aff}(S_{[1]}(K_0({{\tilde{A}}}))$.
Therefore $K_0({\tilde A})=\Z\cdot [1_{\tilde A}]+K_0(A).$
Denote by $\D_0=\rho_A(K_0(A)).$  Since $A\in {\cal D}^d,$ $\D_0$ is dense in $\R\D_0\subset
 \text{Aff}(\mathrm{S}_{[1]}({K_0}({{\tilde A}})))$. 
 Keep in mind that $J([p_i])=J(\af([p_j]))$ as an integer.
Note that   $\af([p_i])=J([p_i])\cdot [1_{\tilde A}]+\af(o([p_i]))$
and $\rho_A(\af([p_i]))=J([p_i])+\rho_A(\af(o([p_i])),$ $i=1,2,...,l.$

Let  $\{S_j\}$ be the sequence of \SCA s in ${\cal C}_0$ in the construction \ref{C201GLN}.
Fix $k\ge 1.$
Let {{$e_1^{k+j}$,}} $e_i^{{{k+j}}}, e_{i,\pm}^{{{k+j}}}\in K_0({\tilde S_{k+j}}), $ 
 ${{i=2,..., l_{k+j}}}$  and
$R_{{{k+j}}}$ be  the $r_{k+j}\times {{(2l_{k+j}-1)}}$ matrix as described in \ref{Constr-Sn}.
Put  {{$\af([ \imath_n\circ {\tilde h_{k+j, n}}(e_1^{{k+j}})])=g^{(n)}_{{s(j)+1}}(=1_{\tilde A})$ and }}
\beq\label{June14n2-1}
\hspace{-0.2in}\af([ \imath_n\circ {\tilde h_{k+j, n}}(e_{i,+}^{{k+j}})])=g^{(n)}_{{s(j)+2i-2}}{{\in K_0({\tilde A})}},\,\,
\af([\imath_n\circ {{{\tilde h_{k+j,n}}}}(e_{i,-}^{{k+j}})])=g^{(n)}_{{s(j)+2i-1}}{{\in K_0({\tilde A})}},
 \eneq
 $j=0,1,\cdots, n-k$ and {{$i=2,3,...,l_k$}}.
Let   $a_j^{(n)}=\rho_B(g_j^{(n)})$ and $a_j^{(o,n)}=\rho_B(o(g_j^{(n)})),$ $j=1,2,...,$ 
$n=1,2,....$ Note that $a_j^{(n)}\in\mathbb{D}^+\backslash\{0\}.$ It follows from Lemma \ref{traces} that (recall that $\af$ is an order isomorphism) for $j=0,1,\cdots$, {{$a^{(n)}_{s(j)+1}=a_{s(j)+1}=1\in \Aff T(B)$}} and {{for  $i=2,{{...,}} l_{k+j},$}}
$$\lim_{n\to\infty}a^{(n)}_{{s(j)+2i-2}}=a_{{s(j)+2i-2}}=
{\lim_{n\to\infty}}{{\rho_A(g^{(n)}_{{s(j)+2i-2}})}}>0~~\mbox{and}$$
 $$\lim_{n\to\infty}a^{(n)}_{{s(j)+2i-1}}=a_{{s(j)+2i-1}}=
{\lim_{n\to\infty}}\rho_A(g^{(n)}_{{s(j)+2i-1}})>0$$ {{uniformly.}}
{{\Wlog, we may assume $J(a_j^{(n)})=J(a_j^{(n+1)})$ as an integer for all $n\in \N.$}}
Since for $i\ge 2,$  $e^{k+j}_i\in K_0(S_{k+j}),$ we have $J(g^{(n)}_{{s(j)+2i-2}})-J(g^{(n)}_{{s(j)+2i-1}})=0,\,i\ge 2.$
It follows that
\beq\label{1207-21-zero}
J(a_{{s(j)+2i-2}})-J(a_{{s(j)+2i-1}})=0,\, 2\le i\le l_j.
\eneq

 {{Also note that $o(a_{s(j)+1})=0$.}} Moreover, by \ref{rho-uniform},
$\sum_{j=1}^n{{x_{i,j}}}a_j^{(o,n)}\to z_i\in \D_0$ uniformly.  Furthermore, by \ref{Constr-Sn} {{(with $l_i$ replaced by $l_{k+i-1}$)}},
$\{x_{i,j}\},$ $g_j^n,$ and $R_n$ satisfy the condition  of Lemma \ref{solveeq}
with $K_0(A)$ in place of $G$, {{$\rho_A(K_0(A))$ in place of $\D_0,$}}
 {{ $T(A)$ in place of  $S$ and ${{2l_{k+j-1}-1}}$ in place of $l_j$.}}
%
So, Lemma \ref{solveeq} applies. Fix  $\dt,$ $N$ and $K$ obtained from Lemma \ref{solveeq}.

Note that $A\in\mathcal D^d$, and hence the strict order on the projections of ${\tilde A}$ is determined by traces of ${\tilde A}.$
{{Let $\af(G_{0,1})$ be  {{the free abelian}} group generated by $\{d_1,d_2, \cdots d_{l'}\}$.}} Write
\beq\label{05-19-2021-1}
\hspace{-0.3in}\af([p_i])=J([p_i])\cdot [1_{\tilde A}]+\sum_{j=1}^{{l'}}m_{i,j}d_j+s_i\andeqn o(\af([p_i]))=\sum_{j=1}^{{l'}}m_{i,j}d_j+s_i
\eneq
 where
$m_{i,j}\in \Z,$  $d_j\in \af(G_{01})$ and
$s_i\in \af(G_{inf}),$
 $j=1,2,...,{{l'}}$ and {{$i=2,3...,r+1$. (Note that $\af([p_1])=\af([1_{\tilde B}])=[1_{\tilde A}]$.)}}
Denote {{${\cal P}_1=\af({\cal P}_B)$ and}} $G_{A}=\af(G({\cal P}_B)).$
Thus $\af(G_{inf})\subset {\rm ker}\rho_A.$

Fix an integer $k(n)\ge 1.$
Applying Corollary {{5.4 of \cite{GLrange}}} 
 to  $M_{k(n)}(A)$ with any finite subset $\mathcal{G}$, any $\epsilon>0$ and any
$0<r_0<\dt/{{(2k(n))}}<1$, one has a $\mathcal{G}$-$\epsilon/k(n)^2$-multiplicative map $L:\mbox{M}_{k(n)}(A)\to\mbox{M}_{k(n)}(A)$ with the following properties:
\begin{enumerate}
\item $[L]|_{{\mathcal P}_1}$ and $[L]|_{{G_A}}$ are well defined;
\item $[L]$ induces the identity maps on
      $\af(G_{inf}),$
      $G_A\cap\mbox{K}_1(A)$,
      $G_A\cap\mbox{K}_0(A,\mathbb Z/k\mathbb Z)$ and $G_A\cap\mbox{K}_1(A, \mathbb Z/k\mathbb Z)$
      for the $k$ with $G_A\cap\mbox{K}_i(A, \mathbb Z/k\mathbb Z)\neq\{0\}$,
      $i=0,1$;
\item $|\tau\circ[L](g)|\leq r_0|\tau(g)
|$ for all $g\in
      G_A\cap\mbox{K}_0(A)$ and $\tau\in T(A)$;
\item There exist elements $\{f_i\}\subset{K_0}(A)$
      such that for  $i=1,...,{{l'}},$
      \beq\label{05-18-2021-3}
      {{d_i-[L](d_i)}}=MK^3(k_0+1)!f_i.\eneq
      \item $d_\tau(e_0)<r_0\rforal \tau\in T(A),$
\end{enumerate}
where $e_0$ is a strictly positive element of $\overline{L(M_{k(n)}(A))M_{k(n)}(A)L(M_{k(n)}(A))}.$
By the choice of $r_0,$ replacing $L$ by ${\rm Ad}\in U\circ L$ for a
unitary $U\in M_{k(n)}(A)^\sim,$    we may assume that $L(M_{k(n)}(A))\subset A.$
%
%
Note {{that,}} by (2) above, $(\af-[L]\circ \af)(s_i)=0.$  Then we have
\beq\label{extBA-nn1}
&&\hspace{-0.8in}{{\alpha([p_i]-J([p_i])[1_{\tilde A}])-([L](\af([p_i])-J([p_i])[1_{\tilde A}]))}}\\
  &=&
 {{(\sum m_{i,j}d_j)-[L](\sum m_{i,j}d_j)}}\\
 &=& MK^3(k_0+1)!(\sum m_{i,j}f_j)\\\label{exBA-nn2}
 &=& MK^3(k_0+1)!f_i',
\eneq
 where $f_i'=\sum m_{i,j}f_j,$
for {{$i=2,3,...,r+1.$}}
Note that $[p_1]=[1_{\tilde B}]=J([p_1]).$ 
Since $T(A)$ is compact, we may choose $r_0$ such that {{(3) above implies that}}
\beq
{{|\tau([L](o(\af([p_i]))))|}}<\dt/2\rforal \tau\in T(A).
\eneq
 Define {{$\bt([p_1])=[1_{\tilde A}]$, and}}\beq\label{July-31-2021-4}
 {{\bt([p_i])=J([p_i]) [1_{\tilde A}]+K^3(k_0+1)!f_i'
 ~~\mbox{and}~~ o(\bt([p_i])=K^3(k_0+1)!f_i',}}\eneq
 for {{$i=2,3,...,r+1.$}}
Let $\tilde{z_i}'=\beta([p_{{i+1}}])$,  $z'_i=\rho_A(\tilde{z_i}')\in\Aff(S_{[1]}(K_0({{\tilde{A}}})))$ and $o(z_i')=\rho_A(o(\bt([p_{{i+1}}]))=
{{K^3}}(k_0+1)!\rho_A({{f'_{i+1}}}),$ $i=1,2,...,{{r}}.$  In particular, $o(\bt)[p_1])=0.$
Then, {{by (\ref{05-19-2021-1}) and (\ref{05-18-2021-3}), }} we have:
$$\begin{array}{lll}
    ||{Mo^r(z')-o^r(z)}||_\infty & = & \max_i\{||\rho_A(\alpha([p_i]-J([p_i])[1_{\tilde A}] )-[L](\af([p_i])-J([p_i])[1_{\tilde A}]))\|\}\\
     & = & \max_i\{\sup_{\tau\in T(B)}\{\tau\circ[L](o(\af([p_i]))\}
     < \delta/2,
  \end{array}$$
  where $z=(z_1,z_2,...,z_r),$ $z_i=J(\af([p_{{i+1}}]))+o(z_i),$ {{(note that, {{here}} we do not {{use}} $p_1=1_{\tilde A}$ as $o(p_1)=0$)}},
  $o(z_i)\in \rho_A(K_0(A)),$ ${{i=1,2...,r,}}$
  and $z'=(z_1',z_2',...,z_r').$
By Lemma \ref{solveeq},  for sufficiently large $n,$ one obtains
$$\tilde{u}^A=(J(a_1)+Mo(u_{0,1}^A),J(a_2)+Mo(u_{0,2}^A),...,J(a_{s(n)})+Mo(u_{0,s(n)}^A))
\in K_0({\tilde{A}})^{{s(n)}}$$
such that (recall that $\af$ is an order isomorphism)
\beq\label{June-16-n1}
\hspace{-0.24in}
\sum_{j=1}^{s(n)} {{x_{i,j}}}o(u_{0,j}^A)={{o(\tilde{z}'_{i})}}.
\eneq
As importantly, we have
\beq\label{June16-3}
{{{\bar R}_n({\tilde u}^A)>0,}}
\eneq
where
${\tilde u}^A=J^{s(n)}({\tilde u}_0^A)+Mo({\tilde u}_0^A).$

 Let  {{$\kappa_0^{k+j}: K_0((S_{k+j})^\sim)\to K_0({\tilde A})$}} be  defined by, for each $j,$
\beq\label{1207-21-111}
&&\hspace{-0.5in}e_{s(j)+1}^{k+j}\mapsto J(a_{s(j)+1})=[1_{\td A}]\andeqn, {\rm for}\,\, 2\le i\le l_{k+j},\\
&&\hspace{-0.5in}e_i^{k+j}\mapsto  J(a_{s(j)+2i-1})+Mo(u^A_{0,{{s(j)+2i-2}}})-J(a_{s(j)})+Mo(u^A_{0,
{{s(j)+2i-1}}}),\,\, {{0}}\leq j\leq n.
\eneq
By  \eqref{1207-21-zero}, for $2\le i\le l_{k+j},$
\beq
e_i^{k+j}\mapsto  M(o(u^A_{0,{{s(j)+2i-2}}})-o(u^A_{0,{{s(j)+2i-1}}})),\,\, {{0}}\leq j\leq n.
\eneq
%
 Then, by  \ref{Constr-Sn} and by {{\eqref{June16-3}}},
 it is strictly positive.
 By {{\eqref{05-20-2021}}} and \eqref{June-16-n1}, for $2\le i\le l_j,$
\beq\nonumber
&&\hspace{-0.6in}\kappa_0^{k+j}([\psi_{k+j}^{k}](o([p_i])))=\sum_{l=1}^{l_{k+j}}
{{m_l^{k+j}}}([\psi_{k+j}^k](o([p_i])))\kappa_0^{k+j}(e_l^{k+j})\\\nonumber
&&=\sum_{l=1}^{l_{k+j}}m_i^{k+j}([\psi_{k+j}^k](o([p_i])))
M({{o((u^A_{0,{{s(j)+2l-2}}})-o(u^A_{0,{{s(j)+2l-1}}})))}}\\
&=&
\sum_{l=1}^{{2l_{k+j}-1}} {{x_{i,s(j)+l}}}{Mo({\tilde{u}^A_{0,{{s(j)}}+l}}).}\\\nonumber
\eneq
{{Put $D_0=S_{k}\oplus S_{k+1}\oplus\cdots\oplus S_{k+n-1}$
and $D=S_{k}^\sim \oplus S_{k+1}^\sim \oplus\cdots\oplus S_{k+n-1}^\sim.$
Define $\kappa_0': K_0(D)\to K_0(\td A)$  by $\kappa_0'|_{K_0(S_{k+j}^\sim)}=\kappa_0^{k+j},$
$0\le j\le n-1.$ As mentioned above, $\kappa_0'$ is strictly positive.
Moreover, by \eqref{June-16-n1},
\beq\label{1206-21-n1}
\kappa_0^{k+j}(\psi^k([p_1]))
=[1_{\td A}].
\eneq
It follows Theorem 5.7 of  \cite{GLrange}, there exists a \hm\, $h': D\to M_n(\td A)$ (for
such that $(h')_{*0}=\kappa_0'$ and $h'|_{D_0}\to M_K(A).$
Let $\psi^k: \td B\to D$ by
$$\psi^k(b)=({\psi^k_{k}}^\sim(b)\oplus {\psi^k_{k+1}}^\sim(b){\oplus}\cdots\oplus {\psi^k_{k+n-1}}^\sim(b))
\rforal b\in \td B.$$}}
It follows  that (see also \eqref{June-16-n1})
\beq
\kappa_0'(\psi^k(o([p_i])))&=&(x_{ij})_{s(n)\times 1} \td u^A=
Mo(\bt([p_i])),\,\,\, \hspace{0.4in} {{i=2,...,r+1}}.
\eneq
This also implies that $\kappa_0'(\psi^k(J([p_i])))=J(\af([p_i])),$ {{$i=1,2,...,r+1.$}}
%

%
{{By \eqref{June-16-n1}, one has, keeping the notation in the construction at the beginning of this section,
\beq\label{June-16-nn2}
J([p_i])+\kappa_0'(o([\psi^{{{k}}}_{k}(p_i)]), o([{{\psi^k_{k+1}}}(p_i)]), ...,o([{{\psi^k_{n+k-1}}}(p_i)]))=J([p_i])+M\beta(o([p_i])),
\eneq
$i=1,...,l.$}}

Now, define $h''^\sim: \td B\to  {{D \to}}  M_n(\td A)$
by $h''^\sim=h'\circ ({\psi^k_{k}}^\sim\oplus {\psi^k}^\sim_{k+1}{\oplus}\cdots\oplus {\psi^k}^\sim_{k+n-1}).$
Then $h''^\sim$ is $\mathcal F$-$\delta$-multiplicative. Moreover $h''=h''^\sim|_B$ maps
$B$ into $M_n(A).$  However, by \eqref{1207-8-im-1},
\beq
d_\tau(h''(e_0))<{{d_{\tau}(h'({\psi_k^k}^\sim(e_0)))+1<2}}\rforal \tau\in T(A).
\eneq
Since $M_n(A)$ has strict comparison and has stable rank one, we may well assume that $h''$ maps
$B$ into {{$M_2(A).$}}

For any $x\in {\rm ker}\td \rho,$  by  Corollary \ref{rho}, $x\in {\rm ker}\rho_B.$  Since $\af$ preserves {{the}}
order on $K_0(\td A),$
$\af(x)\in \af(G_{inf}).$    By (2) above, $\af(x)-[L\circ {{\Theta_n}}](x)=0.$
Define  $\Phi: B\to M_{{3}}(A)$ by, for some sufficiently large $n,$
$\Phi(a)=L\circ {{\Theta_n}}(a)\oplus h''(a)$ for all $a\in A.$
Define $\Phi^\sim: {\tilde B}\to {\widetilde{M_{{3}}(A)}}.$
By (2) above, \eqref{June-16-nn2},{{\eqref{July-31-2021-4}}} \eqref{extBA-nn1} and \eqref{extBA-nn1},
\beq\label{ExtBA-n1}
[\Phi^{\sim}]([p_i])-{{3}}J([p_i])&=&[\Phi](o([p_i])=
[L\circ \af](o[p_i])+[h''](o[p_i])\\
&=&[L\circ \af](o[p_i])+M\bt(o([p_i])\\
 &=&[L\circ \af](o[p_i])+ MK^3(k_0+1)!f_i'\\
 &=& [L\circ \af](o[p_i])
 +\alpha([p_i]-J([p_i])[1_{\tilde A}])\\
 &&\hspace{0.2in}-([L](\af([p_i])-J([p_i])[1_{\tilde A}]))\\
 &=& \alpha([p_i]-J([p_i])[1_{\tilde A}])=\af(o([p_i])),\,\,\,i=1,2,...,l.
\eneq
Thus,
\beq
[\Phi]|_{G^{B0}}=\af|_{G^{B0}}.
\eneq
Since $K_0(D)$ is free and $K_1(D)=\{0\},$ by (2) again,
\beq\label{exBA-n5}
&&\hspace{-0.7in}[\Phi]|_{G_B\cap K_1(B)}=[L\circ \af]|_{G_B\cap K_1(A)}=\af|_{G_A\cap K_1(A)},\\
&&\hspace{-0.7in}{[}\Phi{]}|_{G_B\cap K_0(B)/jK_0(B)}=[L\circ \af]|_{G_B\cap K_0(B)/jK_0(B)}=\af|_{G_B\cap K_0(B)/jK_0(B)},\,\,1\le j\le k_0\andeqn\\
&&\hspace{-0.7in}{[}\Phi{]}|_{G_B\cap K_1(B)/jK_1(B)}=[L\circ \af]|_{G_B\cap K_1(B)/jK_1(B)}=\af|_{G_B\cap K_1(B)/jK_1(B)},
\,\,1\le j\le k_0.
\eneq
It follows that
\beq\label{exBA-n6}
[\Phi]|_{\cal P}=\af|_{\cal P}.
\eneq
{{The lemma}} follows.

\end{proof}



\begin{lem}\label{11Ext2n}
Let $A$   be a simple
separable  amenable {{\CA}} in ${\cal D}^d$
 with continuous scale
which {{satisfies}} the UCT.  Suppose also {{that}} $B_1\in{\cal M}_1$ {{is}} as {{(the algebra $A$)}} described in Theorem \ref{Text1}, which is the form of  (3) of Remark 4.32 of \cite{GLrange}  and $B=B_1\otimes U,$
where $U$ is an infinite dimensional UHF-algebra.
Suppose that
$$
\Gamma: ((K_0(A), T(A), \rho_A), K_1(A))\to ((K_0(B), T(B), \rho_B), K_1(B))
$$
is an isomorphism and
 that $\kappa\in KL(A,B)$ carries $\Gamma|_{K_i(A)},$ $i=0,1,$
$\Gamma$ also gives an affine {{homeomorphism}} $\kappa_T: T(B)\to T(A).$
Suppose also there is
a continuous \hm\, ${{\kappa_{cu}}}: U({\tilde A})/CU({\tilde A})\to U({\widetilde{B}})/CU({\widetilde{B}})$ such that $(\kappa, \kappa_T, {{\kappa_{cu}}})$ is compatible {{(see Definition \ref{Dcompatible})}}.

Then there exists a  sequence of approximate multiplicative  \cpc s $\phi_n: A\to B$ such that
\beq\label{11Ext2n-1}
&&[\{\phi_n\}]=\kappa,\\\label{11Ext2n-2}
&&\lim_{n\to\infty}\sup \{|\tau\circ \phi_n(a)-\kappa_T(\tau)(a)|\}=0\rforal a\in A_{s.a.}\tand\\
&&\lim_{n\to\infty}{\rm dist}({{\kappa_{cu}}}(z), \phi_n^{\dag}(z))=0\rforal z\in U({\tilde A})/CU({\tilde A}).
\eneq
\end{lem}

\begin{proof}
Denote by ${{\Pi_{cu}^{\td A}}}: U({\tilde A})/CU({\tilde A})\to K_1(A)$ the quotient map and fix
a splitting map $J_{cu}^A: K_1(A)\to U({\tilde A})/CU({\tilde A}).$
Since $(\kappa, \kappa_T, {{\kappa_{cu}}})$ is compatible, it suffices to show
that there are $\{\phi_n\}$ which satisfies \eqref{11Ext2n-1} and \eqref{11Ext2n-2} and
\beq
\lim_{n\to\infty}{\rm dist}({{\kappa_{cu}}}(J_{cu}^A(\zeta)), \phi_n^{\dag}(J_{cu}^A(\zeta)))=0\rforal \zeta\in K_1(A).
\eneq

It follows from 
{{Lemma}} \ref{ExtBA}
and {{Lemma}} \ref{ExtTBA}  that there exists a sequence $\{\phi_n\}$  which  satisfies \eqref{11Ext2n-1} and \eqref{11Ext2n-2}.
Let $G_1\subset K_1(A)$ be a finitely generated subgroup.
{{Then, for all sufficiently large $n,$ $\Pi_{cu}^{\td B}\circ \phi_n^\dag\circ J_{cu}^{\td A}|_{G_1}=[\phi_n]|_{G_1}$}}.
Note that, since $(\kappa, \kappa_T, {{\kappa_{cu}}})$ is compatible,
for all sufficiently large $n,$
\beq
{{\Pi_{cu}^{\td B}(\kappa_{cu}\circ J_{cu}^{\td A}|_{G_1}-\phi_n^\dag\circ J_{cu}^{\td A}|_{G_1})=0.}}
\eneq

Put
\vspace{-0.15in}\beq\label{11ExtTn2-1}
\lambda_0={{(\kappa_{cu}\circ J_{cu}^{\td A}-\phi_n^\dag\circ J_{cu}^{\td A})|_{G_1}.}}
\eneq
%
%
%
{{Then}} $\lambda_0$ maps from  {{$G_1$}}
to $\Aff(T({\tilde B}))/\overline{\rho_{B}(K_1({\tilde B}))}.$
However, $\Aff(T({\tilde B}))/\overline{\rho_{B}(K_1({\tilde B}))}$ is divisible.
{{Note that $\kappa_1: K_1(A) \to K_1(B)$ is isomorphism, therefore}} there is a \hm\, $\lambda_1: K_1(B)\to \Aff(T({\tilde B}))/\overline{\rho_{B}(K_1({\tilde B}))}$
such that
\beq
({{\lambda_1\circ \kappa_1}})|_{G_1}
=\lambda_0.
\eneq
Now {{define}} $\Lambda: U({\tilde B})/CU({\tilde B})\to U({\tilde B})/CU({\tilde B})$  as follows{{:}}
\beq
\Lambda|_{\Aff(T({\tilde B}))/\overline{\rho_{B}(K_1({\tilde B}))}}={\rm id}_{\Aff(T({\tilde B}))/\overline{\rho_{B}(K_1({\tilde B}))}},\\
\Lambda|_{J_{cu}^B(K_1(B))}=\lambda_1\circ {{\Pi_{cu}^{\td B}}}+({\rm id}_{B})^{\dag}.
\eneq

Note that $([{\rm id}_{B}], ({\rm id}_{B})_T, \Lambda)$ is compatible.
It follows from {{Theorem}} \ref{Text1}
that there exists a \hm\, $\psi_n: B\to B$
such that
\beq
[\psi_n]=[{\rm id}_{B}], \,\,\, (\psi_n)_T=({\rm id}_{B})_T\andeqn \psi_n^{\dag}=\Lambda.
\eneq
Now let $\Phi_n=\psi_n\circ \phi_n.$
Then, for $z\in J_{cu}^A(G_1),$   by \eqref{11ExtTn2-1},
\beq
\Phi_n^{\dag}(z)=\psi_n^{\dag}\circ \phi_n^{\dag}(z)&=&\lambda_1\circ {{\Pi_{cu}^{\td B}}}\circ \phi_n^{\dag}(z)+\phi_n^{\dag}(z)\\
{{=\ld_1\circ\kappa_1\circ\Pi_{cu}^{\td A}(z)+\phi_n^{\dag}(z)}}&=&\lambda_0\circ {{\Pi_{cu}^{\td A}}}(z)+\phi_n^{\dag}(z)={{\kappa_{cu}}}(z).
\eneq
The lemma  follows immediately from the construction of $\Phi_n.$


\end{proof}

\begin{lem}\label{11Ext2n2}
Suppose that  $A$ and $B$ satisfy  {{exactly}} the same conditions in Lemma \ref{11Ext2n},
and suppose that $\kappa, {{\kappa_T}}$ and $\kappa_{cu}$ are  as stated  in \ref{11Ext2n}.
%
%
%
%
Then there exists a  \hm\, $\phi: A\to B$ such that
\beq\label{11Ext3n-1}
[\phi]=\kappa,\,\,\, \phi_T=\kappa_T\andeqn \phi^{\dag}={{\kappa_{cu}}}.
\eneq
\end{lem}

\begin{proof}
The proof is exactly the same as that of \ref{Text1} but applying \ref{11Ext2n} instead of \ref{Lextcu}.
\end{proof}


%
%
%
%
%

%
%
\section{The Isomorphism Theorem for \CA s in ${\cal D}^d$}


\begin{lem}\label{MUN2}
Let $C\in D^d$ be a non-unital separable simple amenable \CA\, with continuous scale which satisfies the UCT
and let $A\in {\cal D}$ {{be}} with continuous scale.
Suppose that
$\phi_1, \phi_2: C\to A$ are two  monomorphisms which maps strictly positive elements
to strictly positive elements.
Suppose also that
\beq\label{MUN2-1}
[\phi_1]=[\phi_2]\,\,\,{\text in}\,\,\, KL(C,A),\\
(\phi_1)_T=(\phi_2)_T\tand \phi_1^{\dag}=\phi_2^{\dag}.
\eneq
Then $\phi_1$ and $\phi_2$ are approximately unitarily equivalent.
\end{lem}

\begin{proof}
Note, by {{A.10 of \cite{GLII},}} 
{{that}}  both $A$ and $B$ are ${\cal Z}$-stable.
{{Also $\phi_1$  and $\phi_2$ are  given and full. So there is $T: A_+\setminus \{0\}\to \N\times \R_+\setminus \{0\}$
such that $\phi_1$ and $\phi_2$ are exactly $T$-$(A_+\setminus \{0\})$-full (see 5.7 of \cite{eglnp}). }}
Then the theorem follows from Theorem 5.3 of \cite{GLII} (see also 5.2 and {{Proposition 5.5}}  of \cite{GLII}).
\end{proof}

\begin{thm}\label{Misothm}
Let $A_1, B_1\in {\cal D}$  be two separable amenable simple \CA s with continuous {{scales}}
which satisfy the UCT.  Let $A=A_1\otimes U_1$ and $B=B_1\otimes U_2,$
where $U_1$ and $U_2$ are infinite dimensional UHF-algebras.
Then $A\cong B$ if and only if  there is an isomorphism
\beq
\Gamma: ((K_0(A),T(A), r_A), K_1(A))\cong ((K_0(B), T(B), r_B), K_1(B)).
\eneq
Moreover,
let $\kappa_i: K_i(A)\to K_i(B)$ be an isomorphism as abelian groups ($i=0,1$) and
let $\kappa_T: T(B)\to T(A)$ {{be}} an affine homeomorphism.
Suppose that $\kappa\in KL(A,B)$ which gives $\kappa_i$ and
$\kappa_{cu}: U({\tilde A})/CU({\tilde A})\to U({\tilde B})/CU({\tilde B})$ is
{{an
isomorphism}} so that
$(\kappa, \kappa_T, \kappa_{cu})$ is compatible {{(see Definition \ref{Dcompatible})}}.
Then there is an isomorphism $\phi: A\to B$ such that
\beq\label{Misothm-1-1}
[\phi]=\kappa\,\,\,
\phi_T=\kappa_T\andeqn \phi^{\dag}=\kappa_{cu}.
\eneq
\end{thm}

\begin{proof}
{{By (3) of Remark 4.32 of \cite{GLrange},}} there is a non-unital simple \CA\, $B_0\in {\cal M}_1$ {{
which has the form}} {{ as $A$ in  Theorem \ref{Text1} 
 }} 
such that  there is an isomorphism
\beq\label{Misothm-1-0}
\Gamma': ((K_0(B_1), T(B_1), \rho_{B_1}), K_1(B_1))\cong ((K_0(B_0),T(B_0), \rho_{B_0}), K_1(B_0)).
\eneq
It follows that there is an isomorphism
\beq\label{Misothm-1}
\hspace{-0.2in}\Gamma'': ((K_0(B), T(B), \rho_{B}), K_1(B))\cong ((K_0(B_0\otimes U_2),T(B_0\otimes U_2), \rho_{B_0\otimes U_2}), K_1(B_0\otimes U_2)).
\eneq
Thus, if we can show the theorem holds for the case $B=B_0\otimes U_2,$ then the general
case follows.  Therefore, \wilog, we may assume that $B_1\in {\cal M}_1$ {{be as in Theorem \ref{Text1}.}}

We also note that the ``only if " part  of the first part of the theorem is obvious.
So we will prove the ``if" part.
Suppose that $\Gamma$ exists. Then, by the UCT, $\Gamma$ induces a compatible
pair $(\kappa, \kappa_T)$ {{for some $\kappa\in KL(A, B).$}}
Choose any $\kappa_{cu}: {{U({\td A})/CU({\td A}) \to U({\td B})/CU({\td B})}}$ so that $(\kappa, \kappa_T, \kappa_{cu})$
is compatible.
Note that there is always at least {{one}} {{$\kappa_{cu}$ such that}}  $\kappa_{cu}|_{J_{cu}^A(K_1(A))}={{{J_{cu}^B\circ \kappa|_{K_1(A)}\circ \Pi^A_{cu}}|_{J_{cu}^A(K_1(A))}}},$
where $\Pi^A_{cu}: U({\tilde A})/CU({\tilde A})\to K_1(A)$ is the quotient map
and $\kappa_{cu}|_{\Aff(T(\td A))/\overline{\rho_{\td A}(K_0(\td A))}}$ {{is}}  induced by $\kappa_T.$
Therefore, to prove the theorem, it  suffices to show the case
that $(\kappa, \kappa_T)$ is given by $\Gamma.$
{{Hence}}  it remains  to show that
there is  an isomorphism $\phi: A\to B$ such that
\eqref{Misothm-1-1} holds. We will use the Elliott intertwining argument.

Let $\{{\cal F}_{a,n}\}$ be an increasing sequence of finite subsets of $A$ such that
$\cup_{n=1}^{\infty}{\cal F}_{a,n}$ is dense in $A,$ let $\{{\cal F}_{b,n}\}$ be an increasing
sequence of finite subsets of $B$ such that $\cup_{n=1}^{\infty} {\cal F}_{b,n}$ is dense in
$B.$
Let $\{\ep_n\}$ be a sequence of  decreasing positive numbers  such that
$\sum_{n=1}^{\infty}\ep_n<1.$

Let $e_a\in A$ and $e_b\in B$  be  strictly positive elements of $A$  and $B,$ respectively, with $\|e_a\|=1$
and with $\|e_b\|=1.$
Note that $d_\tau(e_a)=1$ for all $\tau\in T(A)$ and $d_\tau(e_b)=1$ for
all $\tau\in T(B).$

It follows from \ref{11Ext2n2} that there is a \hm\, $\phi_1: A\to B$ such that
\beq\label{Misothm-5}
[\phi_1]=\kappa,\,\,\, (\phi_1)_T=\kappa_T\andeqn \phi_1^{\dag}=\kappa_{cu}.
\eneq
Note that $d_\tau(\phi_1(e_a))=1$ {{($\tau\in T(B)$).}}  {{Then
$\phi_1(e_a)\sim e_b.$ Since $B$ has almost stable rank one, by {{Theorem 1.2 of \cite{Rlz},}}
there is a partial isometry $V\in B^{**}$
such that $V^*cV\in B$ for $c\in \overline{\phi_1(e_a)B\phi_1(e_a)},$ $VV^*c=cVV^*=c$
for $c\in \overline{\phi_1(e_a)B\phi_1(e_a)},$
and $V^*\phi_1(e_a)V$
is a strictly positive element of $B.$ {{Define $\Phi_1: A\to B$ by
$\Phi_1(a)=V^*\phi_1(a)V$($a\in A$).}}
{{Put $z=V^*\phi_1(e_a).$ By Theorem 5 of \cite{Pedjot87}, since $B$ has almost stable rank one,
for each $n,$ there is a unitary $U_n\in \td B$ such that
$U_nf_{1/n}(\phi_1(e_a))=V^*f_{1/n}(\phi_1(e_a)).$ It follows
that $U_n\phi_1(a)U_n^*\to  \Phi_1(a)$ for all $a\in A.$ Thus
$[\Phi_1]=[\phi_1],$ $(\Phi_1)_T=(\phi_1)_T$ and $\Phi_1^\dag=\phi_1^\dag.$}}
Replacing $\phi_1$ 
 by  {{$\Phi_1,$}}
we may assume that}}
$\phi_1$ maps $e_a$ to a strictly positive element of $B.$
It follows from \ref{Text1}
that there is a \hm\,
$\psi_1': B\to A$ such that
\beq\label{Misothm-6}
[\psi_1']=\kappa^{-1},\,\, (\psi_1')_T=\kappa_T^{-1}\andeqn (\psi_1')^{\dag}=
 ( \phi_1^{\dag})^{-1}.
\eneq
Thus
\vspace{-0.1in}\beq
[\psi_1'\circ \phi_1]=[{\rm id}_A],\,\,\, (\psi_1'\circ \phi_1)_T={\rm id}_{T(A)}\andeqn
(\psi_1'\circ \phi_1)^{\dag}={\rm id}_{U({\tilde A})/CU({\tilde A})}.
\eneq
It follows from  {{Lemma \ref{MUN2}}}
that there exists a unitary $u_{1,a}\in {\tilde A}$ such that
\beq\label{Misothm-7}
{\rm Ad}\, u_{1,a}\circ \psi_1'\circ \phi_1\approx_{\ep_1} {\rm id}_A\,\,\,{\rm on}\,\,\, {\cal F}_{a,1}.
\eneq
Put $\psi_1={\rm Ad}\, u_{1,a}\circ \psi_1'.$
Then we obtain the following diagram
\vspace{-0.12in} \begin{displaymath}
\xymatrix{
A \ar[r]^{{\rm id}_A} \ar[d]_{\phi_1} & A\\
B \ar[ur]_{\psi_1}
}
\end{displaymath}
which is approximately commutative on the subset ${\cal F}_{a,1}$ within $\ep_1.$

By applying \ref{11Ext2n2},  there exists a \hm\, $\phi_2': A\to B$ such that
\beq\label{Misothm-8}
[\phi_2']=\kappa,\,\,\, (\phi_2')_T=\kappa_T\andeqn  (\phi_2')^{\dag}=
(\psi_1^{\dag})^{-1}=\kappa_{cu}.
\eneq
Then,
\vspace{-0.1in}\beq
[\phi_2'\circ \psi_1]=[{\rm id}_B],\,\,\, (\phi_2'\circ \psi_1)_T\andeqn (\phi_2'\circ \psi_1)^{\dag}={\rm id}_{U({\tilde B})/CU({\tilde B})}.
\eneq
It follows from  Lemma \ref{MUN2} again
 that there exists a unitary $u_{2,b}\in {\tilde B}$ such that
\beq\label{Misothm-7+}
{\rm Ad}\, u_{2,b}\circ {{\phi_2'}}\circ \psi_1\approx_{\ep_2} {\rm id}_B\,\,\,{\rm on}\,\,\, {\cal F}_{b,2}\cup \phi_1({\cal F}_{a,1}).
\eneq
Put $\phi_2={\rm Ad}\, u_{2,b}\circ {{\phi_2'}}.$ Then we obtain the following diagram:
 \begin{displaymath}
\xymatrix{
A \ar[r]^{{\rm id}_A} \ar[d]_{\phi_1} & A \ar[d]^{\phi_2}\\
B\ar[ur]_{\psi_1}\ar[r]_{{\rm id}_B} & B
}
\end{displaymath}
with the upper triangle approximately commutes on $\mathcal F_{a,1}$ within  $\ep_1$ and the lower triangle approximately commutes on ${\cal F}_{b,2}\cup \phi_1({\cal F}_{a,1})$ within $\ep_2.$
Note also
\beq\label{Misothm-10}
[\phi_2]=\kappa,\,\,\, (\phi_2)_T=\kappa_T\andeqn  (\phi_2)^{\dag}=\kappa_{cu}.
\eneq

We then continue this process, and, by the induction,  {{obtain}} an approximate intertwining:
 \begin{displaymath}
 \xymatrix{
A  \ar[r]^{{\rm id}_A}\ar[d]_{\phi_1}  &   A  \ar[r]^{{\rm id}_A}\ar[d]_{\phi_2}  &   A  \ar[r]^{{\rm id}_A}\ar[d]_{\phi_3} &   \cdots \cdots A \\
  B  \ar[r]_{{\rm id}_B}\ar[ru]^{\psi_1}&    B  \ar[r]_{{\rm id}_B}     \ar[ru]^{\psi_2}&   B\ar[r]_{{\rm id}_B}&  \cdots \cdots B  \\
 }
\end{displaymath}

By the Elliott approximate intertwining argument, this implies
that $A\cong B$ and the isomorphism $\phi$  produced by the above diagram meets the requirements
of \eqref{Misothm-1-1}.
\end{proof}

%
%
%
%

%
%


\section{ Sub-homogeneous \CA s and traces}

Most of the materials of this section are  taken directly from Section 16 of \cite{GLN}.
Here we need non-unital versions of Section 16 of \cite{GLN}.  Most of the proof  {{presented here}} is a modification
of the proof in Section 16 of \cite{GLN}. 

In this section, if $A$ is a non-unital {{$C^*$-algebra, }}
 denote by
${\tau_\C^A}$ the tracial state of $\td A$ vanishing on $A.$

\begin{df}\label{DDo}
Let $D$ be a non-unital \CA. Denote by
$C(\T, {\tilde D})^o$ the \SCA\, of $C(\T, {\tilde D})$ generated by
$C_0(\T\setminus \{1\})\otimes 1_{\tilde D}$ and $1_{C(\T)}\otimes D.$
The unitization of $C(\T, {\tilde D})^o$ is $C(\T, {\tilde D})=C(\T)\otimes {\tilde D}.$
Let $C$ be another non-unital \CA, $L: C(\T, {\tilde D})^o\to C$ {{a}} \cpc\,
and $L^{\sim}: C(\T)\otimes {\tilde D}\to {\tilde C}$ {{the}} unitization.
Suppose that $D$ is amenable.
Denote by $z$ the standard unitary generator of $C(\T).$  For any finite subset ${\cal F}\subset C(\T)\otimes {\tilde D},$
any finite subset ${\cal F}_d\subset {\tilde D},$ and $\ep>0,$
there {{exist}} a finite subset ${\cal G}\subset D$ and $\dt>0$ such that, whenever
$\phi: D\to  C$ is a ${\cal G}$-$\dt$ -multiplicative \cpc \, (for any \CA\, $C$) and
$u=\lambda \cdot 1_{\td C}+c\in \td C$ ($|\lambda|=1,$ $c\in C$) {{is a unitary}} such that
$\|[u, \phi(g)]\|<\dt$ for all $g\in {\cal G},$ there exists an ${\cal F}$-$\ep$-multiplicative
\cpc\, $L': C(\T)\otimes {\tilde D}\to {\tilde C}$ such that (see, for example, 2.8 of \cite
{LnHomtp-c}))
\beq
\|L'(z\otimes 1)-{\bar \lambda} \cdot u\|<\ep/4,\, \|L'(f(z)\otimes1)-{\bar \lambda}u-1\|<\ep/4\\
\andeqn \|L'(1\otimes d)-\phi(d)\|<\ep/4\tforal d\in {\cal F}_d,
\eneq
where {{$f(z)=z-1$}}
 and $f\in C_0(\T\setminus \{1\}).$
{{Note that $\|\pi(L'(f(z)\otimes 1))\|<\ep/4$ and $\|\pi(L'(1\otimes d))\|<\ep/4$ (for all $d\in {\cal F}_d$),
where $\pi: \td C\to \C$ is the quotient map.}}
{{Choose an element $c\in C_+,$ such that
$$
\|cL'(g)c-L'(g)\|<\ep/2\rforal g\in \{{{f\otimes 1_{\td D},~}} 1\otimes d: {{f(z)=z-1,~}} d\in {\cal F}_d, 
\}.
$$
Define $L: C(\T)\otimes \td D\to \td C$ by
$L(g)=cL'(g)c$ for all $g\in C(\T, \td D)^o$ and $L(1_{C(\T, \td D)})=1_{\td C}.$
Then,
\beq\nonumber
L(z\otimes 1)=L(f(z)\otimes 1)+1_{\td C}=cL'(f(z)\otimes 1)c+1_{\td C}\approx_{\ep/2} L'(f(z)\otimes 1)+1_{\td C}
=L'(z\otimes 1)\approx_{\ep/4} {\bar \lambda}u.
\eneq
Similarly,
\beq\label{DPhiv}
L(1\otimes d)\approx_{3\ep/4} \phi(d)\rforal d\in {\cal F}_d.
\eneq}}
We will denote such $L$ by $\Phi_{\phi, u}.$ {{In particular, we have $L(C(\T, \td D)^o)\subset C$.}}

Conversely,  there exists a finite subset ${\cal G}'\subset C(\T, {\tilde D})^o$ and $\dt'>0,$
if $L: C(\T, {\tilde D})^o\to C$ is ${\cal G}'$-$\dt'$-multiplicative \cpc, {{then}}
there is a unitary $u\in {\tilde C}$ such that
\beq
\|L^\sim(z\otimes 1)-u\|<\ep
\eneq
and $\phi=L^{\sim}|_{1\otimes D}$ is a \cpc.
{{It is worth reminding {{that,}}  if $u\in CU(\td C),$ then $u=1_{\td C}+c$ for some $c\in C.$}}

{{Suppose}} that $D$ is a simple \CA. Then $T_f(C(\T)\otimes D)=T_f(C(\T))\otimes T(D)$
{{(see \ref{DTtilde} {{for notation of $T_f$}}).}} {{Note that $C(\T)\otimes D$ is an  essential ideal of $C(\T, {\tilde D})^o.$}}
Therefore
$
T_f(C(\T, {\tilde D})^o)=T_f(C(\T))\otimes T(D).
$
{{Suppose that $D$ is a separable simple \CA\, with continuous scale and $e\in D$
($0\le e\le 1$) is a
strictly positive element. Choose a strictly positive function $f_0\in C_0(\T\setminus \{1\})$
($\|f_0\|=1$).
Then $a:=(1/2)(1_{C(\T)}\otimes e+f_0\otimes 1_{\td D})$ (with $\|a\|\le 1$)
 is a strictly positive element of $C(\T, {\tilde D})^o$ and
 $\tau((a)^{1/n})$ {{converges to}}  $1$
 uniformly
on  $T(C(\T, {\tilde D})^o).$}}
\end{df}

\begin{df}\label{Dbeta}
{\rm Let $A$ be a
\CA. Consider the tensor product
$A\otimes C(\T).$
By the K\"{u}nneth {{Formula}} {{(note that $K_*(C(\T))$ is finitely generated)}}, the tensor product induces two canonical injective \hm s
\begin{equation}\label{Dbeta-1}
\bt^{(0)}: K_0(A)\to K_1(A\otimes C(\T))\quad\mathrm{and}\quad
\bt^{(1)}: K_1(A)\to K_0(A\otimes C(\T)).
\end{equation}\index{$\bt^{(i)}$}
In this way (with further application of the  K\"{u}nneth
Formula), one may write
\begin{equation}\label{Dbeta-2}
K_i(A\otimes C(\T))=K_i(A)\oplus \bt^{(i-1)}(K_{i-1}(A)),\,\,\,i=0,1.
\end{equation}
For each $k\ge 2,$ one also obtains the following injective \hm s
\begin{equation}\label{Dbeta-3}
\bt_k^{(i)}: K_i(A, \Z/k\Z)\to K_{i-1}(A\otimes C(\T),\Z/k\Z),\,\,\,i=0, 1.
\end{equation}
Moreover, one may write
\begin{equation}\label{Dbeta-4}
K_i(A\otimes C(\T),\Z/k\Z)=K_i(A,\Z/k\Z)\oplus
\bt_k^{(i-1)}(K_{i-1}(A,\Z/k\Z)),\,\,\,i=0,1.
\end{equation}


If $x\in \underline{K}(A),$ {{let us write}}  ${\boldsymbol{\bt}}(x)$ for $\bt^{(i)}(x)$
if $x\in K_i(A)$ and for
$\bt_k^{(i)}(x)$ if $x\in K_i(A,\Z/k\Z).$ So {{we have}} an injective \hm\,
\beq\label{Dbeta-6}
{\boldsymbol{\bt}}: \underline{K}(A)\to \underline{K}(A\otimes C(\T)),\andeqn
\\\label{Dbeta-7}
\underline{K}(A\otimes C(\T))=\underline{K}(A)\oplus {\boldsymbol{\bt}}(\underline{K}(A)).
\eneq

Let $h: A\otimes C(\T)\to B$  be a
\hm. Then
$h$ induces a \hm\, $h_{*i,k}: K_i(A\otimes C(\T),\Z/k\Z)\to K_i(B,\Z/k\Z),$
$k=0,2,3,...$ and $i=0,1.$
Suppose that $\phi: A\to B$ is a
\hm\,  and
$v\in  U(B)$ (or $v\in U(\td B)$ if $B$ is not unital) 
is a unitary
such that $\phi(a)v=v\phi(a)$ for all $a\in A.$
Then $\phi$ and $v$ determine a
\hm\,
$h: A\otimes C(\T)\to B$ by $h(a\otimes z)=\phi(a)v$ for all $a\in A,$ where
$z\in C(\T)$ is the identity function on the unit circle $\T.$

We use {${\rm Bott}(\phi,\, v): \underline{K}(A) \to \underline{K}(B)$ to denote the collection of}
all \hm s $h_{*i-1,k}\circ \bt_k^{(i)},$ {where $h: A\otimes C(\T)\to B$ is the homomorphism determined by $(\phi, v)$,}
and we write\index{${\rm Bott}$}
\beq\label{Dbeta-8}
{\rm Bott}(\phi, \,v)=0
\eneq
if $h_{{{*i-1}},k}\circ \bt_k^{(i)}=0$ for all $k$ and $i.$  
In particular, {{if}} $A$ is unital, \eqref{Dbeta-8} implies that $[v]=0$ in $K_1(B)$.
We also {{write}} ${\rm bott}_i(\phi, \, v)$ for
$h_{*i-1}\circ \bt^{(i)},$ $i=0,1.$

Suppose that  $A$ is a
separable amenable \CA.
Let ${\cal Q}\subset \underline{K}(A\otimes C(\T)),$
{{${\cal F}_0\subset A$
and ${\cal F}_1\subset {\td A}\otimes C(\T)$  be finite subsets.}}  Suppose that $(\ep, {\cal F}_{{0}}, {\cal Q})$ is a $KL$-triple (see, 2.1.16 of \cite{Lncbms}, for example).
Then,  by Lemma 2.13  of \cite{GLN} (see also Lemma 2.8 of \cite
{LnHomtp-c})),
 there {{exist}} a finite subset ${\cal G}\subset A$ and $\dt>0$ satisfying the following condition:
For any
${\cal G}$-$\dt$-multiplicative \cpc\, $\phi: A\to B$
and
{{any}}  unitary
$v\in B$ (or $v\in \td B$)
such that
\beq\label{Dbeta-9}
\|[\phi(g), \, v]\|<\dt\rforal g\in {\cal G},
\eneq
there exists a
unital
${\cal F}_1$-$\ep$-multiplicative
\cpc\,  $L: {\td A}\otimes C(\T)\to B$ (or $\td B$) such that
\beq\label{Dbeta-10}
\|L({{f\otimes 1}})-\phi(f)\|<\ep\rforal f\in {\cal F}_0\andeqn
\|L(1\otimes z)-v\|<\ep.
\eneq
In particular, $[L]|_{\cal Q}$ is well defined.
Let ${\cal P}\subset \underline{K}(A)$ be a finite subset.
There are $\dt_{\cal P}>0$ and a finite subset ${\cal F}_{\cal P}$ satisfying the following condition: if
$\phi: A\to B$ is a
${\cal F}_{\cal P}$-$\dt_{\cal P}$-multiplicative {{\cpc\,}} and $v\in B$ (or $v\in {\tilde B}$) is a unitary such that
(\ref{Dbeta-9}) holds for $\dt_{\cal P}$ (in place of $\dt$) and ${\cal F}_{\cal P}$ (in place of ${\cal G}$),
then there exists a unital \cpc\,   $L: A\otimes C(\T)\to B$ which satisfies \eqref{Dbeta-10} {{so}}  that
$[L]|_{{\boldsymbol{\bt}}({\cal P})}$ is well defined,
{{and $[L']|_{{\boldsymbol{\bt}}({\cal P})}=[L]|_{{\boldsymbol{\bt}}({\cal P})}$
if $L'$ also satisfies \eqref{Dbeta-10} {{(for the same $\phi$ and $v$)}} (see  {{2.12 of \cite{GLN}}}).}}
In this case, we will write
\beq\label{Dbeta-11}
{\rm Bott}(\phi,\, v)|_{\cal P}{{=[L]\circ \bt|_{\cal P}}}
\eneq
for all $x\in {\cal P}.$
In particular,
when
$
[L]|_{{\boldsymbol{\bt}}({\cal P})}=0,
$
we will write
\beq\label{Dbeta-12}
{\rm Bott}(\phi,\, v)|_{\cal P}=0.
\eneq
When $K_*(A)$ is finitely generated,
$\mathrm{Hom}_{\Lambda}(\underline{K}(A), \underline{K}(B))$ is determined
by a finitely generated subgroup of $\underline{K}(A)$ (see \cite{DL}).
Let ${\cal P}$ be a finite subset which generates this subgroup.
Then, in this case, instead of (\ref{Dbeta-12}), we may write
\beq\label{Dbeta-13}
{\rm Bott}(\phi, \, v)=0.
\eneq
In general, if ${\cal P}\subset K_0(A),$
we will write
${\rm bott}_0(\phi, \,v)|_{\cal P}={\rm Bott}(\phi,\, v)|_{\cal P},$
and if ${\cal P}\subset K_1(A),$
we will write
${\rm bott}_1(\phi,\, v)|_{\cal P}={\rm Bott}(\phi, \, v)|_{\cal P}.$
}

\end{df}

\begin{NN}\label{Rkt2}
Let $D$ be a non-unital but $\sigma$-unital \CA.
Then  $C(\T, {\tilde D})^o$ is an ideal of $C(\T)\otimes {\tilde D}.$
From the short exact sequence $0\to C(\T, {\tilde D})^o\to C(\T)\otimes {\tilde D}\to \C\to 0$
and the six term exact sequence {{in $K$-theory,}}
one easily computes that
$$
K_0(C(\T, {\tilde D})^o)=K_0(D)\oplus {\boldsymbol{\bt}}(K_1(D))~~{\mbox{and}~~K_1(C(\T, {\tilde D})^o)=K_1(D)\oplus {\boldsymbol{\bt}}(K_0({\tilde D})).}
$$
We {{may}} write $\underline{K}(C(\T, {\tilde D})^o)=\underline{K}(D)\oplus {\boldsymbol{\bt}}(\underline{K}({\tilde D}))$ or,  \beq\label{june-19-2021}
\underline{K}(C(\T, {\tilde D})^o)=\underline{K}(D)\oplus {\boldsymbol{\bt}}(\underline{K}(D))\oplus {\boldsymbol{\bt}}(\underline{K}(\C\cdot1_{\tilde D})),\eneq
as ${\boldsymbol{\bt}}(\underline{K}({\tilde D}))={\boldsymbol{\bt}}(\underline{K}(D))\oplus {\boldsymbol{\bt}}(\underline{K}(\C\cdot1_{\tilde D}))$.


{\it For convenience, we may {{write}} $D^\T$ for {{$C(\T, {\tilde D})^o.$}}}

{{Let $\pi_1: C(\T,\td D)\to \td D$ be the point-evaluation at $1\in \T$
and also denote by $\pi_1$ the restriction ${\pi_1}|_{C(\T, \td D)^o}.$
Then there is a splitting short exact sequence\\
\,\,\,\,\,$0\to I:=C_0(\T\setminus \{1\}, \td D)
\to D^\T\stackrel{\pi_1}{\to}D\to 0.$ If $K_1(D)=\{0\},$
then ${\pi_1}_{*0}: K_0(D^\T)\to K_0(D)$ is an isomorphism.  Moreover $\td\pi_1: \widetilde{D^\T}\to \td D$
induces an  order isomorphism from $K_0(\widetilde{D^\T})$ onto $K_0(\td D).$}}

\end{NN}

The following is {{well known}} and follows, for example, from 10.10.4 of  \cite{Bla-Ktheory}.

\begin{lem}\label{LKT2}
Let $E$ be a non-unital and $\sigma$-unital \CA.
Then
$({\rm ker}\rho_{E}(K_0(E)))\oplus {\boldsymbol{\bt}}(K_1(E))= {\rm ker}\rho_{C(\T, {\tilde E})^o}.$
\end{lem}

\begin{proof}
It suffices to prove that ${\boldsymbol{\bt}}(K_1(E))\subset {\rm ker}\rho_{C(\T)\otimes {\tilde E}}.$
Fix a unitary  $u\in M_n({\tilde E})$ {{with $[u]\not=0$ in $K_1(E)$. Then the spectrum of $u$ is the full circle $\T$.}}
Let $C$ be the \SCA\, generated by $z\otimes 1_{\tilde E}$ and $1_{C(\T)}\otimes u.$
Therefore $C\cong C(\T^2).$
Denote by $e=\diag(1, 0)$ and $p\in M_2(C(\T^2)),$ a non-trivial rank one projection.
Let $j: C\to C(\T, M_n({\tilde E}))$ be the embedding.
Then $j(e), j(p)\in C(\T, M_n({\tilde E})).$
We identify $e$ with the obvious projection in $C(\T, M_n({\tilde E}))$ and
$p$ with another projection in $C(\T, M_n({\tilde E})).$
Then $\boldsymbol{\bt}([u])$ may be represented by
$\pm([e]-[p]).$ Note that, for every trace $t\in T(C),$
$t([e]-[p])=0.$ It follows that, for any $\tau\in T(C(\T)\otimes {\tilde E}),$
$\tau(j([e])-j([p]))=0.$ This implies that ${\boldsymbol{\bt}}([u])\subset {\rm ker}\rho_{C(\T)\otimes {\tilde E}}.$
The lemma follows.
\end{proof}

\begin{lem}[16.8 of \cite{GLN}]\label{APPextnn}
Let  $A=C$ or
$A=C(\T, {\td C})^o$ for some {{$C\in {\cal C}_0.$}} 
Let $\Delta: A_+^{q, {{\bf1}}}\setminus\{0\}\to(0, 1)$ be an order preserving map. Let $\mathcal H\subset A$ be a finite subset and let $\sigma>0$. Then there are  a finite subset $\mathcal H_1\subset A_+^{\bf 1}\setminus \{0\},$
$\dt>0,$ a finite subset ${\cal P}\subset K_0(A),$ and  an integer $K\ge 1$ such that, for any
$\tau\in {{\R_+}}\cdot T(A)$  {{(with $0<\|\tau\|\le 1$)}}  which satisfies
\begin{equation*}
\tau(h) > \Delta(\hat{h}) \tforal h\in \mathcal H_{1}
\end{equation*}
 and any { {positive homomorphism} } $\kappa: K_0(\td A)\to K_0(M_s)=\Z$
with $s=\kappa({{[1_{\td A}]}})$
such that
\beq\label{2020-715-n1}
|\rho_{\td A}(x)(\tau^\sim)-(1/s)(\kappa(x))|<\dt
\eneq
for all $x\in {\cal P},$ where $\tau^\sim=\tau+(1-\|\tau\|) \tau_\C^A,$
there is a homomorphism $\phi: A\to M_{sK}$ such that
$\phi_{*0}=K\kappa$ and $$|{\rm tr}\circ\phi(h)-\tau(h)|<\sigma \tforal h\in\mathcal H,$$
where ${\rm tr}$ is the tracial state on $M_{sK}.$
\end{lem}

\begin{proof}
We only consider the case that $A=C(\T, {{\td C}})^o.$
The case
$A=C$ will follow from a simplification of the proof.
Let $\DT,$ ${\cal H}$ and $\sigma$ be given as in the lemma.
Let $A_1:=\td C\otimes C(\T)={\tilde A}$ and $\pi: A_1\to \C$ be the quotient map.
For each $x\in A_{s.a.},$ denote by $\hat x$ the image of $x$ in $A^q_{s.a.}$ (see \ref{Dlambdas}).
%
We will apply  the result of the unital case of 16.8  of \cite{GLN}.
For $x\in {A_1}_+^{\bf 1}\setminus \{0\},$
define
\beq\label{164gln-0}
\Delta_1(\hat{x})=\sup\{\Delta(\widehat{e_n^{1/2}xe_n^{1/2}}): n\ge 1\},
\eneq
where $\{e_n\}$ is an approximate identity for $A.$
Therefore $\Delta_1(\hat{x})>0$ as $x\not=0.$
It is also clear that $\Delta_1$ is order preserving since $\Delta$ is.
Furthermore
\beq
\widehat{e_n^{1/2}xe_n^{1/2}}=\widehat{x^{1/2}e_nx^{1/2}}\le \hat{x}.
\eneq
It follows that $\Delta_1(x)\le \Delta(x)$ of $x\in A_+^{\bf 1}.$
Moreover, $\tau(x)\ge \tau(e_n^{1/2}xe_n^{1/2})$ for all $n$ and $\tau\in T(A).$
For each $x\in {A_1}_+^{\bf 1}\setminus 
\{0\},$ there is $n_x\ge 1$ such that
\beq
\Delta(\widehat{e_{n_x}^{1/2}xe_{n_x}^{1/2}})\ge (3/4)\sup\{\Delta(\widehat{e_n^{1/2}xe_n^{1/2}}): n\ge 1\}.
\eneq
Define
\beq
\Delta_2(\hat{x})=(3/4)\Delta_1(\hat{x})\rforal x\in {A_1}^{\bf 1}_+\setminus \{0\}.
\eneq
%
%
%
%
%
%
%
%
Let ${\cal H}_0\subset (A_1)_+^{\bf 1}\setminus \{0\}$
(in place of ${\cal H}_1$),  $\dt_1>0$ (in place of $\dt$), ${\cal P}_1$
(in place of ${\cal P}$)  and $K$ be given by Lemma 16.8 of \cite{GLN} for $A_1$ (in place $A$),  $\DT_2$ (in place of $\Delta$), $\sigma/2$ (in place of $\sigma$)
and ${\cal H}.$
{{\Wlog, we may assume that ${\cal P}_1=({\cal P}_1\cap (K_0(A)))\sqcup ({\cal P}_1\cap \Z\cdot [1_{\td A}])$. Define ${\cal P}:={\cal P}_1\cap (K_0(A))$.}}
Let $p_1, {\bar p}_1, p_2, {\bar p}_2,  ...,p_m, {\bar p}_m\in M_R(\td A)$ be projections
for some integer $R\ge 1$ such that
$\{[p_1]-[{\bar p}_1], [p_2]-[{\bar  p}_2],..., [p_m]-[{\bar p}_m]\}={\cal P}.$

Choose $\dt=\dt_1/2.$
%
Define
\beq
{\cal H}_1:=\{\widehat{e_{n_x}^{1/2}xe_{n_x}^{1/2}}: x\in {\cal H}_0\}.
\eneq
Now assume that $\tau\in {{\R_+}}\cdot T(A)$ and $\kappa: K_0({\tilde A})\to \Z=K_0(M_s)$ be given as described in
the lemma for ${\cal H}_1,$ $\dt,$ and $K$ above.
Note {{that}} $\kappa_1$ is order preserving and $\kappa_1([1_{\td C}])=s.$
Recall {{that}}
\beq\label{2020-715-n2}
\tau(f)\ge \Delta(\hat f)\rforal f\in {\cal H}_1\andeqn |\rho_{\td A}(x)(\tau^\sim)-(1/s)\kappa(x)|<\dt\rforal x\in {\cal P}.
\eneq
It follows that, for all $f\in {\cal H}_0,$
\beq
\tau(f)\ge \tau(f^{1/2}e_{n_f}f^{1/2})=\tau(e_{n_f}^{1/2}fe_{n_f}^{1/2})
\ge \Delta(\widehat{e_{n_f}^{1/2}fe_{n_f}^{1/2}})\ge \Delta_2(\hat{f}).
\eneq
Since $\kappa_1([1_{\td C}])=s,$  $\kappa_1|_{{\cal P}_3}=0,$  and $\tau\in {{\R_+}} \cdot T(A)$ with $0<\|\tau\|\le 1,$
  by \eqref{2020-715-n2},
\beq
|\rho_{A_1}(x)(\tau^\sim)-(1/s)\kappa(x)|<\dt\rforal x\in {\cal P}_1.
\eneq
Note that we may view $\tau^\sim$ as a tracial state
of $A_1$ by $\tau^\sim(a+\lambda 1_{A_1})=\tau(a)+(1-\|\tau\|)\lambda$
for all $a\in A$ {{and $\lambda\in \C.$}}
%
It follows from 16.8 of \cite{GLN} {{that}} there is {{a}} \hm\, $\psi: A_1\to M_{sK}$
such that  $\psi_{*0}=K\kappa_1$ and
\beq\label{162--2020}
|{\rm tr}(\psi(f)-\tau(f)|<\sigma\rforal f\in {\cal H}.
\eneq
where ${\rm tr}$ is the tracial state of $M_{sK}.$
Let $\phi: A\to M_{sK}$ be defined by $\phi:=\psi|_A.$
{{Lemma}} follows.

\end{proof}

\begin{NN}\label{Rextendable}
Let $A$ be a non-unital \CA\, {{and}} $T_0(A)=\{ r\cdot \tau: r\in [0,1], \tau\in T(A)\}.$
Let $B$ be a unital \CA\, {{and}} $\gamma: T(B)\to T_0(A)$  be a continuous map
such that $\|\gamma(\tau)\|>0$ for all $\tau\in T(B).$
Suppose that there is an affine continuous map $\gamma^\sim: T(B)\to T(\td A)$
such that $\gamma^\sim(\tau)|_A=\gamma(\tau)$ for all $\tau\in T(B).$
Then we say $\gamma$ is extended to an affine continuous map $\gamma^\sim.$
If $B$ is not unital, then we say $\gamma$ is extended to an affine continuous map
$\gamma^\sim: T(\td B)\to T(\td A)$ if there is an affine continuous map
$\gamma^\sim: T(\td B)\to T(\td A)$ such that $\gamma^\sim(t^B_\C)=t^A_\C,$
where $t^A_\C$ (or $t^B_\C$) is the tracial state of $\td A$ which vanishes on $A$ (or on $B$),
and $\gamma^\sim(\tau)|_A=\gamma(\tau)$
for all $\tau\in T(B)(\subset T(\td B)).$
\end{NN}



\begin{lem}{\rm (16.9 of \cite{GLN})}\label{ExtTraceI-D}
Let $A=C$ or $A=C(\T, \td C)^o$
for some $C\in {\cal C}_0.$
Let $\Delta: {{A_+^{q, {\bf 1}}}}\setminus\{0\}\to(0, 1)$ be an order preserving map. Let $\mathcal F, \mathcal H\subseteq {{A}}$ be  finite subsets, and let $\epsilon>0, \sigma>0$.
Then there are  {{finite subsets $\mathcal H_1\subseteq {{A}}_+^{\bf 1}\setminus \{0\}$,
${\cal P}\subset K_0(A),$}}  $\dt>0,$ and a positive integer $K$ such that, for any continuous  map
$\gamma: T(C([0, 1]))\to T_0(A)$ with $\inf\{\|\gamma(\tau)\|: \tau\in T(C([0,1]))\}>0$
which is extended to an affine continuous map $\gamma^\sim: T{{(C([0,1]))}}\to T(\td A)$
satisfying
$$
\gamma(\tau)(h) > \Delta(\hat{h}) \tforal h\in \mathcal H_{1} \tand \tforal \tau\in T(C([0, 1]))
$$
and  any positive homomorphism $\kappa: K_0({\tilde A})\to K_0(M_s(C([0, 1])))$ with $\kappa([1_{\tilde A}])=s$
such that
$$
|\rho_A(x)(\gamma^\sim(\tau))-(1/s)\tau(\kappa(x))|<\delta  \tforal \tau\in T(C([0, 1]))
$$
for all $x\in {\cal P},$
there is an 
$\mathcal F$-$\epsilon$-multiplicative {{\cpc\,}}
$\phi: A\to M_{sK}(C([0, 1]))$ such that
$\phi_0=K\kappa$ and
$$|\tau\circ\phi(h)-\gamma'(\tau)(h)|<\sigma\tforal h\in\mathcal H,$$
where $\gamma': T(M_{sK}(C([0,1])))\to T_0(A)$ is induced by $\gamma$, {{by identifying $T(M_{sK}(C([0,1])))$ with $T(C([0,1]))$ .}} { {Furthermore $\phi_0=\pi_0\circ \phi$ and $\phi_1=\pi_1\circ \phi$ are {{genuine}} homomorphisms.}}

{{In}} the case that $A\in\mathcal C_0$, the map $\phi$ can be chosen to be a homomorphism.
\end{lem}

\begin{proof}
Since any
\CA s in $\mathcal C_0$ are semi-projective,
the second part of the statement follows directly from the first part of the statement. Thus, let us only show the first part of the statement.
We will consider the case $A=C(\T, \td C)^o.$
Let $\DT,$ ${\cal H}$ and $\ep$ and $\sigma$ be given as in the lemma.
Let $A_1:=\td C\otimes C(\T)=\td A$ and $\pi: A_1\to \C$
be the quotient map.
Without loss of generality, one may assume that $\mathcal F\subseteq \mathcal H$.

Let  $\Delta_2: (A_1)_+^{q, \bf 1}\setminus \{0\}\to (0,1)$  be defined as in the proof of \ref{APPextnn} associated with $\DT$
above.
The rest of the proof  is to reduce to the unital case and via a same route as
in the proof of  \ref{APPextnn}.

Let ${\cal H}_0\subset (A_1)_+^{\bf 1}\setminus \{0\}$
(in place of ${\cal H}_1$),  $\dt_1>0$ (in place of $\dt$), ${\cal P}_1$
(in place of ${\cal P}$)  and $K$ be given by Lemma 16.9  of \cite{GLN} for $A_1$ (in place $A$),  $\DT_2$ (in place of $\Delta$), $\sigma/2$ (in place of $\sigma$)
and ${\cal H}.$

{{\Wlog, we may assume that ${\cal P}_1=({\cal P}_1\cap (K_0(A)))\sqcup ({\cal P}_1\cap \Z\cdot [1_{\td A}])$. Define ${\cal P}:={\cal P}_1\cap (K_0(A))$.}}
Let $p_1, {\bar p}_1, p_2, {\bar p}_2,  ...,p_m, {\bar p}_m\in M_R(\td A)$ be projections
for some integer $R\ge 1$ such that
$\{[p_1]-[{\bar p}_1], [p_2]-[{\bar  p}_2],..., [p_m]-[{\bar p}_m]\}={\cal P}.$

Choose $\dt=\dt_1/2.$
Define  (see the proof of \ref{APPextnn})
\beq
{\cal H}_1:=\{\widehat{e_{n_x}^{1/2}xe_{n_x}^{1/2}}: x\in {\cal H}_0\}.
\eneq

Now assume that $\gamma$ and $\kappa: K_0({\tilde A})\to \Z=K_0{{(M_s(C([0,1])))}}$ are given as described in
the lemma for ${\cal H}_1,$ $\dt$ and $K$ above.
Note that, for all $\tau\in T(C([0,1])$ {{and}}  all $f\in {\cal H}_0,$
\beq
\gamma(\tau)(f)\ge  \gamma(\tau)(f^{1/2}e_{n_f}f^{1/2})
\ge \Delta(\widehat{e_{n_f}^{1/2}fe_{n_f}^{1/2}})\ge \Delta_2(\hat{f}).
\eneq
As in the proof of \ref{APPextnn}, we also have, for all $x\in {\cal P}_1,$
\beq
|\rho_A(x)(\gamma^\sim(\tau))-(1/s)\tau(\kappa(x))|<\delta  \tforal \tau\in T(C([0, 1])).
\eneq
By applying 16.9 of \cite{GLN}, we obtain an
$\mathcal F$-$\epsilon$-multiplicative completely positive linear map
$\psi: A_1\to M_{sK}(C([0, 1]))$ such that
$\phi_0=K\kappa$ and
$$|\tau\circ\psi(h)-\gamma'(\tau)(h)|<\sigma\tforal h\in\mathcal H,$$
where $\gamma': T(M_{sK}(C([0,1])))\to T_0(A)$ is induced by $\gamma.$ Furthermore $\phi_0=\pi_0\circ \psi$ and $\phi_1=\pi_1\circ \psi$ are true homomorphisms.

Define $\phi:=\psi|_{A}.$ Then {{the}}  lemma follows.
\end{proof}

\begin{thm}[4.18 of \cite{GLN}]\label{UniqAtoM}
Let $A$ be a non-unital subhomogeneous \CA\, such that  ${\tilde A}\in {\cal {\bar D}}_s$ (see 4.8 of \cite{GLN}) 
{{and}}
let $\Delta: A_+^{q,{\bf 1}}\setminus \{0\}\to (0,1)$ be an order preserving map.

For any $\ep>0$ and finite subset ${\cal F}\subset A,$ there exists $\dt>0,$ a finite subset
${\cal P}\subset \underline{K}(A),$ a finite subset ${\cal H}_1\subset A_+^{\bf 1}\setminus \{0\}$ and
a finite subset ${\cal H}_2\subset A_{s.a.}$ satisfying the following:

If $\phi_1, \phi_2: A\to M_n$ are {{two  \hm s}} such that
\beq\label{UnAM-1}
&&[\phi_1]|_{\cal P}=[\phi_2]|_{\cal P},\\\label{UNAM-1+}
&&\tau\circ \phi_1(g)\ge \Delta(\hat{g})\tforal g\in {\cal H}_1\tand\\
&&|\tau\circ \phi_1(h)-\tau\circ \phi_2(h)|<\dt\tforal h\in {\cal H}_2,
\eneq
where $\tau\in T(M_n),$ then there exists a unitary $u\in M_n$ such that
\beq\label{UnAM-2}
\|{\rm Ad}\, u\circ \phi_1(f)-\phi_2(f)\|<\ep\tforal f\in {\cal F}.
\eneq
{{Furthermore, if both $K_0(A)$ and $K_1(A)$ are torsion free groups, then the subset ${\cal P}\subset \underline{K}(A)$  can be chosen such that ${\cal P}\subset K_0(A)\subset  \underline{K}(A)$.}}
\end{thm}

\begin{proof}
We will apply  the result of the unital case  as 4.18 of \cite{GLN}.

Define $\Delta_2: {\tilde A}_+^{\bf 1}\setminus \{0\}\to (0,1)$ be defined  as in the proof of \ref{APPextnn}.
%
%
%
%
Let $\ep>0$ and  ${\cal F}\subset A$  be given.
Let  $\dt>0,$ ${\cal P}_1\subset {{\underline{K}(\td A)}}$  (in place of ${\cal P}$) be a finite subset,
${\cal H}_1'\subset {\tilde A}_+^{\bf 1}\setminus \{0\}$ (in place of ${\cal H}_1$) and
${\cal H}_2'\subset {\tilde A}_{s.a.}$ be finite subsets required by 4.18 of \cite{GLN}
for  $\ep,$ ${\cal F},$ ${\tilde A}$ and $\Delta_2.$

Consider the short exact sequence:
\beq
0\longrightarrow  A\,\,  \stackrel{j_A}{\longrightarrow}\,\, {\tilde A} \,\, \stackrel{\pi_A}{\rightleftharpoons}_{s_A}\, \C \to 0.
\eneq
Let
${\cal P}=\{x-[s_A\circ \pi_A](x): x\in {\cal P}_1\}\subset \underline{K}(A)$ {{and}}
${\cal H}_1=\{e_{n_x}^{1/2}xe_{n_x}^{1/2}: x\in {\cal H}_1'\}\cup {\cal H}_1'\cap A.$
For each $x\in {\cal H}_2',$ one may write
$x=\lambda + h_x,$ where $\lambda\in \R$ and $h_x\in A_{s.a.}.$
Set ${\cal H}_2=\{h_x: x\in {\cal H}_2'\}.$

Now suppose that $\phi_1, \phi_2: A\to M_n$ are two \hm s which {{are extended}} to
two unital \hm s  $\phi_1^{\sim}, \phi_2^{\sim}: {\tilde A} {{\to}} M_n$ such that
they satisfy the condition for the {{above-mentioned}} ${\cal P},$ ${\cal H}_1$ and ${\cal H}_2.$
It is immediate that
\beq\label{417gln-10}
[\phi_1^{\sim}]|_{{\cal P}_1}=[\phi_2^{\sim}]|_{{\cal P}_1}.
\eneq
We also have that, for all $x\in {\cal H}_1',$
\beq\nonumber
\tau\circ \phi_1^{\sim}(x) &\ge & \tau(\phi_1(e_{n_x}^{1/2}xe_{n_x}^{1/2}))\\\label{417gln-11}
&\ge& \Delta(\widetilde{e_{n_x}^{1/2}xe_{n_x}^{1/2}})\ge (3/4)\sup\{\Delta(\widetilde{e_n^{1/2}xe_n^{1/2}}): n\in \N\}
=\Delta_2(\hat{x}).
\eneq
We further estimate that, for all $x=\lambda+h_x\in {\cal H}_2',$
\beq\label{417gln-12}
|\tau(\phi_1^{\sim}(\lambda+h_x))-\tau(\phi_2^{\sim}(\lambda+h_x))|=|\tau(\phi_1(h_x))-\tau(\phi_2(h_x))|<\dt
\eneq
for $\tau\in T(M_n).$
Combining \eqref{417gln-10}, \eqref{417gln-11} and \eqref{417gln-12} and applying Theorem 4.18 of \cite{GLN},
we obtain  a unitary $u\in M_n$
such that
\beq
\|{\rm Ad}\, u\circ \phi_1(f)-\phi_2(f)\|<\ep\rforal f\in {\cal F}.
\eneq
For the  {{``Furthermore"}} part, one only needs to notice that $K_1(M_n)=0$,  $K_1(M_n, \Z/k\Z)=0$, and an element $\xi\in Hom_{\LD}(\underline{K}(A),{{\underline{K}(M_n)}})$ is completely decided by $\xi|_{K_0(A)}$.

\end{proof}

\begin{thm}{\rm (cf. 16.10 of \cite{GLN})}\label{ExtTraceC-D}
 Let $C=C_1$ or $C=C(\T, \td C_1)^o$ for some $C_1\in {\cal C}_0.$
 Let $\Delta: C_+^{q, {\bf1}}\setminus\{0\}\to (0, 1)$ be an order preserving map. Let $\mathcal F, \mathcal H\subseteq C$ be finite subsets, and let $1>\sigma,\, \epsilon>0$.
There exist  a finite subset $\mathcal H_1\subseteq C_+^{\bf 1}\setminus \{0\}$, $\delta>0$,
a finite subset ${\cal P}\subset K_0(C)$ and a positive integer $K$ such that for any continuous {{map $\gamma: T(D)\to T_0(C)$  with
$\inf\{\|\gamma(\tau)\|: \tau\in T(D)\}>0$ which is
extended to an affine continuous map $\gamma^\sim: T(\td D)\to T_0(\td C)$}}
satisfying
$$\gamma(\tau)(h) > \Delta(\hat{h})\rforal h\in \mathcal H_1\tforal \tau\in T(\td D),$$
where $D$ is a non-unital
 C*-algebra
 in $\mathcal C_0$,
any positive \hm\, $\kappa: K_0(\td C)\to K_0(\td D)$ {with $\kappa([1_{\tilde C}])=s[1_{\tilde D}]$}
 for some integer $s \ge 1$ satisfying
$$|\rho_C(x)(\gamma^\sim(\tau))-{(1/s)}\tau(\kappa(x))|<\delta\tforal \tau\in T(D)$$
and for all $x\in {\cal P},$
there is an $\mathcal F$-$\epsilon$-multiplicative positive linear map $\phi: C\to M_{{sK}}(D)$ such that
\beq
&&\phi_{*0}={K}\kappa\tand\\
&&|(1/(sK))\tau\circ\phi(h)-\gamma(\tau)(h)|<\sigma\tforal h\in\mathcal H\tand \tau\in T(D).
\eneq

In the case that $C\in\mathcal C_0$, the map $\phi$ can be chosen to be a homomorphism.
\end{thm}
\begin{proof}
As in the proof of \ref{ExtTraceI-D}, since \CA s in ${\cal C}_0$ are semi-projective, we will only prove the first part
of the statement.

Without loss of generality, one may {also} assume that $\mathcal F\subseteq \mathcal H$.



Let $\mathcal H_{1, 1}\subset C_+^1\setminus \{0\}$ (in place of $\mathcal H_{1}$),
 $\sigma_1>0$ (in place of ${{\delta}}$)  and ${\cal P}_1\subset \underline{K}(C)$ (in place of ${\cal P}$)  be finite subset
required by
{{Theorem \ref{UniqAtoM}}} with respect to $C$ (in the place of ${{A}}$), $\min\{\sigma/{4}, \epsilon/{2}\}$ (in the place of $\epsilon$), $\mathcal H$ (in the place of $\mathcal F$) and $\Delta$.
Note that
{{$K_i(C)$ is {{a free abelian}} group, $i=0,1.$}}   {{Therefore, by the ``Furthermore" part of  Theorem \ref{UniqAtoM},}}
we may assume
that {{${\cal P}_1\subset K_0(C).$}}


Let $\mathcal H_{1, 2}\subseteq C$ (in place of $\mathcal H_1$) {be a finite subset}, let $\sigma_2$ (in place of $\delta$)
be a positive  number,  ${\cal P}_2\subset K_0(C)$ (in place of ${\cal P}$) be a finite subset,
{and $K_1$ (in place of $K$) be an integer required by  Lemma \ref{APPextnn}} with respect to $\mathcal H\cup\mathcal H_{1, 1}$
and ${\frac{1}{2}\min\{\sigma/16, \sigma_1/4, \min\{\Delta(\hat{h})/2: \ h\in\mathcal H_{1, 1}\}\}}$ (in the place of {$\sigma$}) and {{$\Delta$.}}   

Let $\mathcal H_{1, 3}$ (in place of $\mathcal H_1$), $\sigma_3>0$ (in place of $\delta$),
${\cal P}_3\subset K_0(C)$ (in place of ${\cal P}$) be a finite subset and $K_2$
(in place of $K$) be  finite subset and constants required by  Lemma \ref{ExtTraceI-D}  with respect to $C$, $\mathcal H\cup\mathcal H_{1, 1}$
(in the place of $\mathcal H$), {$\min\{\sigma/16, \sigma_1/4, \min\{\Delta(\hat{h})/2:\ h\in\mathcal H_{1, 1}\}\}$} (in the place of $\sigma$), $\ep/4$ (in place of $\ep$), ${\cal H}$ (in place of
${\cal F}$)  and $\Delta$. 
{{Let ${\cal P}\supset {\cal P}_1\cup {\cal P}_2 \cup {\cal P}_3$ be a finite set generating $K_0(C)$.}}

Put $\mathcal H_1=\mathcal H_{1, 1}\cup\mathcal H_{1, 2}\cup\mathcal H_{1, 3},$
$\delta=\min\{{\sigma_1/2}, \sigma_2, {1/4}\}$ and $K=MK_1K_2$.
Let
\beq\label{04-28-2021}
\hspace{-0.2in}D={{A}}(F_1, F_2, \psi_0, \psi_1)=\{(f,a)\in C([0,1], F_2)\oplus F_1: f(0)=\psi_0(a)\andeqn f(1)=\psi_1(a)\}
\eneq
be a
\CA\, in $\mathcal C_0$, and let $\gamma^\sim: T(\td D)\to T(\td C)$ be a given continuous affine map satisfying $$\gamma(\tau)(h) > \Delta(\hat{h})\rforal h\in \mathcal H_{1} \rforal \tau\in T(D).$$
{L}et $\kappa: K_0({\tilde C})\to K_0(M_s({\tilde D}))$ be any positive  \hm\,  {with $s[1_{\tilde D}]=\kappa([1_{\tilde C}])$} satisfying
$$
|\rho_C(x)(\gamma^\sim(\tau))-{(1/s)}\tau(\kappa(x))|<\delta\rforal \tau\in T(D)
$$
{and for all $x\in {\cal P}.$}
Write $C([0, 1], F_2)=I_1\oplus I_2\oplus \cdots \oplus I_{{k}}$ with $I_i=C([0, 1], M_{r_i})$, $i=1, ..., {{k}}$.
Note that $\gamma^\sim$ induces a continuous affine map $\gamma_i^\sim: T(I_i)\to T(\td C)$
{{(which extends $\gamma_i: T(I_i)\to T_0(C)$ with
$\inf\{\|\gamma_i(t)\|: t\in T(I_i)\}>0$)}}
by
$\gamma_i: T(I_i) \to T_0(C)$ defined by $\gamma_i(\tau)=\gamma(\tau\circ\pi_i)$ for each $1\leq i\leq {k},$   where $\pi_i$ is the restriction map $D\to I_i$
defined by $(f,a)\to f|_{[0,1]_i}$ {{(the restriction on the $i$-th interval).}}  It is clear that
for any $1\leq i\leq {k}$,
one has that
\beq
\label{intv-dense}
\gamma_i(\tau)(h)>\Delta(\hat{h}) \rforal h\in\mathcal H_{1, 3}{\andeqn} \rforal \tau\in T(I_i)
\andeqn\\\label{intv-cpt}
|\rho_C(x)(\gamma_i^\sim(\tau))-\tau((\pi_i)_{*0}\circ\kappa(x))|<\delta\leq \sigma_3\rforal  \tau\in T(M_s(I_i)),
\eneq
{and for all $x\in {\cal P} {{(\supset {\cal P}_3)}}$ and for any $1\le i\le k.$}
Also write $F_1=M_{{R_1}}\oplus\cdots\oplus M_{{R_l}}$ and denote by $\pi'_j: D\to M_{{R_j}}$ the corresponding
quotient map of $D$.
Since
\beq\nonumber
&&\gamma(\tau)(h)>\Delta(\hat{h})\rforal h\in\mathcal { {H}}_{1, 2}\,\,\, {\rm and}\rforal \tau\in T(D),
\andeqn\\
&&|\rho_C(x)(\gamma^\sim(\tau))-{(1/s)}\tau(\kappa(x))|<\delta  \rforal \tau\in T(D)
\eneq
{and for all  $x\in {\cal P}{{(\supset {\cal P}_2)}}, $} one has that, for each $j,$
\beq\nonumber
\gamma\circ(\pi_j')^*(\mathrm{tr})(h)>\Delta(\hat{h})\rforal h\in\mathcal H_{1, 2}\andeqn\\
|\rho_C(x)(\gamma^\sim \circ(\pi_j')^*(\mathrm{tr'}))-\mathrm{tr}([\pi_j']\circ\kappa(x))|<\delta\leq \sigma_2,
\eneq
where ${\rm tr}$ is the tracial state on $M_{sR_j}$ and ${\rm tr}'$ is the tracial
state on $M_{R_j},$ for all  $x\in K_0(C)$ and
where $\gamma\circ (\pi_j')^*({\rm tr})=\gamma({\rm tr}\circ \pi_j).$

It follows from \ref{APPextnn} that there is a homomorphism $\phi'_j: C\to M_{{R}_j}\otimes M_{sK_1K_2}$ such that
\beq\label{C-D-nnn1}
&&{(\phi'_j)_{*0}=(\pi_j')_{*0}\circ {K_1K_2}\kappa}\andeqn\\
\label{pt-pre}
&&\hspace{-0.2in}|\mathrm{tr}\circ \phi'_j(h)- (\gamma \circ (\pi_j')^*)(\mathrm{tr}')(h) | <{\min\{\sigma/16, \sigma_1/4, \min\{\Delta(\hat{h})/2: \ h\in\mathcal H_{1, 1}\}\}}
\eneq
for all $h\in \mathcal H\cup\mathcal H_{1, 1},$
where
${\rm tr}$ is the tracial
state on $M_{R_j}\otimes M_{sK}$ and
where ${\rm tr}'$ is the tracial state on $M_{R_j}.$
{{Put}} $$\phi'=\bigoplus_{j=1}^l \phi'_j: C\to F_1\otimes\mathrm{M}_{sK_1K_2}(\mathbb C).$$
Applying Lemma \ref{ExtTraceI-D} to \eqref{intv-dense} and \eqref{intv-cpt}, one obtains, for any $1\leq i\leq {k}$,
an ${\cal H}$-$\ep/4$-multiplicative \morp\,
$\phi_i: C\to I_i\otimes M_{{sK_1K_2}}$ such that ${(\phi_i){*_0}=(\pi_i)_{*0}\circ K_1K_2\kappa}$ and
\begin{equation}\label{int-pre}
|(1/sK_1K_2)\tau \circ \phi_i(h)- ((\gamma\circ(\pi_i)^*)(\tau))(h) |< {\min\{\sigma/16, \sigma/4, \min\{\Delta(\hat{h})/2:\ h\in\mathcal H_{1, 1}\}\}}
\end{equation}
for all $h\in \mathcal H\cup\mathcal H_{1, 1}\cup\mathcal G_1,$ where $\tau\in T({I_i}).$ {{Furthermore, as
{{in the}}  conclusion of Lemma \ref{ExtTraceI-D}, the restrictions of $\phi_i$ to both boundaries are  homomorphisms.}}


For each $1\leq i\leq {k}$, denote by $\pi_{i, 0}$ and $\pi_{i, 1}$ the evaluations of $I_i\otimes M_s$ at the point $0$ and $1$ respectively. Then one has
\begin{equation}\label{C-D-n1}
{\psi_{0,i}\circ \pi_e=\pi_{i,0}\circ \pi_i,}~~{{\mbox{and}~~\psi_{1,i}\circ \pi_e=\pi_{i,1}\circ \pi_i,}}
\end{equation}
{{where {{$\psi_{0,i}, \psi_{1,i}: F_1\to F_2^i $}} are partial {{maps}} of $\psi_0$ and $\psi_1$ appeared in the definition of $D$ (see (\ref{04-28-2021})).}}
It follows that
\begin{eqnarray}\label{C-Dn-2}
(\psi_{0, i}\circ\phi')_{*0} & = & (\psi_{0,i})_{*0}\circ (\sum_{j=1}^l(\pi_j')_{*0})\circ K_1K_2\kappa\\
&=&(\psi_{0,i})_{*0}\circ (\pi_e)_{*0}\circ K_1K_2\kappa=(\pi_{i,0}\circ \pi_i)_{*0}\circ K_1K_2\kappa\\
&=&(\pi_{i, 0})_{*0}\circ (\phi_i)_{*0}.
\end{eqnarray}
Moreover, note that by \eqref{pt-pre},
\begin{equation}
\mathrm{tr}\circ (\psi_{0, i}\circ\phi')(h)\geq \Delta(\hat{h})/2\rforal h\in \mathcal H_{1, 1},
\end{equation}
and by \eqref{int-pre},
\begin{equation}
\mathrm{tr}\circ (\pi_{i, 0}\circ \phi_i)(h)\geq \Delta(\hat{h})/2\rforal h\in \mathcal H_{1, 1}.
\end{equation}
It also follows from \eqref{pt-pre} and \eqref{int-pre} that
\begin{equation}
|\mathrm{tr}\circ (\psi_{0, i}\circ\phi')(h) - \mathrm{tr}\circ (\pi_{i, 0}\circ \phi_i)(h) | < \sigma_1/2\rforal h\in \mathcal
H_{1,1}.
\end{equation}
Consider {{amplifications}}
\beq\nonumber
\phi'_i:&=&\phi_i\otimes 1_{\mathrm{M}_M(\mathbb C)}: C\to I_i\otimes M_{{sK}}\andeqn\\
\phi'':&=&\phi'\otimes 1_{\mathrm{M}_M(\mathbb C)}: C\to F_1 \otimes M_{{sK}}.
\eneq
{{From (\ref{C-Dn-2}), one has}}
$$[\psi_{0, i}\circ\phi'']=[(\pi_{i, 0})_{*0}\circ (\phi'_i)]: K_0(C) \to K_0(M_{r_isK}). 
$$
In particular $[\psi_{0, i}\circ\phi'']|_{\cal P}=[(\pi_{i, 0})_{*0}\circ (\phi'_i)]|_{\cal P}$.
Therefore, by {{Theorem \ref{UniqAtoM}}}, there is a unitary $u_{i, 0}\in M_{r_i} \otimes M_{{sK}}$ such that
$$\|\mathrm{Ad}u_{i, 0}\circ \pi_{i, 0}\circ \phi'_i(f) -  \psi_{0, i}\circ\phi''(f)\|<\min\{\sigma/{4}, \epsilon/{2}\}\rforal f\in\mathcal H.$$
Exactly the same argument shows that
there is a unitary $u_{i, 1}\in M_{r_i} \otimes M_{{sK}}$ such that
$$\|\mathrm{Ad}u_{i, 1}\circ {{\pi}}_{i, 1}\circ \phi_i{{'}}(f) -  \psi_{1, i}\circ\phi''(f)\|<\min\{\sigma/4, \epsilon/2\}\rforal f\in\mathcal H.$$

Choose {two}  path{s} of {unitaries}  ${\{u_{i, 0}(t):t\in [0,1/2]\}\subset} M_{r_i} \otimes M_{{sK}}$
such that $u_{i, 0}(0)=u_{i, 0}$ and {$u_{i, 0}(1/2)=1_{M_{r_i} \otimes M_{sK}},$} {and
$\{u_{i,1}(t):t\in [1/2, 1]\}\subset M_{r_i}\otimes M_{{sK}}$ such that
$u_{i,1}(1/2)=1_{M_{r_i}\otimes M_{sK}}$ and $u_{i,1}(1)={{u_{i,1}.}}$ }
{Put $u_i(t)=u_{i,0}(t)$ if $t\in [0,1/2)$ and $u_i(t)=u_{i,1}(t)$ if $t\in [1/2,1].$
Define ${\tilde \phi}_i: C\to I_i\otimes M_{sK}$ by
\begin{equation*}
\pi_t\circ {\tilde \phi_i}={\rm Ad}\, u_i(t)\circ \pi_t\circ \phi'_i,
\end{equation*}
where $\pi_t: I_i\otimes M_{sK}\to M_{r_i}\otimes M_{sK}$ is the {{point evaluation}} at $t\in [0,1].$}

{One has that, for each $i,$
\beq\nonumber
&&\| \pi_{i, 0}\circ {\tilde \phi}_i(f) -  \psi_{0, i}\circ\phi''(f)\|<\min\{\sigma/4, \epsilon/2\}
\andeqn\\
&&\| \pi_{i, 1}\circ {\tilde \phi_i}(f) -  \psi_{1, i}\circ\phi''(f)\|<\min\{\sigma/4, \epsilon/2\}\rforal f\in\mathcal H.
\eneq
}
For each $1\leq i\leq {k}$, let $\epsilon_i<1/2$ be a positive number such that
\beq\nonumber
&&\hspace{-0.2in}\|{\tilde \phi}_i(f)(t) - \psi_{0, i}\circ\phi''(f)\|<\min\{\sigma/{4}, \epsilon/{2}\} \rforal f\in\mathcal H \rforal t\in [0, \epsilon_i]\andeqn\\
&&\hspace{-0.2in}\|{\tilde \phi}_i(f)(t) -  \psi_{1, i}\circ\phi''(f)\|<\min\{\sigma/{4}, \epsilon/{2}\}\rforal f\in\mathcal H \rforal t\in [1-\epsilon_i, 1].
\eneq
Define ${\Phi}_i: C\to I_i\otimes M_{{sK}}$ to be
$${\Phi}_i(t)=\left\{
\begin{array}{ll}
\frac{(\epsilon_i-t)}{\epsilon_i}(\psi_{0, i}\circ\phi'') + \frac{t}{\epsilon_i}{{\tilde \phi}}_i(f)(\epsilon_i), & \textrm{if $t\in [0, \epsilon_i]$},\\
{\tilde \phi}_i(f)(t), & \textrm{if $t\in [\epsilon_i, 1-\epsilon_i]$ },\\
\frac{(t-1+\epsilon_i)}{\epsilon_i}(\psi_{1, i}\circ\phi'') + \frac{1-t}{\epsilon_i}{{\tilde\phi}}_i(f)(\epsilon_i), & \textrm{if $t\in [1- \epsilon_i, 1]$}.
\end{array}
 \right.$$
The map ${\Phi}_i$ is not necessarily a homomorphism, but it is $\mathcal H$-$\epsilon$-multiplicative; in particular, it is $\mathcal F$-$\epsilon$-multiplicative. Moreover, it satisfies the relations
\beq\label{C-D-n5}
\hspace{-0.2in}\pi_{i, 0}\circ {\Phi}_i(f) =  \psi_{0, i}\circ\phi''(f)
\andeqn
\pi_{i, 1}\circ {\Phi}_i(f) =  \psi_{1, i}\circ\phi''(f)\rforal f\in\mathcal H, i=1,..., {k}.
\eneq
{Define  $\Phi'(f): C\to C([0,1], F_2)\otimes M_{sK}$ by $\pi_{i,t}\circ \Phi'=\Phi_i,$ where
$\pi_{i,t}: C([0,1], F_2)\otimes M_{sK}\to M_{r_i}\otimes M_{sK}$ {{is}} defined by the point evaluation
at $t\in [0,1]$ (on the $i$-th summand), and define $\Phi'': C\to F_1$ by $\Phi''(f)=\phi'(f)$ for all $f\in C.$  Define
$$\phi(f)=(\Phi'(f),\Phi''(f)).$$
It follows from \eqref{C-D-n5}
that
$\phi$ is {{an}} $\mathcal F$-$\epsilon$-multiplicative positive linear map from $C$ to $D\otimes M_{sK}$}.
It follows from (\ref{C-D-nnn1}) that
\beq\label{C-D-n9}
[\pi_e\circ \phi(p)]=[\phi'(p)]=(\pi_e)_{*0}\circ K\kappa([p]) \rforal p\in {\cal P}.
\eneq
Since $(\pi_e)_{*0}: K_0(D)\to \Z^l$ is injective, one has
\begin{equation}\label{C-D-n11}
\phi_{*0}=K\kappa.
\end{equation}
{{By  \eqref{int-pre} and \eqref{pt-pre},}}
one calculates
that
$$
|(1/sK)\tau\circ \phi(h)-\gamma(\tau)(h)|<\sigma\rforal h\in {\cal H}
$$
and for all $\tau\in T(D).$
\end{proof}

\begin{lem}\label{TtoDeltaN}
{{Let $C$ be a separable \CA\, with $T({{C}})\not=\emptyset$ and
$\lambda_s(C)>0 $(see \ref{Dlambdas}).
For any $\ep>0$ and any finite subset
${\cal H}\subset C_{s.a.},$ there  {{exist}}  a finite
subset of extremal traces ${\cal T}\subseteq T(C)$  and
a continuous affine map $\lambda: T(C)\to \triangle,$
where $\triangle$ is the convex hull {{of}}  ${\cal T}$ such that}}
\beq\label{TtoD-1}
|\lambda(\tau)(h)-\tau(h)|<(\ep+(1-\lambda_s)\|h\|\rforal
h\in\mathcal H\tand \tau\in T(C).
\eneq
\end{lem}

\begin{proof}
Denote $C_1=C$ (if $C$ is unital), or $C_1=\td C.$
By {{Corollary of Theorem 5.2 of \cite{LL}}} (see also Lemma 5.1 of \cite{Tsang}),  there exists a sequence of finite dimensional simplexes $X_n$ such that
$T(C_1)=\lim_{\leftarrow n} (X_n, j_n),$ where $j_n: X_{n+1}\to X_n$ is a {{surjective}} continuous affine
map,  which also  induces the  inductive limit
$\Aff(T(C_1))=\lim_{n\to\infty}(\Aff(X_n), j_n^{\sharp}),$ where each $j_n^{\sharp}: \Aff(X_n)=\R^{r(n)}\to
\Aff(X_{n+1})=\R^{r(n+1)}$ is a unital order preserving map.
Denote by $j_n^\infty: T(C_1)\to X_n$  the continuous affine map induced by the reverse inductive limit system. It follows from the surjectivity of $j_n$ that $j_n^\infty$ is surjective.
{{Write $\partial_e(X_n)=\{x_{1,n}, x_{2,n},...,x_{m(n),n}\}.$}}
Put $T_{i,n}=(j_n^\infty)^{-1}({{x_{i,n}}}){{\not=\emptyset}}.$ Then $T_{i,n}$ is a (close) face of $T(C_1).$
Choose an extremal point $t_{i,n}\in T_{i, n}.$ Then, since $T_{i,n}$ is a face,  $t_{i,n}\in \partial_e(T(C_1)),$
$1\le i\le m(n),$ $n=1,2,....$
If $C_1=\td C,$
let $e\in C$ be a strictly positive element.
Choose $N\ge 1$ such
that
$\inf\{\tau(f_{1/N}(e)):\tau\in T(C)\}>\lambda_s-\ep/8.$

Fix any finite subset ${\cal H}$ and $\ep>0.$
We may assume
that   {{$\|h\|=1$ for all $h\in {\cal H}$}}
and $f_{1/N}(e)\in {\cal H}$ (in the non-unital case).
There is a unital positive linear map  $\gamma: \Aff(T(C_1))\to \Aff(X_{n_1})$  such that
\beq
\|j_{n_1, \infty}^\sharp\circ \gamma(\hat{f})-\hat{f}\|<\ep/8\rforal f\in {\cal H}
\eneq
for some $n_1\in \N.$
Denote by $\gamma_\sharp: X_{n_1}\to T(C_1)$ the continuous affine map
induced by $\gamma.$
It follows that, for all $f\in {\cal H}$  and $\tau\in T(C_1),$
\beq\label{810-n1}
\gamma_\sharp\circ j_{n_1}^\infty(\tau)(f)=j_{n_1}^*(\tau)(\gamma(\hat{f}))
=j_{n_1, \infty}^\sharp\circ \gamma(\hat{f})(\tau)
\approx_{\ep/16} \tau(f).
\eneq
Note that $j_{n_1}^\infty(T(C_1))={{X_{n_1}}}.$  Let $\triangle$  be the convex {{hull}} of $\{t_{1, n_1},...,t_{m(n_1), n_1}\}.$
Define a continuous affine map $\gamma': {{X_{n_1}}}\to \triangle$
by $\gamma'(x_{i,n_1})=t_{i,n_1},$ $1\le i \le m(n_1),$   and
define a continuous affine map $\lambda:=\gamma'\circ j_{n_1}^\infty: T(C_1)\to \triangle.$
For any $\tau\in T(C_1),$ write $j_{n_1}^\infty(\tau)=\sum_{i=1}^{m(n_1)}\af_i(\tau)x_{i,n_1}$
with $\af_j(\tau)\ge 0$ and $\sum_{i=1}^{m(n_1)}\af_i(\tau)=1.$
Then, for all $f\in {\cal H}$ and $\tau\in T(C_1),$ by \eqref{810-n1} twice,
\beq
\lambda(\tau)(f)=
\sum_{i=1}^{m(n_1)}\af_i(\tau) t_{i, n_1}(f)
\approx_{\ep/8} \sum_{i=1}\af_i(\tau)\gamma_\sharp\circ j_{n_1}^\infty(t_{i, n_1})(f)\\\label{810-n2}
=\gamma_{\sharp}(\sum_{i=1}\af_i(\tau)x_{i, n_1})(f)
=\gamma_{\sharp}\circ j_{n_1}^\infty(\tau)(f)
\approx_{\ep/8}\tau(f).
\eneq
This proves the case $C=C_1$ ($\lambda_s=1$).  Now consider the non-unital case. Let $\tau_\C$ be the tracial state of $\td C$ which vanishes on $C.$ If $\tau_\C\not=t_{i, n_1}$
for any $i,$ then $t_{i, n_1}\in T(C)(\subset T(C_1))$ {{for all $i,$}}   and we are done.
Otherwise, suppose that $t_{1, n_1}=\tau_\C.$
For each $\tau\in T(C),$ let $\lambda(\tau)=\sum_{i=1}^{m(n)}\af_i(\tau)t_{i, n_1}.$
Since  (recall that $f_{1/N}(e)\in {\cal H}$)
\beq
|\sum_{i=2}^{m(n)}\af_i(\tau)t_{i, n_1}(f_{1/N}(e))-\tau(f_{1/N}(e))|<\ep/4
\eneq
(by \eqref{810-n2}), we compute that, for any $\tau\in T(C),$
\beq
\sum_{i=2}^{m(n)}\af_i(\tau)\ge \tau(f_{1/N}(e))-\ep/4>\lambda_s-5\ep/8.
\eneq
It follows that $\af_1(\tau)\le 1-{{\lambda_s}}+5\ep/8.$
Let $\triangle_1$ be the convex hull of ${\cal T}:=\{t_{2,n_1},t_{2, n_2},...,t_{m(n), n_1}
\}.$
{{Define a continuous affine map  $\ld': \triangle \to \triangle_1$
by sending both $t_{1,n_1}$ and $t_{2, n_1}$ to $t_{2, n_1}$, and $t_{i, n_1}$ to itself for $i>2$.}} Define $\lambda_1: T(C)\to \triangle_1$ by $\ld_1=\ld'\circ\ld$. {{Hence,}} if $\lambda(\tau)=\sum_{i=1}^{m(n)}\af_i(\tau)t_{i, n_1}$, then
$\lambda_1(\tau)=(\af_1(\tau)+\af_2(\tau))t_{2, n_1}+\sum_{i=3}^{m(n)}\af_i(\tau)t_{i,n_1}$.
Then, since $\lambda(\tau)(f)=\sum_{i=2}^{m(n)} \af_i(\tau)(f)$ for $f\in C,$  we have,  for all $f\in {\cal H}\subset C,$
\beq\nonumber
|\lambda_1(\tau)(f)-\tau(f)|<|\af_1(\tau)t_{2,n_1}(f)|+\ep/4<
(\ep+(1-\lambda_s))\rforal \tau \in T(C).
\eneq
\end{proof}


\begin{prop}{\rm (cf. Lemma 9.4 of \cite{LnTAI})}\label{95tai}
Let  $A$ be a separable non-unital \CA. Let {{$\ep>0$,}} ${\cal F}\subset A$ be a
finite subset {{and let $e_A\in A$ be a strictly  positive element with $\|e_A\|=1$ and   $1>\eta>0,$}}.  There exists $\dt>0$ and a finite subset ${\cal G}\subset A$
satisfying the following:  
if $C$ is a separable \CA\,  with at least one tracial {{state,}}  and  if
$L: A\to C$ is  a  ${\cal G}$-$\dt$-multiplicative \cpc\, such that
$t(L(e_A))\ge 1-\eta$
for all $t\in T(C),$  then, for  any $t\in T(C),$ there exists trace $\tau$ of $A$
with $\|\tau\|\ge 1-\eta$
such that
\beq
|\tau(a)-t\circ L(a)|<\ep \tforal a\in {\cal F}.
\eneq
\end{prop}

\begin{proof}
Otherwise, there would be an $\ep_0>0$ and a finite subset ${\cal F}_0\subset A,$
a sequence of separable {{\CA s}} $C_n,$ a sequence of \cpc s $L_n: A\to C_n$ such that
\beq
\lim_{n\to\infty}\|L_n(a)L_n(b)-L_n(ab)\|=0\rforal a,b\in A,
\eneq
and a sequence $t_n\in T(C_n)$ such that
\beq\label{94tai-3}
&&t_n(L_n(e_A))\ge 1-\eta\andeqn\\
&&\inf\{\max\{|\tau(a)-t_n(L_n(a))|: a\in {\cal F}_0\}: \tau\in T(A)\}\ge \ep_0
\eneq
for all $n.$ Let $s_n=t_n\circ L_n$ be a positive linear functional with $\|s_n\|\le 1.$
Suppose that $\tau$ is a weak *-limit of $\{s_n\}.$  Then
there is a subsequence $\{n_k\}$ such that  $\tau(a)=\lim_{k\to\infty} t_{n_k}\circ L_{n_k}(a)$
for all $a\in A.$   Then one checks that $\tau$ is a trace on $A.$
By \eqref{94tai-3}, $\tau(e_A)\ge 1-\eta.$  Thus $\|\tau\|\ge 1-\eta.$
Therefore, there exists $K>1$ such that
\beq
|\tau(a)-t_{n_k}(L_{n_k}(a))|<\ep_0\rforal a\in  {{{\cal F}_0}}
\eneq
for all $k\ge K.$ A contradiction.
\end{proof}

\begin{lem}{\rm (cf. 16.12 of \cite{GLN})}\label{cut-trace}
Let $C$ be a non-unital stably finite C*-algebra {{and}} $A\in \mathcal D$
with  continuous scale.  Let $\alpha: T(A)\to T(C)$ be a continuous affine map
{{and $\Delta: C_+^{q, {\bf 1}}\setminus \{0\}\to (0,1)$  be an order preserving map.}}
\begin{enumerate}
\item\label{cut-trace-a} For any finite subset $\mathcal H\subseteq \Aff(T(C))$ {{and}}  $\sigma>0$, there is a C*-subalgebra $D\subseteq A$ and a continuous affine map $\gamma: T(D)\to T(C)$ such that $D\in\mathcal C_0$ and
$$| h(\gamma(\imath(\tau)))- h(\alpha(\tau))|<\sigma\tforal\tau\in T(A) \tand \tforal h\in \mathcal H,$$
where $\imath: T(A)\to T(D)$ is the map defined by
$\imath(\tau)=\frac{1}{\|\tau|_D\|}\tau|_D,$  {{and $\|\tau|_D\|>1-\sigma/4$ for all $\tau\in T(A).$}}

\item\label{cut-trace-b} If there are a finite subset $\mathcal H_1\subseteq {{C^{\bf 1}_+\setminus \{0\}}}$
such that
$$\alpha(\tau)(g)>\Delta(\hat{g})
\tforal g\in\mathcal H_1\tforal \tau\in T(A),$$
the affine map $\gamma$ can be chosen so that
$$\gamma(\tau)(g)>(3/4)\Delta(\hat{g})
\tforal g\in\mathcal H_1$$ for any $\tau\in T(D)$.

\item\label{cut-trace-c} If the positive cone of $K_0({\tilde C})$ is generated by a finite subset ${\cal P}$  of {{projections,}} and there is an {{order-unit-preserving}} map $\kappa: K_0({\tilde C}) \to K_0({\tilde A})$ which
maps $K_0(C)$ to $K_0(A)$ and
is compatible with $\alpha$ {{(see Definition \ref{Dcompatible})}} and
$\kappa(K_0({\tilde C})_+\setminus \{0\})\subset K_0({\tilde A})_+\setminus \{0\},$ then, for any $\delta>0$, the C*-subalgebra $D$ and $\gamma$ can be chosen so that  there are
homomorphisms
$\kappa_0: K_0(C)\to K_0(A_0)$ and $\kappa_1: K_0(C)\to K_0(D),$
where $A_0$ is a hereditary \SCA\, of $A$ which is orthogonal to $D,$  such that
 $\kappa_0(K_0({\tilde C})_+\setminus \{0\})\subset K_0({\tilde D})_+\setminus \{0\},$
 $\kappa_1(K_0({\tilde C})_+\setminus \{0\})\subset K_0({\tilde D})_+\setminus \{0\},$ $\kappa=\kappa_0+{\imath}_{*0}\circ \kappa_1,$ {where $\imath: D\to A$ is the embedding,} and
\begin{equation}\label{june2-nn1}
|\gamma(\tau)(p)-\tau(\kappa_1([p])|<\dt\rforal p\in {\cal P}{\tand\,\,\, \tau\in T(D)}.
\end{equation}


\item\label{cut-trace-d} Moreover, in addition to (\ref{cut-trace-c}), if $A\cong A\otimes U$ for some infinite dimensional UHF-algebra, for any given positive integer $K$, the C*-algebra $D$ can be chosen so that  $D=M_K(D_1)$
for some $D_1\in {\cal C}$ ($D_1\in {\cal C}_0$) and $\kappa_1=K\kappa_1',$
where $\kappa_1': K_0(C)\to K_0(D_1)$ is a strictly positive homomorphism. {Furthermore,
$\kappa_0$ can also be chosen to be strictly positive.}
\end{enumerate}

\end{lem}

\begin{proof}
Write $\mathcal H=\{h_1, h_2, ..., h_m\}.$  We may assume that $\|h_i\|\le 1,$ $i=1,2,...,m.$ Choose $f_1, f_2, ..., f_m\in A_{s.a.}$ such {{that}} $\tau(f_i)=h_i(\alpha(\tau))$ for all $\tau\in T(A)$ and $\|f_i\|\le 2,$ $i=1,2,...,m$
(see 9.2 of \cite{LnTAI}). Put ${\cal F} =
\{f_1,f_2,...,f_m\}.$
{{Fix $\sigma>0.$ If ${\cal H}_1\subset C_+^{\bf 1}\setminus \{0\}$ is given,
let $\sigma_0=\inf\{\Delta(\hat{h}): h\in {\cal H}_1\}>0.$
Put $\sigma_1:=\min\{\sigma, \sigma_0\}.$ }}
Since $A\in {\cal D}$ has continuous scale (recall \ref{DD0}), there {{exists}} a strictly
positive element $e\in A$ with $\|e\|=1,$ mutually orthogonal
\SCA s $A_{n,0}$ and $A_{n,1}$ of $A$ with $A_{n,1}\in {\cal C}_0$ and two
\cpc s $\phi_{n,0}: A\to A_{n,0}$ and $\phi_{n,1}: A\to A_{n,1}$ such that
\beq\label{1612gln-c1}
\lim_{n\to\infty}\|x-\diag(\phi_{n,0}(x),\phi_{n,1}(x))\|=0\rforal x\in A,\\\label{1612gln-c2}
\lim_{n\to\infty}\|\phi_{n,i}(ab)-\phi_{n,i}(a)\phi_{n,i}(b)\|=0\rforal a, b\in A,\,\,\, i=0,1,\\\label{1612gln-c3}
\lim_{n\to\infty}\sup\{d_\tau(e_{n,0}): \tau\in T(A)\}=0\andeqn\\\label{1612gln-c4}
t(f_{1/4}(\phi_{n,1}(e)))\ge 1-\sigma_1/32\rforal t\in T(A_{n,1}),
\eneq
where $e_{n,0}\in A_{n,0}$ is a strictly positive element.
Note that \eqref{1612gln-c4} implies also
$\lambda_s(A_{n,1})\ge 1-\sigma_1/32.$  Put $D:=A_{n,1}.$
{{By \eqref{1612gln-c3}, we may assume that $\|\tau|_D\|>1-\sigma_1/4$ for all
$\tau\in T(A).$}}
It follows from \ref{95tai} that, for some large $n,$ for each $\tau\in  T(A_{n,1}),$ there is {{a}}
trace $\gamma'(\tau)$ with $1\ge \|\gamma'(\tau)\|\ge 1-\sigma_1/32$
such that
\begin{equation}\label{june2-n2}
|\tau(\phi_{n,1}(a))-\gamma'(\tau)(a)|<\sigma_1/32\rforal a\in {\cal F}.
\end{equation}
%
Put $\gamma''(\tau)=(1/\|\gamma'(\tau)\|)\gamma'(\tau).$
Then $\gamma''(\tau)\in T(A)$ and
\beq\label{june2-n2+}
|\tau(\phi_{n,1}(a))-\gamma''(\tau)(a)|<\sigma_1/16\rforal  a\in {\cal F}.
\eneq

Applying  \ref{TtoDeltaN}
one obtains $t_1, t_2,...,t_n\in \partial_e{T(D)}$ and
a continuous affine map $\lambda: T(D)\to \triangle,$
the convex hull of $\{t_1,t_2,...,t_n\},$
such that
\begin{equation}\label{june6-n1}
|\tau(\phi_{n,1}(a))-\lambda(\tau)(\phi_{n,1}(a))|<\sigma_1/16\rforal \tau\in T(D)
\end{equation}
and $a\in {\cal F}.$
Define $\lambda_1:  \triangle\to T(A)$ by
$\lambda_1(t_i)=\gamma''(t_i),\,\,\, i=1,2,...,m.$
Define $\gamma=\af\circ \lambda_1\circ \lambda.$
 Then
\begin{eqnarray*}
h_j(\gamma(\imath(\tau)))&=&h_j(\alpha\circ \lambda_1\circ\lambda(\imath(\tau)))
=\lambda_1\circ\lambda(\imath(\tau))(f_j)\\
&\approx_{\sigma_1/16}& \lambda(\imath(\tau))(\phi_{n,1}(f_j))
\approx_{\sigma_1/16}\imath(\tau)(\phi_{n,1}(f_j))\\
&\approx_{\sigma_1/8}& \tau(f_j)= h_j(\alpha(\tau)),
\end{eqnarray*}
and this proves (\ref{cut-trace-a}) by
letting $L=\phi_{n,1}.$
Note that it follows from the construction that $\gamma(\tau)\in\alpha(T(A))$, and hence  (\ref{cut-trace-b})
can also be arranged.

Part  (\ref{cut-trace-c}) follows easily by choosing $A_0=A_{n,0},$ $D=A_{n,1},$ $\kappa_0=[\phi_{n,0}]\circ [\kappa]$
and $\kappa_1=[\phi_{n,1}]\circ [\kappa]$ with sufficiently large $n$ as demonstrated above
{{(recall that $K_0(\td C)_+$  is assumed to be finitely generated).}}

Then  (\ref{cut-trace-d}) follows straightforwardly as part (\ref{cut-trace-c}) except
the ``Furthermore" part.

{To see the ``Furthermore" part, we note that we may choose $D\subset A\otimes 1_U.$ Choose a projection
$e\in U$ such that
$$
0<t_0(e)<\dt_0<\dt-\max\{|\gamma(\tau)(p)-\tau(\kappa_1([p])|: p\in {\cal P} \tand\, \tau\in T(D)\},
$$
where $t_0$ is the unique tracial state of $U.$ We then replace $\kappa_1$ by  $\kappa_2: K_0(A)\to K_0(D_2),$
where $D_2=D\otimes (1-e)$ and $\kappa_2([p])=\kappa_1([p])\otimes [1-e].$
Define $\kappa_3([p])=\kappa_1([p])\otimes [e].$
Define $A_0'$ to be the hereditary \SCA\, generated by $A_0\oplus D\otimes e.$ Then let $\kappa_4: K_0(C)\to K_0(A_0')$ be defined by
$\kappa_4=\kappa_0+[\imath]\circ \kappa_3,$ where $\imath: D\otimes e\to A\otimes U\cong A$ is the embedding.
We then replace $\kappa_0$ by $\kappa_4.$ Note that, now, $\kappa_4$ is strictly positive.}
\end{proof}

\begin{cor}\label{CExtTraceC-D}
 {{Let $C=C_0$ or $C=C(\T, \td C_0)^o$ for some $C_0\in {\cal C}_0.$}}
 Let $\Delta: C_+^{q, {\bf 1}}\setminus\{0\}\to (0, 1)$ be an order preserving map. Let $\mathcal F, \mathcal H\subseteq C$ be finite subsets, and let $1>\sigma, \epsilon>0$.
There exist  a finite subset $\mathcal H_1\subseteq C_+^{\bf 1}\setminus \{0\}$, $\delta>0$,
{a finite subset ${\cal P}\subset K_0(C)$}and a positive integer $K$ such that for any continuous affine map $\gamma: T(B)\to T(C)$ satisfying
$$\gamma(\tau)(h) > \Delta(\hat{h})\rforal h\in \mathcal H_1\tforal \tau\in T(B),$$
where $B$ is a non-unital  C*-algebra in $\mathcal {\cal D}$ with continuous scale,
any positive \hm\, $\kappa: {{K_0(\td C)\to K_0(\td B)}}$ {{with $\kappa(K_0(C))\subset K_0(B),$}}
$\kappa([1_{\tilde C}])=s[1_{\tilde B}]$
for some integer $s \ge 1$ satisfying
$$|\rho_C(x)(\gamma(\tau))-{(1/s)}\tau(\kappa(x))|<\delta\tforal \tau\in T(B)$$
and for all $x\in {\cal P},$
there is an $\mathcal F$-$\epsilon$-multiplicative positive linear map $\phi: C\to M_{{sK}}(B)$ such that
\beq
&&\phi_{*0}={K}\kappa\tand\\
&&|(1/(sK))\tau\circ\phi(h)-\gamma(\tau)(h)|<\sigma\tforal h\in\mathcal H\tand \tau\in T(B).\eneq

In the case that $C\in\mathcal C_0$, the map $\phi$ can be chosen to be a homomorphism.
\end{cor}

\begin{proof}
This is the combination of  (1),  (2) and (3) of \ref{cut-trace}
 and   \ref{ExtTraceC-D}
 as well as Theorem 5.7 of \cite{GLrange}.

 Let ${\cal F},$  ${\cal H},$ $\sigma$ and  $\ep>0$   be given as in the statement above.
 Choose
 $\mathcal H_1\subseteq C_+^{\bf 1}\setminus \{0\}$, $\delta>0$,
{a finite subset {{${\cal P}\subset K_0(C)$} and}} a positive integer $K$ be required by \ref{ExtTraceC-D}
for ${\cal F},$ ${\cal H},$ $\sigma/2,$  $\ep/2$ and $(3/4)\Delta$ (instead of $\Delta$).
By applying the combination of (1), (2) and (3) of \ref{cut-trace} for $\af=\gamma,$
$\sigma/2sK$ (in place of $\sigma$) and $\dt/2sK$ (in place $\dt$), we obtain
$A_0,$ $D,$ $\gamma_1$ (as $\gamma$), $\kappa_0$ and $\kappa_1$ which satisfy
the conclusions of (1), (2) and (3) of \ref{cut-trace}.
Then, applying \ref{ExtTraceC-D},  we obtain an $\mathcal F$-$\epsilon$-multiplicative positive linear map $\phi_1: C\to M_{{sK}}(D)
\subset M_{sK}(B)$ such that
\beq
&&(\phi_1)_{*0}={K}\kappa_1\tand\\
&&|(1/(sK))\tau\circ\phi(h)-\gamma_1(\tau)(h)|<\sigma/2sK\tforal h\in\mathcal H\tand \tau\in T(B).\eneq
Note also {{that}} the conclusion (1) mentioned above implies that $\|\tau|
_D\|>1-\sigma/4sK$ for all $\tau\in T(A),$ and
\beq
|h(\gamma_1(\imath(\tau)))- h(\gamma(\tau))|<\sigma/2\tforal\tau\in T(A) \tand \tforal h\in \mathcal H.
\eneq
Note that $K_1(C_0)=\{0\},$ by \ref{Rkt2}, $K_0(C)=K_0(C_0).$
Moreover, the map $\pi_1:  C\to C_0$ (evaluating at $1$) induces an isomorphism $(\pi_1)_*: K_0( C)\to K_0(C_0)$. 
Thus there is a \hm\, $\kappa_0': K_0(C_0)\to {{K_0(A_0)=K_0(A)}}$
such that the induced map $\td\kappa_0': K_0(\td C_0)\to K_0({{\widetilde{A_0}}})$ is strictly positive and
$\kappa_0'\circ {{(\pi_1)_*}}=\kappa_0.$
If $C=C_0,$ then we may view that $\pi_1$ is the identity map.
By applying Theorem 5.7 of \cite{GLrange}, there is also an $\mathcal F$-$\epsilon$-multiplicative positive linear map $\phi_0':
C_0\to M_{sK}(A_0)$ such that
$[\phi_0]=K\kappa_0'.$  Recall that $A_0\perp D.$
So $\|\tau|_{A_0}\|<\sigma/4sK$ for all $\tau\in T(A).$  Define
$\phi:=\phi_0\circ \pi_1\oplus \phi_1.$  It is ready to check that $\phi$ meets the requirements.


\end{proof}
%
%
%

\section{Homotopy lemmas}
We will retain notation introduced in \ref{DDo}.

\begin{lem}\label{Hext1}
Let $A\in {\cal M}_1$ {{be}} as (3) in Remark 4.32 of \cite{GLrange}
with continuous scale and
let $B$   be  a  separable  simple \CA s
which has the form $B=B_0\otimes U$ for some $B_0\in {\cal M}_1$ with continuous scale
{{for an infinite dimensional UHF-algebra $U.$}}
Suppose that there are $\kappa\in KL(C(\T, {\tilde A})^o,B)$ and
an affine continuous map
$\kappa_T: T(B)\to  T(C(\T, {\tilde A})^o)$ such that $\kappa$ and $\kappa_T$ are compatible (see Definition \ref{Dcompatible}),
where $\kappa_T(T(B))$ lies in a compact convex subset of $T_f(C(\T, {\tilde{D}})^o)$ {{(see \ref{DTtilde} for notation of $T_f$).}}
Then,
there exists a  sequence of approximate multiplicative  \cpc s $\phi_n: C(\T,{\tilde A})^o\to B$ such that
\beq
&&[\{\phi_n\}]=\kappa\tand\\
&&\lim_{n\to\infty}\sup \{|\tau\circ \phi_n(a)-\kappa_T(\tau)(a)|: \tau\in T(B)\}=0\tforal a\in C(\T,{\tilde A})^o_{s.a.}.
\eneq
\end{lem}

\begin{proof}
The proof is similar to that of \ref{ExtAB}  combining with  {{that of}}  \ref{ExtTBA}.
Let $1>\ep>0, 1>\sigma>0,$ ${\cal F}\subset
A^\T= C(\T,{\tilde A})^o$  and ${\cal H}\subset
 C(\T,{\tilde A})^o_{s.a}$ be
 finite subsets. \Wlog, we may assume that
$${\cal F}=\{1_{C(\T)}\otimes a: a\in {\cal F}_A\}\cup
\{1_{C(\T)}-z, 1_{C(\T)}-z^*\},$$
where ${\cal F}_A\subset A^{\bf 1}$ is a finite
subset and $z\in C(\T)$ is the identity function on the unit circle,
${\cal H}
=\{1_{C(\T)}\otimes a: a\in {\cal H}_A\}\cup \{ f\in {\cal H}_\T\},$ where ${\cal H}_A\subset A^{\bf 1}_{s.a.}$ and
${\cal H}_\T\subset C_0(\T\setminus \{1\})_{s.a.}$ are finite subsets.

In what follows, we will identify $A$ with $1_{C(\T)}\otimes A$ {{and $C_0(\T\setminus \{1\})$ with $C_0(\T\setminus \{1\})\otimes 1_{\tilde A}$, both as {{\SCA s}}
of $A^\T=C(\T, {\tilde A})^o$.}}


Fix a finite subset ${\cal P}\subset \underline{K}(A^\T).$

Choose $\dt>0$ and finite subset ${\cal G}\subset
A^\T$ so that $[L]|_{{\cal P}}$ is well defined
for any
${\cal G}$-$\dt$-multiplicative \cpc\, $L$ from $
A^\T.$
We may assume that $\dt<\ep$ and ${\cal F}\cup {\cal H}\subset {\cal G}.$
We may further assume that
$$
{\cal G}={\cal G}_A\cup {\cal G}_\T,
$$
where ${\cal G}_A\subset A^{\bf 1}$ and ${\cal G}_\T\subset C_0(\T\setminus\{1\})^{\bf 1}$ are finite subsets.
Since both $A$ and $B$ have continuous scales,
$T(A)$ and $T(B)$ are compact (5.3 of \cite{eglnp}).

Choose  $a_0\in A_+$ such that $\|a_0\|=1$ and
\beq\label{11Ext1-n2}
d_\tau(a_0)<\min\{\sigma, \dt\}/4\rforal \tau\in T(A).
\eneq
 Let  $e_0\in A^\T$ be a strictly positive element of $A$ with $\|e_0\|=1$
such that $\tau(e_0)>15/16$ for all $\tau \in T(A)$ (see the end of \ref{DDo}).
\Wlog,  
 we may assume {{that,}} for some large $n\ge 1,$
\beq
{\cal G}_A\subset  E_n\oplus C_n\oplus D_n,
\eneq
where $C_n, D_n\in {\cal C}_0$ and $E_n$ {{are}} as in {{Theorem 4.34 of \cite{GLrange}.}} 
We write $E_n'=E_n\oplus D_n.$
As in {{4.34 of \cite{GLrange},}} 
we assume that ${\rm ker}\rho_{A}\cap (\iota_n)_{*0}(K_0(C_n))=\{0\}$
and ${\rm ker}\rho_{E_n'}=K_0(E_n')$ and $(\iota_n)_{*0}(K_0(E_n'))\subset {\rm ker}\rho_{A},$
where $\iota_n: {{E_n'}}
\oplus C_n\to A$ is the embedding. Moreover, we may also assume $e_{00}\oplus e_{01}\le e_0$
and
\beq\label{2020-718-n1}
0\not=e_{00}\le e_0-e_{01}\lesssim a_0\andeqn d_\tau(e_{01})>1-\min\{\sigma,\dt\}/64\rforal \tau\in T(A),
\eneq
where $e_{00}\in (E_n')_+$ and $e_{0,1}\in (C_n)_+$ are strictly positive elements,
and assume that ${\cal G}_A={\cal G}_{eA}'\cup {\cal G}_{cA}',$
where ${\cal G}_{eA}'\subset E_n'$ and ${\cal G}_{cA}'\subset C_n$ are finite subsets.
Let ${E_n'}^\T$ be the \SCA\, generated by $C_0(\T\setminus \{1\})\otimes 1_{\td A}$ and
$E_n'$, and let $C_n^\T$ be the \SCA\, generated by $C_0(\T\setminus \{1\})\otimes 1_{\td A}$
and $C_n,$ respectively. Then {{${E_n'}^\T\cong C(\T, {\td E}'_n)^o$ and $C_n^\T\cong C(\T, \td C_n)^o$, when we identify $1_{{\td E}'_n}=1_{\td C_n}=1_{\td A}$.}}
We may further assume, \wilog, that there are finite
{{subsets}} ${\cal P}_0\subset \underline{K}(C(\T, {\tilde E_n'})^o)$ and
${\cal P}_1\subset\underline{K}(C(\T, {\tilde C_n})^o)$ such that
${\cal P}\subset [\imath_n]({\cal P}_0\cup {\cal P}_1),$
where $\imath_n: C(\T, {\tilde E_n'})^o+ C(\T, {\tilde C_n})^o\to C(\T, {\tilde A})^o$
is the embedding.
Since $K_0(C(\T, {\tilde C_n})^o)$ is finitely generated (see \ref{Rkt2}), we may assume that
${\cal P}_1\cap K_0(C(\T, {\tilde C_n})^o)$ generates $K_0(C(\T, {\tilde C_n})^o).$
Let  $e_c\in C_n$  be a strictly positive element of $C_n$ with $\|{{e_c}}\|=1.$
From \ref{Rkt2}, {{we have}} 
$$
K_0(C(\T, {\tilde E_n'})^o)=K_0(E_n')\oplus {\boldsymbol{\bt}}(K_1(E_n')).
$$

Note here we assume {{that,}} as constructed in {{4.34 of \cite{GLrange},}} 
$K_0(E_n')={\rm ker}\rho_{E_n'}.$  It follows  from
\ref{LKT2}
that $K_0(C(\T, {\tilde E_n'})^o)={\rm ker}\rho_{C(\T, {\tilde E_n'})^o}.$
We may also assume that ${\cal P}_1\cap K_0(C_n)={\cal P}_{1,0}\sqcup {\cal P}_{1,1},$
where ${\cal P}_{1,0}\subset   {\boldsymbol{\bt}}(K_1(C_n))$ and
${\cal P}_{1,1}\subset K_0(C_n).$ Note that ${\cal P}_{1,0}\subset {\rm ker}\rho_{C_n}.$

Denote by $\Psi_0: A^\T\to C(\T, {\tilde E_n'})^o $ and $\Psi_1: {{A^\T}}\to C(\T, {\tilde C_n})^o$ two \cpc s
which are
${\cal G}$-$\dt/16$-multiplicative  and {{such that}}
\beq\label{ExtAB-10+}
\|\Psi_0(b)-b\|<\dt/16\rforal b\in {\cal G}_{eA}'\andeqn \|\Psi_1(c)-c\|<\dt/16\rforal c\in {\cal G}_{cA}'.
\eneq
Let
{{$
{\cal P}_2=[\imath_n]({\cal P}_0\cup {\cal P}_1)
$.}}
Recall  {{that}} we assume  
that ${\cal P}_1$ contains a generating set of $K_0(C(\T,{\tilde C_n})^o).$

We may assume that, for some $m_0\ge 1,$
$
{\cal P}_2\subset \underline{G}\oplus {\boldsymbol{\bt}}(\underline{G}),
$
where
$$
\underline{G}=
K_0(A){{\oplus}} K_1(A){{\oplus}}\bigoplus_{j=1}^{m_0} {{(K_0(A, \Z/j\Z)\oplus K_1(A, \Z/j\Z))}}.
$$
Moreover,   we may assume that, for some {{$m\ge m_0,$}}
$m!x=0$ for all $x\in {\rm Tor}(K_0(A^\T))\cap {\cal P}.$
Let $G_{0, {\cal P}}$ be the subgroup generated by $K_0(A^\T)\cap {\cal P}_2.$
We may write $G_{0, {\cal P}}:=F_0\oplus F_{00}\oplus G_0,$ where $F_{0}$
and $F_{00}$ are   free,   $F_0$
is generated by $[\imath_n]({\cal P}_{1,1}),$
{{$F_{00}\oplus G_0$ is   generated by $[\imath_n]({\cal P}_0\cup {\cal P}_{1,0})\cap K_0(A^\T)$ }}and
$G_0$ is a finite group.
In particular, $m!x=0$ for all $x\in G_0.$
Moreover, $F_0\subset (\imath_n)_{*0}(K_0(C(\T, {\tilde C_n})^o).$


Choose $0<\dt_1<\dt$  and finite subset ${\cal G}_3\subset {{A^{\T}}}$ such that
$[L']|_{{\cal P}_2}$ is well defined for
any ${\cal G}_3$-$\dt_1$-multiplicative \cpc\, {{$L'$}} from ${{A^{\T}}}$.
We assume that {{${\cal G}\subset {\cal G}_3.$}}

Note that, by  {{Theorem 3.3 of \cite{GLrange}}} 
and 14.7 of \cite{GLII}, $A^\T$ satisfies the assumption of 9.8 of \cite{GLII}.
 It follows from {{Theorem 3.4 of \cite{GLrange}}} 
 that
there exists a ${{{\cal G}_3}}$-$\dt_1/4$-multiplicative \cpc\,
$L: A^\T\to B\otimes M_N$
 for some
integer $N\ge 1$ such that
\beq\label{11Ext1-2-Gong}
[L]|_{{\cal P}_2}=
{{\kappa|_{{\cal P}_2}.}}
\eneq
\Wlog, we may assume  ${\cal G}_3\subset (A^\T)^{\bf 1}.$
%
%
Define $\Delta:
{{(C_n^\T)^{q, {\bf 1}}_+}}\setminus \{0\}\to (0,1)$ by
\beq\label{2020-718-n18}
\Delta(\hat{h})=(1/2)\inf\{{\widehat{\imath_n(h)}}(\kappa_T(\tau)):\tau\in T(B)\}
\eneq
for $h\in
{{(C_n^\T)^{\bf 1}_+}}\setminus \{0\}.$
Note that since $\imath_n(h)\in (A^\T)^{\bf 1}_+\setminus \{0\},$ and
$\kappa_T(T(B))$ lies in a compact subset of $T_f(A^\T),$
$\Delta(\hat{h})>0$ for all $h\in
(C_n^\T)^{\bf 1})_+\setminus \{0\}.$ Put ${\cal H}'=\{1\otimes \Psi_1(h): h\in {\cal H}_A\}\cup{\cal H}_\T.$

Let ${\cal H}_1\subset (C^\T)^{\bf 1}_+,$
$\dt_0$ (in place of $\dt$), ${{\cal P}'_{1}}\subset K_0(C(\T,{\tilde C}_n)^o)$
(in place of ${\cal P}$) and $K$ be given by {{Corollary \ref{CExtTraceC-D}}} for
$C=C_n^\T,$ ${\cal G}_3$ (in place of ${\cal F}$), ${\cal H}',$
$s=1$ and  $\sigma/4$ (in place of $\sigma$)  and $\dt_1/4$ (in place of $\ep$) above.
%



Let ${{\cal Q}}\subset \underline{K}(B)$ be a finite subset
which contains $[L]({\cal P}_2).$
We assume that
\beq\label{11Ext1-3+2}
{\cal Q}\subset K_0(B){{\oplus}} K_1(B){{\oplus}}\bigoplus_{i=0,1}\bigoplus_{j=1}^{m_1}K_i(B, \Z/j\Z)
\eneq
for some $m_1\ge 2.$ Moreover, we may assume
that ${{m_1!}}x=0$ for all $x\in {\rm Tor}(G_{0,b}),$ where $G_{0,b}$ is the subgroup
generated by ${\cal Q}\cap K_0(B).$ \Wlog, we may assume that
$m|m_1.$ Choose an integer $m_2$ such that
$m_1K|m_2$ and
\beq
m/m_2<\sigma\dt_0/16.
\eneq
Let $p_1, p_2,...,p_l\in M_r(\td C_n)$ be projections which generate
$K_0(\td C_n)_+.$ Let ${{\bar p_i}}\in M_r(\C\cdot 1_{\td C_n})$ be scalar projections  with rank
$R_i\ge 1$
such that $o([p_i])=[p_i]-[{\bar p}_i]\in K_0(C_n)$ {{for}} $i=1,2,...,k.$
Since $(\kappa, \kappa_T)$ is compatible,
{{$R_i+\rho_{B}(\kappa(o([p_i])))(s)>0$
for all $s\in T(B).$ }}

{{We further assume ${\cal P}'_1\subset \{[p_i]-[{\bar p}_i]: i=1,2,\cdots,l\}$ (note that $K_0(C(\T,{\tilde C}_n)^o)=K_0(C_n)$, since $K_1(C_n)=0$).}}
 Set $R:=\max\{R_i: 1\le i\le l\}$ and
\beq\label{2020-713-11n1}
\eta_1:=\min\{\inf |R_i+{{\rho_{B}}}(\kappa(o([p_i])))(s)|:
s\in {{T(B)}}\}: 1\le i\le l\}.
\eneq
Let $b_0\in {{B_+}}$ with $\|b_0\|=1$ such that
\beq\label{11Ext1-3}
d_\tau(b_0)<\min\{{{\sigma\dt_0/(rR)}}, \dt_1, \eta_1/R\}/16(N+1)(K+1)m_2!\rforal \tau\in T(B).
\eneq
Let $e_b\in   B\otimes M_N$
be a strictly positive element of $B\otimes M_N$
such that
\beq\label{11Ext1-3++}
\tau(e_b)>7/8\rforal \tau\in T(B\otimes M_N).
\eneq


Let ${\cal G}_b\subset B\otimes  M_N$ be a  finite subset and
$1/2>\dt_2>0$ be such that
$[\Phi]|_{\cal Q}$ is well defined for any ${\cal G}_b$-$\dt_2$-multiplicative \cpc\,  $\Phi$ from
$B\otimes M_N.$
Note that $B_0\in {\cal M}_1$ and $B=B_0\otimes U.$ By {{Theorem 4.31 of \cite{GLrange}}}
(see (2) of Remark {{4.32 of \cite{GLrange}}}), 
{{and Theorem \ref{Misothm}, we may assume that}} there are simple \CA s $E_b'$ and $D_b$ with continuous scales
such that  $K_0(E_b')={\rm ker}\rho_{E_b'}$ and
$K_0(D_b)$ is torsion free, ${\rm ker}\rho_{D_b}=\{0\},$ ${\rm ker}\rho_B\cap j_{*0}(K_0(D_b))=\{0\},$
where $j: D_b\to B$ is the embedding, and
%
there are  ${\cal G}_b$-$\dt_2$-multiplicative \cpc s $\phi_{0,b}: B\otimes M_N\to
E_b\subset \overline{\phi_{0,b}(B\otimes M_N)( B\otimes M_N)\phi_{0,b}(B\otimes M_N)}$
and $\psi_{0,b}: B\otimes M_N\to D_b\subset B\otimes M_N$
with $E_b:=(E_b'\otimes U)$ being orthogonal  to  $E_b'':={{M_{N(m_2!)}(D_b)}}$ 
such that
\beq\label{11Ext1-4.1}
&&\hspace{-0.8in}\|b-\diag(\phi_{0,b}(b),\overbrace{\psi_{0,b}(b), \psi_{0,b}(b),..., \psi_{0,b}(b)}^{{N(m_2!)}})\|<\min\{\dt_2, {\ep\over{16}}, {\eta_1 \sigma\dt_0\over{16R^2}}\}\rforal b\in {\cal G}_b\\\label{11Ext1-4nn}
&&\hspace{-0.2in}\andeqn \phi_{0,b}(e_b)\lesssim b_0\andeqn d_\tau(\psi_{0,b}(e_b))>(N-d_\tau(b_0))/Nm_2!\rforal \tau\in T(B).
\eneq
Note that $K_1(D_b)=\{0\}.$ Moreover, {{with $\Psi_{0,b}$ being  direct sum of $N(m_2!)$ copies
of ${{\psi_{0,b}}}$}},
we may also assume
that
\beq\label{2020-718-100}
&&[{\rm id}_B]|_{\cal Q}=[\phi_{0,b}]|_{\cal Q}+[\Psi_{0,b}]|_{\cal Q}
\\\label{11Ext1-4+G}
&&(m_2)![\psi_{0,b}]|_{{\rm Tor}(G_{0,b})}=0\andeqn
(m_2)![\psi_{0,b}]|_{{\cal Q}\cap K_i(B, \Z/j\Z)}=0,\,\,\, j=2,3,...,m_1,\\
&&{[}\psi_{0,b}{]}|_{{\cal Q}\cap {\rm ker}\rho_B}=0.
\eneq

Therefore
\beq\label{11Ext1-5n}
&&[\phi_{0,b}{]}|_{{\cal Q}\cap{\rm ker}\rho_B}=[{\rm id}_B{]}|_{{\cal Q}\cap {\rm ker}\rho_B},\,\,
{[}\phi_{0,b}{]}|_{{\cal Q}\cap K_1(B)}=[{\rm id}_{B}{]}|_{{\cal Q}\cap K_1(B)}\andeqn\\\label{11Ext1-5n+2}
&&{[}\phi_{0,b}{]}|_{{\cal Q}\cap K_i(B, \Z/j\Z)}=[{\rm id}_B{]}|_{{\cal Q}\cap K_i(B, \Z/j\Z)},\,\,\, j=2,3,...,m_1.
\eneq
Since $(\kappa, \kappa_T)$ is compatible, by \eqref{11Ext1-2-Gong}
we may assume
that $[L]({\cal P}_0\cup {\cal P}_{1,0})\subset  {\rm ker}\rho_{B}.$
Therefore we may further  assume that $[L]({\cal P}_0\cup {\cal P}_{1,0})\subset j_{*0}(K_0(E_b)).$

Let $G_{\cal P}$ be the subgroup
generated by ${\cal P}$ and let
$\kappa'=\kappa- \phi_{0,b}\circ [L{]}$
be defined on $G_{\cal P}.$
Then, by {{\eqref{11Ext1-2-Gong}}}, \eqref{11Ext1-5n} and \eqref{11Ext1-5n+2},  since $\kappa$ preserves
the order,
we compute that
\beq\label{11Ext1-101}
&&\kappa'|_{G_0\oplus F_{00}}=0,
\kappa'|_{{\cal P}\cap K_1(A)}=0\andeqn\\
&&\kappa'|_{{\cal P}\cap K_i(A, \Z/j\Z)}=0,\,\,\, j=2,3,...,m.
\eneq
Let $\zeta:=\kappa'\circ {\iota_n}_{*0}: {{K_0(C_n^\T)\to K_0(B)}}.$  By \eqref{11Ext1-5n},
$(\kappa'\circ \iota_{*0})|_{{\rm ker}\rho_{C_n^\T}}=0.$
So $\zeta$ factors through $K_0(C_n).$
Then, by {{\eqref{11Ext1-2-Gong},}} \eqref{11Ext1-3} and \eqref{2020-713-11n1},  for all {{$s\in T(B),$}}
\beq
R_i+{{\rho_B}}(\zeta(o([p_i]))(s)&=&R_i+{{\rho_B}}(\kappa(o([p_i])(s)-
{{\rho_B}}([\phi_{0,b}]\circ [L](o([p_i]))(s)\\
&=&R_i+{{\rho_B}}(\kappa(o([p_i])(s)-
{{\rho_B}}([\phi_{0,b}](\kappa(o([p_i])))(s)\\
&\ge& R_i+{{\rho_B}}(\kappa(o([p_i])(s)- (\eta_1/R16) (2R)\\
&>& R_i+{{\rho_B}}(\kappa(o([p_i])(s)- \eta_1>0
\eneq
for $1\le i\le l.$
Since $\{p_1,p_2,..., p_l\}$ generates $K_0(\td C_n)_+,$
this implies that the unital extension {{$\zeta^\sim : K_0(\td C_n^\T)\to K_0(\td B)$}}  
is strictly positive.

 Note that $K|m_2.$
 Note {{also,}} by  \eqref{2020-718-100},  {{that}} $\zeta=[\Psi_{0,b}]\circ [L].$
Thus $\zeta_0:=(1/K)\zeta$ is a \hm\, from $K_0(C_n)\to K_0(B).$
Note that $B$ has continuous scale. Let $b^K\in B_+$
with $d_\tau(b^K)=1/K$ for all $\tau\in T(B).$ Put
 $B_1=\overline{b^KBb^K}.$
We may view {{$\zeta_0:K_0(C_n^{\T})\to K_0(B_1)$}} (recall {{$C_n^{\T}=C(\T, C_n)^o$).}}
We then extend {{$\zeta_0^\sim: K_0({{\td C_n^{\T}}})\to K_0(\td B_1)$ with $\zeta_0^\sim([1_{{\td C_n^{\T}}}])=[1_{\td B_1}].$}}
Then
$\zeta_0^\sim$ is strictly positive.


Let ${\iota_n}_T: T_f(A^\T)\to \R\cdot T_f(C_n^\T)$
be defined by
${\iota_n}_T(\tau)(c)=\tau(c)$ for all $c\in  C_n^\T.$
Define $\af:={\iota_n}_T\circ \kappa_T: T(B)\to T_f(C_n^\T).$
By \eqref{2020-718-n1}, $\inf\{\|{\iota_n}_T(\kappa_T(\tau))\|: \tau\in T(B)\}\ge 1-\sigma/64.$
Recall, for all $\tau\in T(B),$
\beq\label{2020-718-19}
\rho_{B}(\kappa({\iota_n}_{*0}(o([p_i])))(\tau)=\rho_{{A^\T}}({\iota_n}_{*0}(o([p_i]))(\kappa_T(\tau))
=\rho_{{C^\T}}(o([p_i])({\iota_n}_T(\kappa_T(\tau)).
\eneq
Let $\bt_1:T(B_1)\to T(B)$ be the   affine homeomorphism so that $\bt_1^{-1}(\tau)(b)=K\tau(b)$
for all $b\in B_1$ and
for all $\tau\in T(B).$  Recall that $\|\tau|_{B_1}\|=d_\tau(b^K)=1/K$
for all $\tau\in T(B).$
It follows from  \eqref{2020-718-19}
  {{\eqref{11Ext1-3} and \eqref{11Ext1-4nn}}}
that (see \ref{Rextendable})
\beq
|\rho_{{C^\T}}(o([p_i])(\af^\sim(\bt_1(\tau)))-\rho_{B_1}({{\zeta_{0}}}^\sim(o([p_i]))(\tau)|<{{\dt_0/2}}\rforal \tau\in T(B_1).
\eneq
Also, by \eqref{2020-718-n18}, for $\tau\in T(B_1),$
\beq
\af^\sim(\bt_1(\tau))(h)={\iota_n}_T(\kappa(\bt_1(\tau)))(h)\ge \Delta(\hat{h})\rforal h\in {\cal H}_1.
\eneq
Applying Corollary \ref{CExtTraceC-D} (note that ${\cal P}'_1\subset \{[p_i]-[{\bar p}_i]: i=1,2,\cdots,l\}$), we obtain a ${\cal G}_3$-$\dt_1$-multiplicative \cpc\,
$L_1: C_n^\T\cong C(\T, {\tilde C}_n)^o\to M_K(B_1)\cong B$ such that
$[L_1]=\kappa'\circ [\imath_n]$ and, for all $\tau\in T(B_1),$
\beq
|(1/K)\tau\circ L_1(h)-\kappa_T(\bt_1(\tau))({\imath}_n(h))|<\sigma/4\rforal h\in {\cal H}.
\eneq
By  the definition of $\bt_1,$
for all $t\in T(B),$
\beq\label{2020-718-21}
|t\circ L_1(h)-{{\kappa_T(t)}}({\imath}_n(h))|<\sigma/2\rforal h\in {\cal H}.
\eneq
Define  $\phi: {{A^\T}} \to M_2(B)$ by
$\phi(a)=\phi_{0,b}\circ L(a)\oplus L_1({{\Psi_1}}(a))$  for all $a\in A^\T.$
Then (by choosing  sufficiently
large ${\cal G}_b$)
 $\phi$ is ${\cal G}_3$-$\dt_1$-multiplicative. One then checks that  {{(see also \eqref{11Ext1-3} and \eqref{11Ext1-4nn})}}
\beq
[\phi]|_{\cal P}=\kappa|_{\cal P}\andeqn
|t\circ \phi(a)-\kappa_T(t)(a)|<\sigma\rforal a\in {\cal H}{{\andeqn t\in T(B).}}
\eneq
Therefore, we obtain a  sequence of approximate multiplicative  \cpc s $\phi_n: C(\T,{\tilde A})^o\to M_2(B)$ such that
\beq
&&[\{\phi_n\}]=\kappa\tand\\
&&\lim_{n\to\infty}\sup \{|\tau\circ \phi_n(a)-\kappa_T(\tau)(a)|: \tau\in T(B)\}=0\tforal a\in C(\T,{\tilde A})^o_{s.a.}
\eneq
To modify the $\phi_n$ so it maps to $B$ instead of $M_2(B),$  using a strictly positive element of $A^\T,$
we deploy the same argument {{used}} in the end of the proof of Theorem \ref{ExtTBA}.
\end{proof}

 \begin{df}\label{Dz20}
 Denote by ${\cal Z}_{2,o}$ the \CA\, in ${\cal M}_1$ with continuous scale which is an inductive limit of \CA s in ${\cal C}_0$
with $K_0({{{\cal Z}_{2,o}}})=\Z,$ $K_1({\cal Z}_{2,o})=\{0\}$ and with
 two extremal traces $t_{o,+}$ and $t_{o,-}.$
 Let $1_\Z\in K_0({{{\cal Z}_{2,o}}})$ be the generator of $\Z$ which is represented by
 ${{[1_{{\cal Z}^{\sim}_{2,o}}]}}-[p],$ where $p\in {{M_2({\cal Z}^{\sim}_{2,o})}}$ is a projection. We assume
 $t_{o,+}(1_\Z)=1$ and $t_{o,-}(1_\Z)=-1.$

Denote by ${\cal Z}_{+,o}$ the \CA\, in ${\cal M}_1$ with continuous scale which can be written as  {{an}} inductive limit of
of \CA s in ${\cal C}_0,$
with $K_0({\cal Z}_{+,o})=\Z\cong \Z g_1,$  and ${\cal Z}_{+,o}$ has two extremal tracial states
$\tau_{1,+}$ and $\tau_{1,o}$ such that $\tau_{1,+}(g_1)=1$ and $\tau_{1,o}(g_1)=0.$
%

Let $n>1$ be an integer. Denote by ${\cal Z}_o^n$  {{the}} \CA\, in ${\cal M}_1$ with a unique
tracial state $\tau_o^n$ and with $K_0({\cal Z}_o^n)=\Z/n\Z$ and $K_1({\cal Z}_o^n)=\{0\}.$
Note that ${\cal Z}_o^n$ is in the class of simple \CA s classified in \cite{GLII}. Note that, by Theorem 4.31 of \cite{GLrange}, ${\cal Z}_{2,o},$ ${\cal Z}_{+,o}$  and ${\cal Z}_o^n$ exist. {{Note that the inductive limit  of these algebras (in the construction of Theorem 4.31 of \cite{GLrange}) also can be chosen in the form of (3) of Remark 4.32 of \cite{GLrange}.}}

 \end{df}

 \begin{NN}\label{tvalues}
 Suppose that $A\in {\cal M}_1$ and $x\in K_0(A)\setminus {\rm ker}\rho_A$ such that $\Z x\cong \Z.$
 Then there exists a $t_o\in T(A)$ such that $\rho_A(y)(t_o)=0$ for all $y\in K_0(A)$ {{(see Lemma 3.1 of \cite{GLrange}).}}
 Suppose that $t_1(x)=\theta\not=0.$  There is an integer $m$ such that $m\theta>0.$
 Put $\tau_A=(1-(1/m\theta))t_o+(1/m\theta)t_1.$
 Then $\tau_A(x)=1/m.$  Let {{$U$}} be a UHF-algebra of infinite type.
 We can always choose $m$ so that $1/m\in K_0({{U}}).$

 \end{NN}

\begin{lem}\label{Lr=0}
Let $A,\, B\in {\cal M}_1$ be two simple \CA s with continuous scale.
Let ${\cal P}_0\subset K_0(A)$ be a finite subset, {{${\cal V}\subset  U(M_m(C(\T)\otimes \td A))$
(for some $m\in  \N$)
be a  finite subset such that $\{[v]: v\in {\cal V}\} =\boldsymbol{\bt}({\cal P}_0),$ and
  $0<\ep_0<1/2.$}}
{{Then}} there {{exist}}  $\ep>0$ and a finite subset
${\cal F}$ satisfying the following:
Suppose that $\phi: A\to B$ is a \hm\, and $u\in {\tilde B}$ is a unitary
with $\pi_\C(u)=1,$  where $\pi_\C: {\tilde B}\to \C$ is the quotient map,
such that (see \ref{DDo} for {{the}} notation $\Phi_{\phi, u}$)
$$
\|[\phi(a),\, u]\|<\ep\rforal a\in {\cal F}\andeqn [\Phi_{\phi, u}]|_{\boldsymbol{\bt}({\cal P}_0)}=0,
$$
then $\Phi_{\phi, \,u}$ induces a \hm\, ${{\lambda:=\Phi_{\phi, \,u}^\dag\circ J_{cu}^{C(\T)\otimes \td A}\circ \boldsymbol{\bt}}}: G_0\to \Aff(T({\tilde B}))/\overline{\rho_{\tilde B}(K_0({\tilde B}))}$
such that the image of $\lambda(G_0)$ lies in
$\Aff(T(B))^\iota/\overline{\rho_{\tilde B}(K_0({\tilde B}))},$ where $G_0\subset K_0(A)$ is the subgroup generated by ${\cal P}_0,$
{{and ${\rm dist}({{\overline{\lceil \Phi_{\phi, u}(v)\rceil}}}, \Phi_{\phi, \,u}^\dag( J_{cu}^{C(\T)\otimes \td A}([v])))<\ep_0$ for all $[v]\in {\boldsymbol{\bt}}({\cal P}_0)$}} {{(see \ref{DLddag} for notation).}}
Moreover, if $B\in {\cal D}^d,$ $\lambda({\rm Tor}(G_0))=0.$
\end{lem}

\begin{proof}
Suppose
that $G_0$ is generated by $[p_i]-[q_i],$ where {{$p_i, q_i\in M_m({\tilde A})$}} are  projections, $i=1,2,...,k.$
We may also assume that ${\rm Tor}(G_0)$ is generated by $[p_i]-[q_i],$
$i=k_0, k_0+1,...,k$ (for some $k_0\ge 1$).
We choose $\ep$ sufficiently small and ${\cal F}$ sufficiently large
so that
$[\Phi_{\phi, u}]|_{{\cal P}_0}$ and
$[\Phi_{\phi, u}]|_{\boldsymbol{\bt}({\cal P}_0)}$
are {{well defined,}} and $\Phi_{\phi, u}$ induces a \hm\,
$\lambda: G_0\cong {\boldsymbol{\bt}}(G_0)\to {{U({\tilde B})/CU({\tilde B})}}$ {{such that ${{\dist({{\overline{\lceil \Phi_{\phi, u}(v)\rceil}}}, \lambda\circ \boldsymbol{\bt}^{-1}([v]))<\ep_0<1/2}}$ for all {{$v\in {\cal V}$}}}}
(see 14.5 of \cite{LinLAH}).
{{Since $[\Phi_{\phi, u}]|_{\boldsymbol{\bt}({\cal P}_0)}=0$,}} {{${{\overline{\lceil \Phi_{\phi, u}(v)\rceil}}}\in U_0({\tilde B})/CU({\tilde B})$ for any $v\in {\cal V}$. Consequently,   $\lambda(G_0)\subset U_0({\tilde B})/CU({\tilde B})\cong \Aff(T({\tilde B}))/\overline{\rho_{\tilde B}(K_0({\tilde B}))}$. }}
Moreover,
let
$$
S_i=((1_m-\phi^\sim(p_i))+
\phi^\sim(p_i)\underline{u})((1_m-\phi^\sim(q_i))+
\phi^\sim(q_i)\underline{u^*}),\,\,\,i=1,2,...,k,
$$
where $\underline{u}=\diag(\overbrace{u,u,...,u}^m).$
We may assume that (see \ref{DLddag})
$$
\overline{\lceil S_i\rceil}=\lambda([p_i]-[q_i]).
$$
 Let $\pi_\C: {\tilde B}\to \C$ {{be the quotient map}} and {{continue to denote by  $\pi_\C$}} its extension on $M_m({\tilde B}).$
Then
\beq
\pi_\C(S_i)=\pi_\C((((1_m-\phi^\sim(p_i))+\phi^\sim(p_i))(1_m-\phi^\sim(q_i))+\phi^\sim(q_i))={{\pi_\C(1_m).}}
\eneq
It follows that
$\pi_\C(\lceil S_i\rceil)=\pi_\C(1_m).$
By {{assumption,}} we may write $\lceil S_i\rceil=\sum_{j=1}^n\exp(ih_j),$ where $h_j\in M_m(B_{s.a.}),$
$j=1,2,...,n$ (see Lemma 6.1 of \cite{GLII}). In particular,
$$
{{\tau^B_\C}}(\sum_{j=1}^nh_j)=0.
$$
It follows that $\lambda([p_i]-[q_i])\in \Aff(T(B))^\iota/\overline{\rho_{\tilde B}(K_0({\tilde B}))},$
$i=1,2,...,k.$

It remains to show $\lambda([p_i]-[q_i])=0$ for $i=k_0, k_0+1,...,k,$ when $B\in {\cal D}^d.$
But this follows from the fact that $\Aff(T(B))^\iota/\overline{\rho_{\tilde B}(K_0({\tilde B}))}$
is torsion free (by \ref{PK0divisible}).

\end{proof}


\begin{lem}\label{Lfinitelimit}
Let $A=A_1\otimes U$ be a simple \CA\, with continuous scale which satisfies the UCT,
where  $U$ is an infinite dimensional UHF algebra and $A_1\in {\cal D}.$
Then $A=\lim_{n\to\infty}(A_n'\otimes U, j_n),$
where
$A_n'$ is a simple $C^*$-algebra {{as in (3) of Remark 4.32  of \cite{GLrange}}} 
 with finitely generated $K_i(A_n')$
($i=0,1$){{(so it also satisfies {{condition}} in Theorem 4.31 of \cite{GLrange})}},  {{and}} each $j_n$ maps strictly positive elements to strictly positive elements.

Moreover, there is a sequence of \SCA s  $A_n\subset A$
such that $A=\overline{\cup_{n=1}^\infty A_n}$,
 \beq\label{june-25-2021}
 \lim_{n\to\infty}{\rm dist} (x, A_n)=0~~\mbox{for any}~~x\in A,
 \eneq
 and
each $A_n$ has the form $A_n'\otimes M_{k(n)}$ for some integer $k(n)$

Furthermore, for any finite subset ${\cal P}\subset \underline{K}(A),$
there is an integer $N\ge 1$ such that ${\cal P}\subset [\iota_n](\underline{K}(A_n)),$
where $\iota_n:={j_{n,\infty}}|_{A_n}$ is the embedding.
%
%
\end{lem}

 \begin{proof}
 {{Recall that, since $A_1\in {\cal D},$  we have ${\rm Ped}(A_1)=A_1$ (see \ref{DD0}).
 Since $A$ has continuous scale, $T(A)$ is compact (see Theorem 5.3 of \cite{eglnp}).
 It follows that $T(A_1)=T(A)$ is compact.  Hence $A_1$ also has continuous scale (see also Theorem 5.3 of \cite{eglnp}).
 By Theorem 4.31 of \cite{GLrange}, there is  $A_0\in {\cal M}_1$ such that $A_0$ and $A_1$ have same Elliott invariants.}}
 By Theorem \ref{Misothm}, 
 $A\cong A_0\otimes U$. 
 There exists a sequence of finitely generated {{subgroups}}
 $G_n\subset G_{n+1}\subset K_0(A)$ and $F_n\subset F_{n+1}\subset K_1(A_0)$
 such that
 $\cup_{n=1}^{\infty} G_n=K_0(A_0)$ and $\cup_{n=1}^{\infty} F_n=K_1(A_0).$
 Fix $J_{cu}^{A_0}: K_1(A_0)\to U({\tilde A_0})/CU({\tilde A_0}),$ a splitting map of the following exact sequence
 $$
 0\to \Aff(T({\tilde A_0}))/\overline{\rho_{\tilde A_0}(K_0({\tilde A}_0))}\to  U({\tilde A}_0)/CU({\tilde A}_0)\to
 K_1(A_0)\to 0.
 $$

 It follows from {{(3) of Remark 4.32 of \cite{GLrange}}} 
 that there are {{\CA s}} $A_n'\in {\cal M}_1$  with continuous {{scales}}
 such that $K_0(A_n')=G_n,$  $K_1(A_n')=F_n$ and $T(A_n')=T(A_0),$ $n=1,2,....$  Moreover,
 $\rho_{A_n'}=\rho_{A_0}|_{G_n}.$
 Define $\imath_{0,n}: K_0(A_n')=G_n\to G_{n+1}=K_0(A_{n+1})$ and
 $\imath_{1,n}: K_1(A_n')=F_n\to F_{n+1}=K_1(A_{n+1}')$ to be the embeddings.
 Define $\imath_{T,n}: T(A_{n+1}')=T(A)\to {{T(A)=T(A_{n}')}}$ to be the identity map.
 Let $J_{cu}^{A_n'}: K_1(A_n')=G_n\to U({\tilde A_n'})/CU({\tilde A_n'})$  {{be}} a splitting map of the following exact sequence
 $$
 0\to \Aff(T({\tilde A_n'}))/\overline{\rho_{\tilde A_n'}(K_0({\tilde A_n'}))}\to  U({\tilde A_n'})/CU({\tilde A_n'})\to
 K_1(A_n')\to 0.
 $$
We may identify
 $T(A_n')$ with $T(A_0)=T(A).$ Therefore, we may further identify $T({\tilde A}_0)$ with $T({\tilde A_n'}).$
 Thus $\imath_{T,n}$ induces an affine continuous map
 $$
 \overline{\imath_{T,n}}: \Aff(T({\tilde A_n'}))/\overline{\rho_{\tilde A_n'}(K_0({\tilde A_n'}))}
 \to \Aff(T({\tilde A_{n+1}}))/\overline{\rho_{\tilde A_{n+1}'}(K_0({\tilde A_{n+1}'}))}.
 $$
 Define ${{\gamma_n}}: U({\tilde A_n'})/CU({\widetilde{A_{n}'}})\to U({\widetilde{A_{n+1}'}})/CU({\widetilde{A_{n+1}'}})$ by
 $({{\gamma_n}})|_{ \Aff(T({\widetilde{A_n'}}))/\overline{\rho_{\tilde A_n'}(K_0({\tilde A_n'}))}}=\overline{\imath_{T,n}}$ and
 $({{\gamma_n}})|_{J_{cu}^{A_n'}(K_1(A_n'))}=J_{cu}^{A_{n+1}'}\circ \imath_{1,n}\circ {J_{cu}^{A_n'}}^{-1}.$
 Let ${{\kappa_n}}\in KL(A_n', {{A_{n+1}'}})$ be an element such that
 $({{\kappa_n}})|_{K_i(A_n')}=\imath_{i,n},$ $i=0,1,$ given by the UCT.
 Then $(\kappa_n, \imath_{T,n}, \gamma_n)$ is compatible.
 Write $U=U_1\otimes U_2$, where both $U_1$ and $U_2$ are infinite dimensional UHF algebras.
 Since $A_n'$ is as in (3) of Remark 4.32 {{of \cite{GLrange}}}, it follows from Theorem \ref{Text1},
 that there exists a \hm\, $\phi_n: A_n'\otimes U_1\to A_{n+1}'\otimes U_1$ such that
 $$
 [\phi_n]=\kappa_n\otimes [\id_{U_1}],\,\,\, (\phi_n)_T=\imath_{T,n}\andeqn \phi_n^{\dag}=\gamma_n\otimes \id_{{U_1}}^\dag.
 $$
 Note that since each $A_n'$ has continuous scale and $\imath_n$ is {{the}}  identity map when we identify
 $T(A_n')$ with $T(A),$  ${\phi_n}$ 
 maps strictly positive elements to strictly positive elements.
 Let $B_1=\lim_{n\to\infty}(A_n'\otimes {{U_1}}, \phi_n).$ 
 and $B=B_1\otimes U_2$. Then, one checks easily from the construction that
 ${\rm Ell}(B)={\rm Ell}(A).$  The first part of lemma then follows from the isomorphism theorem \ref{Misothm} by
 {{setting}} $A=B$ and $j_n=\phi_n\otimes \id_{U_2}.$


 {{Let ${\cal G}_1\subset {\cal G}_2\subset \cdots\subset {\cal G}_n\cdots\subset A$ be a sequence of subsets with $A=\overline{\cup_{n=1}^{\infty} {\cal G}_n}$ and $\ep_1>\ep_2>\cdots>\ep_n>\cdots >0$ be a sequence of positive numbers with $\sum_{n=0}^{\infty}\ep_n<\infty$. For each $n$, there is an $l(n)$ such that $\dist(f, A'_{l(n)}\otimes U)<\ep_n/2$ for all $f\in {\cal G}_n$. Hence there is a finite subset ${\cal G}'_n\subset A'_{l(n)}\otimes U$ such that $\dist(f, {\cal G}'_n)<\ep_n/2$ for all $f\in {\cal G}_n$.  Since $U=\lim(M_{k(m)}, \iota_{m, m+1})$ with $k(1)|k(2)|\cdots$, there is a $k(n)$ such that
 $\dist(f, A'_{l(n)}\otimes M_{k(n)})<\ep_n/2$ for all $f\in {\cal G}'_n$. Set $A_n=A'_{l(n)}\otimes M_{k(n)}$. Evidently (\ref{june-25-2021}) holds.}}

 To see the last part of the statement, one notes that, since ${\cal P}$ is a finite subset,
 there exists $N\ge 1$ such that ${\cal P}\subset [j_{n, \infty}](\underline{K}(A_n'\otimes U))$
 for all $n\ge N.$   It {{is}} then clear, that, by passing to a subsequence,
 we may also assume that ${\cal P}\subset [\iota_n](\underline{K}(A_n)).$


 \end{proof}

\begin{lem}\label{Lpurterbation}
For any $1/2>\ep>0,$ there exists $\dt>0$ satisfying the following:
For any pair of positive elements $a,\, b$ in a \CA\, $A$ with $\|a\|=\|b\|=1,$
let $p_1$ be the spectral projection of $a$ corresponding to the close subset $\{1\}$ in $A^{**};$
if $\|a-b\|<\dt,$ then,
\beq
\|p_1q_{[0, 1-\ep)}\|<\ep/2 \andeqn \|p_1q_{[1-\ep, 1]}-p_1\|<{{\ep/2,}}
\eneq
where $q_{S}$ is the spectral projection of $b$ corresponding to the {{subset}} $S$ in $A^{**}.$
\end{lem}

\begin{proof}
Let $1/2>\ep>0$ be given.  Let $g_\ep\in C_0((0,1])$ such that $0\le g_\ep\le 1,$
$g_\ep(t)=0$ if $t\in [0,1-\ep]$ and  $g_\ep(t)=1$ if $t\in [1-\ep/2, 1].$

There is a universal constant $\dt>0$ independent of $a$ and $b$ such that
\beq
\|g_\ep(a)-g_\ep(b)\|<\ep/2
\eneq
wherever $a, b$ satisfy the assumption and $\|a-b\|<\dt.$
Therefore
\beq
\|g_\ep(a)p_1-g_{\ep}(b)p_1\|<\ep/2.
\eneq
Note that $g_\ep(a)p_1=p_1.$ Thus,
\beq
\|q_{[0,1-\ep)}p_1-q_{[0,1-\ep)}g_\ep(b) p_1\|<\ep/2
\eneq
Since $q_{[0,1-\ep)}g_\ep(b)=0,$
This implies
\beq
\|q_{[0,1-\ep)}p_1\|<\ep/2, \,\,\, \andeqn \|p_1q_{[0,1-\ep)}\|<{{\ep/2.}}
\eneq
Hence
\beq
 \|p_1q_{[1-\ep,1]}-p_1\|=\|p_1q_{[1-\ep,1]}-p_1(q_{[1-\ep,1]}+q_{[0,1-\ep)})\|<\ep/2.
\eneq
\end{proof}

 \begin{NN}\label{Djcubt}
   Let $A$ be a {{stably}} projectionless simple \CA\, with stable rank one.  Let $g_1, g_2,...,g_k\in K_0(A)$
 such that  $[p_i]-[q_i]=g_i,$ where $p_i, q_i\in M_r({\tilde A})$ are projections.
 In what follows, we will set
 $$
 J_{cu}^{C(\T)\otimes \td A}\circ {\boldsymbol{\bt}}(g_i)=\overline{((1_r-p_i)+p_i\otimes z)((1_r-q_i)+q_i\otimes z^*)}
 $$
 as elements in $U({\widetilde{A^{\T}}})/CU({\widetilde{A^\T}}),$ {{where $z$ is the standard unitary generator of $C(\T).$}}
   \end{NN}


\begin{lem}\label{LAtensorC}
{{Let $A$ be a separable stably projectionless simple \CA\, and $B$ be a separable simple \CA\, with stable rank one,
and }}
$\phi:A\otimes C(\T)\to B$
{{be a
\hm\,
which maps}}
strictly positive elements to strictly positive elements.
Suppose that $[\phi]|_{K_1(A\otimes C(\T))}=0.$
For any $\ep>0,$ any finite {{subsets}} ${\cal F}\subset A$
{{and}}
${\cal P}_0\subset K_0(A),$
any $\ep_0>0,$
there exists a unitary $v\in CU({\tilde B})$
such that
\beq\label{LAtensorC-1}
\|\phi_A(a),\, v]\|<\ep\tforal a\in {\cal F},
\eneq
where $\phi_A=\phi|_{A\otimes 1_{C(\T)}},$
and, for all $g\in {\cal P}_0,$
\beq\label{LAtensorC-2}
{\rm dist}((\Phi_{\phi_A,\, v})^{\dag}(J_{cu}^{C(\T)\otimes \td A}\circ{\boldsymbol{\bt}}(g)), \phi^{\dag}(J_{cu}^{C(\T)\otimes \td A}\circ {\boldsymbol{\bt}}(g)))<\ep_0.
\eneq
%
%
\end{lem}

\begin{proof}
We will {{write}} $\phi^\sim: {\widetilde{A\otimes C(\T)}}\to {\tilde B}$ for the unital extension.
We may assume that ${\cal F}\subset  A^{\bf 1}.$

Let $G_0$ be the subgroup generated by ${\cal P}_0.$
We may assume that $G_0$ has a set of generators
$g_1, g_2,...,g_k.$
We may further assume that there are  positive integers, $m,$ $m_1,m_2,...,m_k,$
with $m_i\le m,$ $i=1,2,...,k,$
and projections $p_1, p_2,...,p_k\in M_m({\tilde A})$ such that
$1_{m_i}-p_i\in M_m(A)$ and $[1_{m_i}]-[p_i]=g_i,$  $i=1,2,...,k.$
Let $x_i:=p_i-1_{m_i}\in M_m(A).$
Let $\{e_n\}$ be an approximate identity for $A$ such that
$e_{n+1}e_n=e_n,$ $n=1,2,....$
Write $E_n=\diag(\overbrace{e_n,e_n,..,e_n}^m)\in M_m(A),$ $n=1,2,....$
Then $E_n1_{m_i}=1_{m_i}E_n,$ $i=1,2,....$
Note that $p_i$ is close to $1_{m_i}+E_nx_iE_n$ for sufficiently large $n,$
therefore there is a projection $p_i'\in 1_{m_i}+E_{n+1}M_m(A)E_{n+1}$
which is close to $p_i.$ Therefore, \wilog, we may assume
that $x_iE_1=E_1x_i=x_i,$ $i=1,2,...,k.$ {{Note that $1_{m_i}E_1=E_1 1_{m_i}.$}}
Hence,  we  may also assume, \wilog,
that
\beq\label{LAtensorC-10}
&&(1_{m_i}-p_i)E_1=E_1({{1_{m_i}}}-p_i)=1_{m_i}-p_i,\\\label{LAtensorC-11}
&&(p_i-1_{m_i}p_i)E_1=E_1(p_i-1_{m_i}p_i)=p_i-1_{m_i}p_i\andeqn\\\label{LAtensorC-12}
&&(1_{m_i}-1_{m_i}p_i)E_1=E_1(1_{m_i}-1_{m_i}p_i)=1_{m_i}-1_{m_i}p_i
,\,\,\,i=1,2,...,k.
\eneq
Moreover, \wilog, we may also assume that
$e_1a=ae_1=a$ for all $a\in {\cal F}.$

Let $z\in C(\T)$ be the identity function on $\T.$
Put $Z=\diag(\overbrace{z,z,...,z}^m).$
Consider
\beq
W_i&=&((1_m-1_{m_i})\otimes 1_{C(\T)}+1_{m_i}\otimes Z)((1_m-p_i)\otimes 1_{C(\T)}+p_i\otimes Z^*)\\
  &=& (1_m-1_{m_i})(1_m-p_i)\otimes 1_{C(\T)}+(p_i-1_{m_i}p_i)\otimes Z^*\\
  &&+(1_{m_i}-1_{m_i}p_i)\otimes Z +1_{m_i}p_i\otimes 1_{C(\T)}.
  \eneq
Write $z=1+x,$ where $x\in C(\T)$ such that $x(1)=0$ {{and}}  $\T$ is identified with the unit circle.
Then $x^*x=xx^*$ and $x+x^*+x^*x=0.$
%
Let $y=1+e_2\otimes x$ and $Y=1_m+{{E_2\otimes x}}.$
Note that
\beq\label{y^*y=yy^*}
y^*y=1+(e_2^2-e_2)\otimes x^*x=yy^*.
\eneq
{{Note that $sp(x)=\{\ld\in \C: |\ld+1|=1\}$ which implies $sp(y)\subset \{\ld\in \C: |\ld|\leq 1\}$. Hence
\beq\label{August-10-2021}
\|y\|=\|Y\|\leq 1.
\eneq}}
By \eqref{LAtensorC-10}, \eqref{LAtensorC-11}, {{\eqref{LAtensorC-12} and \eqref{y^*y=yy^*},}}
\beq
&&\hspace{-0.5in}(p_i-1_{m_i}p_i)\otimes Z^*={{(p_i-1_{m_i}p_i)\otimes1_{C(\T)}+}} (p_i-1_{m_i}p_i)E_2\otimes x^*\in
M_m({\widetilde{A\otimes C(\T)}}),\\
&&\hspace{-0.5in} (1_{m_i}-1_{m_i}p_i)\otimes Z={{(1_{m_i}-1_{m_i}p_i)\otimes 1_{C(\T)}+}} ({{1_{m_i}}}-1_{m_i}p_i)E_2\otimes x\in M_m({\widetilde{A\otimes C(\T)}}),\\\label{91-2}
&&((p_i-1_{m_i}p_i)\otimes 1_{C(\T)})Y^*Y=(p_i-1_{m_i}p_i)\otimes 1_{C(\T)}\andeqn\\\label{91-2+}
&&((1_{m_i}-1_{m_i}p_i)\otimes 1_{C(\T)})Y^*Y=((1_{m_i}-1_{m_i}p_i)\otimes 1_{C(\T)}).
\eneq
It follows that $W_i\in M_m({\widetilde{A\otimes C(\T)}}).$
So $\phi^\sim(W_i)$ is defined {{for}} $i=1,2,...,k.$ Let $P_1$ be the spectral projection of $\phi(y^*y)$ in $B^{**}$ corresponding to $\{1\}.$
{{Then, by \eqref{91-2},}}
\beq\label{91-3}
{{\phi((p_i-1_{m_i}p_i)\otimes 1_{C(\T)}){\bar P_1}=\phi((p_i-1_{m_i}p_i)\otimes 1_{C(\T)}),\,\,({\bar P}_1=\diag({{\overbrace{P_1, P_1,...,P_1}^m}})).}}
\eneq


{{Fix}}  ${{\ep>0}},\ep_0>0.$  {{Let $\eta<\min\{\ep/8,\ep_0/8\}$ be a sufficient small positive number (to be determined later).}}
Choose $\dt>0$ as in \ref{Lpurterbation} for {{$\eta$}} (instead of $\ep$).
 Since ${\tilde B}$ has stable rank one, one obtains an invertible
element $y_1\in 1_{\tilde B}+B$ such that
\beq\label{LAtensorC-n0}
&&\|\phi^\sim(y)-y_1\|<\min\{\ep_0/2, {{\eta/2,}} \dt\},\,\,\,\|\phi(y^*y)-y_1^*y_1\|<\min\{{{\eta}}, \dt\}\andeqn\\
&&\|\phi(yy^*)-y_1y_1^*\|<\min\{{{\eta}}/2, \dt\}.
\eneq
Let
$Q_{S}$ be the spectral projection of $y_1^*y_1$ in $B^{**}$ corresponding to the subset $S.$
Then, {{by}} \ref{Lpurterbation},
\beq\label{LAtensorC-n1}
\|P_1Q_{[0, 1-{{\eta}})}\|<{{\eta}}/2\andeqn \|P_1Q_{[1-{{\eta}},1]}-P_1\|<{{\eta}}/2.
\eneq
Note that
\beq\label{LAtensorC-n2}
\|Q_{[1-{{\eta}}, 1]}(y_1^*y_1)^{-1/2}-Q_{[1-{{\eta}},1]}\|<|1-(1-{{\eta}})^{1/2}|.
\eneq
Put
$
{{\eta_0}}=|1-(1-{{\eta}})^{1/2}|.
$
Let $v_1=y_1(y_1^*y_1)^{-1/2}\in 1_{\tilde B}+B.$ Note {{also}} that ${{v_1}}$ is a unitary.
Put $V_1=\diag(\overbrace{v_1,v_1,...,v_1}^m),$ $Y_1=\diag(\overbrace{y_1,y_1,...,y_1}^m)$
and
${\bar Q}_{S}=\diag(\overbrace{Q_S, Q_S,...,Q_S}^m).$
Then, by {{\eqref{91-3},}} \eqref{LAtensorC-n1}, \eqref{LAtensorC-n2},  {{\eqref{LAtensorC-n0}, and \eqref{August-10-2021},}}
\beq
\hspace{-0.3in}\phi((p_i-1_{m_i}p_i)\otimes 1_{C(\T)})V_1^*&&=\phi((p_i-1_{m_i}p_i)\otimes 1_{C(\T)})(Y_1^*Y_1)^{-1/2}Y_1^*\\
&&=\phi((p_i-1_{m_i}p_i)\otimes 1_{C(\T)}){\bar P}_1{\bar Q}_{[1-{{\eta}}, 1]}(Y_1^*Y_1)^{-1/2}Y_1^*\\
&&+\phi((p_i-1_{m_i}p_i)\otimes 1_{C(\T)}){\bar P}_1{\bar Q}_{[0, 1-{{\eta}})}(Y_1^*Y_1)^{-1/2}Y_1^*\\
&&\approx_{{2\eta_0+\eta}}\phi((p_i-1_{m_i}p_i)\otimes 1_{C(\T)}){\bar P}_1{\bar Q}_{[1-{{\eta}}, 1]}Y_1^*\\
&&\approx_{{{\eta}}/2}\phi((p_i-1_{m_i}p_i)\otimes 1_{C(\T)}){\bar P}_1Y_1^*\\
&&\approx_{{{\eta}/2}} \phi((p_i-1_{m_i}p_i)\otimes 1_{C(\T)})\phi^\sim(Y)^*\\
&&=\phi((p_i-1_{m_i}p_i)\otimes  Z^*).
\eneq
Note that $V_1^*=Y_1^*(Y_1Y_1^*)^{1/2}.$
Similarly,  using  {{\eqref{91-3},}}
\beq
V_1^*\phi((p_i-1_{m_i}p_i)\otimes 1_{C(\T)})&&=Y_1^*(Y_1Y_1^*)^{1/2}\phi((p_i-1_{m_i}p_i)\otimes 1_{C(\T)})\\&&\approx_{{2\eta+2\eta_0}} \phi((p_i-1_{m_i}p_i)\otimes  Z^*).
\eneq
Also  {{(see \eqref{91-2+})}}
\beq
\hspace{-0.4in}V_1\phi((1_{m_i}-1_{m_i}p_i)\otimes 1_{C(\T)})&&=Y_1(Y_1^*Y_1)^{-1/2}\phi((1_{m_i}-1_{m_i}p_i)\otimes 1_{C(\T)})\\
&&\approx_{{2\eta_0+\eta}}Y_1{\bar Q_{[1-{{\eta}},1]}}{\bar P_1}\phi((1_{m_i}-1_{m_i}p_i){{\otimes 1_{C(\T)})}}\\
&&\approx_{{{\eta}/2}} Y_1{\bar P_1}\phi((1_{m_i}-1_{m_i}p_i){{\otimes 1_{C(\T)})}}\\
&&=Y_1\phi((1_{m_i}-1_{m_i}p_i){{\otimes 1_{C(\T)})}}\\
&&\approx_{{{\eta}}/2}\phi^\sim(Y)\phi((1_{m_i}-1_{m_i}p_i){{\otimes 1_{C(\T)})}}\\
&&=\phi(((1_{m_i}-1_{m_i}p_i)\otimes Z).
\eneq
Moreover,
\beq
\phi((1_{m_i}-1_{m_i}p_i)\otimes 1_{C(\T)})V_1&=&\phi((1_{m_i}-1_{m_i}p_i)\otimes 1_{C(\T)})(Y_1Y_1^*)^{1/2}Y_1\\
&&\approx_{{2\eta+2\eta_0}} \phi(((1_{m_i}-1_{m_i}p_i)\otimes Z).
\eneq

Furthermore,
\beq\label{LAtensorC-30-1}
&&V_1\phi(1_{m_i}-p_i)\approx_{2\eta+2\eta_0} \phi(1_{m_i}-p_i)V_1,\,\,\,i=1,2,...,k\andeqn
\eneq
{{In fact, using $(a\otimes 1_{C(\T)})y=y(a\otimes 1_{C(\T)}),$  the same argument as above  shows
that}}
\beq\label{LAtensorC-30}\
&&v_1{{\phi_A}}(a)\approx_{{2\eta+2\eta_0}} {{\phi_A}}(a)v_1\rforal a\in {\cal F}.
\eneq
We {{further}} note that
\beq
V_1\phi^\sim(1_{m_i})=\phi^\sim (1_{m_i})V_1.
\eneq
It follows {{(\ref{LAtensorC-30-1})}} that
\beq
V_1\phi^\sim (p_i)\approx_{2\eta+{{2\eta_0}}} \phi^\sim(p_i)V_1.
\eneq
Put
\beq\nonumber
W'_i=\left(\phi^\sim((1_m-1_{m_i})\otimes 1_{C(\T)})+(1_{m_i}\otimes 1_{C(\T)}) {{V_1}}\right)
\left(\phi^{\sim}((1_m-p_i)\otimes 1_{C(\T)})+\phi^{\sim}(p_i\otimes 1_{C(\T)}) V_1^*\right).
\eneq
Combining these estimates, one obtains
\beq\label{LAtensorC-30+}
W_i'\approx_{{4\eta+4\eta_0}} \phi^\sim(W_i).
\eneq
We may write $v_1=1+b$ for some normal element $b\in \overline{{{(\phi(e_1\otimes 1_{C(\T)}))}}B(\phi(e_1\otimes 1_{C(\T)}))}.$
Note $(e_3-e_2)b=b(e_3-e_2)=0.$ Put
$C_1=\overline{\phi((e_3-e_2)\otimes 1_{C(\T)})B\phi((e_3-e_2)\otimes 1_{C(\T)})}.$

Since $B$ {{is separable  and simple, and}}
has stable rank one,
there exists $b_1\in C_1$ such that
$\zeta_1=1_{{\tilde C}_1}+b_1\in U({\tilde C_1})$ such that $[\zeta]+[v_1]=0$ in $K_1({\tilde B}).$
Put $v_2=1_{{\tilde A}}+b_1.$ Then we may even assume that ${{v_2v_1}}\in CU({\tilde B})$
(see \ref{Unitary}).

Set $v={{v_2v_1}}.$
Note that  $v_2$ commutes with $\phi_A(a)$ for all $a\in {\cal F}.$
By \eqref{LAtensorC-30},
$$
\|[\phi_A(a), \, v]\|<{{2\eta+2\eta_0}}\rforal a\in {\cal F}.
$$
By choosing sufficiently small $\eta_0$ and $\dt,$
as in  \ref{DDo}, we may assume that $\Phi_{\phi_A,v}$ is {{well}} defined.
Set  $V=\diag(\overbrace{v,v,...,v}^m)$ and
$V_2=\diag(\overbrace{v_2,v_2,...,v_2}^m).$ Then $V_2$ commutes with  $V_1$,~~$\phi^{\sim}((1_m-p_i)\otimes 1_{C(\T)}$ and $\phi^{\sim}(p_i\otimes 1_{C(\T)}))$.
Thus $W'_i$ equals
\beq\nonumber
&&\hspace{-0.7in}{{\left((\phi^\sim((1_m-1_{m_i})\otimes 1_{C(\T)}))+(1_{m_i}\otimes 1_{C(\T)}) {{V}}\right)\left(\phi^{\sim}((1_m-p_i)\otimes 1_{C(\T)})+\phi^{\sim}(p_i\otimes 1_{C(\T)}) {{V^*}}\right)}}
\eneq
for
$i=1,2,...,k.$ By Lemma {{\ref{Lr=0},}} 
${\rm dist}(\overline{\lceil W'_i \rceil}, (\Phi_{\phi_A,\, v})^{\dag}(J_{cu}^{C(\T)\otimes \td A}\circ{\boldsymbol{\bt}}(g_i)))<\ep_0/4,$ if $\eta$ is small enough (recall that  $[\phi]|_{K_1(A\otimes C(\T))}=0$).
Finally,  combining  \eqref{LAtensorC-30+},
condition \eqref{LAtensorC-2} holds.  {{The lemma}} follows.
\end{proof}
%

\begin{lem}\label{Ldeterm2}
Let $A,\, B_0\in {\cal M}_1$ be such that $A\in {\cal D}^d$ 
and $B=B_0\otimes {{ U}}$ for some infinite dimensional UHF-algebra {{$U.$}}
Suppose that both $A$ and $B$ have continuous scale.
Suppose that $G_1\subset  K_0(A)$ is a
{{subgroup}} generated by $g_1, g_2,...,g_k,$
and
$\lambda: G_1\to \Aff(T(B))^\iota/\overline{{{\rho_{\tilde B}}}(K_0({\tilde B}))}$ is a \hm.
Then, for any $\ep>0,$ and any finite {{subsets ${\cal F}\subset A$  and}}
${\cal P}\subset \underline{K}(A),$ there is {{a}}
\hm\, $\phi: A\to B$
and a unitary $u\in CU({\tilde B})$
such that
\beq
&&\|[\phi(a),\, u]\|<\ep\tforal a\in  {\cal F}_1,\\
&&{[}\Phi_{\phi, u}{]}|_{\cal P\cup {\boldsymbol{\bt}}(\cal P)}=0,\\
&&\phi^{\dag}|_{J_{cu}^A(K_1(A))}=0,\andeqn\\
&&{\rm dist}( (\Phi_{\phi, u})^{\dag}\circ J_{cu}^{C(\T)\otimes {\td A}}(\boldsymbol{\bt}(g_i)), \lambda(g_i))<\ep,\,\,\,i=1,2,...,k,
\eneq
where $\Phi_{\phi, u}: A^\T\to B$ is the \cpc\,  induced by $\phi$ and $u.$
\end{lem}

 \begin{proof}

 First we note that $\Aff(T(B))^\iota/\overline{\rho_B(K_0({\tilde B}))}$ is a divisible group.
 Therefore we may assume that
 $\lambda$ is defined on $K_0(A).$

 Let $C_A\in {\cal M}_1$ {{be}} with continuous scale and with a unique tracial state $t_c$ such that
 $K_0(C_A)=K_0(A)$ as an abelian group,  $K_0(C_A)={\rm ker}\rho_{C_A}$ and $K_1(C_A)=\{0\}.$
 Therefore $C_A\in {\cal D}_0.$  Note {{that,}} by Corollary 13.4 of \cite{GLII},
 $C_A\otimes \zo\cong C_A.$ {{Note that $C_A$ can be written as the form of (3) of Remark 4.32 of \cite{GLrange}.}}

  Define $\af\in KL(A, C_A)$ such that $\af|_{K_0(A)}: K_0(A)\to K_0(C_A)$ is an isomorphism.
  It follows from {{Corollary 3.2 of \cite{GLrange}}} 
  that there exists  $t_{o,A}\in T(A)$ such that
  $\rho_A(x)(t_{o,A})=0$ for all $x\in K_0(A).$ Define $\af_T: T(C_A)\to T(A)$ by
  $\af_T(t_c)=t_{o,A}.$ One then checks that $(\af, \af_T)$
   is compatible.  Define $\af_{cu}: U({\tilde A})/CU({\tilde A})\to U({\tilde C_A})/CU({\tilde C_A})$
   as follows: {{on}} $\Aff(T({\tilde A}))/\overline{\rho_{{\tilde A}}(K_0(A))},$ define $\af_{cu}$ to be
   the map induced by $\af_T$ and $\af|_{K_0(A)};$ {{on}} $J_{cu}^A(K_1({\tilde A})),$ define $\af_{cu}$  to be zero.
   Then $(\af, \af_T, \af_{cu})$ is compatible. So, by \ref{11Ext2n2},
    there exists a \hm\,
   $h: A\to C_A$ such that
   $([h],h_T, h^{\dag})=(\af, \af_T, \af_{cu}).$

   Consider \CA\, $D=C_A\otimes F,$ where $F$ is a unital classifiable simple separable
   \CA s  with finite nuclear dimension and satisfies the UCT such that
   $K_0(F)=\Z,$ $K_1(F)=\Z$ and $T(F)=T(B).$  Note that, by  Theorem 15.8 of \cite{GLII}
   (see also Theorem {{4.27 of \cite{GLrange}}}), 
   $D\in {\cal M}_1\cap {\cal D}_0.$ {{Furthermore, we can assume $D$ is as the algebra $A$ in Theorem \ref{Text1},  by (3) of Remark 4.32 of \cite{GLrange}.}}
   Let $C_W={\cal W}\otimes F.$

   Put $\kappa\in KL(D, C_W)=0.$ Let $\kappa_T: T(C_W)=T(B)\to T(D)$ be the identity map.
   Define $\kappa_{cu}: U({\tilde D})/CU({\tilde D})\to \Aff(T({{\tilde{C}_W}}))/\Z$
   as follows:
   $\kappa_{cu}|_{\Aff(T({\tilde D}))/\overline{\rho_{{\tilde D}}(K_0({\tilde D})}}$ is defined
   to be the map induced by $\kappa_T$ and $\kappa$ ($=0$).
  Write $G_1=\Z^k \oplus {\rm Tor}(G_1).$
  Since $\overline{\iota_B^\sharp(\rho_B(K_0(B)))+\Z}/\Z$ is divisible (by
  Lemma \ref{PK0divisible}),
 {{we}} obtain a \hm\, $\Lambda: G_1\to
  \Aff(T(B))^{\iota}/\Z$ such that $\pi \circ \Lambda=\lambda|_{G_1},$
  where $\pi: \Aff(T(B))^\iota/\Z\to \Aff(T(B))^\iota/\overline{\rho_{\td B}(K_0(\td B))}$
  is the quotient map.  Since $\Aff(T(B))^\iota/\Z$ is divisible,
  we may assume that $\Lambda$ is defined on $K_0(A).$
  We may view $\Lambda$ maps $K_0(A)$ into ${{\Aff(T(B))^{\iota}/\Z=}}\Aff(T(C_W))^\iota/\Z=\Aff(T(C_W))^\iota/\overline{\rho_{\tilde C_W}(K_0({\tilde C_W}))}$.
Recall that $K_1(D)=K_1(C_A)\otimes \Z\oplus K_0(C_A)\otimes {{\Z}}=K_0(C_A\otimes C(\T))=K_0(A).$
     Now define
   $\kappa_{cu}|_{J_{cu}^D(K_1(C_A)\otimes \Z)}=0$ and $\kappa_{cu}|_{J_{cu}^D(K_0(C_A)\otimes \Z)}
   =\Lambda\circ {{\Pi_{cu}^{\td D}}}|_{J_{cu}^D(K_0(C_A)\otimes \Z)}$ (see \ref{DkappaJ} for ${{\Pi_{cu}^{\td D}}}$),
    where we identify {{$K_1(D)$ with}} $K_0(C_A)$ {{and}} $K_0(A)$, {{and also identify $\Aff(T(B))^{\iota}/\Z$ with $\Aff(T(C_W))^\iota/\Z$}}  .
   By {{Theorem}} \ref{Text1},
   there is  {{a}} \hm\, $h_1: D\to C_W$ which induces $(\kappa, \kappa_T, \kappa_{cu}).$
   Let $w\in U(F)$ {{be}} such that {{$\Z[w]=K_1(F).$}}
   Let $\iota: A\otimes C(\T)\to C_A\otimes F$ be defined by
   $\iota(a\otimes f)=h(a)\otimes f(w)$ for all $a\in A$ and $f\in C(\T).$
   Define $\psi_1:=h_1\circ \iota: A\otimes C(\T)\to C_W$
   and $\phi_1: A\to C_W$ by $\phi_1:={\psi_1}|_{A\otimes 1_{C(\T)}}.$  Then,
   (note that $K_i(C_W)=\{0\},$ $i=0,1$), by \ref{LAtensorC},
   there exists a unitary $v\in  CU({{\tilde{C}_W}})$
   such that
   \beq
   \|[\phi_1(a),\, v]\|<\ep\rforal a\in {\cal F},\\
   {\rm dist}(\Phi_{\phi_1, \,v}^\dag(\boldsymbol{\bt}(g_i)), \psi_1^{\dag}(\boldsymbol{\bt}(g_i)))<\ep,\,
   \,i=1,2,...,k.
   \eneq
   We may also assume that $[\Phi_{\phi_1, v}]|_{{\cal P}\cup{\boldsymbol{\bt}}({\cal P})}$ is {{well defined}}
   ({{if $\ep$}} is sufficiently small and ${\cal F}$ is sufficiently large).
    Note {{that}} since $K_i(C_W)=\{0\}$ ($i=0,1$), $[\Phi_{\phi_1, v}]|_{{\cal P}\cup{\boldsymbol{\bt}}({\cal P})}=\{0\}.$
    Recall also $\psi_1^\dag|_{G_1}= \Lambda.$
    There is also a \hm\, $h_2: C_W\to B$ such that $[h_2]=0$ {{and}} $(h_2)_T={\rm id}_{T(B)}.$
    Define $\phi=h_2\circ \phi_1.$
    %
  Then,
  for  $u:=h_2^\sim(v)\in  CU({\tilde B}),$ we have
   \beq
   \|[\phi(a),\, u]\|<\ep\rforal a\in {\cal F},\\
   {[}\Phi_{\phi, u}{]}|_{{\cal P}\cup{\boldsymbol{\bt}}({\cal P})}=0\andeqn\\
   {\rm dist}(\Phi_{\phi, \,u}^\dag(\boldsymbol{\bt}(g_i)), \psi^{\dag}(\boldsymbol{\bt}(g_i)))<\ep,\,
   \,i=1,2,...,k.
   \eneq
  Note that   $K_0(\td B)=K_0(B)\oplus \Z.$ Hence
  $\psi^\dag\circ \boldsymbol{\bt}|_{G_1}=h_2^{\dag}\circ \pi\circ \Lambda|_{G_1}=\lambda.$
  Note also $\phi^\dag|_{J_{cu}^A(K_1(A))}=0.$
  The lemma follows.

 \end{proof}

\begin{df}\label{DlocalKL}
Let $A$ be a separable \CA,  ${\cal P}\subset \underline{K}(A)$ be a {{finite}} subset and
$G({\cal P})$ be the subgroup generated by ${\cal P}.$
Suppose that $A=\overline{\cup_{n=1}^\infty} A_n$, $\lim_{n\to\infty}{\rm dist} (x, A_n)=0$ for any $x\in A$, and $K_*(A_n)$ are finitely generated.
Let  $\iota_n: A_n\to A$ be the embedding.
Suppose that $G({\cal P})\subset [\iota_{n}](\underline{K}(A_n).$
Put $F=[\iota_{n, \infty}]^{-1}(G({\cal P})).$
Let $B$ be another \CA\, and
let $\Gamma: G({\cal P})\to \underline{K}(B)$ be a \hm.
If there is  {{an}} $\af\in {\rm Hom}_\Lambda(\underline{K}(A_n), \underline{K}(B))$
{{such that}}  $\af(x)=\Gamma([\iota_{n}](x))$ for all $x\in F,$ then
we write $\Gamma\in KL_{loc}^A(G({\cal P}), \underline{K}(B)).$
In fact, when {{$K_i(A)$ ($i=0,1$)}} 
is  finitely generated,  $KL_{loc}^A(G({\cal P}), \underline{K}(B))$ can be defined without assuming
$A=\overline{\cup_{n=1}^\infty}A_n$ (see the end of 2.1.16 of \cite{Lncbms}).
\end{df}

Let $m>1$ be an integer.  In what follows
we may write
\beq\label{2020-920-1}
K_0(A, \Z/m\Z)=K_0(A)/mK_0(A)\oplus {\rm Tor}(K_1(A),{{\Z/m\Z}}),
\eneq
where ${\rm Tor}(K_1(A),{{\Z/m\Z}})$ is identified with those elements $x\in K_1(A)$ such that $mx=0.$
It should be noted that the direct sum is not natural.

 \begin{lem}\label{Lbotz2}
 Let $A$ be a finite direct sum of \CA s in ${\cal M}_1$ as constructed in Theorem 4.31 of \cite{GLrange}  and $B={\cal Z}_{2,o}\otimes U$  for some   UHF-algebra of infinite type.
 Then, for any $\ep>0,$  any finite subset ${\cal F}\subset A$ and any finite subset
 ${\cal P}\subset \underline{K}(A),$
there exists $\eta>0$  and a finite subset ${\cal Q}\subset \underline{K}(A)$
such that ${\cal P}\subset {\cal Q}$ and $G({\cal Q})\cap K_1(A)=\Z^k\oplus ({\rm Tor}(K_1(A))\cap G({\cal Q})),$
where $G({\cal Q})$ is the subgroup generated by ${\cal Q},$
which has a free generating set $\{g_1, g_2,...,g_k\}$
for $\Z^k$
satisfying the following:
if $1_\Z\in K_0({\cal Z}_{2,o})$ is a generator and $r\in K_0(U)_+\subset \R_+$
with
 \beq\label{Lbotz2-1}
0<r<\eta,
 \eneq
 then there exists  $\kappa\in KL^{A^\T}_{loc}(G({\cal Q})\cup \boldsymbol{\bt}(G({\cal Q})), B),$
 ${\cal F}$-$\ep$-multiplicative \cpc\, $\phi: A\to B$ and a unitary
 $u\in U({\tilde B})$ such that
 \beq\label{Lbotz0-2}
 &&\|[\phi(a),\, u]\|<\ep\rforal a\in {\cal F},\\
&& {\rm bott}(\phi,\, u)(g_1)=1_\Z\otimes r, {\rm bott}(\phi, \, u)(g_j)=0,\,\,\, j\not=1, {\rm bott}(\phi,\, u)|_{{\rm Tor}(G({\cal Q})}=0,\\
&&{[}\phi{]}|_{\cal Q}=0,
\kappa|_{{\boldsymbol{\bt}}(K_0(A)\cap {\cal Q})}=0, \kappa|_{{\boldsymbol{\bt}}(K_0(A,\Z/m\Z)\cap {\cal Q})}=0,\,m=2,3,...,\tand\\
&&\kappa|_{{\cal Q}\cup {\boldsymbol{\bt}}({\cal Q})}=[\Phi_{\phi, u}]|_{{\cal Q}\cup {\boldsymbol{\bt}}({\cal Q})}.
 \eneq
(Note that $K_1(B)=0$ and $K_1(B, \Z/m\Z)=0$ for all $m$.)

Moreover, if ${\rm Tor}(G({\cal Q})\cap K_1(A))\cong (\Z/m_1\Z)\cdot g_{1,t}\oplus (\Z/m_2\Z)\cdot g_{2,t}\oplus \cdots \oplus
(\Z/m_N\Z)g_{N,t},$ where $g_{1, t}, ...,g_{N,t}$ are generators,
$m_i\ge 2$ are integers, and $B_i={\cal Z}_o^{m_i}\otimes U$ ($1\le i\le N$) with $1_{\Z/m_i\Z}\in K_0({\cal Z}_o^{m_i})$
as a generator, then, for any $\ep>0$ and finite subset ${\cal F}\subset A,$ there exists $\eta>0$  {{with the following property:}}
 for any $r\in K_0(U)_+$ with $r<\eta,$
 {{there}} exists  $\kappa\in KL^{A^\T}_{loc}(G({\cal Q})\cup \boldsymbol{\bt}(G({\cal Q})),B_i),$
 {{an}} ${\cal F}$-$\ep$-multiplicative \cpc\, $\phi: A\to B_i$ and a unitary
 $u\in U({\tilde B})$ such that
 \beq\label{Lbotz0-2.1}
 &&\|[\phi(a),\, u]\|<\ep\rforal a\in {\cal F},\\
&&\hspace{-0.2in} {\rm bott}(\phi,\, u)(g_{i,t})=1_{\Z/m_i\Z}\otimes r, {\rm bott}(\phi, \, u)(g_{j,t})=0,\,\,\, j\not=i, {\rm bott}(\phi,\, u)|_{\Z^k}=0,\\\label{Lbotz0-2+}
&&{[}\phi{]}|_{\cal Q}=0,
\kappa|_{{\boldsymbol{\bt}}(K_0(A)\cap {\cal Q})}=0, \kappa|_{{\boldsymbol{\bt}}((K_0(A)/mK_0(A))\cap {\cal Q})}=0,\,m=2,3,...,\tand\\
&&\kappa|_{{\cal Q}\cup {\boldsymbol{\bt}}({\cal Q})}=[\Phi_{\phi, u}]|_{{\cal Q}\cup {\boldsymbol{\bt}}({\cal Q})}.
 \eneq

 \end{lem}

 \begin{proof}
 We may  write $A=\lim_{n\to\infty} (A_n, \iota_n),$ where $K_i(A_n)$
 is finitely generated ($i=0,1$)   and $A_n$ is a finite direct sum of simple \CA s in ${\cal M}_1.$
 We may also assume
 that  ${\cal F}\subset A_1$ and ${\cal P}\subset [\iota_{1, \infty}](\underline{K}(A_1)).$
 Let $[\iota_{1, \infty}](K_1(A))=\Z^k\oplus {\rm Tor}([\iota_{1, \infty}](K_1(A_1)).$
 Let $\{g_1, g_2,...,g_k\}$ be a free generator set for $\Z^k.$
 Let ${\cal Q}'\subset \underline{K}(A_1)$ be a finite
 subset such that $[\iota_{1, \infty}](G({\cal Q}'))\supset G({\cal P})$
 and $K_1(A_1)\subset G({\cal Q}).$ Put ${\cal Q}=[\iota_{1, \infty}]({\cal Q}').$
 In what follows we also  {{view}}
 $A_n$ as \SCA s of $A$ and
 $A_n^\T$ as \SCA\, of $A^\T.$

 Consider $\af\in Hom_{\Lambda}(\underline{K}({\tilde A_1}\otimes C(\T)), \underline{K}(B))$
 such that $\af|_{\underline{K}(A_1)}=0,$ and by the UCT, one may also require
 that $\af|_{\boldsymbol{\bt}(K_0(A_1))}=0,$  $\af|_{{\boldsymbol{\bt}}(K_0(A_1,\Z/m\Z))}=0,\,m=2,3,...,$
 $\af({\boldsymbol{\bt}}(g_1))=1_\Z,$
 $\af({\boldsymbol{\bt}}(g_j))=0$ if  $j\not=1,$ and $\af({\boldsymbol{\bt}}({\rm Tor}(K_1(A_1))))=0.$
 Note that ${\tilde{A^\T}}={\tilde A}\otimes C(\T).$
 Let $\kappa'\in Hom_{\Lambda}(\underline{K}(A_1^\T), \underline{K}(B))$ be given by $\af$
 which also gives
 $\kappa''\in KL_{loc}^{A^\T}(G({\cal Q})\cup \boldsymbol{\bt}(G({\cal Q}), \underline{K}(B)).$

 Let ${\cal G}\subset A_1^\T$ be a finite subset and $\dt>0$ such that
 ${\cal F}\subset {\cal G}$ as we identify $A_1$ with the corresponding \SCA\, of $A_1^\T$ and
 such that
 $[L]|_{{\cal Q}\cup{\boldsymbol{\bt}}({\cal Q})}$ is {{well defined}} for any ${\cal G}$-$\dt/2$-multiplicative
 \cpc s {{$L$}} from $A_1^\T.$
 We also assume
 that $\dt<\ep/2.$

 Note that, by {{Theorem 3.3 of \cite{GLrange}}} 
  and {{14.7}} of \cite{GLII}, $A_1^\T$ satisfies the assumption of 9.8 of \cite{GLII}.

 It follows from {{Theorem 3.4 of \cite{GLrange}}} 
 that
there exists a ${\cal G}$-$\dt/2$-multiplicative \cpc\,
$\Phi': A_1^\T\to B\otimes M_N$
 for some
integer $N\ge 1$ such that
\beq\label{11Ext1-2}
[\Phi']|_{{\cal Q'}\cup{\boldsymbol{\bt}({\cal Q'})}}
=\kappa'|_{{\cal Q'}\cup{\boldsymbol{\bt}({\cal Q'})}}.
\eneq
Since $A_1^\T$ is amenable, \wilog, there is a \cpc\, $j: A^\T\to A_1^\T$
such that  $\|j(a)-a\|<\dt/2$
for all $a\in {\cal G}.$ Put $\Phi=\Phi'\circ j: A^\T\to B\otimes M_N.$
Then $\Phi$ is ${\cal G}$-$\dt$-multiplicative and $[\Phi]|_{{\cal Q} \cup{\boldsymbol{\bt}({\cal Q})}} =\kappa''|_{{\cal Q} \cup{\boldsymbol{\bt}({\cal Q})}}.$

Let $\{k(n)\}$ be the sequence of integers such that $U=\lim_{n\to\infty}(M_{k(n)}, j_n).$
We may assume that $k(n_1)\ge N$ and view $\Phi$ maps $A^\T$ into $B\otimes M_{k(n)}$
for all {{$k(n)\ge k(n_1)\ge N.$}} Let $\eta=1/k(n_1)>0.$
Note that $1/k(n_1)\in K_0(U)_+.$     For any $0<r<\eta,$ choose a projection $p\in U$ such that $\tau_U(p)=rk(n_1)$
($<1$).  Recall {{that}} $U\otimes U\cong U.$
Let $\Psi: B\otimes M_{k(n_1)}\to
B\otimes U\otimes U$ be the \hm\, defined by $\Psi(b\otimes c)=b\otimes j_{n_1}(c)\otimes p$ for all $c\in U.$
Define $L=\Psi\circ \Phi$ and $\phi=L|_A.$ Let $u=1+(L((z-1)\otimes 1_{\td A}),$
where $z\in C(\T)$ is the standard unitary generator of $C(\T).$
Put $\kappa:=[\Psi]\circ \kappa''$ and
 write $\Phi_{\phi, u}=L.$  Note that $\kappa,$ $\phi$ and $u$ meet the requirements.

The proof for the ``Moreover" part is exactly the same, but replacing
$B$ by $B_i$ with the obvious modification.  In particular,
 we let  $\af\in Hom_{\Lambda}(\underline{K}({\tilde A_1}\otimes C(\T)), \underline{K}(B_i))$  {{be}}
 such that $\af|_{\underline{K}(A_1)}=0$ and by the UCT, one may also require
 that $\af|_{\boldsymbol{\bt}(K_0(A_1))}=0,$  $\af|_{{\boldsymbol{\bt}}(K_0(A_1)/mK_1(A_1))}=0,\,m=2,3,...,$
 $\af({\boldsymbol{\bt}}(g_1))=1_{\Z/m_i\Z},$
 $\af({\boldsymbol{\bt}}(g_j))=0$ if  $j\not=i,$ and $\af({\boldsymbol{\bt}}({\rm Tor}(K_1(A_1))))=0.$
 The rest of the proof remains the same.
So we only have \eqref{Lbotz0-2+}.

 \end{proof}

 \begin{lem}\label{LZ2toB}
 Let $B\in {\cal M}_1$ be a simple \CA\, with continuous scale.
 Suppose that $y\in K_0(B)\setminus \{0\}$ with
 $|\rho_B(y)(\tau)|<r$ for all $\tau\in T(B)$ and for some $0<r<1.$
 Then there exists a  nonzero \hm\, $\phi: {\cal Z}_{2,o}\to B$ such that
 $\phi_{*0}(1_\Z)=y.$

 Suppose that $B=B_1\otimes U,$  where $U$ is an infinite dimensional UHF-algebra, {{and}}
 $m\ge 2$ is an integer {{with}} $my=0.$
 Then there exists a nonzero \hm\, $\phi: {\cal Z}_o^m\to B$
 such that $\phi_{*0}(1_{\Z/m\Z})=y.$
  \end{lem}

 \begin{proof}
 Note that $t_{o,+}(1_\Z)=1$ and $t_{o,-}(1_\Z)=-1$ and $\Aff(T({\cal Z}_{2,o}))\cong \R^2.$
 Define $\Lambda: \Aff(T({\cal Z}_{2,o}))\to \Aff(T(B))$ as follows:
 {{We identify $\Aff(T({\cal Z}_{2,o}))$ with $\R^2;$ define}} $\Lambda((1,-1))=\rho_B(y)$
 and $\Lambda((1,1))=e\in \Aff(T(B)),$ where $e(\tau)=1$ for all $\tau\in T(B).$
 Note that
 $$
 \Lambda((1,0))(\tau)=(1/2)\Lambda((1,1)+(1,-1))(\tau)=(1/2)(e+\rho_B(y))(\tau)>0
 $$
 for all $\tau\in T(B).$
 Also
 $$
 \Lambda((0,1))=(1/2)(\Lambda((1,1)-(1,-1))(\tau)=(1/2)(e-\rho_B(y))(\tau)>0
 $$
 for all $\tau\in T(B).$ This $\Lambda$ induces a \hm\,
 $\af: Cu^{\sim}({\cal Z}_{2,o})\to Cu^{\sim}(B)$ which maps
 $1_\Z$ to $y.$  Since $B$ has stable rank one (see \ref{DD0}), it follows from \cite{Rl} that
 there is a \hm\, $\phi: {\cal Z}_{2,0}\to B$ such that $Cu^\sim(\phi)=\af.$
 One checks that $\phi$ meets the requirements.

For the second part of the lemma, note $\rho_B(y)=0.$
Let $\kappa\in KL({\cal Z}_o^m, B)$ {{be}} induced
by $\kappa([1_{\Z/m\Z}])=y.$  Define $\kappa_T: T(B)\to T({\cal Z}_0^m)$
by $\kappa_T(\tau)=\tau_o$ for all $\tau\in T(B).$
Define $\kappa_{cu}: U(\td {\cal Z}_o^m)/CU(\td {\cal Z}_o^m)=\R/\Z\to U(\td B)/CU(\td B)$
to be the map {{induced}} by $\kappa_T$ from $\R/\Z\to \Aff(T(\td B))/\overline{\rho_{\dt B}(K_0(B))}.$
Then $(\kappa, \kappa_T, \kappa_{cu})$ is compatible. By {{Theorem}} \ref{Text1} (note that ${\cal Z}_o^m$ is in the form of (3) of Remark 4.32 of \cite{GLrange}), there is a \hm\,
$\phi: {\cal Z}_o^m\to B$ such that $[\phi]=\kappa,$ $\phi_T={{\kappa_T}}$ and
$\phi^\dag=\kappa_{cu}.$  In particular, $\phi_{*0}(1_{\Z/m\Z})=y.$

 \end{proof}

\begin{lem}\label{Lbotextfinite}
 Let $A\in {\cal M}_1$ {{be a simple $C^*$-algebra as constructed in Theorem 4.31 of \cite{GLrange}}} and let {{$U$}} be an infinite dimensional  {{UHF-algebra}}.
Then, for any $\ep>0,$ any finite {{subsets}} ${\cal F}\subset A$ and
 ${\cal P}\subset \underline{K}(A),$
there exists $\eta>0$  and a finite subset ${\cal Q}\subset \underline{K}(A)$
with  ${\cal P}\subset {\cal Q}$ and
$G({\cal Q})\cap K_1(A)\cong \Z^k\oplus {\rm Tor}(K_1(A))\cap G({\cal Q})$
{{that}} satisfy the following:
 if $\af_0: K_1(A)\cap G({\cal Q})\to K_0(C\otimes U)$ is a \hm,
 where $C\in {\cal D}$ is a separable simple \CA\, with {{continuous scale,}} 
 such that
 \beq\label{Lbotextfinite-1}
 \|\rho_{C\otimes U}(\af_0(g_j))\|<\eta,\,\,\,j=1,2,...,k,
 \eneq
 where $\{g_1, g_2,...,g_k\}$ is a free generating set of $\Z^k,$ {{then}}
there exists an element $\af\in  KL_{loc}^{A^T}(G({\cal Q})\cup\boldsymbol{\bt}(G({\cal Q})), \underline{K}(C\otimes U)),$
an ${\cal F}$-$\ep$-multiplicative \cpc\,
$\phi: A\to C\otimes U$ and a unitary $u\in U({\widetilde{B_1\otimes U}})$ such that
\beq\label{Lbotextfinite-2}
&&\|[\phi(a),\, u]\|<\ep\rforal a\in {\cal F},\\\label{Lbotextfinite-3}
&&\hspace{-0.8in}\af|_{{\boldsymbol{\bt}}(K_1(A)\cap {\cal Q})}=\af_0\circ {\boldsymbol{\bt}}^{-1},\,\, \af|_{{\boldsymbol{\bt}}(K_0(A)\cap {\cal Q})}=0,\,\,
\af|_{{\boldsymbol{\bt}}((K_0(A)/mK_0(A))\cap {\cal Q})}=0,\,\,\,m=2,3,...,\\
&&\af|_{{\cal Q}\cup {\boldsymbol{\bt}}({\cal Q})}=[\Phi_{\phi, u}]|_{{\cal Q}\cup {\boldsymbol{\bt}}({\cal Q})},\\
&&{{{\rm bott}_1}}(\phi, u)=\af_0, \,\,\, [\Phi_{\phi, u}]|_{{\boldsymbol{\bt}}({\cal Q}\cap K_0(A))}=0\andeqn\\
&&{[}\phi{]}|_{\cal Q}=\af|_{\cal Q}\andeqn \af|_{\underline{K}(A)}=0.\,\,\,
\eneq

\end{lem}

\begin{proof}
Let ${\rm Tor}(G({\cal Q})\cap K_1(A))=(\Z/m_1\Z)g_{k+1}\bigoplus (\Z/m_2\Z)g_{k+2}\bigoplus\cdots \bigoplus (\Z/m_{k+N}\Z)g_{k+N},$ where $\{g_{k+1}, g_{k+2},...,g_{k+N}\}$ forms a  basis for ${\rm Tor}(G({\cal Q})\cap K_1(A))$
and $g_{k+j}$ has order ${{m_j}},$
$m_i=p_i^{n_i},$ $p_i$ is a prime number and $n_i\ge 1$ is an integer, $i=1,2,...,N.$
Let $1_{\Z/m_i\Z}$ be a generator of order $m_i,$  $i=1,2,...,N.$
Let $p_i^{n_i'}$  be the order of $1_{\Z/m_i\Z}\otimes [1_U]$ in $(\Z/m_i\Z)\otimes K_0(U),$
$0\le n_i'\le n_i,$ $i=1,2,...,N.$   We may assume that $n_i'>0,$ for
$i=1,2,...,N_0<N$ and $n_i'=0$ if $N_0<i\le N.$
Write  $U=M_{\p}$. Note that for each $p_i$, there is a nonzero element of $(\Z/m_i\Z)\otimes K_0(U)$ of order $p_i$. Consequently each prime factor $p_i$ appears in $\p$ at most finitely many times. On the other hand, since  $\p$ is a super-nature number ($M_{\p}$ is infinite dimensional), $\p$ has
infinitely many {{(possible repeating)}} prime factors $\{q_1,q_2,...,\},$ each of which is relatively prime
to all $p_1^{n_1'},p_2^{n_2'},...,p_{N_0}^{n_{N_0}'}.$
Choose a product ${\bar q}$ of   finitely many $\{q_1, q_2,...,\}$ which is minimum
among those products such that
$
{\bar q}\ge 6(k+N).
$
Put $B=C\otimes U.$ Let $e_1, e_2,...,e_{k+N_0}\in U$ be mutually
orthogonal non-zero projections such that $\sum_{j=1}^{k+N_0}e_i\le 1_U$
and $\tau_U(e_j)=1/{\bar q}.$
Note that $1_{\Z/m_i\Z}\otimes e_{k+i}$
has {{order}} $p_i^{n_i'},$ $1\le i\le N_0,$
where $\tau_U$ is the unique tracial state of $U.$

  Define
$B_i=(1\otimes e_i)B(1\otimes e_i),$ $i=1,2,...,k+N_0.$
Put $D={\cal Z}_{2,o}\otimes U,$  $D_i={\cal Z}_o^{m_i}\otimes U,$
$i=1,2,...,N_0.$
It follows from \ref{Lbotz2} that, there exists $\eta_0>0,$ and
a finite subset ${\cal Q}\subset \underline{K}(A)$ with ${\cal P}\subset {\cal Q},$
$\kappa^{(i)}\in KL_{loc}^{A^\T}(G({\cal Q})\cup \boldsymbol{\bt}(G({\cal Q})), \underline{K}(D))$
($1\le i\le k$),
${{\kappa^{(k+j)}}}\in KL_{loc}^{A^\T}(G({\cal Q}),\cup \boldsymbol{\bt}(G({\cal Q})), \underline{K}(D_j))$
($1\le j\le N$),
{{${\cal F}$-$\ep$}}-multiplicative \cpc s
$\psi_i: A\to D$  and unitaries $u_i\in D$ ($1\le i\le k$),
$\psi_{k+j}: A\to D_j$ and unitaries {{$u_{k+j}\in D_j$}} ($1\le j\le {{N_0}}$),
 such that $0<r_i<\eta_0<1,$
\beq\label{Lbotextfinite-10}
&&\kappa^{(i)}({\boldsymbol{\bt}}(g_i))=1_\Z\otimes r_i,\,\,\,\kappa^{(i)}({\boldsymbol{\bt}}(g_j))=0,\,\,\,j\not=i
\text{(for}\,\,1\le i\le k),
\\
&&\kappa^{(k+j)}({\boldsymbol{\bt}}({{g_{k+j}}}))=1_{\Z/m_j\Z}\otimes r_{k+j},\,\,\, \kappa^{(i)}({\boldsymbol{\bt}}({{g_{k+j}}}))=0,
\,\,i\not=k+j,\andeqn,\\
&&\kappa^{(i)}|_{{\boldsymbol{\bt}}(K_0(A))}=0,\,\,
\kappa^{(i)}|_{{\boldsymbol{\bt}}(K_0(A)/mK_0(A))}=0,\\
&&{[}\Phi_{\psi_i, u_i}{]}|_{{\cal P}\cup {\boldsymbol{\bt}}({\cal P})}=\kappa^{(i)}|_{{\cal P}\cup {\boldsymbol{\bt}}({\cal P})},\\
&&{[}\psi_i{]}|_{\cal P}=0\andeqn\\
&&\|[\phi(a),\, u_i]\|<\ep\rforal a\in {\cal F},
\eneq
$i=1,2,...,k+N_0.$
Here we   choose $r_{j}=1/{\bar q}_j,$ where ${\bar q}_j$ is a finite product
of prime factors $\{q_1,q_2,...,\}.$ As a consequence,
 $1_{\Z/m_j\Z}\otimes r_{k+j}$ has order $p_j^{n_j'},$ $1\le j\le N_0.$
\Wlog, we may assume  that $u_i=1_{\td D}+x_i$ for some $x_i\in D,$ $1\le i\le k,$ and
$u_{k+j}=1_{\td D_j}+x_{k+j}$ for some
$x_{k+j}\in D_j,$  $j=1,2,...,{{N_0}}.$

Let $\eta=\min\{{{r_{j}\eta/{\bar q}}}: 1\le j\le k+N_0\}$ and $\af_0$ satisfy the condition \eqref{Lbotextfinite-1}.
Note that $K_0(B)=K_0(B)\otimes K_0(U)$ and $\af_0(g_{k+j})=0$ for $N_0<j\le N.$
 Note {{also}} that
 \beq{{
 \frac{1}{r_i\tau_U(e_i)}={\bar q}_i {\bar q}\andeqn
 (\frac{1}{r_i\tau_U(e_i)})\|\rho_{C\otimes U}\af_0(g_i)\|\leq \eta_0<1.}}
 \eneq
For each $i,$ by \ref{LZ2toB},  there exists a \hm\, $\phi_i: D\to B_i\otimes U\cong B_i$
(with the form $\phi_i'\otimes \id_U$ where $\phi_i': {\cal Z}_{2,0}\to B_i$)
such that  $(\phi_i)_{*0}(1_\Z)={\bar q}_i{\bar q}\af_0(g_i)\otimes e_i,$ $i=1,2,...,k,$
and  a \hm\, $\phi_{k+i}: D_i\to B_{k+i}$
such that $(\phi_{k+i})_{*0}(1_{\Z/m_i\Z})={\bar q}_{k+i}{\bar q}\af_0(g_{k+i})\otimes e_{k+i},$ $j=1,2,...,N_0.$
Thus $(\phi_j\circ \psi_j)_{*0}(g_j)=\af_0(g_{j})$ and $(\phi_j\circ \psi_j)_{*0}(g_i)=0,$
if $i\not=j,$  $1\le j\le N_0.$

Let $\af\in KL_{loc}^{A^\T}(G({\cal Q})\cup\boldsymbol{\bt}(G({\cal Q})),
\underline{K}(B))$
{{satisfy}}
 \eqref{Lbotextfinite-3}
and
$\af|_{\underline{K}(A)}=0.$
Put $u=1_{\td B}+\sum_{i=1}^{{k+N_0}} \phi_i(x_i).$ One checks that $u$ is a unitary in $\td B.$
Define $\phi: A\to B$ by
$$
\phi(a)=\diag(\phi_1\circ \psi_1(a), \phi_2\circ \psi_2(a),...,\phi_{{k+N_0}}\circ \psi_{k+N_0}(a))
$$
for all $a\in A.$ One then checks
that $\af,$ $\phi$  and $u$ satisfy the requirements.

\end{proof}
%

\begin{lem}\label{Tbotext}
 Let $A\in {\cal M}_1$ be {{as in (3) of Remark 4.32 of \cite{GLrange} (see Theorem \ref{Text1} also)}} with continuous scale.
 Suppose that {{$K_i(A)$  is finitely generated ($i=0,1$)}}  and $K_1(A)=\Z^k\oplus {\rm Tor}(K_1(A)),$ where $\Z^k$ is generated by
 $g_1, g_2,...,g_k.$  Then, for any $\ep>0,$ any finite subset ${\cal F}\subset A$ and
 any finite subset ${\cal P}\subset \underline{K}(A)$ with
 $\{g_1, g_2,...,g_k\}\subset {\cal P}\cap K_1(A),$ there exists $\eta>0$
 {{satisfying}} the following:
if {{$\af\in KL(A^\T, B),$}} where $B=C\otimes U,$
$C\in {\cal D}$  is an amenable simple \CA\, with continuous scale  and
$U$ is a UHF-algebra of infinite type
such that
\beq\label{Tbotext-1.1}
|\rho_B({{\af({\boldsymbol{\bt}}(g_i))}}(\tau)|<\eta\, \tforal \tau\in T(B),\,\,\, i=1,2,...,k,
\eneq
 and if $\phi: A\to B$ is  a \hm\, which maps strictly positive elements to
 strictly positive elements, {{then}}
there exists a unitary $u\in CU({\widetilde{B}})$ such that
\beq\label{Tbotex-2}
&&\|[\phi(a),\, u]\|<\ep\rforal a\in {\cal F},\\\label{Tbotex-3}
&&{\rm Bott}(\phi, \, u)|_{\cal P}=\af({\boldsymbol{\bt}})|_{\cal P}.
\eneq

\end{lem}

  \begin{proof}
  Let $e_a\in A$ be a strictly positive element of $A$ with $\|e_a\|=1.$
Since $A$ has continuous scale, \wilog, we may assume that
\beq\label{11ExtT1-5.1}
\min\{\inf\{\tau(e_a): \tau\in T(A)\},\inf\{\tau(f_{1/2}(e_a)):\tau\in T(A)\}\}>3/4.
\eneq
Let $T: A_+\setminus \{0\}\to \N\times \R_+\setminus \{0\}$ be given by
Theorem 5.7 of \cite{eglnp} {{corresponding to $e_a$ (in place of $e$) and $3/8$ (in place of $d$).}}

Note that, by {{A.10 of}} \cite{GLII}, both $A$ and $B$ are ${\cal Z}$-stable.
Thus, by 5.5 of \cite{GLII},  $K_0({\tilde B})$ is weakly unperforated.
Let ${\bf T}(n,k)=n$ for all $(n,k)\in \N\times \N.$ Then, as explained in 5.2 of \cite{GLII},
${\cal D}={\cal D}_{{\bf{T}}(n,k)}.$

We now apply 5.3 of \cite{GLII}.
Let $\ep>0$ and  a finite subset ${\cal F}$ be given.
We may assume that
\beq\label{Tbote-n1}
{\rm Bott}(\Phi, \, u)|_{\cal P}={\rm Bott}(\Phi',\, u)|_{\cal P},
\eneq
if $\|[\Phi(a),\, u]\|<\ep,$  $\|[\Phi'(a),\, u]\|<\ep$  and
$$
\|\Phi(a)-\Phi'(a)\|<\ep\rforal a\in {\cal F}
$$
for any ${\cal F}$-$\ep$-multiplicative \cpc s $\Phi, \Phi': A\to B.$
Note that $U$ is assumed to {{be of}}  infinite type and hence is strongly absorbing.
We identify $A$ with $A\otimes 1_U$ as a \SCA\, of $A\otimes U$ and $B$ with $B\otimes 1_U$
as a \SCA\, of {{$B\otimes U.$}}

We may assume that ${\cal F}\subset A^{\bf 1}.$
Let $\dt>0, \gamma>0,$ $\eta_0>0$ (in place of $\eta$),  ${\cal G}'\subset A\otimes U,$  ${\cal H}_1'\subset (A\otimes U)_+\setminus \{0\},$
${\cal P}_a'\subset {\underline{K}}(A\otimes U)$ (in place of ${\cal P}$), $\{v'_1, v'_2, ..., v'_{m_0'}\}\subset  U({\widetilde{A\otimes U}})$ such that
${\cal P}_a'\cap K_1(A\otimes U)=\{[v_1'], [v_2'],...,[v_{m_0'}']\}$ and ${\cal H}_2'\subset (A\otimes U)_{s.a.}$ be finite subset
required  by 5.3 of \cite{GLII} for $\ep/4$ (in place of $\ep$) and ${\cal F}$  {{and $T$ above.}}
\Wlog, we may assume that ${\cal H}_2'$ is in the unit ball of $(A\otimes U)_{s.a.}.$
We may also assume that ${\cal H}_1'={\cal H}_1\otimes {\cal H}_1^U$ and
${\cal H}_2'={\cal H}_2\otimes {\cal H}_2^U,$ where
${\cal H}_1\subset A_+^{\bf 1}\setminus \{0\},$ ${\cal H}_1^U\subset U_+^{\bf 1}\setminus \{0\},$
and ${\cal H}_2\subset A_{s.a.}^{\bf 1}$ and ${\cal H}_2^U\subset U_{s.a.}^{\bf 1}$ are finite subsets.

We may {{further}} assume that ${\cal P}_a'\subset \{z\otimes y: z\in {\cal P}_a, y\in {\cal P}^U\},$
where ${\cal P}_a\subset \underline{K}(A)$ is a finite subset and
${\cal P}^U\subset K_0(U)$ is also a finite subset.

Without loss of generality, we may assume
that ${\cal G}'=\{g\otimes f: g\in {\cal G} \andeqn f\in {\cal G}^U\},$ where
${\cal G}\subset A$ is a finite subset and $1_U\in {\cal G}^U\subset U$ is also a finite subset.
{{Since $K_1(A\otimes U)=K_1(A)\otimes K_0(U)$, without loss of generality, we may assume $v'_i=1_{(A\otimes U)^{\sim}}+y_i\otimes f_i$ for $1\leq i\leq m_0'$, where $y_i\in A$ and $f_i\in U$, such that $f_i$ is a projection and $v_i=1_{\tilde A}+y_i$ is a unitary. }}

Recall that $K_i(A)$ is finitely generated ($i=0,1$). We may assume further {{that  $[L]$ is a well-defined
element  of}} $KL(A, B)$   for any ${\cal F}$-$\ep$-multiplicative \cpc\, $L: A\to B.$
Note also {{that}} $K_i(A^\T)$ is also finitely generated ($i=0,1$).

Put $\eta_{00}=\min\{\eta_0, \gamma/2\}.$
We identify $U$ with $U\otimes U.$  \Wlog, we may assume that $\phi$ maps {{$A$ into
$B\otimes 1_U.$}}

 Let $e_0'\in U$ be  {{a non-zero projection}} such that $t_U(e_0')<\eta_{00}/2,$ where
$t_U$ is the unique tracial state of $U,$  and
let $e_0=1_{\tilde C}\otimes 1_U\otimes e_0'.$
Put $r_0=t_U(e_0)\in K_0(U)$ {{(regarding $K_0(U)$ as a subgroup of $\R$).}}
Let  $B_0=e_0(C\otimes U\otimes U)e_0.$ Note {{that}} $B_0$ also has continuous scale.
%
Let ${\cal P}_1={\cal P}_a\cup {\cal P}.$

Recall that $K_1(A)=\Z^k\oplus {\rm Tor}(K_1(A)).$ Write $\Z^k=\Z g_1 \oplus \Z g_2\oplus \cdots \Z g_k.$
{{Therefore,}} \wilog, we may assume that $g_i=[v_i],$ $i=1,2,...,k,$ and $m_0=k.$
By 2.11 of \cite{DL},  since $K_i(A)$ is finitely generated ($i=0,1$), there exists $K\ge 1$ such that
\beq\label{72-3+1}
Hom_{\Lambda}(F_{K}\underline{K}(A),\, F_{K}\underline{K}(B))=Hom_{\Lambda}(\underline{K}(A), \underline{K}(B))
\eneq
(see also \cite{DL} for the notation $F_K$ there).  Thus, there is a finite subset
${\cal Q}_1\subset \underline{K}(A)$ such that $G({\cal Q}_1)=F_K\underline{K}(A)$
and ${\boldsymbol{\bt}}(G({\cal Q}_1)={\boldsymbol{\bt}}(F_K\underline{K}(A)).$

{{Since $K_*(A)$ is finitely generated, the inductive limit algebra in (3) of Remark 4.32 of \cite{GLrange} also in the class of Theorem 4.31 of \cite{GLrange}. Hence we can apply \ref{Lbotextfinite}.}} Define $\eta_1$ to be the number $\eta$
and ${\cal Q}$ in \ref{Lbotextfinite} associated with $\ep/4,$ ${\cal F}$
and ${\cal P}_1$ (as well as $A$ and  $U$).  As mentioned above, we may assume  ${\cal Q}={\cal Q}_1{{\supset {\cal P}_1}}.$
Put $\eta=t_U(e_0')\cdot \eta_1/2.$ We may assume that $\eta_1<1/4.$

Suppose ${{\af\in KL(A^\T, B)}}$ satisfies the assumption of the {{lemma}} for the above $\eta.$
{{Let $\phi: A\to B$ be also  {{a}} given homomorphism satisfying the assumption of the {{lemma}}.}}  
Then we obtain an element $\af_1\in Hom_{\Lambda}(\underline{K}(A^\T), \underline{K}(B))$
such that ${\af_1}|_{\underline{K}(A)}=0$ and
${{\af_1\circ {\boldsymbol{\bt}}}}|_{F_K\underline{K}(A)}=\af\circ {\boldsymbol{\bt}}.$
Let
$
{{\af_0={\af_1}\circ{\boldsymbol{\bt}}|_{(K_1(A))}:}} K_1(A)\to K_0(B).$ Note $K_0(B_0)=K_0(B).$
Then
\beq
|\rho_{{\tilde B}_0}(\af_0(g_j))(\tau)|<\eta_1\rforal \tau\in T(B_0),\,\,\, j=1,2,...,k.
\eneq

It follows from \ref{Lbotextfinite} that there {{exist}}
${{\af'}}\in KL(A^\T, B_0),$
{{an ${\cal F}$-$\ep/4$-multiplicative \cpc\,}}
$\phi_0: A\to B_0$ and a unitary $u_0\in {\widetilde{B_0}}$ such that
\beq\label{Tbotext-10}
&&\|[\phi_0(a),\, u_0]\|<\ep/4\rforal a\in {\cal F},\\\label{Tbote-10+1}
&&{\rm Bott}(\phi_0,\,u_0)|_{{\boldsymbol{\bt}({{\cal Q}})}}={{\af'}}|_{{\boldsymbol{\bt}({{\cal Q}})}},\\
&&{{\af'}}|_{\boldsymbol{\bt}(K_1(A))}=
{{\af_0\circ\boldsymbol{\bt}^{-1}}}|_{\boldsymbol{\bt}(K_1(A))},\\\label{Tbotext-10+n1}
&&[\phi_0]|_{{\cal Q}}={{\af'}}|_{{\cal Q}}\andeqn {{\af'|}}_{\underline{K}(A)}=0.
\eneq
We may write $u_0=1_{\tilde B_0}+\zeta_0,$ where $\zeta_0\in B_0.$


Let $e_1'\in {{(1-e_0')U(1-e_0')}}$
be a non-zero projection with
$t_U(e_1')<\eta_{00}/3$
and let $e_1=1_{\tilde C}\otimes 1_U\otimes e_1'.$
Put $B_1={{(1-e_0-e_1)B(1-e_0-e_1)}}.$
{{Let  {{$r_1:=r_0+t_U(e_1')<\eta_{00}$}} and $s_1: B\to {{B_1}}$  be
defined by $s_1(c\otimes a)=c\otimes a\otimes (1-e_1{{-e_0}})$
for all
$c\in C$ and $a\in U.$   Then $[s_1](x)=(1-r_1)x$ for all $x\in \underline{K}(B).$}}
Define
$\kappa|_{\underline{K}(A)}= ([\phi]-
[\phi_0])\times [s_1]=[s_1\circ \phi]$ (as $[\phi_0]=0$)
and $\kappa|_{\boldsymbol{\bt}(\underline{K}(A))}=(\af-{{\af'}})|_{\boldsymbol{\bt}(\underline{K}(A))}.$

Define $\kappa_T: T(B_1)\to T_f(A^\T)$ as follows:
$\kappa_T(\tau)(1_{C(\T)}\otimes a)=\tau(\phi(a)\otimes (1-e_1-e_0))$ for all $\tau\in T(B_1)$ and
$a\in A,$ and
$$\kappa_T(\tau)(f\otimes 1_{{\tilde A}})=\int_\T f(t) dm(t)$$
for all $f\in C(\T),$ where $m$ is the normalized Lebesgue measure on {{$\T.$}}

Note that $K_0(A^\T)=K_0(A)\oplus {\boldsymbol{\bt}}(K_1(A)).$ Using \ref{LKT2}, one checks
that   $\kappa_T$ and $\kappa$ {{are}} compatible.  Since $B_1$ and $A$  have continuous {{scales,}}
$T(B_1)$ and $T(A)$  are compact. Therefore $\kappa_T(T(B_1))$ lies in a compact subset of $T_f(A^\T).$
It follows from \ref{Hext1} that there exists a sequence of \cpc s $\Psi_n: A^\T\to B_1$
such that
\beq\label{Tbote-11}
&&\lim_{n\to\infty}\|\Psi_n(ab)-\Psi_n(a)\Psi_n(b)\|=0\rforal a, b\in A^\T,\\\label{Tbote-11+}
&&[\{\Psi_n\}]=\kappa,\\\label{Tbote-11++}
&&\lim_{n\to\infty}\sup\{|\tau\circ \Psi_n(a)-\kappa_T(\tau)(a)|: \tau\in T(B)\}=0\rforal a\in A^\T.
\eneq
Define $\psi_n: A\to B_1$ by $\psi_n(a)=\Psi_n(a)$ for all $a\in A\subset A^\T.$  One checks that,
for all large $n,$
\beq\label{2020-907-1}
[\psi_n]|_{{\cal P}_1}=[s_1\circ \phi]|_{{\cal P}_1}.
\eneq
There is also a unitary $w_n\in {\tilde B_1}$ with $w_n=1_{{\tilde B}_1}+x_n,$
where $x_n\in B$ such that
\beq\label{Tbote-12}
\lim_{n\to\infty}\|\Psi_n((1-z)\otimes 1_{\tilde A})-x_n\|=0,
\eneq
where $z\in C(\T)$ is the identity function on $\T.$ It follows
that
\beq\label{Tbote-13}
\lim_{n\to\infty}\|[\psi_n(a), \, w_n]\|=0\rforal a\in A.
\eneq

Since we assume that $\phi$ maps strictly positive elements to strictly positive elements,  by \eqref{Tbote-11++}  and
by passing to a subsequence, we may assume that
\beq\label{Tbote-14}
\tau(f_{1/2}(\psi_n(e_a)))>3/8\rforal \tau\in T(B).
\eneq

View $\phi_0$ and $\psi_n$ as maps from $A$ {{into}} $B.$
Let $G_u=J_{cu}(K_1(A)),$
where $J_{cu}:=J_{cu}^{\td A}: K_1(A)\to CU({\tilde A})/CU({\tilde A})$ is the splitting map in \ref{DkappaJ}.

Let $B_2={{(e_1+e_0)B(e_1+e_0)}}$ and
$s_2: B\to B_2$ by $s_2(c\otimes a)=c\otimes a\otimes (e_1+e_0)$  for all $c\in  C$ and $a\in U.$
{{Let $\iota_e: B_2\to B$ be the embedding. Then
${\bar \iota_e}: U(\td B_2)/CU(\td B_2)\to
 U(\td B)/CU(\td B)$ induced by $\iota_e$}}
  is an isomorphism (see also {{Proposition}} 4.5 of \cite{GLII}).
{{Denote by
 ${\iota_e}_{*1}: K_1({{B_2}})\to K_1(B)$ the  isomorphism
  induced by $\iota_e.$
}}  Then
 ${\iota_e}_{*1}\circ {{\Pi_{cu}^{\td B_2}}}=\Pi_{cu}^{\td B}\circ {\bar \iota_e}.$
Let $\lambda_n': G_u\to U({\tilde B})/CU({\tilde B})$ be defined
by $\lambda_n'={{(\phi^{\dag}-\psi_n^{\dag}-\phi_0^{\dag})|_{G_u}}}.$
Define $\gamma_{n,0}: U({\tilde A})/CU({\tilde A})\to  {{U(\td B_2)/CU(\td B_2)}}$
as follows:  $\gamma_{n,0}|_{G_u}=({\bar \iota_e})^{-1}\circ \lambda_n'$ and
$\gamma_{n,0}$ on $\Aff(T({{\tilde A}}))/\overline{\rho_{{\tilde A}}(K_0({{\tilde A}}))}$ 
is induced by ${{(s_2\circ \phi)_T.}}$  Let $\kappa^0\in KL(A, B_2)$ be defined by $\kappa^0|_{{\underline{K}(A)}}=[s_2\circ \phi].$
Recall that $[\phi_0]=0${{,}} $[\psi_n]|_{{\cal P}_1}=[s_1\circ \phi]|_{{\cal P}_1}$ {{and $[s_1\circ \phi]+[s_2\circ \phi]=[\phi]$.}}
Then $(\kappa^0, (s_2\circ\phi)_T, \gamma_{n,0})$ is compatible.

Therefore, by Theorem \ref{Text1}, there is a \hm\, $\psi_{n,0}: A\to B_2$
such that
\beq\label{Tbote-15}
[\psi_{n,0}]=\kappa^0\,\,\,{\rm in}\,\,\, KL(A, B_2),\\
(\psi_{n,0})_T=
(s_2\circ\phi)_T
\andeqn (\psi_{n,0})^{\dag}={{\gamma_{n,0}}}.
\eneq
Define $\theta_n: A\to M_2(B)$ by
$\theta_n(a)=\diag(\psi_n(a), \psi_{n,0}(a), \phi_0(a))\rforal a\in A.$
{{Note that
\beq\label{2020-907-2}
\theta_n^\dag|_{G_u}|=(\psi_n^\dag+\psi_{n,0}^\dag+\phi_0^\dag)|_{G_u}=\phi^\dag|_{G_u}.
\eneq}}
Define $u_n=1_{\tilde B}+\diag(x_n, \zeta_0).$
Then $\Phi_{\theta_n, u_n}$ is ${\cal F}$-$\ep$-multiplicative.
Note {{that,}}
by \eqref{Tbote-14}, we have
\beq\label{Tbote-21}
\tau(f_{1/2}(\theta_n(e_a)))>3/8\tforal \tau\in T(B).
\eneq
{{Since $T$ is chosen as in Theorem 5.7 of \cite{eglnp}}},  for all sufficiently large $n,$
$\theta_n$ is exactly $T$-${\cal H}_1$-full.
We also have {{(by \eqref{Tbotext-10+n1}, \eqref{2020-907-1}, \eqref{Tbote-15}, \eqref{2020-907-2} and $t_U(e_0')<\eta_{00}/2$)}}
\beq\label{Tbote-22}
&&[\theta_n]|_{{\cal P}_a}=[\phi]|_{{\cal P}_a},\\
&&\lim_{n\to\infty}{\rm dist}(\theta_n^{\dag}(f),\phi^{\dag}(f))<\eta_{00}\rforal f\in U({\tilde A})/CU({\tilde{A}})\andeqn\\
&&\sup\{|\tau\circ \theta_n(a)-\tau\circ \phi(a)|: \tau\in T(B)\}<\eta_{00}\tforal a\in {\cal H}_2.
\eneq
{{Define ${\bar \kappa}_T: T(B)\to T_f(A^\T)$ by
${\bar \kappa}_T(\tau))(1_{C(\T)}\otimes a)=\tau(\phi(a))$
for all $a\in A_{s.a.},$ and ${\bar \kappa}_T(\tau)(f\otimes 1_{\td A})=\int_\T f(t) dm(t)$
for all $f\in C(\T),$ where $m$ is the normalized Lebesgue measure on $\T.$}}
Then,  in fact, by  \eqref{Tbote-11++} and {{choices}} of $e_0$ and $e_1,$   we also have
\beq\label{2020-720-1}
\sup\{|\tau\circ \Phi_{{\theta_n, u_n}}(a)-{\bar \kappa}_T(\tau)(a)|: \tau\in T(B)\}<\eta_0\tforal a\in {\cal H}_2.
\eneq
Choosing a sequence of $\eta_{0,n}<\eta_0/n.$ {{Let $\theta_n'$ be of  the same form as $\theta_n$, but with $\eta_{0,n}$ in place of $\eta_0$. }} {{Then there is}} a sequence
of approximately multiplicative maps $L_n: A^\T\to M_2(B)$
such that  ($z\in C(\T)$ is the identity function on $\T$)
\beq\label{2020-720-2}
&&\lim_{n\to\infty}\|L_n(a)-\theta_n'(a)\|=0\rforal a\in A,\\\label{Tbote-12.1}
&& \lim_{n\to\infty}\|{{L_n((1-z)\otimes 1_{\tilde A})}}-\diag(x_n, \zeta_0)\|=0,\\
&&[\{L_n\}]|_{{\cal P}\cup{\cal P}_a}=[\{\Phi_{\theta_n', u_n}\}]|_{{\cal P}\cup {{{\cal P}_a},}}\\\label{Tbote-12+1}
&&\lim_{n\to\infty}{\rm dist}(L_n^{\dag}(f),\phi^{\dag}(f))=0\rforal f\in U({\tilde A})/CU({\tilde{A}})\andeqn\\\label{Tbote-12-2}
&&\sup\{|\tau\circ L_n(a)-{\bar \kappa}_T(a)|: \tau\in T(B)\}=0\rforal a\in A^\T.
\eneq
{{As mentioned in the end of the proof of \ref{Hext1},
applying the argument in the end of the proof \ref{ExtTBA},  \wilog, we may assume
that $L_n$ maps $A^\T$ to $B.$}}

Consider $L_n'': A\otimes U\to B\otimes U$ {{defined}} by $L_n''=L_n\otimes {\rm id}_U$ and
$\phi': A\otimes U\to B\otimes U.$  Then $\phi'$  is a \hm\, and $L_n''$ is exactly $T$-${\cal H}_1'$-full.
{{By \eqref{Tbote-22} and \eqref{2020-720-2},
$[L_n'']|_{{\cal P}_a'}=[\phi']|_{{\cal P}_a'}.$ {{By}}  \eqref{Tbote-12+1},
$\lim_{n\to\infty}{\rm dist}(L_n^{\dag}(v_j'),\phi^{\dag}(v_j'))=0,$ $1\le j\le m_0',$ as
$v_j'$ has the form $1_{(A\otimes U)^\sim}+{{y_j\otimes f_{j}}}$ given {{earlier.}}  Also, by \eqref{Tbote-12-2},
$\sup\{|\tau\circ L_n''(a)-\phi'(a)|: \tau\in T(B)\}=0\rforal a\in A\otimes U.$
}}

It follows from {{Theorem}} 5.3 of \cite{GLII} (see also 5.2 and {{Proposition}} 5.5 of \cite{GLII}) that, for all large $n,$
there exists  a unitary $W_n\in {\tilde B}$ such that
\beq\label{Tbote-23}
\|W_n^*L_n''(a\otimes 1_U)W_n-\phi'(a\otimes 1_U)\|<\ep/2\rforal a\in {\cal F}.
\eneq
Note that ${{L_n''}}(a\otimes 1_U)=L_n(a)\otimes 1_U$ and $\phi'(a\otimes 1_U)=\phi(a)\otimes 1_U.$
(Recall that we identify $\phi(a)$ with $\phi(a)\otimes 1_U$ for all {{$a\in A.$}})
Put $V_n=W_n^*u_nW_n.$ Then, for all large $n,$
\beq\label{Tbote-24}
\|[\phi(a), \,V_n]\|<\ep\rforal a\in {\cal F}.
\eneq
By the definition of $\kappa,$ $\theta_n,$
\eqref{Tbote-10+1}, \eqref{Tbote-11+} and \eqref{Tbote-n1},
and \eqref {2020-720-2},
we compute that (for all large $n$)
\beq\label{Tbote-25}
{\rm Bott}(\phi, \, V_n)|_{\cal P}=\af({\boldsymbol{\bt}}({\cal P})).
\eneq

Fix a sufficiently large $n$ and any sufficiently large finite subset ${\cal G}_A.$
Since $B$ also has stable rank one, by {{Lemma}} \ref{Unitary},  for any $\eta>0,$
there exists  $v_0\in U({\tilde B})$ such that
$\|v_0a-av_0\|<\eta$ for all $a\in {\cal G}_A$ such that $[v_0]=[V_n]$ and
$v_0^*V_n\in CU({\tilde B}).$ Put $u=v_0^*V_n.$
With sufficiently large ${\cal G}_A$ and sufficiently small $\eta,$ $u$ meets all requirements.
  \end{proof}

%
%

 \begin{thm}\label{TTbote}
  Let ${{A'}}\in {\cal M}_1$
  and ${{A:=A'}}\otimes U_1=\overline{\cup_{n=1}^\infty A_n},$ {{with $\lim_{n\to\infty}{\rm dist} (x, A_n)=0~~\mbox{for any}~~x\in A$,}}
  where $U_1$ is a UHF-algebra with infinite {{type and}} $A_n\in {\cal M}_1$ with
  finitely generated $K_i(A_n)$ ($i=0,1$) as in Lemma \ref{Lfinitelimit}.
Then, for any $\ep>0,$ any finite subset ${\cal F}\subset A$ and
 any finite subset ${\cal P}\subset \underline{K}(A),$
 there exists $\eta>0$  and   an integer $N\ge 1$
 satisfying the following: ${\cal P}\subset [\iota_n](\underline{K}(A_n))$
 (where $\iota_n: A_n\to A$ is the embedding, for all $n\ge N$),
if $\af\in  Hom_{\Lambda}(\underline{K}({{A^\T_n}}), \underline{K}(B))$
for some $n\ge N,$  and $K_1(A_n)=\Z^k\oplus {\rm Tor}(K_1(A_n)),$
where $B=C\otimes U_2$ which satisfies the UCT,  $C\in {\cal D}$  is an amenable simple \CA\, with continuous scale  and
$U_2$ is a UHF-algebra of infinite type
such that
\beq\label{Tbotext-1.2}
|{{\af({\boldsymbol{\bt}}(g_i))}}(\tau)|<\eta\, \tforal \tau\in T(B),\,\,\, i=1,2,...,k,
\eneq
 where ${\cal P}^{(0)}:=\{g_1, g_2,...,g_k\}$ is a set of free generators for $\Z^k{{\subset K_1(A_n)}},$
 and if $\phi: A\to B$ is {{a \hm\,}} which maps strictly positive elements to
 strictly positive elements,  {{then}}
there exists a unitary $u\in CU({\widetilde{B}})$ such that
\beq\label{Tbotex-2.1}
&&\|[\phi(a),\, u]\|<\ep\rforal a\in {\cal F},\\\label{Tbotex-3.1}
&&{\rm Bott}(\phi\circ \iota_n, \, u)|=\af({\boldsymbol{\bt}}).
\eneq

 \end{thm}

 \begin{proof}
 Fix $\ep>0,$ {{and finite subsets ${\cal F}$ and}} ${\cal P}.$
 We may choose $N\ge 1$ such
 that ${\cal P}\subset [\iota_n](\underline{K}(A_n))$ for all $n\ge N$ (see {{Lemma}} \ref{Lfinitelimit}). We may also assume
 that ${\cal F}\subset A_n.$
 {{Let}} $\phi_n: A_n\to B$ {{be}} defined by $\phi_n:=\phi|_{A_n}.$
 {{Note that the algebra $A_n$ in Lemma \ref{Lfinitelimit} is as (3) of Remark 4.32 of \cite{GLrange} with finiely generated K-theory.}} {{Hence Lemma \ref{Tbotext} can be applied}} to $\ep,$ ${\cal F},$  ${\cal P}$ and $\phi_n$ above.
 {{The theorem}} then follows.
 %
 %
 \end{proof}

\begin{df}\label{DB2} {{(Definition 3.6 of \cite{LN} and 7.2.6 of \cite{Lncbms})}}
Let $C$ be a separable \CA. Let $1/4>\Delta_c(t,{\cal F}, {\cal P}_0, {\cal P}_1,h)>0$ be a function defined on $t\in [0,1],$ the family of all finite subsets ${\cal F}\subset C,$ the family of all finite subsets ${\cal P}_0\subset K_0(C),$ and family of all finite subsets ${\cal P}_1\subset K_1(C),$  and the set of all \hm s
$h: C\to A$ for some \CA\, $A.$  We say that $A$ has Property (B2) associated with $C$ and $\Delta_c$ if the following holds:

For any
\hm\, $h: C\to A,$ any $\ep>0,$ any  finite subset ${\cal F}\subset C,$ any finite subset ${\cal P}_0\subset K_0(C),$ and any finite subset ${\cal P}_1\subset K_1(C),$ there are finitely generated subgroups $G_0\subset K_0(C)$  with ${\cal P}_0\subset G_0$ and $G_1\subset K_1(C)$  and a finite subset ${\cal Q}\subset G_1$ which generat{e}s $G_1,$
satisfying the following: for any homomorphisms $b_0: G_0\to K_1(A)$ and $b_1: G_1\to K_0(A)$ such that
\beq\label{Pbott-2}
|\rho_A(b_1(g))(\tau)|<\Delta_c(\ep, {\cal F}, {\cal P}_0, {\cal P}_1, h)
\eneq
for any $g\in {\cal Q}$ and any $\tau\in\mathrm{T}(A)$, there exists a unitary $u\in U(\td A)$ with $[u]=0$ in
$K_1(A)$ such that
\beq\label{Pbott-3}
\text{bott}_0(h, \, u)|_{{\cal P}_0}=b_0|_{{\cal P}_0},\,\,\, \text{bott}_1(h,\, u)|_{{\cal P}_1}=b_1|_{{\cal P}_1}\andeqn\\\nonumber
\|[h(c),\, u]\|<\ep\tforal c\in {\cal F}.
\eneq

If $C$ is unital, one may assume that $[1_C]\in {\cal P}_0.$ Then the first part of \eqref{Pbott-3} implies that
$[u]=0$ in $K_1(A),$ {{if we also assume that
$b_0([1_C])=0$.}}

Let $C$ be $A$ in Lemma \ref{TTbote},  $\ep>0,$  finite subsets ${\cal F}\subset A,$ ${\cal P}_0\subset K_0(A)$
and ${\cal P}_1\subset K_1(A)$ be given.   Set ${\cal P}={\cal P}_0\cup {\cal P}_1.$
Let $\eta>0$ and $N$ be given by Lemma \ref{TTbote}.   Note ${\cal P}\subset [\iota_n](\underline{K}(A_n)).$
Let $G_i=[\iota_n](K_i(A_n))$ ($i=0,1$).    One may write {{that}} $G_1=\Z^k\oplus {\rm Tor}(K_1(A_n))$
and {{assume that}} ${\cal P}_0=\{g_1, g_2,...,g_k\}\subset \Z^k\subset G_1$ is a set of free generators.
Put ${\cal Q}={\cal P}_0\cup G_{1,t},$ where $G_{1,t}$ is a finite subset generating ${\rm Tor}(K_1(A_n)).$
Let $B$ be as in Lemma \ref{TTbote}.
Suppose that $b_i: G_i\to K_i(B)$ is a \hm\, ($i=0,1$) such that
$|\rho_B(b_1(g))|<\eta$ for all $g\in {\cal Q}.$
By the UCT, there exists $\af'\in {\rm Hom}_{\Lambda}(\underline{K}(A_n), \underline{K}(B))$
such that $\af'\circ [\iota_n]|_{G_i}=b_i,$ $i=0,1.$
Define $\af\in {\rm Hom}_{\Lambda}(\underline{K}(A^\T), \underline{K}(B))$ as follows (see \eqref{june-19-2021}):
$\af|_{\underline{K}(A)}=0,$ $\af|_{{\boldsymbol{\bt}}(\underline{K}(A))}=\af'\circ {\boldsymbol{\bt}}$
and $\af|_{{\boldsymbol{\bt}}(\underline{K}(\C\cdot1_{\tilde A}))}=0.$
By Lemma \ref{TTbote}, there exists a unitary $u\in U_0(\td B)$ such that
\beq
\|[h(a),\, u]\|<\ep\rforal a\in {\cal F},\,
{\rm bott}_0(h,\, u)|_{{\cal P}_0}=b_0\andeqn  {\rm bott}_1(h,\,u)|_{{\cal P}_1}=b_1.
\eneq
In other words,
(with $\Delta_B=\eta$), $A$ satisfies Property (B2) associated with $B$ and $\Delta_B.$

\end{df}

  \begin{thm}\label{TTboteplus}
  Let $A$ be a simple \CA\, in {{ \ref{Misothm}}} 
  and $A=\overline{\cup_{n=1}^\infty A_n}$ {{with $\lim_{n\to\infty}{\rm dist} (x, A_n)=0~~\mbox{for any}~~x\in A$,}} {{as}} in Lemma \ref{Lfinitelimit}
  so that $K_i(A_n)$ is finitely generated ($i=0,1$).
Then, for any $\ep>0,$ any finite {{subsets}} ${\cal F}\subset A$  {{and}}
 ${\cal P}\subset \underline{K}(A),$  and
 $\{s_1,s_2,...,s_m\}\subset {\cal P}\cap K_0(A),$  there exists $\eta>0$  and
 an integer $N\ge 1$
{{satisfying}} the following:
${\cal P}\subset [\iota_n](\underline{K}({{A^\T_n}}))$ ($\iota_n: A_n\to A$ is the embedding)
for all $n\ge N,$
if $\af\in {\rm Hom}_{\Lambda}(\underline{K}({{A_n^{\T}}}), \underline{K}(B))$
and $K_1(A_n)=\Z^k\oplus {\rm Tor}(K_1(A_n)),$
where $B=C\otimes U,$ $C\in {\cal D}$  is amenable  \CA\, with continuous scale  and
$U$ is a UHF-algebra of infinite type
such that
\beq\label{Tbotext-1}
|{{\af({\boldsymbol{\bt}}(g_i))}}(\tau)|<\eta\, \tforal \tau\in T(B),\,\,\, i=1,2,...,k,
\eneq
where ${\cal P}^{(0)}:=\{g_1,g_2,...,g_k\}$ is a free generator set for $\Z^k{{\subset K_1(A_n)}},$
$\lambda: G_0\to U({\tilde B})/CU({\tilde B}),$
where $G_0$
is the subgroup generated by $\{{\bar s_1},{\bar s_2},...,{\bar s_m}\},$
${\bar s_i}\in K_0(A_n)$ and
$[\iota_n]({\bar s}_i)=s_i,$ such that
${{\Pi_{cu}^{\td B}}}\circ \lambda=\af\circ {\boldsymbol{\bt}}|_{G_0},$  $\sigma>0,$ and if $\phi: A\to B$ is  {{a}} \hm\, which maps strictly positive elements to
 strictly positive elements, {{then}}
there exists a unitary $u\in CU({\widetilde{B}})$ such that
\beq\label{Tbotex-2.2}
&&\|[\phi(a),\, u]\|<\ep\rforal a\in {\cal F},\\\label{Tbotexplus-3}
&&{\rm Bott}(\phi\circ \iota_n, \, u)=\af({\boldsymbol{\bt}}),\\
&& {\rm dist}(\Phi_{\phi,\, u}^{\dag}\circ J_{cu}^{C(\T)\otimes \td A}\circ
\boldsymbol{\bt}(s_j), \lambda(s_j))<\sigma,\,\,\,j=1,2,...,m.
\eneq

 \end{thm}

\begin{proof}
Let $\ep>0,$ ${\cal F}\subset A,$  ${\cal P}\subset \underline{K}(A),$
and $\{s_1,s_2,...,s_m\}$ be {{as}}  given.
We may assume that ${\cal F}$ is in the unit ball of $A.$
By {{Theorem}} \ref{TTbote}, there exists $\eta>0$ and   an integer $N\ge 1$
 {{satisfying}}
the conclusion of \ref{TTbote} for $\ep/2,$ ${\cal F}$ and ${\cal P}.$
We may assume that ${\cal F}\subset A_N.$

We may assume that
$
[\Phi_{\phi', u'}]|_{\cal P}
$ is {{well defined}} and
\beq
[\Phi_{\phi', u'}]|_{\cal P}=[\Phi_{\phi'', u''}]|_{\cal P}{{\andeqn {\rm Bott}(\phi', u')|_{\cal P}={\rm Bott}(\phi'',u'')|_{\cal P},}}
\eneq
if $\|[\phi'(a),\, u']\|<\ep,$
$\|\phi'(a)-\phi''(a)\|<\ep$ and $\|u'-u''\|<\ep,$ for any
\hm s $\phi',\phi'': A_N\to B$ and any unitaries $u', u''\in {\tilde B}.$
%
Choose ${\bar s}_i\in K_0(A_n)$ such that  $[\iota_n]({\bar s}_i)=s_i,$ $1,2,...,k.$
{{Let $G_1$ be the subgroup of $K_0(A)$ generated by $s_i$ ($1\le i\le k$).}}
Let $\af\in {\rm Hom}_{\Lambda}(\underline{K}({{A^\T_n}}),\underline{K}(B))$  ($n=N$) and  $\phi: A\to  B$ be as
in Theorem \ref{TTbote}  so that \eqref{Tbotext-1} holds.

  By Theorem \ref{TTbote}, there is a unitary $u_1\in CU({\tilde B})$ which satisfies the conclusion
of Theorem \ref{TTbote}
for $\ep/2,$ ${\cal F}$ and ${\cal P}.$    In particular,
\beq\label{2020-824-n1}
{{{\rm Bott}(\phi\circ \iota_n, u_1)=\af({\boldsymbol{\bt}}).}}
\eneq
\Wlog,
we may assume that $\Phi_{\phi,u_1}^\dag$
is well defined on $J_{cu}^{C(\T)\otimes \td A}\circ {\boldsymbol{\bt}}(G_0).$
Thus,   by \eqref{2020-824-n1}, for $j=1,2,...,m,$
\beq\label{2020-824-n2}
\Pi_{cu}^{\td B}\circ \Phi_{\phi,u_1}^\dag\circ J_{cu}^{C(\T)\otimes A}\circ [\iota_n]\circ {\boldsymbol{\bt}}({\bar s}_j)
= \af({\boldsymbol{\bt}}({\bar s}_j).
\eneq
  Let $e_a\in A$ be a strictly positive element of $A$ with $\|e_a\|=1.$
Since $A$ has continuous scale, \wilog, we may assume that
\beq\label{11ExtT1-5.2}
\min\{\inf\{\tau(e_a): \tau\in T(A)\},\,\inf\{\tau(f_{1/2}(e_a)):\tau\in T(A)\}\}>3/4.
\eneq
Let $T: A_+\setminus \{0\}\to \N\times \R_+\setminus \{0\}$ be given by
Theorem 5.7 of \cite{eglnp} {{corresponding to $e_a$ (in place of $e$) and $3/8$ (in place of $d$).}}.

Note that $B$ has continuous scale and $\phi$ maps strictly positive elements to
strictly positive elements. Since  $\tau\circ \phi\in T(A)$ for all $\tau\in T(B),$
we have
\beq
\tau\circ \phi(f_{1/2}(e_a))> 3/4\rforal \tau\in T(B).
\eneq
By   5.7 of \cite{eglnp},
$\phi$ is exactly $T$-${\cal H}'$-full for any finite subset
${\cal H}'\subset A_+^{\bf 1}\setminus \{0\}.$

{{We will apply 5.3 of \cite{GLII} (see the earlier part of the proof of \ref{Tbotext}).}}
Fix {{finite subsets}} ${\cal G}\subset A,$
${\cal H}_0\subset A_+^{\bf 1}\setminus \{0\},$
${\cal H}_1\subset A_{s.a.},$
${\cal Q}_0\subset \underline{K}(A),$
${\cal U}\subset U({\tilde A}),$ {{and fix}}
$\ep_0>0,$ $\eta_1>0$ and $\eta_2>0.$
We may assume that ${\cal U}=\{\af_i\cdot 1_{\tilde A}+y_i: \af_i\in \T, y_i\in A,\,\, 1\le i\le k_1\}.$
We {{also}} assume that ${\cal P}\subset {\cal Q}_0.$
We may write that  $u_1=1+x$ for some normal element $x\in A.$

Let $e_B\in B_+$ be a strictly positive element. Choose $n_0\ge 1$ such
that
\beq
&&\|f_{1/n_0}(e_B) \phi(a)f_{1/n_0}(e_B)-\phi(a)\|<\min\{\ep_0/32, \eta_2,\sigma/4\}\\
&&\rforal a\in {\cal G}\cup\{y_1, y_2,..,y_{k_1}\}\andeqn\\
&&\|(1+f_{1/n_0}(e_B)xf_{1/n_0}(e_B))-u_1\|<\min\{\ep/16, \ep_0/16, \sigma/4\}.
\eneq
Put $B_1=\overline{f_{1/n_0}(e_B)Bf_{1/n_0}(e_B)}.$ There exists
$x_1\in B_1$ such that
\beq
1+x_1\in U({\tilde B})\andeqn \|(1+x_1)-u_1\|<\min\{\ep/16, \ep_0/16,\sigma/4\}.
\eneq
Choose $b\in \overline{(1_{\tilde B}-f_{1/(n_0+2)}(e_B))B(1_{\tilde B}-f_{1/(n_0+2)}(e_B))}_+$ with
$\|b\|=1$ and let $B_0=\overline{bBb}.$  We assume that $d_\tau(b)<\eta_1$ for all $\tau\in T(B).$
Moreover, we may assume that $B_0$ has continuous scale {{(see Lemma 6.8 of \cite{GLrange}).}}
Note that $B_0B_1=0.$
Put $\sigma_0=\inf\{d_\tau(b): \tau\in T(B)\}>0.$

Define $\lambda_0=\lambda-\Phi_{\phi, u_1}^{\dag}\circ J_{cu}^{C(\T)\otimes \td A}\circ {\boldsymbol{\bt}}|_{G_0}.$
{{By the assumption $\Pi_{cu}^{{\td B}}\circ \lambda=\af({\boldsymbol{\bt}})|_{G_0}$ and
by \eqref{2020-824-n2}, $\Pi_{cu}^{{\td B}}\circ \lambda_0=0.$}}
Therefore
$\lambda_0$ maps $G_0$ into $\Aff({{T(B)}})^\iota/\overline{\rho_{\tilde B}(K_0({\tilde B}))}.$
By applying Lemma \ref{Ldeterm2}, we obtain a \hm\, $\phi_0: A\to B_0$ and
a unitary {{$u_0=1+x_0\in CU({\tilde B}_0)$ with $x_0\in B_0$}} such that
\beq
&&\|[\phi_0(a), \, u_0]\|<\ep/2\rforal a\in {\cal F},\\\label{81017-3}
&&[\Phi_{\phi_0, u_0}]|_{{\cal Q}_0\cup {\boldsymbol{\bt}}({\cal Q}_0)}=0\,\,\, {\rm in}\,\, KL(A^\T, B),\\
&&\phi_0^{\dag}|_{J_{cu}^A(K_1(A))}=0\andeqn\\\label{81017-1}
&&{\rm dist}((\Phi_{\phi_0, u_0})^{\dag}({\boldsymbol{\bt}}(s_i)), \lambda_0(s_i)){{<}}\sigma\cdot \sigma_0/4.
\eneq
Define  $\phi_1: A\to B$ by  $\phi_1(a)=f_{1/n_0}(e_B)\phi(a)f_{1/n_0}(e_B)+ \phi_0(a)$ for all $a\in A$ and
$u_2=1_{\td B}+x_1+x_0\in {\tilde B}.$  Note that, since $B_0B_1=0,$  $u_2$ is a unitary.
As in the end of the proof of
\ref{Tbotext}, we may choose {{$u_2$}} so that it is in $CU({\tilde B}).$
{{We compute that, for $1\le j\le m,$
\beq\label{2020-824-1}
\hspace{-0.2in}\Phi_{\phi_1, u_2}^{\dag}\circ J_{cu}^{C(\T)\otimes \td A}\circ {\boldsymbol{\bt}}(s_j)\approx_{\sigma/4}
\Phi_{\phi, u_1}^{\dag}\circ J_{cu}^{C(\T)\otimes \td A}\circ {\boldsymbol{\bt}}(s_j)+
\Phi_{\phi_0,u_0}^\dag\circ J_{cu}^{C(\T)\otimes \td A}\circ {\boldsymbol{\bt}}(s_j).
\eneq}}
{{Also,}}  $\phi_1$ is ${\cal G}$-$\ep_0$-multiplicative.
\Wlog, we may assume
that
\beq
\tau\circ f_{1/2}(\phi_1(e_a))>3/8\rforal \tau\in T(B).
\eneq
Therefore, with sufficiently small $\ep_0$ and large ${\cal G},$ by {{5.7  of \cite{eglnp},}}
$\phi_1$ is exactly $T$-${\cal H}_0$-full.

Note that
\beq
[\phi_1]|_{{\cal Q}_0}=[\phi]|_{{\cal Q}_0},\\
{\rm dist}(\phi^{\dag}(x), \phi_1^{\dag}(x))<\eta_2\andeqn\\
|\tau\circ \phi_1(c)-\tau\circ \phi(c)|<\eta_1\rforal c\in {\cal H}_1.
\eneq
Since ${\cal Q}_0,$ ${\cal H}_0,$ ${\cal U},$ $\ep_0,$ $\eta_1,$ and $\eta_2$ are arbitrarily chosen,
{{we may}} choose in such a way so that 5.3 of \cite{GLII} can be applied  for $\ep/2$ and ${\cal F}.$
By applying {{5.3 of \cite{GLII},}} we obtain a unitary $w\in {\tilde B}$ such that
\beq
\|w^*\phi_1(a)w-\phi(a)\|<\ep/2\rforal a\in {\cal F}.
\eneq
We then choose $u=w^*u_2w.$
Moreover, we may assume that ${\cal P}\subset {\cal Q}_0.$
{{Note that $\Phi_{\phi, u}={\rm Ad}\, w\circ \Phi_{\phi_1,u_1}.$}}
Note that, by \eqref{81017-3},
\beq\label{81017-2}
[\Phi_{\phi, u}]|_{\cal P}=[\Phi_{\phi_1, u_1}]|_{\cal P}.
\eneq
It follows from \eqref{81017-1} and \eqref{2020-824-1} that
\beq
{\rm dist}(\Phi_{\phi,\, u}^{\dag}\circ J_{cu}^{C(\T)\otimes \td A}\circ \boldsymbol{\bt}(s_j), \lambda(s_j))<\sigma,\,\, j=1,2,...,m.
\eneq
By \eqref{81017-2} and the choice of $u_1,$
\beq
{\rm Bott}(\phi,\, u)|_{\cal P}={\rm Bott}(\phi, u_1)|_{\cal P}=\af({\boldsymbol{\bt}})|_{\cal P}.
\eneq
{{The lemma}} follows.
\end{proof}
%

%

%
%
%
%

%
%
%
%
%
%
\section{Maps to the mapping tori}

%
%

%
%


\begin{lem}\label{141-0}
Let $A$ be a \CA, $a_1, a_2\in A_+$ such that $a_1+a_2=p$ is a projection,
${{b_1,}} b_2\in A_+^{\bf 1}$ with $b_1b_2=0$ {{and}} $1/4>\dt>0.$
Suppose that, for any $\ep>0,$ there are $x_1, x_2\in A^{\bf 1}$ such
that $x_1x_1^*\approx_{\ep}a_1,$ $x_2x_2^*\approx_{\ep} a_2,$
$x_1^*x_1f_{\dt}(b_1)=f_{\dt}(b_1)x_1^*x_1=x_1^*x_1$ and $x_2^*x_2f_{\dt}(b_2)=f_{\dt}(b_2)x_2^*x_2=x_2^*{{x_2.}}$
Then there is $v\in A$ such that $v^*v=p$ and $vv^*\le f_{\dt/2}(b_1+b_2).$

\end{lem}

\begin{proof}
Let $1/4>\ep>0$ and $x_1, x_2$ be as described.
Put $z=(x_1+x_2).$ Then
\beq\nonumber
&&\hspace{-0.3in}z^*z=x_1^*x_1+x_2^*x_2+x_1^*x_2+x_2^*x_1\\\nonumber
&&\hspace{-0.2in}\le f_{\dt}(b_1)+f_{\dt}(b_2)+f_{\dt}(b_1)x_1^*x_2f_{\dt}(b_2)+f_{\dt}(b_2)x_2^*x_1f_{\dt}(b_1)\in {\rm Her}(f_{\dt}(b_1)+f_{\dt}(b_2)).
\eneq
Also
\beq
zz^*=(x_1+x_2)(x_1+x_2)^*=x_1x_1^*+x_2x_2^*+x_1x_2^*+x_2x_1^*=x_1x_1^*+x_2x_2^*{{\approx_{2\ep}}} p.
\eneq
There is $y\in pAp$ such that $yzz^*y^*=p.$
Then $q:=z^*y^*yz\le \|y\|^2z^*z\in {\rm Her}(f_{\dt}(b_1)+f_{\dt}(b_2)).$
Note that $q$ is a projection and $q\le f_{\dt/2}(b_1+b_2).$

\end{proof}

One should note that $A$ is not unital in the following statement.

\begin{lem}\label{Lcont}
Let $A$ be a non-unital and $\sigma$-unital \CA\, and $B$  {{a}} non-unital separable
amenable \CA.
Let
$\{e_n\}$ be an approximate identity for $A\otimes B$
with the property that
$$
e_{n+1}e_{n}=e_{n}e_{n+1}=e_{n}
$$
and {{$\|e_{n+1}-e_{n}\|=1$ 
for}} {{all $n,$}}
and let $\{a_k\}$ be a sequence of orthogonal
elements in $(A\otimes B)_+$
such that $0<\dt_0\le  \|a_k\|\le M$ for some $M>1>\dt_0>0,$
one of $a_k$ is full in $A\otimes B,$ and such that
$\sum_{k=1}^\infty a_k$ converges strictly to an element in $M(A\otimes B).$
Suppose also, {{for any $i\ge 1$ and
$0<d\leq\dt_0/2,$}} there exists $N\ge 1$ such that, when
{{$m\ge N$ and $k\in\N$,}}
\beq
e_{m+k}-e_m\lesssim {{f_d(a_i)}}\,\,\,{\rm in}\,\,\, A\otimes B.
\eneq
Then, $\sum_{k=1}^\infty a_k$ is full in $M(A\otimes B).$
Moreover, there is a projection\\
 $p\in \pi({{\overline{f_{\dt_0/8}(\sum_{k=1}^\infty a_k)M(A\otimes B)f_{\dt_0/8}(\sum_{k=1}^\infty a_k)}}})$
and $v\in M(A\otimes B)/(A\otimes B)$
such that $v^*v=p$ and $vv^*=1_{M(A\otimes B)/(A\otimes B)},$ where $\pi: M(A\otimes B)\to M(A\otimes B)/(A\otimes B)$
is the quotient map.

Furthermore,  suppose, in  addition,  that $B$ is simple, and $\{e_{a, n}\}$ is an approximate identity
for $A$ with $e_{a,n+1}e_{a, n}=e_{a,n}=e_{a,n}e_{a,n+1}$ and $\|e_{a,n+1}-e_{a,n}\|=1,$
%
%
and suppose, for any $j, j_1\in \N$ (with $j<j_1$)
and $0<d<1/4,$ {{and}}
any $b\in B_+\setminus \{0\}$
with $\|b\|=1,$   there exists {{$N\ge 1$}} such that, when $m\ge N,$ for any $k\in \N,$
\beq
e_{m+k}-e_m \lesssim f_d((e_{a, j_1}-e_{a, j})\otimes  b)\,\,\,{\rm in}\,\, A\otimes B,
\eneq
{{and, for any  $j_1>j,$ $(e_{a, j_1}-e_{a, j})\otimes  b$ is full.}}
Then, for any
$c\in B_+\setminus\{0\},$ $1_{\tilde A}\otimes c$ is full
in $M(A\otimes B).$

\end{lem}

\begin{proof}
Let us first show that the ``Furthermore" part follows from the first part of the statement.
Fix $c\in B_+\setminus \{0\}.$ We may assume that $\|c\|=1.$
By {{assumption,}} we have that $\|e_{n+1}-e_n\|=1$ and $\|e_{a, j+1}-e_{a,j}\|=1$
for all $n$ and $j.$
%
%
Put $g_n=(e_{n+1}-e_{n})$
and  $g_{a,n}=e_{a, n+1}-e_{a, n},$
 $n=1,2,....$
 Then, if $|j-i|\ge 2,$
 \beq
g_{a,j}g_{a,i}=0\andeqn {{g_{j}g_{i}=0.}} 
 \eneq

Let ${{J\subset M(A\otimes B)}}$ be the closed ideal generated by $x:=1_{\tilde A}\otimes c.$
Note that $x\not\in A\otimes B.$
We will show that $1\in J.$

Note that $x=\sum_{k=1}^{\infty}g_{a,2k-1}\otimes c+\sum_{k=1}^{\infty}g_{a,2k}\otimes c\in J_+\setminus A\otimes B$. By choosing either $n(k)=2k$ or $n(k)=2k-1$, we have that
{{$\sum_{k=1}^{\infty} g_{a, n(k)}\otimes c\in J_+\setminus A\otimes B,$}}
$g_{a, n(k)}\otimes c\not=0$ and $0\le {g_{a,n(k)}}\otimes c\le 1,$
where $n(k+1)>n(k)+1.$
Since $\|g_{a,n(k)}\|=1,$
we may further assume that there is $0<\dt_0<1/2$ such that
$\|g_{a, {{n(k)}}}\otimes c\|\ge \dt_0.$

Put $a_k=g_{a,n(k)}\otimes c.$  {{Note that, for some $k,$ $a_k$ is full in $A\otimes B.$}}
By the first part of the statement,
$\sum_{k=1}^\infty g_{a,2k}\otimes c$ is full in $M(A\otimes B).$ Similarly, $\sum_{k=1}^\infty g_{a,2k-1}\otimes c$
is also full.  It follows that $x:=1_{\td A}\otimes c=(\sum_{k=1}^\infty g_{a, 2k}\otimes c)+(\sum_{k=1}^\infty g_{a, 2k-1}\otimes c)$ is full in $M(A\otimes B).$

It remains to show the first part of the statement.

Let $y_1=\sum_{k=1}^\infty a_{2k-1},$  $y_2=\sum_{k=1}^\infty a_{2k}$ and $y=y_1+y_2.$
Note {{that}} $y_1, y_2\le y.$
Now let $J$ be the closed ideal generated by $y.$ By the assumption, there exists $n_0\ge 1$ such that
$$
\sum_{n=n_0+1}^m g_{2k}\lesssim  f_{\dt_0/2}(a_2)
$$
for all $m\ge n_0.$ 
By the induction,  there is an integer $n_k>n_{k-1}$ such that $$
\sum_{n=n_k+1}^m g_{2k}\lesssim  f_{\dt_0/2}(a_{2(k+1)})
$$ for all $m\ge n_k.$ 
Let  {{$N_1=\{n_0+1,n_0+2,\ldots, n_1\}$, $N_2=\{n_1+1,n_1+2,\ldots, n_2\},\ldots$, $N_k=\{n_{k-1}+1,n_{k-1}+2,\ldots, n_k\}$.}} In such a way, we get partition $\{n_0+1, n_0+2,...\}$ into finite subsets $N_1, N_2,....$ (of consecutive integers)
such that
\beq
\sum_{j\in N_k}g_{2j}\lesssim f_{\dt_0/2}(a_{2k}),\,\,\, k=1,2,....
\eneq
Fix  $1>\ep>0.$ There are $x_k$ of the form {{$x_k=r_k^*f_{\dt_0/2}(a_{2k})r_k$}}
such that
\beq
\|x_k-\sum_{j\in N_k} g_{2j}\|<\ep/2^{k+1}\tand\|x_k^{1/2}-\sum_{j\in N_k} g_{2j}^{1/2}\|<\ep/2^{k+1}.
\eneq
We may assume that $0\le x_k\le 1.$  
{{Set $z_k=r_k^*f^{1/2}_{\dt_0/2}(a_{2k})$. Then}} $z_kz_k^*=x_k,$ $z_k^*z_k\le {{\|r_k\|^2\cdot}}f_{\dt_0/2}({{a_{2k}}})$ and
\beq\label{141-nn1}
z_k^*z_kf_{\dt_0/4}(a_{2k})=f_{\dt_0/4}(a_{2k})z_k^*z_k
=z_k^*z_k.
\eneq
{{Note that $\{a_k\}$ are  mutually orthogonal. Hence}} $z_iz_j^*=0,$ if ${{i\not=j}}.$
Therefore
\beq
(\sum_{k=1}^nz_k)(\sum_{k=1}^n z_k)^*=\sum_{k=1}^n z_kz_k^*
\eneq
and $\{\|\sum_{k=1}^n z_kz_k^*\|\}$ is bounded.
It follows that $\{\|\sum_{k=1}^n z_k\|\}$ is bounded.
It is then easy to see that $\sum_{k=1}^n z_k$
converges in the left strict topology to the element
$z=\sum_{k=1}^{\infty}z_k$ in the left multiplier $LM(A\otimes B).$
To show that $\sum_{k=1}^n z_k$ also
converges strictly to $z,$ it suffices to show that, for each $m,$
$g_m\sum_{k={{N}}}^{\infty} z_k$ converges in norm to zero as
$N\to\infty.$ Write
$z_k=(z_kz_k^*)^{1/2}u_k.$ Then
\beq
\|g_m \sum_{k=N}^{\infty} z_k\| &\le &\|g_m(\sum_{k=N}^{\infty}z_k-
\sum_{k=N}^{\infty}(\sum_{j\in N_k}g_{2j}^{1/2})u_k)\|+\|g_m\sum_{k=N}^{\infty}(\sum_{j\in N_k} g_{2j}^{1/2})u_k\|\\
&\le & \|\sum_{k=N}^{\infty} \|x_k^{1/2}-\sum_{j\in N_k}g_{2j}^{1/2}\|+
\|g_m\sum_{k=N}^{\infty} (\sum_{j\in N_k}g_{2j}^{1/2})u_k\|\\
&<& \sum_{k=N}^{\infty}\ep/2^{k+1}+\|g_m\sum_{k=N}^{\infty}(\sum_{j\in N_k}g_{2j}^{1/2})u_k\|.
\eneq
However,
$$
g_m(\sum_{k=N}^{\infty}(\sum_{j\in N_k}( g_{2k}^{1/2})u_k)=0\rforal N>m+1.
$$
One concludes that $\lim_{N\to\infty}\|g_n\sum_{k=N}^{\infty} z_k\|=0.$ Therefore
$z\in M(A\otimes B).$ On the other hand,
\beq
zf_{\dt_0/4}(\sum_{k=1}^\infty a_{2k})=(\sum_{k=1}^{\infty}z_k)(\sum_{k=1}^{\infty} f_{\dt_0/4}(a_{2k}))
=\sum_{k=1}^{\infty}z_k=z.
\eneq
It follows that $z\in J.$ But
$$
\|zz^*-\sum_{k\ge n_0}g_{2k}\|<\ep.
$$
It follows that $\sum_{k=1}^{\infty} g_{2k}\in J+A\otimes B,$ as $y_2\in J.$
Similarly, $\sum_{k=1}^{\infty} g_{2k-1}\in J+A\otimes B.$ Therefore
$1_{M(A\otimes B)}\in J+A\otimes B.$
Suppose that $a_{k_0}$ is full in $A\otimes B.$
Since $a_{k_0}\le \sum_{k=1}^\infty a_k,$
the closed ideal generated by $y$
contains $A\otimes B.$ Therefore $\sum_{k=1}^\infty a_k$ is full in $M(A\otimes B).$

Note {{that,}}  by \eqref{141-nn1}, $z^*zf_{\dt_0/4}(y_2)=f_{\dt_0/4}(y_2)z^*z=z^*z.$
Symmetrically,  we may assume that, for some $n_0'\ge 1,$ and, for any $\ep>0,$    there is $z'\in M(A\otimes B)$ such that
$(z')^*z'f_{\dt_0/4}(y_1)=f_{\dt_0/4}(y_1)(z')^*z'=(z')^*z'$ and $\|z'(z')^*-\sum_{n\ge n_0}g_{2k-1}\|<\ep.$
Note that $\pi(\sum_{n\ge n_0}g_{2k}+\sum_{n\ge n_0'} g_{2k-1})=1_{M(A\otimes B)/A\otimes B}.$
Then the ``Moreover" part follows from Lemma \ref{141-0}.
\end{proof}

\begin{lem}\label{LMtorusW}
Suppose that $C$ and $B_0$ are two separable \CA s such that
$B_0$  has  property (W). Let $\phi_0, \phi_1: C\to B_0$ be two monomorphisms.
Then $M_{\phi_0, \phi_1}$ has property (W).
\end{lem}

\begin{proof}
 Let $T: (B_0)_+\setminus\{0\}\to  \N\times {{\R_+}}\setminus \{0\}$
be a map {{such}} that there exists a sequence of approximately multiplicative \cpc s $\psi_n: B_0\to {\cal W}$ such that, for any
finite subset ${\cal H}\subset (B_0)_+^{\bf1} \setminus \{0\},$
$\psi_n$ are exactly $T$-${\cal H}$-full for all $n\ge n_0$ {{(for certain $n_0$ depending on ${\cal H}$).}}
Assume $T(b)=(N(b), M(b)),$ where $N: (B_0)_+\setminus \{0\}{{\to \N}}$ and $M: (B_0)_+\setminus \{0\}\to \R_+\setminus \{0\}.$
Note that ${\cal W}\cong {\cal W}\otimes Q,$ where $Q$ is the universal UHF-algebra (see \cite{eglnkk0}, for example).
Therefore, for any $k(n)\in \N,$ there is a \hm\,
$\phi^{(k)}: {\cal W}\otimes M_{k(n)}\to {\cal W}$ which maps strictly positive elements to strictly positive elements.
Let $e_W\in {\cal W}$ with $\|e_W\|=1$ be a strictly positive element.
We assume that $N^w\times M^w:{\cal W}_+\setminus \{0\}\to \N\times {{\R_+}}\setminus \{0\} $ {{is}} a map {{such that}} $f_{1/2}(e_W)$ is $N^w\times M^w$-full in ${\cal W},$ i.e.,  for any $w\in {\cal W}_+\setminus \{0\}$ {{with $\|w\|\leq 1$,}}
there are $w_1,w_2,...,w_{N^w}\in {\cal W}$ with $\|w_i\|\le M^w$ ($1\le i\le N^w$)
such that $\sum_{i=1}^{N_w} w_i^*f_{1/2}(e_W)w_i=w.$
Note {{that,}} by Lemma 3.3 of \cite{eglnp},  if $\|w'-f_{1/16}(e_W)\|<1/32,$ there is $r\in {\cal W}$ with $\|r\|\le 2$ such that
\beq
f_{1/2}(e_W)=r^*w'r.
\eneq
Let $a\in {M_{\phi_0, \phi_1}}_+\setminus\{0\}$ {{with $\|a\|\leq 1$.}}
There is $t_a\in [0,1]$ such that $a(t_a)\in {B_0}_+\setminus \{0\}.$
Therefore, for all sufficiently large $n,$ there are $y_1,y_2,...,y_{N(a(t_a))}{{\in {\cal W}}}$ with $\|y_i\|\le {{M(a(t_a))}}$ ($1\le i\le N(a(t_a))$)
such that
\beq
\sum_{i=1}^{N(a(t_a))}y_i^*{{\psi_n(a(t_a))}}y_i=f_{1/16}(e_W).
\eneq
Choose $\dt(a)>0$ such that $\|a(t_a)-a(t)\|<1/{{(65\|N(a(t_a))\| \max\{\|M(a(t_a))\|\})}},$ whenever  $|t_a-t|\le \dt(a)$ ($t\in [0,1]$).
It follows that
\beq
\|\sum_{i=1}^{N(a(t_a))} y_i^*{{\psi_n(a(t))}}y_i-f_{1/16}(e_W)\|<1/32.
\eneq
So, for all $t\in [t_a-\dt(a), t_a+\dt(a)],$
 there is $r(t)\in {\cal W}$ with $\|r(t)\|\le 2$ such that
\beq
\sum_{i=1}^{N(a(t_a))} r(t)^*y_i^* {{\psi_n(a(t))}}y_ir(t)=f_{1/2}(e_W).
\eneq
Define $N_1: {M_{\phi_0, \phi_1}}_+\setminus \{0\}\to \N$  and
$M_1:{M_{\phi_0, \phi_1}}_+\setminus \{0\}\to\R_+\setminus\{0\}$  by
$N_1(a):=[N^wN(a(t_a))/\dt(a)]+1$ and $M_1(a){{:=2}}M^wN(a(t_a)).$
Define $T_M: {M_{\phi_0, \phi_1}}_+\setminus \{0\}$ by  $T_M=(N_1, M_1).$

Let $\{{\cal F}_n\}$ be an increasing sequence of {{finite subsets of}} $A$ such that $\cup_{n=1}^{\infty} {\cal F}_n$ is dense in $A.$
Let  $\{\ep_n\}$  be a decreasing sequence of positive numbers such that $\sum_{n=1}^{\infty}\ep_n<\infty.$
There are $t_0=0<t_{1,n}<\cdots {{t_{k(n),n}}}=1$ such that  ${{t_{i,n}-t_{i-1,n}}}=1/k(n)$ and
\beq
\|g(t)-g(t_j)\|<\ep_n/2\rforal t\in [t_{i-1,n}, t_{i+1,n}]\andeqn \rforal g\in {\cal F}_n.
\eneq
Define $\Psi_n: M_{\phi_0, \phi_1}\to M_{k(n)}({\cal W})\stackrel{\phi^{k(n)}}{\to} {\cal W}$ by
$$
\Psi_n((a, g))=\phi^{k(n)}(\diag(\psi_n\circ  \pi_0(a), \psi_n\circ \pi_1(a), \psi_n\circ \pi_{t_{1,n}}(g), ...,\psi_n\circ \pi_{t_{k(n)-1},n}(g))),
$$
where $\pi_t: M_{\phi_0, \phi_1}\to B_0$ is the {{point evaluation}} at $t\in [0,1].$
Fix a finite subset ${\cal H}_M\subset M_{\phi_0, \phi_1}.$
It follows that $\{\Psi_n\}$ is a sequence of approximately multiplicative \cpc s from $M_{\phi_0, \phi_1},$ which are
eventually exactly $T_M$-${\cal H}_M$-full.  In other words, $M_{\phi_0, \phi_1}$ has property (W).

\end{proof}

\begin{NN}\label{ConstphiW}
Let $B_0$ and $C$ be  non-unital
separable {{stably}} projectionless simple \CA s with  stable rank one, {{with continuous scales and with
$T(B_0)\not=\emptyset$}}. Let $\phi_0, \phi_1: C\to B_0$ be \hm s which
send strictly positive elements to strictly positive elements.
 It follows from \cite{Rl} that  there is {{an}} embedding $j_w: {\cal W}\to C$ which maps strictly positive elements to strictly positive elements.

Suppose that $\tau\circ \phi_0=\tau\circ \phi_1$ for all $\tau\in T(B_0).$
We also assume that
$B_0$ is
${\cal Z}$-stable.
Put $B=M_{\phi_0, \phi_1}.$
Note that $\tau\circ \phi_0\circ j_w=\tau\circ \phi_1\circ j_w.$  It follows that
$Cu^\sim(\phi_0\circ j_w)=Cu^\sim(\phi_1\circ j_w)$ (see \cite{Rl}).
Then, by \cite{Rl}, $\phi_0\circ j_w$ and $\phi_1\circ j_w$ are approximately
unitarily equivalent.  Let $\{e_{w,n}\}$ be an approximate identity for ${\cal W}$
such that $e_{w,n+1}e_{w,n}=e_{w,n}=e_{w,n}e_{w, n+1}.$
{{By passing to a subsequence, we may assume that there is $a_{w,n}\in {\cal W}_+^{\bf 1}\setminus \{0\}$ such
that $(e_{w, n+1}-e_{w,n})a_{w,n}=a_{w,n}$ for $n\in \N.$}}
Then $\{\phi_i\circ j_w(e_{w,n})\}$ is an approximate identity for $B_0.$  It follows  from Theorem 5.7 of \cite{eglnp}
that there exists $T_0: {\cal W}_+\setminus \{0\}\to \N\times \R_+\setminus \{0\}$  such that
both $\phi_0\circ j_w$ and $\phi_1\circ j_w$ are  exactly $T_0$-${\cal H}_1$-full for every
finite subset ${\cal H}_1\subset {\cal W}_+^{\bf 1}\setminus \{0\}.$
Put $T=2T_0.$ Let $\{{\cal F}_n\}\subset {\cal W}$ be an increasing sequence of finite subsets
{{such that  $\cup_{n=1}^\infty {\cal F}_n$ is dense in ${\cal W}$}} and
$\{\ep_n\}$ is a decreasing sequence of positive numbers such that $\sum_{n=1}^{\infty}\ep_n<\infty.$
Let ${\cal G}_n=\{{\cal F}_n\}\cup \{ab: a, b\in {\cal F}_n\}.$
There is, by {{Theorem 1.0.1 of}} \cite{Rl},  a sequence of unitaries $u_n\in U_0({\tilde{B_0}})$ (see \ref{Lunitary}) such that
\beq
\|u_n^*\phi_1\circ j_w(a)u_n-\phi_0\circ j_w(a)\|<\ep_n/4\rforal a\in {\cal G}_n.
\eneq
Let $\{u_n(t): t\in [1/2,1]\}$ be a continuous path of unitaries in ${\tilde{B_0}}$ such
that $u_n(1/2)=u_n$ and $u_n(1)=1.$
Define $\Phi_n: {\cal W}\to B=M_{\phi_1, \phi_2}$ by $\Phi_n(a)=(\Psi_n(a), j_w(a)),$ where
\beq\label{2020-721-12-1}
\Psi_n(a)(t)=\begin{cases} 2(1/2-t)\phi_0\circ j_w(a)+2t{\rm Ad}\, u_n \circ \phi_1\circ j_w(a), &\,\,\, t\in [0,1/2];\\
                                       {\rm Ad}\, u_n(t)\circ \phi_1\circ j_w(a), &\,\,\,t\in (1/2,1].
                                       \end{cases}
\eneq
Note that $\Psi_n(a)(0)=\phi_0\circ j_w(a)$ and $\Psi_n(a)(1)=\phi_1\circ j_w(a)$ for all $a\in {\cal W}.$
It follows that $\Psi_n$ is ${{{\cal F}_n}}$-$\ep_n$-multiplicative.
Therefore $\Phi_n$ is ${{{\cal F}_n}}$-$2\ep_n$-multiplicative. Moreover, if $a\in {\cal W}$ is a strictly positive element, so is
$\Phi_n(a)$ (in $B$). Fix a finite subset ${\cal H}\subset {\cal W}_+^{\bf 1}\setminus \{0\}.$
Then, for all sufficiently large $n,$ $\Phi_n$ is $T$-${\cal H}$-full.

Let $\{f_n\}$ be an approximate identity for $B$ such that
$f_{n+1}f_n=f_n$ for all $n.$
For each $n\ge 1,$ there is an integer $k(n)$ such that
$f_{k(n)}\Phi_n(a)\approx_{\ep_n/8} f_{k(n)}^{1/2}\Phi_n(a)f_{k(n)}^{1/2}\approx_{{\ep_n/8}}
\Phi_n(a)$ for all $a\in {\cal F}_n.$
Passing to a subsequence of $\{f_n\},$ we may assume
that ${\bar\Phi_n}$ defined by ${\bar\Phi_n}(a):=f_{n}^{1/2}\Phi_n(a)f_{n}^{1/2}$ is ${\cal F}_n$-$\ep_n$-multiplicative.


Consider $B\otimes \zo.$
Let $\phi_0', \phi_1': C\otimes \zo\to B_0\otimes \zo$ be defined by
$\phi_i'=\phi_i\otimes {\rm id}_{\zo},$ $i=0,1.$
Then $B\otimes \zo=M_{\phi_0', \phi_1'}.$
Let $\Phi_n': {\cal W}\to M_{\phi_0', \phi_1'}$ be as defined above (for $\phi_i'$  instead of
$\phi_i$). {{Let ${{{\bar\Phi'}_n(a)}}=(f_n^{1/2}\otimes 1_{\tilde \zo})\Phi_n'(a)(f_n^{1/2}\otimes 1_{\tilde \zo})$
for all $a\in {\cal W}.$ As mentioned above we may assume that ${\bar\Phi'}_n(a)$  is ${\cal F}_n$-$\ep_n$-multiplicative.}}


Let $j_{w,z}: {\cal W}\to \zo$ and $j_{z,w}: \zo\to {\cal W}$ be embeddings which map strictly positive elements
to strictly positive elements (see, for example,
\cite{Rl}  {{and 4.33 of \cite{GLrange}}}) and let
$e_n':=j_{w,z}(e_{w,n})$ {{and $a_n':=j_{w,z}(a_{w,n}).$}}
{{Let}} $W_n:={\rm Her}(a_{w,n})$ and
$Z_n:={\rm Her}(a_{w,n}').$
Then $B\otimes W_n$ and $B_n':=B\otimes Z_n$ are hereditary \SCA s of $B\otimes{\cal W}$ and
$B\otimes \zo,$ respectively.
There is an isomorphism $j_n: {\cal W}\to W_n.$
Define $h_n: B\otimes \zo \to B\otimes  W_n$ by $h_n(b\otimes z)=b\otimes j_n\circ j_{z,w}(z)$
for all $b\in B$ and $z\in \zo$ ($n\in \N$). Define ${{\Lambda}}_{n,k}: {\cal W}\to B_n:=B\otimes Z_n\subset B\otimes \zo$
by ${{\Lambda_{n,k}(a)}}=({{{\rm id}_{B}}}\otimes j_{w,z})\circ h_n\circ{{{\bar\Phi'_{n+k}}(a)}}$
%
%
%
%
%
for all $a\in {\cal W}.$
{{We note that, for all $a\in {\cal W},$}}
%
%
%
%
%
%
\beq\label{09-15-2021}
\Lambda_{n,k}(a)=(f_{n+1}^{1/2}\otimes (e_{n+1}'-e_n'))(({\rm id}_{B}\otimes j_{w,z})\circ h_n\circ {\bar\Phi'_{n+k}}(a))(f_{n+1}^{1/2}\otimes (e'_{n+1}-e'_n)).
\eneq
Let $\phi_{k,W}:
{\cal W}\to M(B\otimes \zo)$ be defined
by
\beq
&&\phi_{k,even}=\sum_{n=1}^{\infty} {{\Lambda_{4n,k}}},\,\,\,
\phi_{k, odd}=\sum_{n=1}^{\infty} {{\Lambda_{4n+2,k}}},\\
&&\phi_{k,W}=\sum_{n=1}^{\infty} {{\Lambda_{2n,k}}}=\diag(\phi_{k, even}, \phi_{k,odd}).
\eneq
Note that, for any finite subset ${\cal F}\subset {\cal W}$ and $\ep>0,$  there exists $k_{{{\cal F},\ep}}\ge 1$ such that
$\phi_{k, W}, \phi_{k, even}, \phi_{k, odd}$ are  ${\cal  F}$-$\ep$-multiplicative for $k\ge k_{{{\cal F},\ep}}.$
Also, for any $a\in  {\cal W}_+\setminus \{0\},$ $\Lambda_{n,k}(a)$ is full in $B\otimes \zo$ for all sufficiently large $n.$
Moreover,  since $\lim_{n\to\infty}\|\Phi'_n(a)\Phi'_n(b)-\Phi'_n(ab)\|=0$
for all $a, b\in {\cal W},$
$\pi\circ \phi_{k, W},$
$\pi\circ \phi_{k, even},$ $\pi\circ \phi_{k, odd}$ are \hm s from ${\cal W}$ to $M(B\otimes \zo)/B\otimes \zo,$
where $\pi: M(B\otimes \zo)\to M(B\otimes \zo)/B\otimes \zo$ is the quotient map.

It follows from the first part of  Lemma \ref{Lcont} that $\pi\circ \phi_{k,even}(a)$ and $\pi\circ \phi_{k, odd}(a)$ are full
in $M(B\otimes \zo)/B\otimes \zo$ for all $a\in {\cal W}_+\setminus \{0\},$ where $\pi: M(B\otimes \zo)\to
M(B\otimes \zo)/B\otimes \zo$ is  the quotient map.  This statement follows from Lemma \ref{Lfullness} below.

We will keep these notation in the next four statements.  
\end{NN}

\begin{lem}\label{L141N}
Let $C,$ $B_0$  and $B=M_{\phi_0, \phi_1}$ be as in \ref{ConstphiW}.
Let $\{e_n\}$ be  an approximate identity of ${{B\otimes \zo}}$ such that
$e_{n+1}e_n=e_ne_{n+1}=e_n$ and $\|e_{n+1}-e_n\|=1.$
Then,
for any $b\in (B\otimes \zo)_+\setminus \{0\}$ such that
$\tau(b)>0$ for all $\tau\in T(B\otimes \zo),$ there is {{an integer}} $N\ge 1$ such that
\beq
e_{m+k}-e_m\lesssim b\tforal m\ge N\tand k\ge 1.
\eneq
{{In particular}}, $b$ is full in $B\otimes \zo.$ Moreover, one may choose an
approximate identity $\{e_{b,n}\}\subset B$ such that $(e_{b,m}-e_{b,n})\otimes c$
is full for any $c\in {\zo}_+\setminus \{0\}$ and $m>n.$
\end{lem}

\begin{proof}
Recall {{that,}} using notation in \ref{ConstphiW}, $B\otimes \zo=M_{\phi_0', \phi_1'},$
and recall
$$
M_{\phi_0', \phi_1'}=\{(f, c)\in C([0,1], B_0\otimes \zo)\oplus C\otimes \zo): \phi_0'(a)=f(0)\andeqn \phi'_1(a)=f(1)\}.
$$
Since both $C\otimes \zo$ and $B_0\otimes \zo$ have continuous {{scales,}}
a compactness argument shows
that
\beq\label{200925-14n1}
\eta:=\inf\{\tau(b): \tau\in T(M_{\phi_0', \phi_1'})\}>0.
\eneq
Moreover,
\beq
t\circ \pi_e'(e_m)\to 1\andeqn \tau\circ \pi_s'(e_m)\to 1, \,\, {\rm as}\,\, n\to\infty
\eneq
for all $t\in T(C\otimes \zo)$ and $\tau\in T(B_0\otimes \zo)$ and for all $s\in (0,1),$
where $\pi_e':  M_{\phi_0', \phi_1'}\to C\otimes \zo$  is the quotient map
and $\pi_s':  M_{\phi_0', \phi_1'}\to B_0\otimes \zo$ is the {{point evaluation}} at $s\in (0,1).$
Since both $B_0\otimes \zo$ and $C\otimes \zo$ have continuous {{scales,}} a standard compactness
argument shows that, for any $\eta>0,$ there exists
$n\ge 1$ such that
\beq\label{17830-1}
t\circ \pi_e'(e_m)>1-\eta/2\andeqn \tau\circ \pi_s'(e_m)\ge 1-\eta/2
\eneq
for all $m\ge n,$ $t\in T(C\otimes \zo),$ $\tau\in T(B_0\otimes \zo)$ and $s\in (0,1).$
Equivalently,   for any $k\in \N,$
\beq\label{17830-1+}
&&{{d_{t\circ \pi_e'}(e_{m+k}-e_{m+1})\le}} t\circ \pi_e'(e_{m+k}-e_m)<\eta\andeqn\\
&&d_\tau(e_{m+k}-e_{m+1})\le  \tau\circ \pi_s'(e_{m+k}-e_m)<\eta\rforal m\ge n.
\eneq
Then the conclusion follows from  \eqref{17830-1+} and {{Theorem}} \ref{CMphi0phi1}.

{{To see that $b$ is full, let {$a\in (B\otimes {\cal Z}_0)_+,$} then, there is an integer $K\ge 1$ such that
$\tau(a)<K\eta$ for all $\tau\in T({{B\otimes {\cal Z}_0}}).$  By Theorem \ref{CMphi0phi1}, $a\lesssim b\otimes 1_K$
(in $M_K({{B\otimes {\cal Z}_0}})$) which implies that $a$ is in the closed ideal generated by $b.$ So $b$ is full in $B\otimes \zo.$}}

{{Let $e_b\in B_+$ be a strictly positive element. Then $\tau_b({{e_b}})>0$ for all $\tau_b\in T(B).$
One then easily choose  an approximate identity such that $\tau_b(e_{b,m}-e_{b,n})>0$
for all $\tau_b\in T(B)$ and $m>n.$  It follows that $\tau((e_{b,m}-e_{b,n})\otimes c)>0$
for any $c\in {\zo}_+\setminus \{0\}$ and $\tau\in T(B\otimes \zo).$}}

\end{proof}

Recall that ${\cal W}\otimes Q\cong {\cal W}$ (see \cite{eglnkk0}).

\begin{lem}\label{Lfullness} {{Let $C$, $B_0$, $\phi_0, \phi_1: C\to B_0$ and the mapping torus   $B=M_{\phi_0, \phi_1}$ be  as in \ref{ConstphiW}.}}
Let $a\in {\cal W}_+$  with $\|a\|=1$ and $e\in  Q$ be non-zero projection.
Then there is a unitary $U\in M_2(M(B\otimes \zo)/B\otimes \zo)$ such that
\beq
U^*\diag(1, 0)U\in {\rm Her}(\pi\circ \phi_{k,odd}(f_{1/8}(a\otimes e))).
\eneq
\end{lem}

\begin{proof}
{{Put $b=a\otimes e\in {\cal W}\otimes Q={\cal W}.$}}
By the construction, $\tau(\Lambda_{n, k}(f_{1/4}({{b}})))>0$ for all $n, k$ and for all $\tau\in T(B\otimes \zo).$
It follows from Lemma \ref{L141N} and \ref{Lcont} that there is
a partial isometry $v\in M(B\otimes \zo)/B\otimes \zo$ such that
\beq
vv^*=1_{M(B\otimes \zo)/B\otimes \zo}\andeqn v^*v\le  f_{1/8}(\pi\circ \phi_{k, odd}({{b}})).
\eneq
Let $V=\diag(v, v^*).$ Then
\beq
V^*\diag(1, 0)V\le f_{1/8}(\pi\circ \phi_{k, odd}({{b}})).
\eneq
Put $v^*v=p.$   Choose $U=\begin{pmatrix} v &0\\ 1-p & v^*\end{pmatrix}.$ Then $U$ is a
unitary in $M_2(M(B\otimes \zo)/B\otimes \zo)$
such that
$$U^*\diag(1, 0)U=V^*\diag(1, 0)V\le f_{1/8}(\pi\circ \phi_{k,odd}({{b}})).$$
\end{proof}

\begin{lem}\label{L1032gl}
Let $C,$ $B_0,$ {{$B,$}} {{$\phi_0,\phi_1: C\to B_0$ and $\phi_{k,odd}:
{\cal W}\to M(B\otimes \zo)$ }} be as in {{\ref{ConstphiW}}}.
Let $\psi_1, \psi_2: {\cal W}\to M(B\otimes \zo)/B\otimes \zo$ be \hm s.
Then, for any {{$\ep>0$}}   and any finite subset ${\cal F}\subset {\cal W},$ there exists a {{unitary}}
$U\in M_8(M(B\otimes \zo))$ such that (for any $k$), for all $a\in {\cal F},$
\beq
\|\pi(U)^*\diag(\psi_1(a), \pi\circ \phi_{k, odd}(a))\pi(U)-\diag(\psi_2(a), \pi\circ \phi_{k, odd}(a))\|<{{\ep.}}
\eneq
(Here an element $a\in M_2(M(B\otimes \zo))$ is  identified with $\diag(a, 0_6)\in M_8(M(B\otimes \zo))$ as $M_2(M(B\otimes \zo))$ is identified as a corner subalgebra of $M_8(M(B\otimes \zo))$, where $0_6$ is the zero element in $M_6(M(B\otimes \zo))$.)
\end{lem}

\begin{proof}
Let $T: {\cal W}_+\setminus\{0\}\to \N\times \R\setminus \{0\}$
be a map such that $\id_{\cal W}$ is exactly $T$-${\cal W}_+\setminus \{0\}$-full.
Let ${\cal G}\subset {\cal W}$ (and $\dt>0$)  {{and}}  ${\cal H}\subset {\cal W}_+\setminus \{0\}),$
and {{the}} integer $K\ge 1$ be given by Corollary {{3.16}} of \cite{eglnkk0} for $\ep$ (and $A={\cal W}$ and ${\cal F}\subset {\cal W}$).

We will {{identify}}  ${\cal W}$ and ${\cal W}\otimes Q.$
Fix  a strictly positive element {{$e_W\in {\cal W}_+^{\bf 1}$}} and let $e_1, e_2,...,e_K\in Q$
be equivalent projections such that $\sum_{i=1}^K e_i=1_Q.$
Let $b_i=\pi\circ \phi_{k, odd}(e_W\otimes e_i),$ $i=1,2,..,K.$
There is, by Lemma \ref{Lfullness},  a unitary $W_1\in M_2(M(B\otimes \zo)/B\otimes \zo)$ such
that
\beq
W_1={{\begin{pmatrix} v & 0 \\ 1-p & v^*\end{pmatrix} }}
\andeqn W_1^*\diag(1, 0)W_1:=p\le f_{1/8}(b_1).
\eneq
Put $D_0:={\rm Her}(b_1), D_1:={\rm Her}(\pi\circ \phi_{k, odd}(e_W\otimes 1_Q))\subset M(B\otimes \zo)/B\otimes \zo.$
There is a unitary $V_i\in {\td{D_1}}$ such that
$V_i^*b_1V_i=b_i,$ $i=1,2,...,K.$ 
Put 
$\Psi_0: {\cal W}\to {{D_0}}$ by $\Psi_0(a)=\pi\circ \phi_{k, odd}(a\otimes e_1)$ for all $a\in {\cal W}.$
To obtain the result, \wilog, we may write
${{\pi\circ \phi_{k, odd}(a)}}=\diag(\Psi_0(a), \Psi_0(a),...,\Psi_0(a))$ (for all $a\in {\cal W}$ {{and}}
where $\Psi_0(a)$ repeats $K$ times)
from ${\cal W}$ to $M_K(D_0),$ and
{{view}} ${\rm Ad}\, W_1\circ \psi_i$ ($i=1,2$) as \hm s from ${\cal W}$ to $D_0.$
Note {{that,}} by Lemma \ref{Lfullness}, $\Psi_0$ is $T$-${\cal W}_+\setminus \{0\}$-full.
Applying Corollary {{3.16}} {{of \cite{eglnkk0},}} we obtain
a unitary $W_2\in M_{K+1}(D_0)^\sim$
such that, for all $a\in {\cal F},$
\beq
W_2^*\diag( {\rm Ad}\, W_1\circ \psi_1(a), \overline{\Psi_0(a)})W_2\approx_{\ep/2}
\diag({\rm Ad}\, W_1\circ \psi_2(a),  \overline{\Psi_0(a)}),
\eneq
where $\overline{\Psi_0(a)}=\diag(\Psi_0(a), \Psi_0(a),...,\Psi_0(a))$ ($\Psi_0(a)$
repeats $K$ times).
{{We view}}  $M_K(D_0)$ as \SCA\, of $M(B\otimes \zo)/B\otimes \zo).$
Then, one obtains a unitary $W_3\in M_2(M(B\otimes \zo)/B\otimes \zo)$ such
that  
\beq
W_3^*\diag( {\rm Ad}\, W_1\circ \psi_1(a), \pi\circ \phi_{k, odd}(a) )W_3\approx_{\ep}
\diag({\rm Ad}\, W_1\circ \psi_2(a), \pi\circ \phi_{k, odd}(a)).
\eneq

Set $W_4=\diag(1,W_0, 1)^*\cdot \diag(W_1, 1_2)\diag(1, W_0,1)\in {{ M_4(M(B\otimes \zo)/B\otimes \zo)}},$ where $W_0=\begin{pmatrix} 0 & 1\\ 1&0\end{pmatrix}.$
Then
$
W_4^*\diag(\psi_i(a), \pi\circ\phi_{k,odd}(a),0,0)W_4=\diag({\rm Ad}\, W_1\circ \psi_i(a), \pi\circ\phi_{k,odd}(a),0,0).
$
Set $W_5=W_4\diag(W_3, 1_2) W_4^*.$
Then $W_5$ is a unitary
{{in $ M_4(M(B\otimes \zo)/B\otimes \zo)$}}
such that (identifying {{$EM_4(M(B\otimes \zo)/B\otimes \zo))E$ with $M_2(M(B\otimes \zo)/B\otimes \zo),$
where $E=\diag(1,1,0,0)$}})
\beq
W_5^*\diag(\psi_1(a), \pi\circ \phi_{k,odd}(a){{,0,0}})W_5\approx_{\ep} \diag(\psi_2(a), \pi\circ \phi_{k, odd}(a){{,0,0}})
\eneq
Let $W_6=\diag(W_5, W_5^*).$ By replacing $W_5$ by $W_6,$
we may assume that there exists a unitary $U\in M_8(M(B\otimes \zo))$ such that, for all $a\in {\cal F},$
\beq
\|\pi(U)^*\diag(\psi_1(a), \pi\circ \phi_{k, odd}(a){{,0_6}})\pi(U)-\diag(\psi_2(a), \pi\circ \phi_{k, odd}(a){{,0_6}})\|<\ep.
\eneq
\end{proof}

\begin{prop}[cf.10.5 of \cite{GLII}]\label{AnuniqW}
Let $B_0$ and $C$ be  non-unital separable  simple stably projectionless \CA s with  stable rank one, {{with continuous scales and with $T(B_0)\not=\emptyset$}}.
Suppose  that $\phi_0, \phi_1: C\to B_0$ are \hm s which map strictly positive elements
to strictly positive elements such that $\tau\circ \phi_0=\tau\circ \phi_1$ for all $\tau\in T(B_0).$
Put $B=M_{\phi_0, \phi_1}.$
Fix an integer $k_0\ge 1.$
Let $j_{w,z}: {\cal W}\to  M_{k_0}(\zo)$ be an embedding  which maps strictly positive elements
to strictly positive elements
and $d: \zo\to \C\cdot 1_{M_{k_0}({\tilde B})}\otimes \zo\subset
M_{k_0}({\tilde B}\otimes \zo)\subset M(M_{k_0}(B\otimes \zo))$ be the embedding defined by $d(z)=1\otimes z$ for all $z\in \zo.$

Let $\ep>0$ and ${\cal F}\subset {{\cal W}}$ be a finite subset.
Then there are  {{integers}} $k\ge 1,$  $K\ge 1$ and a unitary $u\in M_{K+1}(M(M_{k_0}(B\otimes \zo)))$ such that
$$
\|u^*(d_K\circ j_{w,z}(a),0)u-(d_K\circ j_{w,z}(a)\oplus \phi_{k,odd}(a))\|<\ep\tforal a\in {\cal F},
$$
where
$$
d_K(z)=\diag(\overbrace{d(z),d(z),...,d(z)}^K)\tforal z\in \zo.
$$

\end{prop}

\begin{proof}
Keep in mind  that $B_0$ has continuous scale.
So $B_0\otimes \zo$ also has continuous scale.
Put $\phi_i'=\phi_i\otimes {\id}_{\zo}: C\otimes \zo\to B_0\otimes \zo.$
Then $M_{\phi_0, \phi_1}\otimes \zo=M_{\phi_0', \phi_1'}.$ Moreover, $B_0\otimes \zo$ is ${\cal Z}$-stable.
Since ${\cal W}$ is simple and
$d\circ j_{w,z}$ maps a strictly positive element
to that of $\C \cdot 1_{M_{k_0}({\tilde B})}\otimes \zo$  which is not in
$M_{k_0}(B\otimes \zo).$  {{Moreover, }}
by \ref{Lcont} {{(also \ref{L141N}),}}
$d\circ j_{w,z}(a)$ is full in $M(M_{k_0}(B\otimes \zo))$ for every
$a\in {{\cal W}}_+\setminus \{0\}.$
There is a map $T: {{\cal W}}_+\setminus \{0\}\to \N\times \R_+\setminus \{0\}$
such that
$d\circ j_{w,z}$ is $T$-${{\cal W}}_+\setminus \{0\}$-full in $M(M_{k_0}(B\otimes \zo)).$

Let $K\ge 1$ be the integer required by   {{Cor.\,3.16}} of \cite{eglnkk0}
for $\ep/2$ (in place of $\ep$),
${\cal F}$ and $T.$
Applying Cor. 3.16 of \cite{eglnkk0}, one obtains (note that $M(M_{k_0}(B\otimes \zo))$ is unital),
a unitary $v\in M_{K+1}(M(M_{k_0}(B\otimes \zo)))$ such that
$$
\|u^*(d_K\circ j_{w,z}(a),0)u-(d_K\circ j_{w,z}(a)\oplus \phi_{k,odd}(a))\|<\ep\tforal a\in {\cal F}.
$$
\end{proof}

\begin{thm}\label{Tappmul}
Let $A$ be a non-unital separable  amenable \CA.
Let $\ep>0$ and ${\cal F}\subset A$ be finite subset.

There exists $\dt>0$ with $\dt<\ep/2,$  a finite subset ${\cal G}\subset A$
with ${\cal F}\subset {\cal G}$ and an integer $K\ge 1$ satisfying the following:
Suppose that $B$ is a mapping torus  as in \ref{AnuniqW} {{and}}  $\phi: A\to
M_k({\tilde B}\otimes \zo)$ is a ${\cal G}$-$\dt$-multiplicative \cpc\,  ($k\ge 1$ is an integer),
%
$\psi_{z,w}: M_k(\zo)\to {{\cal W}}$ and $\psi_{w,z}: {{\cal W}}\to M_k(\C\cdot 1_{{\tilde B}}\otimes \zo)\cong M_k(\zo)$
{{are \hm s}}
which
map strictly positive elements to strictly positive elements such that
\beq\label{109-e-1}
\|\pi\circ (\phi(a))-(\psi_{w,z}\circ \psi_{z,w}\circ \pi\circ (\phi(a))\|<\dt\rforal a\in {\cal G},
\eneq
where $\pi: M_k({\tilde B}\otimes \zo)\to M_k(\C\cdot 1_{{\tilde B}}\otimes \zo)$
is the quotient {{map.}}
{{Then}} there exists an ${\cal F}$-$\ep$-multiplicative \cpc\, $L_0: A\to M_{{8(K+2)}}(M_k(B\otimes \zo))$
and an ${\cal F}$-$\ep$-multiplicative \cpc\, $L_1: A\to M_{{8(K+2)}}(M_k({\tilde B}\otimes \zo))$
such that  {{$L_0(A)\perp L_1(A)$ and}}
$$
\|{{(L_0(a)+ L_1(a))}}-(\phi(a)\oplus d_{K+1}\circ s\circ \phi^{\pi}(a))\|<\ep\rforal a\in {\cal F},
$$
{{(viewing $(\phi(a)\oplus d_{K+1}\circ s\circ \phi^{\pi}(a))$ as an element
in $M_{8(K+2)}(M_k(\td B\otimes \zo)$),}}
where  $\phi^{\pi}=\psi_{w,z}\circ \psi_{z,w}\circ \pi\circ \phi,$ $s: M_k(\C\cdot 1_{\tilde{B}}\otimes \zo)\to M_k({\tilde B}\otimes \zo)$
is  the nature embedding, and such that $L_0$ and $L_1$ are of the following forms:
$$
L_0(a)=p_m^{1/2}(\phi(a)\oplus {{d_{K+1}}}\circ s\circ \phi^{\pi}(a))p_m^{1/2}
\tforal a\in A
$$
for some $m\ge m_0,$  where  $\{p_m\}$ is an approximate identity for  $M_{{8(K+2)}}(M_k(B\otimes \zo))$ and,
there  are  ${\cal G}$-$\dt$-multiplicative \cpc\, $L_{0,0}: A\to {{\cal W}}$ and    $L_{0,0}({\cal F})$-$\ep/2$-multiplicative  \cpc\,
$L_{w,b}: {{\cal W}}\to {{M_{8(K+2)}}}(M_k({\tilde B}\otimes \zo))$ such that
$L_1=L_{w,b}\circ L_{00}.$
\end{thm}

\begin{proof}
The proof is almost identical to that of 10.7 of \cite{GLII}.
Fix $1/2>\ep>0$ and a finite subset ${\cal F}\subset A.$ We may assume
that ${\cal F}\subset A^{\bf 1}.$
Let
${\cal G}=\{ab: a, b\in {\cal F}\}\cup {\cal F}.$
{{Let $e_n:=f_{n+2}\otimes e_n'\in M_k(B\times \zo)$ {{with $f_{n+2}$ and $e_n'$}} as}}
described  in \ref{ConstphiW}.  Note
that $\{e_n\}$ forms an approximate identity for $M_k(B\otimes \zo)$ and
$e_{n+1}e_n=e_n$ ($n\in \N$).
%
Let $\dt_1>0$ (in place $\dt$) be in 10.6 of \cite{GLII}
for $\ep/64.$

Let $\dt=\min\{\dt_1/2^{12}, \ep/2^{12}\}.$
We view $M_k({\tilde B}\otimes \zo)$ as a  \SCA\, of $M(M_k(B\otimes \zo)).$
Suppose that $\phi: A\to M_k({\tilde B}\otimes \zo)$ is {{${\cal G}$-$\dt$-multiplicative}} \cpc.
Suppose that there are  \hm s
$\psi_{z,w}:  M_k(\zo) \to {{\cal W}}$ and $\psi_{w,z}: {{\cal W}}\to M_k(\C\cdot 1_{{\tilde B}}\otimes \zo)$
such that
\beq\label{Tappmul-9}
\|\pi\circ \phi(a)-(\psi_{w,z}\circ \psi_{z,w}\circ \pi\circ (\phi(a)))\|<
\dt\rforal a\in {\cal G}.
\eneq
Recall that  $\phi^\pi=\psi_{w,z}\circ \psi_{z,w}\circ \pi\circ \phi.$
Put $\phi^W=\psi_{z,w}\circ \pi\circ \phi.$ Thus $\psi_{w,z}\circ \phi^W=\phi^{\pi}.$
Let $K$ be  the integer in \ref{AnuniqW} associated with
$\dt$ (in place of $\ep$) and ${{\phi^W({\cal G})}}$ (in place of ${\cal F}$).

 {{Applying}}
 {{\ref{L1032gl} {{to homomorphisms $\Pi\circ \phi_{k_0,even}$ and $d_{K+1}\circ \psi_{z,w}$}}, we obtain}}
 a unitary
$U_1\in M_{8(K+2)}(M(M_k(B\otimes \zo)))$ such that (for some large $k_0$)
\beq\label{Tappmul-10}
\hspace{-0.3in}\Pi(U_1)^*{{(\Pi\circ \phi_{k_0,W}(\phi^W(a)))}}\Pi(U_1)\approx_{\dt}\diag(d_{K+1}\circ\psi_{w,z}\circ \phi^W(a)), \Pi\circ \phi_{k_0,odd}(\phi^{W}(a))
\eneq
for all $ a\in {\cal G},$  where {{$\Pi: M_{l}(M(M_k(B\otimes \zo)))\to M_{l}(M(M_k(B\otimes \zo)))/M_{l}(M_k(B\otimes \zo))$}} is the quotient map, for $l=1$ or $l=K+2$. 

Recall that $s: M_k(\C\cdot 1_{\tilde{B}}\otimes \zo)\to M_k({\tilde B}\otimes \zo)$
is the {{natural embedding}}
such that
$$
\pi\circ s(a)=a\rforal a\in M_k(\C\cdot 1_{\tilde{B}}\otimes \zo).
$$

Consider $L_{1,1}: A\to M(M_k(B\otimes \zo))$ 
defined by
$L_{1,1}= \phi_{k_0,W}\circ  \phi^{W}$ and $L_{1,0}': A\to M_{K+2}(M(M_k(B\otimes \zo)))$ 
defined by
$$
L_{1,0}'(a)=\diag(d'_{K+1}\circ s\circ\psi_{w,z}\circ \phi^{W}(a)), \phi_{k_0,odd}(\phi^{W}(a))\rforal a\in A,
$$
where
$
d'_{m}(c):=\diag(\overbrace{c,c,...,c}^{m}).
$
By \ref{AnuniqW},   there is another unitary $U_2\in M_{K+2}(M(M_k(B\otimes \zo)))$
such that
\beq\label{Tappmul-n10}
\|U_2^*L_{1,0}'(a)U_2-d'_{K+1}\circ s\circ\psi_{w,z}\circ \phi^W(a)\|<\dt\rforal a\in {\cal G}.
\eneq
%
Define $L_{1,0}:  A\to M_{K+1}(M_k(\td B\otimes \zo))$ by
$$
L_{1,0}(a)=d'_{K+1}\circ s\circ\phi^{\pi}(a)\rforal a\in A.
$$
Put
$\Phi=\phi\oplus {{d'_{{K}}}}\circ s\circ \phi^{\pi}$
and $U=U_1U_2.$  By \eqref{Tappmul-10} and \eqref{Tappmul-n10},
for each $a\in {\cal G},$ there exist $b(a), b'(a)\in M_{{8(K+2)}}(M_k(B\otimes \zo))$
(see also \eqref{109-e-1})
 with
$\|b(a)\|,\,\|b'(a)\|\le 1$
such that
\beq\label{Tappmul-11}
&&\|U^*L_{1,1}(a)U- L_{1,0}(a)+b(a)\|<2\dt\andeqn\\
&& \|U^*L_{1,1}(a)U-\Phi(a) +b'(a)\|<2\dt\rforal a\in {\cal G}.
\eneq

%
\vspace{-0.05in}
Put ${\bar e}_n=\diag(\overbrace{e_n,e_n,...,e_n}^{8(K+2)}),$
${\bar e}_n'=\diag(\overbrace{e_n',e_n',...,e_n'}^{8(K+2)}),$
$n=1,2,....$
Let $p_n=U^*{\bar e_n}U,$ $n=1,2,....$
Then $\{p_n\}$ is an approximate identity for ${{M_{8(K+2)}}}(M_k(B\otimes \zo)).$
Let $S=\N\setminus \{4n,\, 4n-1: n\in \N\}.$
If $m\in S,$ {{then}}
\beq\label{Tappmul-12+1}
&&\hspace{-0.45in}(1-p_m) ({{U^*}}f_{n+1}^{1/2}\otimes ({\bar e}_{4n}'-{\bar e}_{4n-1}'){{U}})
=\begin{cases}  {{U^*}}f_{n+1}^{1/2}\otimes ({\bar e}_{4n}'-{\bar e}_{4n-1}'){{U}}
\,\,\, &\text{if}\,\,\, m<4n-1;\\
   0 \,\,\, &\text{if} \,\,\, m>4n                                                                      
   \end{cases}\\\label{Tappmul-12+2}
 && \hspace{-0.45in}   \andeqn  p_m(1-p_m)({{U^*}}f_{n+1}^{1/2}\otimes ({\bar e}_{4n}'-{\bar e}_{4n-1}'){{U}})
      =0\rforal {{m\in S}}.
                                                                      \eneq
{{(Here $1-p_m:=1_{M_{8(K+2)}}-p_m.$)}}
There is {{an integer}} $N\ge 1$ such that, for any $m\ge N$ and $m\in S,$
\beq\label{Tappmul-13}
&&\|(1-p_m)(U^*L_{1,1}(a)U)-(1-p_m)L_{1,0}(a)\|<4\dt,\\\label{Tappmul-13+1}
&&\|(U^*L_{1,1}(a)U)(1-p_m)-L_{1,0}(a)(1-p_m)\| <4\dt\\\label{Tappmul-13+2}
&&\|(1-p_m)(U^*L_{1,1}(a)U)-(1-p_m)\Phi(a)\|<4\dt\andeqn\\\label{Tappmul-13+3}
&&\|(U^*L_{1,1}(a)U)(1-p_m)-\Phi(a)(1-p_m)\| <4\dt \rforal a\in {\cal G}.
\eneq
Note that, by the construction of $\phi_{k_0,W}$ {{(see (\ref{09-15-2021}) also)}} and \eqref{Tappmul-12+2},
\beq\label{Tappmul-14}
(1-p_m)(U^*L_{1,1}(a)U)&=&(U^*L_{1,1}(a)U)(1-p_m)\\\label{Tappmul-14+1}
&=&(1-p_m)(U^*L_{1,1}(a)U)(1-p_m)\rforal a\in A.
\eneq
It follows from  \eqref{Tappmul-13},
\eqref{Tappmul-13+1}, \eqref{Tappmul-13+2},
\eqref{Tappmul-13+3} and
\eqref{Tappmul-14}, for all $m\ge N$ and $m\in S,$
\beq\label{Tappmul-15}
\|p_m\Phi(a)-\Phi(a)p_m\|<8\dt \andeqn \|(1-p_m)L_{1,0}(a)-L_{1,0}(a)(1-p_m)\|<8\dt
\eneq
for all $a\in {\cal G}.$
By the choice of $\dt_1$ and   10.6 of \cite{GLII},
for all $a\in {\cal G},$
\beq\label{191026-1}
&&\|p_m^{1/2}\Phi(a)p_m^{1/2}-p_m\Phi(a)\|<\ep/64\andeqn\\\label{191026-2}
&&\|(1-p_m)^{1/2}L_{1,0}(a)(1-p_m)^{1/2}-(1-p_m)L_{1,0}(a)\|<\ep/64.
\eneq
Moreover, the map $a\mapsto (1-p_m)(U^*L_{1,1}(a)U)$ is a ${\cal G}$-$\dt$-multiplicative.
By \eqref{191026-2} and \eqref{Tappmul-13},
$a\to (1-p_m)^{1/2}L_{1,0}(a)(1-p_m)^{1/2}$ is ${\cal F}$-$\ep$-multiplicative.
Define  {{(for $m\in S$ and $m\ge N$)}}
$$
L(a)=p_m\Phi(a)+(1-p_m)(U^*L_{1,1}(a)U)\rforal a\in A.
$$
Then, by \eqref{Tappmul-13},
\beq\label{Tappmul-15+}
\|L(a)-\Phi(a)\|<4\dt\rforal a\in {\cal G}.
\eneq
Consequently,
\beq\label{Tappmul-16}
\|L(ab)-L(a)L(b)\|<8\dt \rforal a, b\in {\cal G}.
\eneq
We compute  that
\beq\label{Tappmul-17}
L(ab)=p_m\Phi(ab)+(1-p_m)(U^*L_{1,1}(ab)U)\rforal a, b\in A,
\eneq
and, for all $a, b\in {\cal G},$ by \eqref{Tappmul-12+2}, \eqref{Tappmul-14+1} and \eqref{Tappmul-15},
\beq\nonumber
&&\hspace{-1in}L(a)L(b)=(p_m\Phi(a)+(1-p_m)(U^*L_{1,1}(a)U))(p_m\Phi(b)+(1-p_m)(U^*L_{1,1}(b)U))\\\nonumber
&=& p_m\Phi(a)p_m\Phi(b)+((1-p_m)((U^*L_{1,1}(a)U))(1-p_m)(U^*L_{1,1}(b)U))\\\nonumber
\hspace{0.2in}&\approx_{8\dt+\dt}& p_m\Phi(a)\Phi(b)p_m+(1-p_m)(U^*L_{1,1}(ab)U).
\eneq
Combining this with \eqref{Tappmul-17} and \eqref{Tappmul-16}, {{we obtain}}
\beq\label{Tappmul-18}
\|p_m\Phi(ab)-p_m\Phi(a)\Phi(b)p_m\|<8\dt+8\dt+\dt=17\dt\rforal a, b\in {\cal G}.
\eneq
Therefore (see 10.6 of \cite{GLII})
\beq\label{Tappmul-19}
\|p_m^{1/2}\Phi(ab)p_m^{1/2}-p_m^{1/2}\Phi(a)p_m^{1/2}p_m^{1/2}\Phi(b)p_m^{1/2}\|<
17\dt+3\ep/64<\ep/16.
\eneq
Define $L_0(a)=p_m^{1/2}\Phi(a)p_m^{1/2}$ and $L_1(a)=(1-p_m)^{1/2}(U^*L_{1,1}(a)U)(1-p_m)^{1/2}.$
 {{It follows from (\ref{Tappmul-12+2}) and $m\in S$ that the images of $L_0$ and $L_1$ are mutually orthogonal.}}
By \eqref{Tappmul-15+} and the choice of $\dt_1,$
we finally have
$$
\|(L_0(a)+L_1(a))-\Phi(a)\|<\ep\rforal a\in {\cal F}.
$$
Let $L_{00}=\phi^W: A\to {\cal W}$ and $L_{w,b}: {\cal W}\to M_{{8(K+2)}}(M_k({\tilde B}\otimes \zo))$
be defined by $L_{w,b}(b)=(1-p_m)^{1/2}(d_K'\circ s\circ \psi_{w,z}(b))(1-p_m)^{1/2}$ for $b\in {\cal W}.$
Then $L_1=L_{00}\circ L_{w,b}.$
%
\end{proof}


%
%


\begin{thm}\label{TExistence}
Let $A$ be a non-unital separable  amenable \CA\, which satisfies the UCT
{{and}} has the property (W)
and let $B$ be
a separable  \CA\, as in \ref{AnuniqW}.
For any $\af\in KL(A,{{B\otimes \zo}}),$
there exists an asymptotic sequential morphism $\{\phi_n\}$ from
$A$ into $B\otimes \zo\otimes {\cal K}$ such that
$$
[\{\phi_n\}]=\af.
$$
\end{thm}

\begin{proof}

Let ${\cal P}\subset \underline{K}(A)$ be a finite subset.  Let $\ep>0$ and ${\cal F}\subset A$ be a finite subset.
We assume that, {{for}} any ${\cal F}$-$\ep$-multiplicative \cpc\, $L$ from $A,$
$[L]|_{\cal P}$ is well-defined.

If follows from 10.2 of \cite{GLII}
that there exist  sequences of approximately multiplicative \cpc s
$\Phi_n: A\to B^{\sim}\otimes \zo \otimes {\cal K}$   and $\Psi_n: A\to \C \cdot 1_{B^{\sim}}\otimes \zo\otimes {\cal K}$
such that, for any
finite subset ${\cal Q}\subset \underline{K}(A),$
$$
[\Phi_n]|_{\cal Q}=\af|_{\cal Q}+[\Psi_n]|_{\cal Q}
$$
for all sufficiently large $n,$  where $\Psi_n=s\circ \pi\circ \Phi_n$
(without loss of generality)
and $\pi: B^{\sim}\otimes \zo \otimes {\cal K}\to \C \cdot 1_{B^{\sim}}\otimes \zo\otimes {\cal K}$ is the quotient
map {{and $s: \C \cdot 1_{B^{\sim}}\otimes \zo\otimes {\cal K}\to B^{\sim}\otimes \zo \otimes {\cal K}$
is the splitting map.}}
Fix a sufficiently large $n.$

Let $\{e_{i,j}\}$ be a system of matrix unit for ${\cal K}$ and
let $E$ be the unit of the unitization of  $1_{B^{\sim}}\otimes \zo.$
By considering   maps $a\mapsto (E\otimes \sum_{i=1}^k e_{i,i})\Phi_n(a)(E\otimes \sum_{i=1}^k e_{i,i})$
and $a\mapsto (E\otimes \sum_{i=1}^k e_{i,i})\Psi_n(a)(E\otimes \sum_{i=1}^k e_{i,i}),$
\wilog, we may assume that
the image of $\Phi_n$ is in $M_k(B^{\sim}\otimes \zo)$ and
that of $\Psi_n$ is also in $M_k(\C\cdot 1_{B^{\sim}}\otimes \zo)$ for some
sufficiently large $k.$

Define $\imath^{\circledast}: B^{\sim}\otimes \zo\otimes {\cal K}\to  B^{\sim}\otimes \zo\otimes {\cal K}$
by
$\imath^{\circledast}(b\otimes z\otimes k)=b\otimes j^{\circledast}(z)\otimes k$
for all $b\in B^{\sim},$ $z\in \zo$ and $k\in {\cal K},$  where
$j^{\circledast}: \zo\to \zo$ is an automorphism such that $j_{*0}=-{\rm id}|_{K_0(\zo)}$ as defined
in 8.13 of \cite{GLII}.
Note that
$$
s\circ \pi(\Phi_n\oplus s\circ \pi\circ i^{\circledast}\circ \Phi_n)=\Psi_n\oplus s\circ \pi \circ i^{\circledast}\circ \Phi_n.
$$

Let $\dt>0$ and {{${\cal G}\subset A$}}  be a finite subset.

{{By}} virtue of  8.14 of \cite{GLII},
replacing $\Phi_n$ by $\Phi_n\oplus  s\circ \pi\circ i^{\circledast}\circ \Phi_n$ and
replacing $\Psi_n$ by $\Psi_n\oplus  s\circ \pi\circ i^{\circledast}\circ \Phi_n ,$
and by implementing a unitary in unitization of $M_k(\C\cdot 1_{B^{\sim}}\otimes \zo),$
we may assume that
$$
\|\pi\circ \Phi_n(g)-\phi_{w,z}\circ \phi_{z,w}\circ \pi(\Phi_n(a))\|<\dt\rforal g\in {\cal G}
$$
and  $\Psi_n$ factors through ${{\cal W}},$ in particular,
$[\Psi_n]|_{\cal P}=0.$
In other words,
\vspace{-0.13in}\beq\label{Texist-n1}
[\Phi_n]|_{\cal P}=\af|_{\cal P}.
\eneq
By applying \ref{Tappmul}, we obtains  an integer $K\ge 1,$ ${\cal F}$-$\ep$-multiplicative \cpc s
$L_{0,n}: A\to {{M_{8(K+2)k}}}(B\otimes \zo),$ $L_{1,n}: A\to {{M_{8(K+2)k}}}(B^{\sim}\otimes \zo)$
and $L_{2,n}: A\to M_{(K+1)k}(B^\sim\otimes \zo)$  such that
\beq\label{TE-9}
\|L_{0,n}(a)\oplus L_{1,n}(a)-\Phi_n(a)\oplus L_{2,n}(a)\|<\ep\rforal a\in {\cal F},
\eneq
where $L_{1,n}$ and $L_{2,n}$ factor through ${{\cal W}}.$
In particular,
\beq\label{TE-10}
[L_{1,n}]|_{\cal P}=[L_{2,n}]|_{\cal P}=0.
\eneq
It follows that, using \eqref{Texist-n1} and {{\eqref{TE-9},}}
\beq\label{TE-11}
[L_{0,n}]|_{\cal P}=\af|_{\cal P}.
\eneq

Choose $\phi_n=L_{0,n}$ (for all sufficiently large $n$).

\end{proof}


\begin{cor}\label{CExistence}
Let $A$ be a non-unital separable  amenable \CA\, which satisfies the UCT
{{and}} has the property {{(W),}}
and let $B$ be
a separable  \CA\, as in \ref{AnuniqW} which is ${\cal Z}$-stable.
For any $\af\in KL(A,B),$
there exists  a sequence of \cpc\, maps $\phi_n: A\to B\otimes M_{k(n)}$
(for some increasing sequence $\{k(n)\}$)   such that
$$
[\{\phi_n\}]=\af.
$$
\end{cor}

\begin{proof}
It follows from  \ref{TExistence}
 that
there exists a sequential morphism $\{\psi_n\}$ from $A$ into $B\otimes \zo\otimes {\cal K}$ such that
$[\{\psi_n\}{]}=\kappa_{\zo}^{-1}\circ \af$ {{(see 12.4 and 12.5 of \cite{GLII} with $C=B$).}}
{{Note, by \ref{LMtorusW}, $B$ has the property (W).}}
Let $\{\Psi_n\}$ be a sequential morphism from $B\otimes \zo$ to $B\otimes {\cal K}$
such that $[\{{{\Psi_n}}\}]=\kappa_{\zo}$ given
by 12.5 of \cite{GLII}.  Define $\phi_n=\Psi_{k(n)}\circ \psi_n$ for
a choice of $\{k(n)\}.$ Then $\{\phi_n\}$ meets the requirement.

\end{proof}

  \section{Stable homotopy}

Let us state the following non-unital version {{of}} stable uniqueness theorem.

\begin{thm}\label{TSuniq}
Let $C$ be an amenable separable \CA\, which satisfies the UCT.
For any $\ep>0$ and any finite subset
${\cal F}\subset C,$ there exists a finite subset ${\cal G}\subset C,$ $\dt>0,$
and a finite subset ${\cal P}\subset \underline{K}(C)$ satisfying the following
condition:   Suppose that $M\times N: C_+\setminus \{0\}\to
\R_+\setminus \{0\}\times \N$ is  a map,
{{$A$ is}} a $\sigma$-unital \CA\, with the property
that there is an embedding \, $j_{c,A}: C\to A$
such that $j_{c,A}(c)$ is  $(M(c), N(c))$-uniformly
full in $A$ for all $c\in C_+\setminus \{0\}$
(see Definition 5.5 of \cite{eglnp} {{and Definition  3.11 of \cite{eglnkk0}}}),  and
$\phi_1, \phi_2: C\to A$ are two
${\cal G}$-$\dt$-multiplicative \cpc s such that
\beq\label{2020-721-1}
[\phi_1]|_{\cal P}=[\phi_2]|_{\cal P}.
\eneq
{{Then}} there exist an integer $K\ge 1$ and a unitary $U\in U(M_{K+1}(A))^\sim$
such that, for all $c\in {\cal F},$
\beq
\|{\rm Ad}\, U\circ \diag(\phi_1(c), j_{c,A}(c)\otimes 1_K)-\diag(\phi_2(c), j_{c,A}(c)\otimes 1_K)\|<\ep.
\eneq
\end{thm}

\begin{proof}
The proof is known and a number of versions have appeared.
For this non-unital version, we will apply, for example, Theorem 3.12 of \cite{eglnkk0}.
This application is essentially the same as that of Theorem 9.2  of \cite{Lnamj} and that of Theorem 4.15 of
\cite{GLII}.  {{Here we}} will  {{keep the proof}} brief.

Suppose that the conclusion  is false. We then obtain  a positive number $\ep_0>0,$ a finite
subset {{${\cal F}_0\subset C,$}} a sequence of finite {{subsets}} ${\cal P}_n\subset
\underline{K}(C)$ with ${\cal P}_n\subset {\cal P}_{n+1}$ and
$\cup_n {\cal P}_n=\underline{K}(C),$  a sequence of {{\CA s}} $A_n$ {{with embeddings \, $j_{c,A_n}: C\to A_n$
such that $j_{c,A_n}(c)$ are  $(M(c), N(c))$-uniformly
full in $A_n$ for all $c\in C_+\setminus \{0\}$,}} and
sequences of \cpc s $\{L_n^{(1)}\}$ and $\{L_n^{(2)}\}$ from $C$  to $A_n$  such that
\beq
&&\hspace{-0.3in}\lim_{n\to\infty}\|L_n^{(i)}(ab)-L_n^{(i)}(a)L_n^{(i)}(b)\|=0\rforal a, b\in C,\\
&&\hspace{-0.3in}{[L_n^{(1)}]}|_{{\cal P}_n}=[L_n^{(2)}]|_{{\cal P}_n}\andeqn\\\label{2020-721-0}
&&\hspace{-0.3in}\inf \{\sup\{\|u_n^*\diag(L_n^{(1)}(c), \Psi_{c,a,n}^{(k)}(c))u_n-\diag(L_n^{(2)}(c), \Psi^{(k)}_{c,a,n}(c))\|:c\in {\cal F}_0\}\ge \ep_0,
\eneq
where  $\Psi_{c,a,n}^{(k)}(c)=j_{c,A_n}(c)\otimes 1_k,$ and where {{the}}
infimum is taken among all integers $k,$  and all possible unitaries in $M_{k+1}(A)^\sim.$
We assume that, for each $n,$ $[L_n^{(j)}]|_{{\cal P}_ n}$ is well defined.

We may write ${\cal P}_k=\cup_{i=1}^{m(k)}{\cal P}_{k,i},$
where ${\cal P}_{k,i}\subset K_0(C\otimes D_i)$ (where $D_0=\C,$ {{$D_1=C_0(\T\backslash \{1\}),$}}
and $D_i$ is some commutative \CA\,  with finitely generated {{torsion $K_0$ group or $K_1$ group}}, 
{{$i=2,...,m(n).$}})


For $x\in {\cal P}_{n,i},$ we may write
$x=[p_{x,n,i}]-[q_{x, n,i}],$ where $p_{x, n,i}, q_{x, n,i}\in M_{R(n)}({{(C\otimes D_i)^\sim}})$
are projections.
By \eqref{2020-721-1}, there is an integer $r(n)>R(n)$
such that
\beq\label{2020-721-2}
&&v_{x,n,i}^*\diag({\bar p}_{x, n,i,1}, 1_{r(n)-R(n)})v_{x, n,i}=\diag({\bar p}_{x,n,i,2}, 1_{r(n)-R(n)})\andeqn\\\label{2020-721-3}
&&w_{x,n,i}^*\diag({\bar q}_{x,n,i,1}, 1_{r(n)-R(n)})w_{n,i}=\diag({\bar q}_{x,n,i,2}, 1_{r(n)-R(n)})
\eneq
for all $x\in {\cal P}_{n,i},$ $i=1,2,...,m(n),$
where $v_{x,n,i}, w_{x,n,i}\in M_{r(n)}{{((C\otimes D_i)^\sim)}}$ are partial isometries, {{and}}
${\bar p}_{x,n,i,j}$ and ${\bar q}_{x, n, i,j}$
are projections representing $[(L_n^{(j)}\otimes {\rm id}_{D_i})^\sim(p_{x, n,i})],$ $j=1,2.$
We also assume that $r(n+1)\ge r(n).$

Define $B_n:=A_n\otimes M_{r(n)},$ $B=\prod_{n=1}^\infty B_n,$   and $Q_1=B/\bigoplus_{i=1}^\infty B_n.$
Let $\pi: B\to Q_1$ be the quotient map.   Define $\phi_j: C\to B$
by $\phi_j(c)=\{L_n^{(j)}(c)\}$ for all $c\in C,$ and ${\bar \phi_j}:=\pi\circ \phi_j,$ $j=1,2.$
 It follows from \eqref{2020-721-2} and \eqref{2020-721-3} (recall $B_n= A_n\otimes M_{r(n)}$) that, for all $n,$
 \beq
 [\phi_1]_{{\cal P}_n}={{[\phi_2]|_{{\cal P}_n}.}}
 \eneq
Therefore
\beq
[{\bar \phi}_1]=[{\bar \phi}_2] \,\, {{{\rm in} \,\,\, KL(C, Q_1).}}
\eneq
Define $\Phi_n: C\to B_n$ by $\Phi_n=j_{c,A_n}\otimes {\rm id}_{M_{r(n)}}$
and $\Phi: C\to B$ by $\Phi(c)=\{\Phi_n(c)\}$ for all $c\in C.$
We claim that $\Phi(c)$ is {{$(M(c), N(c))$-full}} in $B$ for all $c\in C_+\setminus \{0\}.$
To see this, let $b=\{b_n\}\in B_+\setminus \{0\}$  with $\|b\|\le 1.$
Let $\{e_m^{(n)}\}$ be an approximate identity for $A_n,$ $n=1,2,....$
Put $E_m^{(n)}=e_m^{(n)}\otimes 1_{M_{r(n)}}.$ Then $\{E_m^{(n)}\}$
is an approximate identity for $B_n,$ $n=1,2,....$
Fix $\eta>0.$ There is $\{m_n\}$ such that
\beq
\|b_n^{1/2}(E_{m(n)}^{(n)})b^{1/2}_n-b_n\|<\ep/2,\,\,n=1,2,....
\eneq
Since $j_{c,A_n}(c)$ is $(M(c), N(c))$ full, there are $y_{n,1}, y_{n,2},...,y_{n, N(c)}\in A_n$
with $\|y_{n,j}\|\le M(c),$ $j=1,2,...,N(c),$ such that
\beq
\|\sum_{j=1}^{N(c)} y_{n,j}^*j_{c, A_n}(c)y_{n,j}-e_{m_n}^{(n)}\|<\ep/2.
\eneq
Let ${\bar y}_{n,j}=y_{n,j}\otimes 1_{r(n)}.$
Then
\beq
\|\sum_{j=1}^{N(c)} {\bar y}_{n,j}^* \Phi_n(c) {\bar y}_{n,j}-E_{m_n}^{(n)}\|<\ep/2.
\eneq
Let $Y_j=\{{\bar y}_{n,j}\}.$ Then $Y_j\in \prod_{n=1}^\infty B_n,$ $j=1,2,...,N(c).$
{{Also,}}
\beq
\|\sum_{j=1}^{N(c)}Y_j^*\Phi(c)Y_j-\{E_{m_n}^{(n)}\}\|<\ep/2.
\eneq
It follows that
\beq
\|\sum_{j=1}^{N(c)}b^{1/2}Y_j^*\Phi(c)Y_j b^{1/2}-b\|<\ep.
\eneq
This proves the claim.  Then $\pi\circ \Phi$ is full in $Q_1.$
Applying Theorem 3.12 of \cite{eglnkk0}, one obtains an integer $K_1\ge 1$ (in place of $n$), and a unitary
$u\in {{(M_{{K_1+1}}(Q_1))}}^\sim$such that
\beq
\|u^*\diag({\bar \phi}_1(c), \Psi(c))u-\diag({\bar \phi}_2(c), \Psi(c))\|<\ep_0/3\rforal c\in {\cal F}_0,
\eneq
where $\Psi(c)=\diag(\overbrace{\pi\circ \Phi(c), \pi\circ \Phi(c),...,\pi\circ \Phi(c)}^{K_1}).$
We may assume that there is a sequence of unitaries $u_n\in {{(M_{{K_1+1}}(B_n))}}^\sim$
such that $\pi(\{u_n\})=u.$
Thus, there exists $N_0\ge 1$  {{such}} that, for all $n\ge N_0,$ for all $c\in {\cal F}_0,$
\beq
\|u_n^*\diag(L_n^{(1)}(c), {\bar \Phi}_n(c))u_n-\diag(L_n^{(2)}(c), {\bar \Phi}_n(c))\|<\ep/2,
\eneq
where   ${\bar \Phi}_n(c)=\diag(\overbrace{\Phi_n(c),\Phi_n(c), ...,\Phi_n(c)}^{K_1}).$
Note that, for a fixed $N\ge N_0,$
we may write ${\bar \Phi}_N(c)=j_{c, A_N}(c)\otimes 1_{K}$ for $K=r(N) K_1.$
This {{contradicts}} \eqref {2020-721-0}.
\end{proof}

\begin{lem}\label{stablehomtp}
Let $C$ be {{a non-unital}}  amenable separable \CA\,  {{which satisfies}} the UCT.
Suppose that there is an embedding $h_{c,w}: C\to {{{\cal W}}}.$
 For any $\ep>0,$ any finite subset ${\cal F}\subset C,$ any finite subset ${\cal P}\subset \underline{K}(C),$ any  \hm\, $h: C\to A,$ where $A$ is any  \CA\, with the property that
 there is a locally uniformly full {{(see 5.5 of \cite{eglnp})}} embedding  ${{j_{w,a}}}:  {\cal W}\to A$ which maps strictly positive elements to
 strictly positive {{elements,}} and any ${\kappa}\in Hom_{\Lambda}(\underline{K}({SC}), \underline{K}(A)),$ there exists an integer $N\ge 1,$ a full \hm\,
$h_0: C\to M_N({{\cal W}})\subset M_N(A)$  and  a unitary 
{{$u\in U_0(M_{N+1}({\tilde A}))$}}
such that
\beq\label{sthomp-1}
\|H(c),\, u]\|<\ep\tforal c\in {\cal F}\andeqn {\rm Bott}(H,\, u)|_{\cal P}=\kappa|_{\cal P},
\eneq
where $H(c)={\rm daig}(h(c), h_0(c))$ for all $c\in C.$
{{If $h_{c,w}^\dag\circ J_{cu}^C=0$ (for a choice of $J_{cu}^C$), we may choose $h_0$ so that $h_0^\dag\circ J_{cu}^C=0$
(see  Definition \ref{DkappaJ}).}}

\end{lem}

\begin{proof}
Define $S=\{z, 1_{C(\T)}\},$ where $z$ is  the identity function on the unit circle. Define $x\in {\rm Hom}_{\Lambda}(\underline{K}(C^\T), \underline{K}(A))$ as follows (see (\ref{june-19-2021}) in \ref{Rkt2}):
\beq\label{sthomp-2}
x|_{\underline{K}(C)}=[h]\andeqn x|_{\boldsymbol{\bt}(\underline{K}(C))}=\kappa,\andeqn  x|_{{\boldsymbol{\bt}}(\underline{K}(\C\cdot1_{\tilde C}))}=0.
\eneq
Fix a finite subset ${\cal P}_1\subset \boldsymbol{\bt}(\underline{K}(C)).$
Choose $\ep_1>0$ and a finite subset ${\cal F}_1\subset C$ satisfying the following:
\beq\label{sthomp-3}
[L']|_{{\cal P}_1}=[L'']|_{{\cal P}_1}
\eneq
for any pair of  ${\cal F}_2$-$\ep_1$-multiplicative \morp s $L',L'':C^\T\to B$ (for any  \CA\, $B$), provided  that
\beq\label{sthomtp-4}
L'\approx_{\ep_1} L''\,\,\,{\rm on} \,\, {\cal F}_2,
\eneq
where 
${\cal F}_2=S\otimes {\cal F}_1\cup{{\{(1-z)\otimes 1_{\tilde C}\}}}.$

Let  {{$\ep>0$ and  finite subsets ${\cal F}$  and ${\cal P}\subset \underline{K}(C)$}} be given.
We may assume, without loss of generality,  that
\beq\label{sthomtp-4+1}
{\rm Bott}(H',\, u')|_{\cal P}={{{\rm Bott}(H',\, u'')|_{\cal P},}}
\eneq
provided $\|u'-u''\|<\ep$ for any unital \hm\, {{$H'$}} from $C$ {{(by choosing ${\cal P}_1=\boldsymbol{\bt}({\cal P})$ in (\ref{sthomp-3})).}}
Put
$\ep_2=\min\{\ep/2, \ep_1/2\}$ and ${\cal F}_3={\cal F}\cup {\cal F}_1$.

Let $j_{c,a}=j_{w,a}\circ h_{c,w}.$
Note that, by 5.6 of \cite{eglnp}, $h_{c,w}$ is strongly locally uniformly full.
One then checks that $j_{c,a}$ is locally uniformly full.
Let $M\times N: C_+\setminus \{0\}\to \R\setminus \{0\}\times \N$ be a map
so that $j_{c,a}(c)$ is {{$(M(c), N(c))$-full}} for all $c\in C_+\setminus \{0\}.$

Let $\dt>0,$ ${\cal G}\subset C$ be a finite {{subset}} and ${\cal P}_0\subset
\underline{K}(C)$ (in place of ${\cal P}$) be as required by {{\ref{TSuniq}}}
for $\ep_2/2$ (in place of $\ep$) and {{${\cal F}_3$}} 
 (in place of ${\cal F}$).
Without loss of generality, we may assume that {{${\cal F}_3$}} 
 and ${\cal G}$ are in the unit ball of $C$ and $\dt<\min\{1/2, \ep_2/16\}.$
Fix another finite subset ${\cal P}_2\subset \underline{K}(C)$ and defined
${\cal P}_3={\cal P}_0\cup {\boldsymbol{\bt}}({\cal P}_2)$ (as a subset of
$\underline{K}(C(\T)\otimes C)$). We may assume that ${\cal P}_1\subset
{\boldsymbol{\bt}}({\cal P}_2).$

Set ${\cal G}_1=S\otimes {\cal G}\cup \{1_{C(\T)}\otimes g: g\in {{{\cal F}_3}}
\}.$
It follows from
{{Theorem 3.5 of \cite{GLrange}}} 
that there are integer $N_1\ge 1$
and a ${\cal G}_1$-$\dt/2$-multiplicative \morp\, $L: C^\T\to M_{N_1}(A)$ such that
\beq\label{sthomtp-5}
[L]|_{{\cal P}_3}=x|_{{\cal P}_3}.
\eneq
We may assume that there is a unitary $v_0\in M_{N_1}({\tilde A})$  {{(by the last identity of \eqref{sthomp-2},
we may even assume that $v_0\in U_0(M_{N_1}(\td A))$ with a larger $N_1$)}} such that
\beq\label{sthomtp-5+}
\|(L(z-1)\otimes 1_{\tilde C})+1_{\tilde A})-v_0\|<\ep_2/2.
\eneq

Define $L_1: C\to M_{N_1}(A)$ by $L_1(c)=L(c\otimes 1)$ for all $c\in C$
and $H_1: C\to M_{N_1}(A)$ by
$H_1(c)={{\diag(h(c), \overbrace{j_{w,a}(h_{c,w}(c)),j_{w,a}(h_{c,w}(c)),...,j_{w,a}(h_{c,w}(c))}^{N_1-1})}}$
for all $c\in C.$  Note that
\beq\label{sthomtp-7}
[L_1]|_{{\cal P}_0}=[h]|_{{\cal P}_0}.
\eneq
It follows from \ref{TSuniq}
that there exists an integer $N_2\ge 1,$ a  \hm\,
and a unitary {{$U\in {{(M_{(N_2+1)N_1}(A))^\sim}}$}} such that
\beq\label{sthomtp-8}
U^*(L_1(c)\oplus h_1(c))U\approx_{\ep/4} H_1(c)\oplus h_1(c)\rforal c\in {\cal F}_3,
\eneq
{{where $h_1=(j_{w,a}\circ h_{c,w})\otimes {\rm id}_{M_{N_2N_1}}:C\to {{j_{w,a}}}(M_{N_2N_1}({{\cal W}})){{\subset M_{N_2N_1}(A)}}.$}}
Put {{$N=(N_2+1)N_1-1=N_2N_1+(N_1-1).$}} 
Now define $h_0: C\to M_N({{\cal W}})$ and $H: C\to M_{N+1}(A)$ by
\beq\label{sthomtp-9}
{{h_0(c)=\diag(h_1(c),\overbrace{h_{c,w}(c),\cdots,h_{c,w}(c)}^{N_1-1})}}\andeqn H(c)=h(c)\oplus h_0(c)
\eneq
for all $c\in C.$  {{If $h_{c,w}^\dag\circ J_{cu}^C=0,$ then $h_0^\dag\circ J_{cu}^C=0.$}}
Define
$u=U^*(v_0\oplus {{1_{M_{N}}}})U.$
Then, by ({{\ref{sthomtp-8}}}) {{and}} {{the fact that $L_1$ is}} ${\cal G}_1$-$\dt/2$- multiplicative, we have
\beq\label{sthomtp-11}
&&\hspace{-0.2in}\|[H(c),\, u]\| \le \|(H(c)-{\rm Ad}\, U\circ (L_1(c)\oplus h_1(c)))u\|\\
&&+\|{\rm Ad}\, U\circ (L_1(c)\oplus h_1(c)),\, u]\|+
\|u(H(c)-{\rm Ad}\, U\circ (L_1(c)\oplus h_1(c)))\|\\
&&<\ep/4+\dt/2+\ep/4<\ep\rforal c\in {\cal F}_2.
\eneq
Define $L_2: C\to M_{N+1}(A)$ by $L_2(c)=L_1(c)\oplus h_1(c)$ for all $c\in C.$ Then, {{by \eqref{sthomtp-8},
the choice of ${\cal F}_2,$  and by \eqref{sthomtp-5}}}
we compute that
\beq\label{sthomtp-12}
{\rm Bott}(H,\, u)|_{\cal P}&=& {\rm Bott}({\rm Ad}\, U\circ L_2,\,u)|_{\cal P}=
 {\rm Bott}(L_2,\, v_0\oplus 1_{M_{N_2(N_1+1)}})|_{\cal P}\\
&=& {\rm Bott}(L_1,\, v_0)|_{\cal P}+{\rm Bott}(h_1,\, 1_{M_{N_2(N_1+1)}})|_{\cal P}\\
&=& [L\circ {\boldsymbol{\bt}}]|_{\cal P}+0=
x\circ{\boldsymbol{\bt}}|_{\cal P}=\kappa|_{\cal P}.
\eneq
%
{{The fact $u\in U_0(M_{N+1}({\tilde A}))$ follows from  $v_0\in U_0(M_{N_1}(\td A))$.}}
%
\end{proof}

\begin{lem}\label{LWhomotopy}
Let $C$ be a separable amenable \CA\, which {{admits}} an embedding {{$\phi_w: C\to {{{\cal W}}}.$}}
Then there exists an embedding $\Phi: C^\T\to {{{\cal W}}}$ {{satisfying}} the following:
there exists a continuous path of unitaries $\{V(t)\in {{{\widetilde {\cal W}}}}: t\in [0,1]\}$ such that
\beq
\|[\Phi(1_{C(\T)}\otimes c), V(t)]\|=0 {{\tforal}} c\in C, \\
V(0)=1_{\widetilde {\cal W}}+\Phi((z-1)\otimes 1_{\tilde C}){{\tand}} V(1)={{1_{\widetilde {\cal W}}}},
\eneq
where $z\in C(\T)$ is the identity function on the unit circle $\T.$   
Moreover,
\beq
{\rm length}(\{V(t)\})\le 2\pi.
\eneq
Furthermore, we can choose $\Phi|_C=\phi_w.$
\end{lem}

\begin{proof}
We refer this to the proof of {{Corollary}} 14.7 of \cite{GLII}.
We use the isomorphism ${{\cal W}}\otimes {\cal U}\cong {{\cal W}}.$
In that proof, we choose {{(recall $\|h\|=1$-- see the proof of {{Corollary}} 14.7 of \cite{GLII})}}
$$
V(t)={{1_{{\widetilde {\cal W}}\otimes {\cal U}}}}+\sum_{n=1}^{\infty} {{(i(1-t)e\otimes h)}^n\over{n!}}=
\exp(i(1-t)e\otimes h)\rforal t\in [0,1].
$$
One easily checks that this $\{V(t): t\in [0,1]\}$ meets the requirements of this lemma
with $\Phi$ {{defined}} in the proof of {{Corollary}} 14.7 of \cite{GLII}.
Note  also  that, as in that proof, $\Phi|_C=\phi_w.$
\end{proof}

 \begin{thm}\label{STHOM}
Let $C$ be a  separable amenable \CA\,  which   satisfies the UCT.
Suppose that there is {{an
embedding}} $h_w: C\to {{\cal W}}.$ For any $\ep>0$ and any finite subset ${\cal
F}\subset C,$ there exists $\dt>0,$ {{finite subsets}} ${\cal G}\subset
C$ {{and}} ${\cal P}\subset \underline{K}(C)$
 satisfying the following:

Suppose that $A$ is a separable \CA\, with the property that
there is a  locally uniformly full embedding  (see 5.5 of \cite{eglnp}  and Definition  3.11 of \cite{eglnkk0}) $j_w: {{\cal W}}\to A$ which maps strictly positive elements to
strictly positive elements.   Suppose that  $h: C\to A$ is a
\hm\, and suppose that $u\in U_0({{\td A}})$ is a unitary such that
\beq\label{Shomp1}
\|[h(a), u]\|<\dt\tforal a\in {\cal G}\tand
{\rm{Bott}}(h,u)|_{\cal P}=0.
\eneq
Then there exists an integer $N\ge 1$, {{a \hm\, $H_0: C^\T\to j_w(M_N({{\cal W}}))$ ($\subset M_N(A)$)}} and a continuous path of
unitaries $\{U(t): t\in [0,1]\}$ in $M_{N+1}({\tilde A})$ such that
\beq\label{Shomp2}
U(0)=u',\,\,\, U(1)=1_{M_{N+1}({\tilde A})}\andeqn \|[h'(a),
U(t)]\|<\ep\tforal a\in {\cal F},
\eneq
where
$$
u'={\rm diag}(u, H_0((z-1)\otimes 1)+1_{M_N({\tilde A})})
$$
{{$h'(c)=h(c)\oplus H_0(c\otimes 1)$ for $c\in C,$}}
and $z\in C(\T)$ is the identity function on the unit circle.
%
Moreover,
\beq\label{Shomp2+}
{{{\rm{Length}}}}(\{U(t)\})\le  {{2 \pi}}+\ep.
\eneq

\end{thm}

\begin{proof}
Let $\ep>0$ and ${\cal F}\subset C$  be given. Without loss of
generality, we may assume that ${\cal F}$ is in the unit ball of
$C.$ {{Let $S=\{z, 1_{C(\T)}\}\subset C(\T),$ where $z$ is  the identity function on the unit circle. For any finite subset ${\cal F}\subset C$, let ${\cal F}_S=S\otimes {\cal F}\cup \{(1-z)\otimes 1_{\tilde C}\}$.}}

Let $\dt_1>0,$  ${\cal G}_1\subset
C^\T,$ ${\cal
P}_1\subset \underline{K}(C^\T)$
 be required by {{Theorem \ref{TSuniq}}}
for $\ep/4$ and ${\cal F}_S.$

Without loss of
generality, we may assume that 
{{${\cal G}_1=({\cal G}'_1)_S=S\otimes {\cal G}_1'\cup \{(1-z)\otimes 1_{\tilde C}\},$ for a finite set
 ${\cal G}_1'$  in the unit ball of $C$.}}
Moreover, without loss of
generality, we may assume that ${\cal P}_1={\cal P}_2\cup{\cal
P}_3,$ where ${\cal P}_2\subset \underline{K}(C)$ and ${\cal
P}_3\subset {\boldsymbol{\bt}}(\underline{K}({\tilde C})).$ Let ${\cal P}= {\cal P}_2 \cup ({\boldsymbol{\beta}}^{-1}({\cal P}_3)\cap \underline{K}(C)) \subset \underline{K}(C)$. Furthermore, we
may assume that any ${\cal G}_1$-$\dt_1$-multiplicative \morp\,
$L'$ from $C^\T$
 to any \CA\, well defines $[L']|_{{\cal P}_1}.$

Let $\dt_2>0$ and ${\cal G}_2\subset C$ be a finite subset
required by 2.8 of \cite
{LnHomtp-c} for $\dt_1/2$ and ${\cal G}_1'$
above.

Let $\dt=\min\{\dt_2/2, \dt_1/2, \ep/2\}$ and ${\cal G}={\cal
F}\cup {\cal G}_2.$

Suppose that $h$ and $u$ satisfy the assumption with above $\dt,$
${\cal G}$ and ${\cal P}.$  {{\Wlog, we may assume $u=1_{\td A}+a$
for some $a\in  A.$}}

Thus, by 2.8 of \cite
{LnHomtp-c},
 there is
${\cal G}_1$-$\dt_1/2$-multiplicative \morp\, $L: {\tilde C}\otimes
C(\T)\to {\tilde A}$ such that
\beq\label{Shomp4}
&&\|L(f\otimes 1)-h(f)\|<\dt_1/2\rforal f\in {\cal G}_1'\andeqn\\
&& \|(L(1\otimes(z-1))+1)-u\|<\dt_1/2.
\eneq

Define $y\in Hom_{\Lambda}(\underline{K}(C^\T),
\underline{K}(A))$ as follows:
$$
y|_{\underline{K}(C)}=[h]|_{\underline{K}(C)}\andeqn
{{y|_{\boldsymbol{\bt}(\underline{K}(\td C))}=0.}}
$$
{It follows from ${\rm{Bott}}(h,u)|_{\cal P}=0$ that $[L]|_{\bt({\cal P})}=0$.}
Then {{(also using the fact $u\in U_0(\td A)$)}}
\beq\label{Shomp4+1}
[L]|_{{\cal P}_{ 1}}=y|_{{\cal P}_{ 1}}.
\eneq

Define $H: C^\T\to A$ by $H(1_{C(\T)}\otimes c)={{h(c)}}$ for all $c\in C,$
$$
H(g\otimes c)=h(c)\cdot g(1)\cdot 1_A\,\,\,{{{\rm for}\,\, g\in C(\T)}} \andeqn  H(f\otimes 1_{{\tilde C}})=f(1)=0
\,\,\,{{{\rm for}\,\, f\in C_0(\T\setminus \{1\}),}}
$$
 where $\T$ is identified with
the unit circle (and $1\in \T$).
It follows that
\beq\label{Shomp3}
[H]|_{{\cal P}_{ 1}}= y|_{{\cal P}_{ 1}}=[L]|_{{\cal P}_{ 1}}.
\eneq

By {{\ref{LWhomotopy},}} there exists
an embedding  $\phi_w: C^\T\to {{\cal W}}$  such that
there exists a continuous path of unitaries $\{V(t)\in {\widetilde {{\cal W}}}: t\in [0,1]\}$ such that
\beq
\|[\phi_w(1_{C(\T)}\otimes c), V(t)]\|=0\rforal c\in C, \\
V(0)=1_{\widetilde {{\cal W}}}+\phi_w((z-1)\otimes 1_{\tilde C})\andeqn V(1)=1_{\widetilde {{\cal W}}},
\eneq
where $z\in C(\T)$ is the usual unitary generator.
Moreover,
\beq
{{{\rm Length}(\{V(t)\})}}\le 2\pi.
\eneq
{{Then,}}  as in the proof
of \ref{stablehomtp}, we may assume that $j_{c,a}:=j_w\circ \phi_w: C^\T\to A$ is
locally uniformly full associated with
the map
$M\times N: {{(C^\T)_+\setminus}} \{0\}\to \R\setminus \{0\}\times \N.$

It follows  {{\ref{TSuniq}}}
that there is an integer $N\ge 1$
 and a unitary $Z\in U(M_{1+N}({\tilde A}))$ such that {{for $H_0: C^\T\to A$ defined by $
H_0(x)=\diag(\overbrace{j_w\circ \phi_w(x),...,j_w\circ \phi_w( x)}^N)\rforal x\in C^\T$, we have}}
\beq\label{Shomp5}
Z^*(H(c)\oplus H_{0}(c))Z\approx_{\ep/4} L(c)\oplus H_{0}(c)\rforal c\in {\cal F}_S.
\eneq
Define
$
V'(t)=\diag(\overbrace{j_w(V(t)), j_w(V(t)),...,j_w(V(t))}^N\}.
$
Then
\beq\label{Shomp6}
&&V'(0)=H_0((z-1)\otimes {{1_{\tilde C}}})+1_{M_N({\tilde A})},\,\,\,
V'(1)=1_{M_{N}({\tilde A})}\andeqn\\
&&H_0(1_{C(\T)}\otimes c)V'(t)=V'(t)H_0(1_{C(\T)}\otimes c)
\eneq
for all $c\in C$ and $t\in [0,1].$  
Moreover,
\beq\label{Shomp6+}
\text{Length}(\{V'(t)\})\le 2\pi.
\eneq

Now define $U(1/4+3t/4)=Z^*{\rm diag}(1, V'(t))Z$ for $t\in [0,1]$
and
$$
u'=u\oplus   V'(0)
\andeqn h'(c)=h(c)\oplus H_0(1_{C(\T)}\otimes c)
$$
for $c\in C$ for $t\in [0,1].$  Then
\beq\label{Shomp7}
\|u'-U(1/4)\|<\ep/4\andeqn \|[U(t),\, h'(a)]\|<\ep/4
\eneq
for all $a\in {\cal F}$ and $t\in [1/4,1].$ The theorem follows by
connecting $U(1/4)$ {{and}} $u'$ with a short path as follows: There
is a self-adjoint element $a\in M_{1+N}({\tilde A})$ with $\|a\|\le {\ep
\pi\over{8}}$ such that
\beq\label{Shomp-8}
\exp(i a)=u'{{U(1/4)^*.}}
\eneq
Then the path of unitaries $U(t)=\exp(i (1-4t) a)U(1/4)$ for $t\in [0,1/4)$ satisfy the above.
\end{proof}

\begin{df}[cf. Definition 3.4 of \cite{Lininv}]\label{Mappingtorus}
Let $C$ and $B$ be   \CA s {{and}} ${{\phi,\psi}}: C\to B$ be two monomorphisms.
Suppose that $[\phi]=[\psi]$ {{is}} in $KL(C,B).$ Then $M_{\phi, \psi}$  (see \eqref{Mtoruses}) corresponds {{to}}
{{the zero}} element of $KL({{C}},B).$ In particular, the corresponding extensions
$$
0\to K_i(SB) \stackrel{\imath_*}{\to} K_i(M_{\phi, \psi})\stackrel{\pi_e}{\to}K_i( C)\to 0\,\,\,\,\,\,\,{\rm (}i=0,1{\rm )}
$$
are pure {(see  Lemma 4.3 of \cite{LnAUCT}).}
{Suppose that $T(B)\not=\emptyset.$ Let $u\in M_l(M_{\phi, \psi}^\sim)$ (for some integer $l\ge 1$) be a unitary
which is a {{piecewise smooth continuous}} function on $[0,1].$
Put
$$
D_B(\{u(t)\})(\tau)={1\over{2\pi} i}\int_0^1 \tau({du(t)\over{dt}}u^*(t))dt\, \tforal \tau\in T(B).
$$
(see \ref{DTtilde}
for the extension of $\tau$ to $M_l(B)$).
}
{Suppose that $\tau\circ \phi=\tau\circ \psi$ for all $\tau\in T(B).$
Then there exists a \hm\,
$$
R_{\phi, \psi}: K_1(M_{\phi, \psi})\to \Aff(T(B)),
$$
defined by $R_{\phi, \psi}([u])(\tau)=D_B(\{u(t)\}){{(\tau)}}$ as above, which is independent  of the choice of the piecewise smooth
path $u$ in $[u].$
We have the following commutative diagram:}
$$
\begin{array}{ccc}
K_0(B) & \stackrel{\imath_*}{\longrightarrow} & K_1(M_{\phi, \psi})\\
 \rho_B\searrow && \swarrow R_{\phi, \psi}\\
  & \Aff(T(B))
  \end{array}.
  $$
{Suppose, in addition,  that $[{{\phi}}]=[{{\psi}}]$ {{is}}  in $KK(C,B).$  Then the following exact sequence splits:
\begin{equation}\label{Aug2-2}
0\to \underline{K}(SB)\to \underline{K}(M_{{\phi, \psi}})\,{\overset{[\pi_e]}{\underset{\theta}{\rightleftharpoons}}} \,\underline{K}(C)\to 0.
\end{equation}
We may assume that $[\pi_0]\circ [\theta]=[{{\phi}}]$ and
$[\pi_1]\circ [\theta]=[{{\psi}}].$
In particular, one may write $K_1(M_{\phi, \psi})=K_0(B)\oplus K_1(C).$ Then we obtain a \hm\,
$$
R_{\phi, \psi}\circ \theta|_{K_1(C)}: K_1(C)\to \Aff(T(B)).
$$
We shall say ``the rotation map vanishes" if there exists  a splitting  map
$\theta,$ {{as above, such}} that $R_{\phi, \psi}\circ \theta|_{K_1(C)}=0.$}

{Denote by ${\cal R}_0$ the set of those elements $\lambda\in {\rm Hom}(K_1(C), \Aff(T(B)))$ for which there is a \hm\,
$h: K_1(C)\to K_0(B)$ such that $\lambda=\rho_B\circ h.$ It is a subgroup of ${\rm Hom}(K_1(C), \Aff(T(B))).$
If $[\phi]=[\psi]$  in $KK(C,B)$ and $\tau\circ \phi=\tau\circ \psi$
for all $\tau\in T(B),$ one has a well-defined element ${\overline{R_{\phi, \psi}}}\in {\rm Hom}(K_1(C), \Aff(T(B)))/{\cal R}_0$
(which is independent of the choice of $\theta$).}

{Under the assumptions that $[\phi] = [\psi]$ in $KK(C, B)$,  $\tau\circ \phi=\tau\circ \psi$ for all $\tau\in T(B)$,  and $C$ satisfies the UCT},  there exists a \hm\, $\theta'_1: K_1(C)\to K_1(M_{\phi, \psi})$ such that
$(\pi_e)_{*1}\circ \theta_1'={\rm id}_{K_1(C)}$ and $R_{\phi, \psi}\circ \theta'_1\in {\cal R}_0$ if, and only
if, there is $\Theta\in {\rm Hom}_{\Lambda}(\underline{K}(C),\underline{K}(M_{\phi, \psi}))$ such that
$$
[\pi_e]\circ \Theta=[{\rm id}_{C}]\,\,\,{\rm in}\,\,\,KK(C,C)\andeqn R_{\phi, \psi}\circ \Theta|_{K_1(C)}=0
$$
{{(see Definition 3.4 of \cite{Lininv}).}}
In other words, $\overline{R_{\phi, \psi}}=0$ if, and only if, there is $\Theta$ as described above
such that $R_{\phi, \psi}\circ \Theta|_{K_1(C)}=0.$  When $\overline{R_{\phi, \psi}}=0,$ one has that
$\theta(K_1(C))\subset {\rm ker}R_{\phi, \psi}$ for some $\theta$ such that (\ref{Aug2-2}) holds. In this case
$\theta$ also gives the following {{decomposition:}}
$$
{\rm ker}R_{\phi, \psi}={\rm ker}\rho_B\oplus K_1(C).
$$
We will consider the case $B$ is {{non-unital}} but ${\rm Ped}(B)=B$
(see 7.11 of \cite{Lncbms} for more details).
\end{df}

\begin{lem}\label{VuV}
Let $C$ be a separable  non-unital {{algebraically}} simple amenable \CA\, which {{admits}} a \hm\,
$\phi_w: C\to {{\cal W}}$  and {{a \hm\,}}
$\phi_{w,c}: {{\cal W}}\to C$ which map strictly positive elements to
strictly positive elements,  and $B$ be a ${\cal Z}$-stable simple separable projectionless  \CA\, with $T(B)\not=\emptyset$ and with an embedding
$j_w: {{\cal W}}\to B.$ Suppose that
both $C$ and $B$ have stable rank one and {{continuous scale.}}

Suppose  {{that}} $\phi_1, \, \phi_2:C\to B$ are two  monomorphisms
such that
\beq\label{2020-721-VuV-1}
[\phi_1]=[\phi_2]\,\,\,{\rm in} \,\,\, KK(C,B)\tand\,
\tau\circ \phi_1=\tau\circ \phi_2.
\eneq
Let $\theta: \underline{K}(C)\to \underline{K}(M_{\phi_1, \phi_2})$ be the splitting map defined {{in \ref{Mappingtorus}.}}

For any $1/2>\ep>0,$ any finite subset ${\cal F}\subset C$ and any
finite subset ${\cal P}\subset \underline{K}(C),$ there are
{{an}} integer $N_1\ge 1,$  an ${\cal F}$-$\ep/2$-multiplicative \morp\,
$L: C\to M_{1+N_1}(M_{\phi_1,\phi_2}),$  a  \hm\,
$h_0: C\to {{\cal W}}\to M_{N_1}(C)$
and a continuous path of unitaries $\{V(t): t\in [0,1-d]\}$ of
$M_{1+N_1}({{\td B}})$ for some $1/2>d>0,$ such that $[L]|_{\cal P}$ is
well defined, $V(0)=1_{M_{1+N_1}({{\td B}})},$
\beq\label{VuV1}
[L]|_{\cal P}&=&\theta|_{\cal P},\\
\label{VuV2}
\pi_t\circ L&\approx_{\ep}&{\rm Ad}\, V(t)\circ
(\phi_1\oplus \phi_1\circ h_0)\,\,\,\text{on}\,\,\,{\cal F} \tforal t\in (0,1-d],\\
\label{VuV3}
\pi_t\circ L&\approx_{\ep}&{\rm Ad}\, V(1-d)\circ
(\phi_1\oplus \phi_1\circ h_0)\,\,\,\text{on}\,\,\,{\cal F} \tforal t\in (1-d, 1]\andeqn\\
\label{VuV33}
\pi_1\circ L&\approx_{\ep}&\phi_2\oplus \phi_2\circ h_0\,\,\,\text{on}\,\,\,{\cal F},
\eneq
where $\pi_t: M_{\phi_1,\phi_2}\to B$ is the {{point evaluation}} at
$t\in (0,1),$ and $h_0$ factors through ${{\cal W}}.$


\end{lem}

\begin{proof}
Let $\ep>0$ and {{${\cal F}\subset C$}}  a finite subset.
Let $\dt_1>0,$ ${\cal G}_1\subset C$ be a finite subset and ${\cal
P}\subset \underline{K}(C)$ be a finite subset required by
\ref{STHOM} for $\ep/4$ and ${\cal F}$ above.
We may further assume that $\dt_1$ is sufficiently small such that
\beq\label{1510-n1}
{\rm Bott}(\Phi, U_1U_2U_3)|_{\cal P}=\sum_{i=1}^3{\rm Bott}(\Phi, U_i)|_{\cal P},
\eneq
provided that
$\|[\Phi(a), U_i]\|<\dt_1\rforal a\in {\cal G}_1,\,\,\,i=1,2,3.$

Let $\ep_1=\min\{\dt_1/2, \ep/16\}$ and ${\cal F}_1={\cal F}\cup
{\cal G}_1.$ We may assume that ${\cal F}_1$ is in the unit ball
of $C.$ We may also assume that $[L']|_{\cal P}$ is well defined
for any ${\cal F}_1$-$\ep_1$-multiplicative \morp\, from $C$ to
any  \CA.

Let $\dt_2>0$ and ${\cal G}\subset C$ be a finite subset and
${\cal P}_1\subset \underline{K}(C)$ be finite subset required by {{\ref{TSuniq}}}
for $\ep_1/2$ and ${\cal F}_1.$ We may assume that
$\dt_2<\ep_1/2,$ ${\cal G}\supset {\cal F}_1$ and ${\cal
P}_1\supset {\cal P}.$ We also assume that ${\cal G}$ is in the
unit ball of $C.$

It follows from {{Corollary}}  \ref{CExistence}
that there exists an integer $K_1\ge
1$  and
a ${\cal G}$-$\dt_2/2$-multiplicative \morp\, $L_1:  C\to
{{M_{K_1}(M_{\phi_1, \phi_2})}}$ such that
\beq\label{VuV-1}
[L_1]|_{{\cal P}_1}=\theta|_{{\cal P}_1}.
\eneq

%
Recall that $\pi_e: M_{\phi_1,\phi_2} \to C$ is the canonical projection.
Note that  $[\pi_e]\circ \theta=[{\rm id}|_C],$  $[\pi_0]\circ \theta=[\phi_1]$ and $[\pi_1]\circ \theta =[\phi_2]$ and, for each $t\in (0,1),$
\beq\label{VuV-4}
[\pi_t]\circ \theta=[\phi_1]{{=[\phi_2]}}.
\eneq
Note that $\phi_{w,c}\circ \phi_w: {{C}} \to C$ factors through ${{\cal W}}.$  Let $h''_0: C\to M_{K_1-1}(C)$ be
a \hm\, factoring through ${\cal W}.$
Applying \ref{TSuniq} to the pair $\pi_e\circ L_1, \id\oplus h''_0: C \to M_{K_1}(C)$,
 we obtain an integer $K_0,$ unitaries   {{$V\in M_{K_1+K_0}(C)^\sim$}}
and a \hm\, $h_0': C\to {{\cal W}}\to
M_{K_0}(C)$ such that (also recall  5.6 of \cite{eglnp}---note that $M_{K_0}(C)$ is algebraically simple)
\beq\label{VuV-5--1}
{{{\rm Ad}\, V \circ (\pi_e\circ  L_1\oplus h_0')\approx_{\ep_1/4}
(\id \oplus h''_0\oplus h_0')\,\,\,\text{on}\,\,\,{\cal F}_1,}}
\eneq
{{where $h_0'=(\phi_{w,c}\circ \phi_w)\otimes {\rm id}_{M_{K_0}}.$}}
%
Note that $\tau\circ \phi_1\circ h_0'=\tau\circ \phi_2\circ h_0'$ for all $\tau\in T(B)$ and
$h_0'$ factors through ${\cal W},$
where we continue to {{write}} $\phi_i$ for $\phi_i\otimes {\rm id}_{M_{K_0}}$ (and
in what follows, we will continue to use this practice).
By   (the {{second}}
paragraph of) \ref{ConstphiW} (see  around \eqref{2020-721-12-1}), we may assume that
there are  an  ${\cal F}_1$-${\ep_1\over{16}}$-multiplicative \cpc\, $L_{00}: C\to M_{K_0}(M_{\phi_1, \phi_2})$ and
 a {{unitary path}} $U_1\in C([0,1], M_{K_0}(B))$ with
$U_1(0)=1_{M_{K_0}}$
such that
\beq
&&{\rm Ad}\, {{U_1(1)}}\circ {{\phi_1}}\circ h_0'\approx_{\ep_1/4} {{\phi_2}}\circ h_0'\,\,\,\text{on}\,\,\,{\cal F}_1,\\ \label{june-30-2021}
&&\pi_0\circ L_{00}=\phi_1\circ h_0', \,\pi_1\circ L_{00}=\phi_2\circ h_0'\andeqn\\
&&\pi_t\circ L_{00}\approx_{\ep_1/16} {{\rm Ad}}\, U_1(t)\circ {{\phi_1}}\circ h_0'\,\,\,\text{on}\,\,\,{\cal F}_1,\,\,\, \rforal t\in (0,1).
\eneq


Write $V_{00}=\phi_1(V)$ and
$V_{00}'=\phi_2(V).$ The assumption that
$[\phi_1]=[\phi_2]$ implies that $[V_{00}]=[V_{00}']$ in $K_1(B).$
Note that, in fact,
$V_{00}$ and $V_{00}'$ are in the same component of
{{$U(M_{K_1+K_0}(B^\sim)).$}}

One obtains a continuous path of unitaries $\{Z(t): t\in [0,1]\}$
in ${{M_{K_1+K_0}(B^\sim)}}$ such that
\beq\label{VuV-5--2}
Z(0)=V_{00}
\andeqn Z(1)={{V_{00}'}}. 
\eneq
It follows that $Z\in M_{K_1+K_0}(M_{\phi_1,\phi_2})^\sim.$ By
replacing $L_1$ by ${{{\rm Ad}\, Z\circ (L_1\oplus   L_{00})}}$ {{and $h''_0\oplus h'_0$ by $h'_0$}}
{{it follows from (\ref{VuV-5--1}) and (\ref{june-30-2021}) that}}
\beq\label{VuV-5--3}
\pi_0\circ L_1\approx_{\ep_1/2}\phi_1\oplus
\phi_1\circ h_0'\,\,\,\text{on}\,\,\,{\cal F}_1\andeqn
\pi_1\circ L_1\approx_{\ep_1/2} \phi_2\oplus \phi_2\circ h_0'
,\,\,\text{on}\,\,\, {\cal F}_1.
\eneq


There is a partition
$0=t_0<t_1<\cdots <t_n=1$
such that
\beq\label{VuV-6}
\pi_{t_{i}}\circ L_1\approx_{\dt_2/8}\pi_t\circ
L_1\,\,\,\text{on}\,\,\,{\cal G} \rforal t_i\le t\le t_{i+1},\,\,\,i=1,2,...,n-1.
\eneq
Applying  Theorem \ref{TSuniq}
again,
and by choosing even larger $K_0,$
 we may assume there is
 a unitary ${{V_{t_i}\in M_{K_1+K_0}(B^\sim)}}$ 
 such that
\beq\label{VuV-7}
{\rm Ad}\, V_{t_i}\circ (\phi_1\oplus \phi_1\circ h_0')
\approx_{\ep_1/2} (\pi_{t_i}\circ L_1)
\,\,\,\text{on}\,\,\, {\cal F}_1.
\eneq

Note that, by (\ref{VuV-6}), (\ref{VuV-7}) and (\ref{VuV-5--3}),
$$
\|[(\phi_1\oplus \phi_1\circ h_0')(a), V_{t_i}V_{t_{i+1}}^*]\|<\dt_2/4+\ep_1\rforal
a\in {\cal F}_1.
$$

Denote by  $\eta_{-1}=0$ and
$$
\eta_k=\sum_{i=0}^k {\rm Bott}(\phi_1\oplus \phi_1\circ h_0',
V_{t_i}V_{t_{i+1}}^*)|_{\cal P},\,\,k=0,1,...,n-1.
$$

Now we will construct, for each $k\leq n-1$ 
a \hm\, ${{F_k:}} 
~C\to M_{{J_{k}}}({{\cal W}})\subset M_{{J_{k}}}(C)$ and
a unitary {{$W_k\in {M_{K_1+K_0+\sum_{i=1}^kJ_i}(B^{\sim})}$}} 
such that
\beq\label{150110-L2}
\|[{{H_k}}(a),\, {{W_k}}]\|<\dt_2/4\tforal a\in {\cal F}_1\andeqn {\rm Bott}({{H_k}},\, {{W_k}})={{\eta_{k-1}}},
\eneq
where ${{H_k}}=\phi_{{1}}\oplus \phi_1\circ h_0'\oplus \phi_1\circ ({{\oplus_{i=1}^k F_i}}),$ ${{k}}=1,2,...,n-1.$

Let ${{W}}_0=1_{M_{{K_1+K_0}}}.$ 
It follows from \ref{stablehomtp} that there is an integer $J_1\ge 1,$
a  \hm\, $F_1: C\to M_{J_1}({{\cal W}})\subset M_{J_1}(C)$ and a unitary ${{W}}_1\in
{{U_0}}({{M_{K_1+K_0+J_1}(B^\sim)}})$
such that
\beq\label{VuV-8}
\|[H_1(a), \, {{W}}_1]\|<\dt_2/4\rforal a\in {\cal F}_1\andeqn
\text{Bott}(H_1, {{W}}_1)=\eta_0,
\eneq
where $H_1=\phi_1\oplus \phi_1\circ h_0'\oplus  \phi_1\circ  F_1$.

Assume that, we have construct {{the}} required $F_i$ and $U_i$ for $i={{1,2,}}...,k<n-1.$
It follows from \ref{stablehomtp} that there is an integer $J_{k+1}\ge 1,$
a unital \hm\, $F_{k+1}: C\to M_{J_{k+1}}({{\cal W}})\subset M_{J_{k+1}}(C)$ and a unitary ${{W}}_{k+1}\in
{{U_0}}(M_{{{K_1}}+K_0+\sum_{i=1}^{k+1}J_i}({{B^\sim}}))$
such that
\beq\label{VuV-9}
\|[H_{k+1}(a), \, {{W}}_{k+1}]\|<\dt_2/{ 4}\rforal a\in {\cal F}_1\andeqn
\text{Bott}(H_{k+1}, {{W}}_{k+1})={{\eta_{k},}}
\eneq
where {{$H_{k+1}=\phi_1\oplus \phi_1\circ h_0'\oplus \phi_1\circ ( \oplus_{i=1}^{k+1}F_i)$.}}

Now define $F_{00}= \oplus^{n-1}_{{i=1}}F_i$
{{and}}
$K={{K_1}}+K_0+\sum_{i=1}^{n-1}J_i.$
Define
$$
v_{t_k}= {\rm diag}({{W}}_k{\rm diag}(V_{t_k}, {\rm id}_{1_{M_{\sum_{i=1}^kJ_i}}}), 1_{M_{\sum_{i=k+1}^{n-1} J_i}}),
$$
$k={{1,2,}}...,n-1$ and
$v_{t_0}=1_{M_{{{K_1}}+K_0+
\sum_{i=1}^{n-1}J_i}}.$
Then, for $i=0,1,2,...,n-2,$
\beq\label{VuV-14}
&&{\rm Ad}\, v_{t_i}\circ (\phi_1\oplus \phi_1\circ h_0'\oplus
\phi_1\circ F_{00})\approx_{\dt_2+\ep_1} \pi_{t_i}\circ (L_1\oplus
F_{00})\,\,\,\text{on}\,\,\,{\cal F}_1,\\
&&\|[\phi_1\oplus  \phi_1\circ h_0'\oplus \phi_1\circ F_{00}(a),\,v_{t_i}v_{t_{i+1}}^*]\|<\dt_2/2+2\ep_1\rforal a\in {\cal F}_1\andeqn\\
\label{VuV-15}
&&\text{Bott}(\phi_1\oplus \phi_1\circ h_0'\oplus \phi_1\circ F_{00}, v_{t_i}v_{t_{i+1}}^*)\\
&=&\text{Bott}(\phi_1', W_i')+
\text{Bott}(\phi_1', V_{t_i}'(V_{t_{i+1}}')^*) +\text{Bott}(\phi_1', (W_{i+1}')^*)\\
&=& \eta_{i-1}+\text{Bott}(\phi_1', V_{t_i}V_{t_{i+1}}^*)-\eta_{i}=0,
\eneq
where $\phi'_{{1}}=\phi_1\oplus \phi_1\circ h_0'\oplus \phi_1\circ F_{00},$   $W_i'={\rm diag}(W_i, 1_{M_{\sum_{j=i+1}^{n-1}J_i}})$
and $V_{t_i}'={\rm diag}(V_{t_i}, 1_{M_{\sum_{i=1}^{n-1}J_i}}).$

It follows from \ref{STHOM} that there is an integer $N_1\ge 1,$
another \hm\, $F_0': C\to M_{N_1}({{\cal W}})\subset M_{N_1}(C)$ and a continuous path of
unitaries $\{w_i(t): t\in [t_{i-1}, t_i]\}$ such that
\beq\label{VuV-17}
&&\hspace{-0.4in}
w_i(t_{i-1})=v_{i-1}'(v_{i}')^*, w_i(t_i)=1, \,\,\, i=1,2,...,n-1\andeqn\\
&&\hspace{-0.4in}\|[\phi_1\oplus \phi_1\circ h_0'\oplus \phi_1\circ F_{00}\oplus  \phi_1\circ F_{0}'(a), \,w_i(t)]\|<\ep/4\rforal a\in {\cal F},
\eneq
$i=1,2,...,n-1,$
where $v_i'={\rm diag}(v_i, 1_{M_{N_1}}(B)),$ $i=1,2,...,n-1.$
Define $V(t)=w_i(t)v_i'$ for $t\in [t_{i-1}, t_i],$ $i=1,2,...,n-1.$
Then $V(t)\in C([0,t_{n-1}], M_{{{K}}+N_1}(B)).$ Moreover,
\beq\label{VuV-18}
{\rm Ad}\, V(t)\circ (\phi_1\oplus \phi_1\circ h_0'\oplus \phi_1\circ (F_{00}\oplus
F_0'))\approx_{\ep} \pi_t\circ L_1\oplus F_{00}\oplus F_0'
\,\,\,\text{on}\,\,\,{\cal F}.
\eneq

Define  $h_0=h_0'\oplus F_{00}\oplus F_0',$ $L=L_1\oplus
F_{00}\oplus F_0'$ and $d=1-t_{n-1}.$ Then, by (\ref{VuV-18}),
(\ref{VuV2}) and (\ref{VuV3}) hold. From (\ref{VuV-5--3}),
(\ref{VuV33}) also holds. 

\end{proof}


 \section{Asymptotically unitary equivalence}

{{Much of this section is taken from  Section 7 of \cite{Lininv}  and Section 27 of \cite{GLN2}.
Some definition {{therein}} will be used (for $\td A,$ $\td B$ and $\td C$ instead of $A,$ $B$ and $C,$
and maps like $\phi^\sim$ will also be used when it is convenient).  Moreover, recall that if $\Phi: A_1\to A_2$ is a map, we will continue
to write $\Phi$ for $\Phi\otimes \id_{M_r}: M_r(A_1)\to M_r(A_2).$ We will continue to use this convention without
further notice.}}
%
%
%
%
%
%

\begin{lem}\label{Lappinductivelimit}
Let $A$ be a  \CA\,  which satisfies the UCT
and let $A_n\subset A$ be  a sequence of separable amenable $C^*$ sub algebras satisfying UCT with finitely generated
$K_i(A_n)$ ($i=0,1$) such that  $A=\overline{\cup_{n=1}^\infty A_n}$,   \beq\label{june-25-2021-1}\lim_{n\to\infty}{\rm dist} (x, A_n)=0~~\mbox{for any}~~x\in A. \eneq
Then there is a subsequence $\{m(n)\}$
and a \hm\, $j_n\in {\rm Hom}_{\Lambda}(\underline{K}(A_{m(n)}), \underline{K}(A_{m(n+1)}))$
such that ${{[\iota_{n+1}]}}\circ j_{n}=[\iota_n]$ 
(where $\iota_n: {{A_{m(n)}\to A}}$ is the embedding),
\beq
\underline{K}(A)=\lim_{n\to\infty}(\underline{K}(A_{m(n)}), j_n) \andeqn [\iota_n](\underline{K}(A_{m(n)}))\subset
[\iota_{n+1}](\underline{K}(A_{m(n+1)})).
\eneq
\end{lem}
\begin{proof}
Fix $k\ge 1.$   Since $K_i(A_n)$ is finitely generated, there is an integer $T(k)\ge 1$
such that $T(k)x=0$ for all $x\in {\rm Tor}(K_i(A_k))$ ($i=0,1$).
It follows from 2.11 of \cite{DL} (note that $A_k$ satisfies UCT) that ${\rm Hom}_{\Lambda}(\underline{K}(A_k), \underline{K}(B))
={\rm Hom}_{\Lambda}(F_{T(k)}\underline{K}(A_k), \underline{K}(B))$
for any $\sigma$-unital \CA\, $B.$
Let
$${\cal P}\subset F_{T(k)}\underline{K}(A_k)=\bigoplus_{i=0,1}(K_i(A_k)\bigoplus_{m|T(k)}K_i(A_k; \Z/m\Z))$$
 be a finite generating set.

 Since each $A_n$ is a separable amenable {{\CA,}}
by (\ref{june-25-2021-1}), there  is a sequence of \cpc s $\Phi_{n}: A\to A_n$  
such that  (see 2.1.13 of \cite{Linbook})
\beq
\lim_{n\to\infty}\|\Phi_{n}(a)-a\|=0\rforal a\in A.
\eneq
Fix $k,$ let $\iota_k: A_k\to A$ be the embedding. 
For  each $n\ge k,$ let $L_{k,n}: A_k\to A_n$ be defined by $L_{k,n}=\Phi_n\circ \iota_k$. Then
\beq\label{2020-825-1}
\lim_{n\to\infty} \|L_{k,n}(x)-x\|=0\rforal x\in A_k,
\eneq
in particular, $\lim_{n\to\infty}\|L_{k,n}(ab)-L_{k, n}(a)L_{k,n}(b)\|=0$ for all $a, b\in A_k.$
Define $L^k: A_k\to \prod_{n=1}^\infty A_n$ by $L^k(a)=\{L_{k,n}(a)\}$ for all $a\in A_k.$
Let $\pi: \prod_{n=1}^\infty A_n\to \prod_{n=1}^\infty A_n/\bigoplus_{n=1}^\infty A_n$ be the quotient map.
Define ${\bar L}^k: =\pi\circ L^k: A_k\to \prod_{n=1}^\infty A_n/\bigoplus_{n=1}^\infty A_n.$
Then ${\bar L}^k$ is a \hm.
Since $K_i(A)$ is finitely generated ($i=0,1$), by 7.2 of \cite{Lin-BDF}
(see also p.99 of \cite{GL}),
as in the last part of 2.1.15 of \cite{Lncbms}, there is an integer $N(k)$
such that, for all $n(k)\ge N(k)$ there exists $j_{n(k)}\in {\rm Hom}_{\Lambda}(\underline{K}(A_k), \underline{K}(A_{n(k)}))$ such
that
\beq
[L_{k, n(k)}]|_{\cal P}=j_{n(k)}|_{\cal P}.
\eneq
Applying the above to $N(k),$ we obtain $N(N(k))\ge N(k)$  and $j_{N(N(k))}\in
{\rm Hom}_{\Lambda}(\underline{K}(A_{N(k)}), \underline{K}(A_{N(n(k))}))$
such that
\beq
[L_{N(k), N(N(k))}]|_{{\cal Q}}=j_{N(N(k))}|_{{\cal Q}}.
\eneq
By \eqref{2020-825-1},  we may assume that,
\beq
[L_{k, NN((k))}]|_{\cal P}=[L_{N(k), N(N(k))}]\circ [L_{k, N(k)}]|_{\cal P}.
\eneq
It follows that $j_{k, N(N(k))}=j_{N(N(k))}\circ j_{N(k)}.$
Moreover, by \eqref{2020-825-1}, we may assume that,  if $x\in {\rm ker}[\iota_k],$ then, $x\in {\rm ker}j_{N(k)}.$
Thus, we obtain a subsequence $\{m(k)\}$ and $j_k\in {\rm Hom}_{\Lambda}(\underline{K}(A_{m(k)}),
\underline{K}(A_{m(k+1)})),$  and an inductive limit
$\lim_{n\to\infty}(\underline{K}(A_{m(k)}), j_k)$
such that,  for each fixed $K$ and some generating set ${\cal P}_k\subset \underline{K}(A_{m(K)}),$
\beq
[L_{m(K), m(k)}]|_{{\cal P}_k}={{j_{K, k}|_{{\cal P}_k}.}}
\eneq
Thus we have the following commutative diagram:
\begin{displaymath}
    \xymatrix{
        \underline{K}(A) \ar[r]^{\id} & \underline{K}(A) \ar[r]^{\id}
        & \underline{K}(A) \ar[r] & \cd \underline{K}(A)\\
        \underline{K}(A_{m(1)}) \ar[r]^{j_1}\ar[u]_{[\iota_1]} &
         \underline{K}(A_{m(2)})\ar[r]^{j_2}\ar[u]_{[\iota_2]}&
         \underline{K}(A_{m(3)})\ar[r] \ar[u]_{[\iota_3]}& \cd
         \lim_{n\to\infty}\underline{K}(A_{m(n)})\ar@{-->}[u]_{\imath}. }
\end{displaymath}
 {{Since $K_i(A, \Z/k\Z)=\cup_{n=1}^\infty [\iota_n](K_i(A_n, \Z/k\Z),$ $k\ge 0,$ $i=0,1,$
the}} \hm\, $\imath$ induced by the diagram
is surjective.
Also, if   $x\in {\rm ker}\, \imath,$ we may assume   that there is
$y\in \underline{K}(A_k)$ such that $j_{k, \infty}(y)=x.$
Then $[\iota_n]\circ j_{k,n}(y)=0$ for some $n\ge k.$
As mentioned above, this implies $j_{n, m}(j_{k,n}(y))=0.$   Hence $x=0.$
 It follows that $\imath$ is an isomorphism.
Moreover, $j_{n, \infty}=[\iota_n].$  Since $j_n(\underline{K}(A_n)\subset \underline{K}(A_{n+1}),$
we also have $[\iota_n](\underline{K}(A_n))\subset j_{n+1, \infty}(\underline{K}(A_{n+1}))=[\iota_{n+1}](\underline{K}(A_{n+1})).$

\end{proof}

\begin{df}\label{D642cbms}
The discussion below is taken from  6.4.2 of \cite{Lncbms} (see 10.4 and 10.5 of \cite{Lnamj}).
For any unital \CA\, $D$ and $u\in U(D),$ {{define}}
$R(u,t)\in C([0,1], M_2(D))$ by
\beq\label{DRU(t)}
R(u,t):=\begin{pmatrix} u & 0\\
                                    0 & 1\end{pmatrix}\begin{pmatrix}\cos({\pi t\over{2}}) &\sin({\pi t\over{2}})\\ -\sin({\pi t\over{2}})
                                    & \cos({\pi t\over{2}})\end{pmatrix} \begin{pmatrix} u^* & 0\\
                                    0 & 1\end{pmatrix} \begin{pmatrix}\cos({\pi t\over{2}}) &-\sin({\pi t\over{2}})\\ \sin({\pi t\over{2}})
                                    & \cos({\pi t\over{2}})\end{pmatrix}.
\eneq
Note $R(u,0)=\diag(1,1)$ and $R(u,1)=\diag(u,u^*).$
Let $C$ {{(with subalgebras $C_n$)}} be as in \ref{Lappinductivelimit} {{in place of $A$ (with subalgebras $A_n$)}}.  Let $B$ be a \CA\, and $\phi_1, \phi_2: C\to B$
be two \hm s. Let ${\cal F}\subset C$ and ${\cal P}\subset \underline{K}(C)$ be  finite subsets, and let
$\dt>0.$ Suppose that $u\in B$ {{(or $u\in {\tilde B}$ if $B$ is not unital)}} is a unitary such that
\beq
{{u^* \phi_1(a)u}}\approx_\dt\phi_2(a)\rforal  {{a\in}} {\cal F}.
\eneq
Define
\beq
L(a):=\begin{cases} R(u,2t)^*\diag(\phi_1(a),0)R(u, 2t), & t\in [0,1/2],\\
                               2(1-t){{\diag(u^*\phi_1(a)u,0)}}
                               +(1-2(1-t))\diag(\phi_2(a),0), &t\in [1/2, 1]\end{cases}.
\eneq
Put $V_1(t)=R(u,2t).$
Note that $L$ maps $C$ to $M_2(M_{\phi_1,\phi_2}).$   If $\dt$ is sufficiently small and ${\cal F}$ is sufficiently large,
there is $n\ge 1$ and $\gamma_n\in KL_{loc}^C(G^{\cal P}, \underline{K}(M_{\phi_1, \phi_2}))$
such that $[L]|_{\cal P}={\gamma_n}|_{\cal P}$ (see the end of 2.1.16 of \cite{Lncbms} {{and Definition \ref{DlocalKL} above}}),
where $G^{\cal P}$ may be viewed as a subgroup of $\underline{K}(C).$
Suppose that ${\cal P}\subset [\iota_n](F_{T(n)}\underline{K}(C_n))$ (see the proof of \ref{Lappinductivelimit}) and ${\cal P}$ generates
$[\iota_n](F_{T(n)}\underline{K}(C_n)).$ We may write $\gamma_n\in KL_{loc}^C([\iota_n](\underline{K}(C_n)), \underline{K}(M_{\phi_1, \phi_2}))$ (see \ref{DlocalKL}).
Let ${\cal Q}\subset F_{T(n)}\underline{K}(C_n)$ be a  finite generating subset such that $[\iota_n]({\cal Q})={\cal P}.$
Note that, assuming $\phi_1$ and $\phi_2$ are injective,  {{we have (see \ref{DfC1})}}
\beq
\pi_e\circ L(x) =\diag(x,0) \rforal x\in {\cal F}.
\eneq
It follows that, with sufficiently small $\ep$ and large ${\cal F},$ we may assume
that
\beq\label{2020-826-1}
[\pi_e]\circ \gamma_n\circ [\iota_n]|_{\cal Q}=[\iota_n]|_{\cal Q}, \,\,\,{\rm or}\,\,\,
[\pi_e]\circ \gamma_n\circ [\iota_n]=[\iota_n].
\eneq

{{Suppose that $\{v(t): t\in [0,1/2]\}$ is a continuous path of unitaries in $B$ {{(or ${\tilde B}$ if $B$ is not unital)}}
such that  $v(0)=1$ and $v(1/2)^*\phi_1(c)v(1/2)\approx_\dt \phi_2(c)$ for all $c\in {\cal F}.$
Put $V_2(t)=\diag(v(t),1).$
Define $L_2: C\to M_{\phi_1, \phi_2}$ by
$L_2(c)(t):=V_2(t)^*\diag(\phi_1(c), 0)V_2(t)$ for $t\in [0,1/2],$  and \\ {{$L_2(c)(t):= 2(1-t)L_2(c)(1/2)
                               +(1-2(1-t))\diag(\phi_2(c),0)$ for $t\in [1/2,1].$}}\\
                               Let $\af_n\in KL_{loc}^C(G^{\cal P}, \underline{K}(M_{\phi_1, \phi_2})$
such that $[L_2]|_{\cal P}={\af_n}|_{\cal P}.$  Note $[\pi_e]\circ [L_2]|_{\cal P}=[{\rm id}_C]|_{\cal P}.$
}}
{{For any \CA\,  $E,$ let $E^1$ be the minimal unitization. Denote $\bar V_i=V_i\otimes 1_{E^1}.$
Let ${\bar L}$ be so defined by replacing $c$ by $c\in C\otimes E^1$ and $V_1$ by $\bar V,$
and ${\bar L}_2$ by replacing $c$ by $c\in C\otimes E^1$ and $V_2(t)$ by $ V_2(t)\otimes 1_{E^1}.$
Let $z\in U(C\otimes E_1).$
Consider $\lambda(z)(t)={\bar L}(z){\bar L_2}(z^*)(t)$ and $\lambda_1(z)(t)={\bar L}(z)(t)$ and $\lambda_2(z)(t)={\bar L_2}(z)(t).$
Note $\lambda(z)\in M_2((S(B\otimes E))^\sim).$
Then (with sufficiently large ${\cal F}$ and small $\dt$) $\|\lambda(z)(t)-1\|<1/6\rforal t\in [1/2,1]$ and
$[\lambda(z)]=[\lambda_1(z)]-[\lambda_2(z)].$
Choose $E=\C,$  $C_0((0,1)),$  or a commutative \CA\, with $K_0(E)=\Z/k\Z$ $(k\ge 2$) and $K_1(E)=0.$
It follows that  (see 10.6 of  \cite{Lnamj} for the definition of $\Gamma({\rm Bott}(-,-)$, also  27.4 of {{\cite{GLN2}}} and
7.3.1 of \cite{Lncbms})}}
$[\lambda(z)]=\Gamma({\rm Bott}(\phi_1, u(1/2)V^*(1/2)))(z)$ and
\beq\label{12-2-Bott}
\Gamma({\rm Bott}(\phi_1, u(1/2)V^*(1/2)))|_{{\cal P}_n}=
[L]|_{{\cal P}_n}-[L_2]|_{{\cal P}_n}={\gamma_n}|_{{\cal P}_n}-[L_2]|_{{\cal P}_n}.
\eneq
\end{df}

 \begin{lem}\label{inv71}
Let $C_0$ be a  simple \CA\,  {{in ${\cal M}_1$}}
with continuous scale,  {{$A_1$}} a  separable simple \CA\, in {${\cal D}$} with continuous scale, and {{$U_1$}}  and $U_2$ {{two}} UHF-algebras of infinite type. Let $C=C_0\otimes U_1$ and $A=A_1\otimes U_2.$ Suppose that $\phi_1, \phi_2: C\to A$ are two  monomorphisms which
map strictly positive elements to strictly positive elements.  Suppose also that
\beq\label{71-1}
&&[\phi_1]=[\phi_2]\,\,\,{\rm in}\,\,\, KL(C, A),\\\label{71-2}
&&\phi_1^{\dag}=\phi_2^{\dag}, \,\,\, (\phi_1)_T=(\phi_2)_T\tand\\\label{71-3}
&&{{R_{\phi_1, \phi_2}}}(K_1(M_{\phi_1,\phi_2})){ \subset} {{\rho_A(K_0( A))}}.
\eneq
Then, for any  sequence of finite subsets $\{{\cal F}_n\}$ of $C$ whose union is dense in $C,$ any increasing sequence of finite
subsets ${\cal P}_n$ of $K_1(C)$ with
$\cup_{n=1}^{\infty} {\cal P}_n=K_1(C)$ and any decreasing sequence of positive numbers $\{\dt_n\}$ with $\sum_{n=1}^{\infty} \dt_n<\infty,$ there exists a sequence of unitaries $\{u_n\}$ in ${{U_0({\tilde A})}}$ such that
\beq\label{71-4}
{\rm Ad}u_n\circ \phi_1\approx_{\dt_n} \phi_2\,\,\,{\rm on}\,\,\,{\cal F}_n\tand\\
\rho_A({\rm bott}_1(\phi_2, u_n^*u_{n+1})(x))=0\rforal x\in {\cal P}_n
\eneq
and for all sufficiently large $n.$

\end{lem}

 \begin{proof}
 The proof is a modification of that of 7.1 of \cite{Lininv}.
 Note that $A\cong A\otimes U_2.$ Moreover, {there is   {{an isomorphism}} $s: A\otimes U_2\to A$ such that
$s\circ \imath$} is approximately unitarily equivalent to the identity map on $A,$ where
 ${\imath}: A\to A\otimes U_2$ defined by $a\to a\otimes 1_{U_2}$ {for all $a\in A$}.
 Therefore we
may assume that $\phi_1(C),\, \phi_2(C)\subset A\otimes 1_{U_2}.$
It follows from \ref{MUN2} that there exists a sequence of unitaries
$\{v_n\}\subset {\tilde A}$ such that
$$
\lim_{n\to\infty} {\rm Ad}\, v_n\circ \phi_1(c)=\phi_2(c)\tforal c\in C.
$$
{{It follows from \ref{Unitary} that we may assume that $v_n\in U_0({\tilde A}).$}}
In what follows in this proof, we will {{write}} $\phi_1^{\sim}$ and $\phi_2^\sim$ for the extension
of $\phi_1$ and $\phi_2,$ respectively on ${\tilde C}.$
We may assume that ${\cal F}_n$ are in the unit ball and $\cup_{n=1}^{\infty} {\cal F}_n$ is dense in the unit ball of $C.$

Put $\ep_n'=\min\{1/2^{n+2}, \dt_n/2\}.$ Let
$C_n\subset C$ be a  {{simple}} \SCA\, in ${\cal M}_1$
(see \ref{Lfinitelimit})
 which
has finitely generated $K_i(C_n)$ ($i=0,1$) such that $C=\overline{\cup_{n=1}^\infty C_n}$ {{and $\lim_{n\to \infty}\dist (x, C_n)=0$ for all $x\in C$.}}
{{We may assume that
${\cal F}_n\subset C_n$.}} {{Denote by $\iota_n: C_n\to C$ the embedding.}}
Let ${\cal P}_1'\subset K_1(C_1)$ be a finite generating set and ${\cal P}_1=[\iota_1](K_1(C_1).$
Let $\dt_1'$ (in place of $\eta$) and $n_1\ge 1$ be as in \ref{TTbote} for $C$ (in place of $A$),
$\ep_1'$ (in place of $\ep$),  ${\cal F}_n,$ and let $Q_1\subset K_1(C_1)$ be a finite generating set.
 By relabeling, we may assume that $n_1=2.$
Let ${\cal Q}_n$ be a finite set of generators of $K_1(C_n),$
 let $\dt_n'>0$ (in place of {{$\eta$}})  and $k(n)\ge n$ be as in
 \ref{TTbote} for $C$ (in place of $A$), $\ep_n'$ (in
place of $\ep$), ${\cal F}_n$ (in place of ${\cal F}$) and $[\imath_n]({\cal Q}_{n-1})$
(in place of ${\cal P}$).
Note that {{(see \ref{Lappinductivelimit})}}, we assume that
\beq\label{71-9}
[\imath_{n+1}]({\cal Q}_{n+1})\supset {\cal P}_{n+1}\cup [\imath_n]({\cal Q}_n).
\eneq

Write $K_1(C_n)=G_{n,f}\oplus {\rm Tor}(K_1(C_n)),$ {where $G_{n,f}$ is a finitely generated free { abelian} group}.
Let $z_{1,n}, z_{2,n},...,z_{f(n),n}$ be the free generators of $G_{n,f}$ and $z'_{1,n},z'_{2,n},...,z_{t(n),n}'$ be generators of ${\rm Tor}(K_1(C_n)).$
We may assume that
$$
{\cal Q}_n=\{z_{1,n}, z_{2,n},...,z_{f(n),n},z'_{1,n},z'_{2,n},...,z_{t(n),n}'\}.
$$

Let $1/2>\ep_n''>0$ so that ${\text{bott}}_1(h', u')|_{K_1(C_n)}$ is  {{a well-defined}} group \hm, ${\text{bott}}_1(h', u')|_{[\iota_n]({\cal Q}_n)}$ is well defined and
$(\text{bott}_1(h',u')|_{[\iota_n](K_1(C_n))})|_{{\cal Q}_n}={\text{bott}}_1(h', u')|_{[\iota_n]({\cal Q}_n)}$ for any  \hm\, $h': C\to A$ and any unitary $u'\in {\tilde A}$ for which
\beq\label{71-10}
\|[h'(c), u']\|<\ep_n''\tforal c\in {\cal G}_n'
\eneq
for some finite subset ${\cal G}_n'\subset C$ which contains ${\cal F}_n.$

 Let $w_{1,n},w_{2,n},...,w_{f(n),n}, w_{1,n}',w_{2,n}',...,w_{t(n),n}'\in  {\tilde C}$ be unitaries (note that  $C$ has stable rank one)
such that $[w_{i,n}]=(\imath_n)_{*1}(z_{i,n})$ and
$[w_{j,n}']=(\imath_n)_{*1}(z_{j,n}'),$ $i=1,2,...,f(n),$
$j=1,2,...,t(n)$ and  $n=1,2,....$
\Wlog, one may write {{$w_{i,n}=1_{\tilde C}+a_{i,n}$}}  and ${{w_{i,n}'}}=1_{\tilde C}+a_{i,n}',$ where $a_{i,n},a_{i,n}'\in C.$
To simplify notation, without loss of generality, we may assume that $a_{i,n}{{, a_{i,n}'}}\in {\cal G}_n',$ $n=1,2,....$

Let $\dt_1''=1/2$ and, for $n\ge 2,$ let $\dt_n''>0$ (in place of $\dt$)
and ${\cal G}_n''$ (in place of ${\cal F}$) be as in 27.2 of \cite{GLN2}
(with $B={\tilde C}$)
associated
with $w_{1,n},w_{2,n},...,w_{f(n),n}, w_{1,n}',w_{2,n}',...,w_{t(n),n}'$
(in place of $u_1, u_2,..,u_n$) and
$$
\{w_{1,{n-1}}, w_{2,n-1},...,w_{f(n-1), n-1}, w_{1,n-1}',w_{2,n-1}',...,w_{t(n-1), n-1}'\}
$$
(in place of $v_1,v_2,...,v_m$). It is clear that, \wilog, we may assume that ${\cal G}_n''\subset C_n.$

Put $\ep_n=\min\{\ep_n''/2, \ep_n'/2,\dt'_n, \dt_n''/2\}$ and
${\cal G}_n={\cal G}_n'\cup {\cal G}_n''.$
We may assume that
\beq\label{71-11}
{\rm Ad}\, v_n\circ \phi_1\approx_{\ep_n} \phi_2\,\,\,{\rm on\,\,\,} {\cal G}_n,\,\,\, n=1,2,....
\eneq
Thus ${\rm bott}_1(\phi_2\circ \imath_n, v_n^*v_{n+1})$ is well defined.
Since $\Aff(T({\tilde A}))$ is torsion free,
\beq\label{71-12-1}
\tau\big({\rm bott}_1(\phi_2\circ \imath_n, v_n^*v_{n+1})|_{{\rm Tor}(K_1(C_n))}\big)=0\tforal \tau\in T(\tilde{A})
\eneq
(recall $K_i(C_n)$ {{is}} finitely generated). We have
\beq\label{71-12}
\|\phi_2^\sim(w_{j,n}){\rm Ad}\, v_n(\phi_1^\sim(w_{j,n})^*)-1\|<(1/4)\sin(2\pi \ep_n)<\ep_n,\,\,\, n=1,2,....
\eneq
Define
\beq\label{71-13}
h_{j,n}={1\over{2\pi i}}\log(\phi_2^\sim(w_{j,n}){\rm Ad}\, v_n(\phi_1^\sim(w_{j,n})^*)),\,\,\,j=1,2,...,f(n), n=1,2,....
\eneq
Moreover, since $\pi_A(w_{j,n})=1,$ by \eqref{71-12}, $\pi_A(h_{j,n})=0,$ where $\pi_A: {\tilde A}\to \C$ is the quotient map.
Then, for any $\tau\in T({\tilde A}),$
\beq\label{09-07-2021}
{{|\tau(h_{j,n})|<\ep_n<\dt_n',\,\,\, j=1,2,...,f(n), n=1,2,\cdots.}}
\eneq 
Since $\Aff(T({\tilde A}))$ is torsion free, it follows from
27.2 of \cite{GLN2}
that
\beq\label{71-13+}
\tau({1\over{2\pi i}}\log(\phi_2^\sim(w_{j,n}'){\rm Ad}\, v_n(\phi_1^\sim({w_{j,n}'}^{*}))))=0,
\eneq
$j=1,2,...,t(n)$ and $n=1,2,....$
It is standard to {{show}} that the inclusion $M_{\phi_1,\phi_2}\hookrightarrow M_{\td{\phi}_1, \td{\phi}_2}$ induces the iosomorphism  $K_1(M_{\phi_1,\phi_2})\cong K_1(M_{\td{\phi}_1, \td{\phi}_2})$. Hence\\  $R_{\td{\phi}_1, \td{\phi}_2}(K_1(M_{\td{\phi}_1, \td{\phi}_2}))|_{T(A)}=R_{\phi_1, \phi_2}({ K_1}(M_{\phi_1,\phi_2}))$. It is easy to see $R_{\td{\phi}_1, \td{\phi}_2}(x)(\tau_\C)=0$, for any $x\in K_1(M_{\td{\phi}_1, \td{\phi}_2})=K_1(M_{\phi_1,\phi_2})$, where $\tau_{\C}$ is the tracial state of ${\tilde A}$ which vanishes on $A$. Consequently,  the assumption (\ref{71-3}) implies that
$R_{\td{\phi}_1, \td{\phi}_2}(K_1(M_{\td{\phi}_1, \td{\phi}_2})){{|_{T(A)}}} \subset \rho_A(K_0(A))$. By
{{the Exel}} formula (see, for example,  \cite{HL}) and by Lemma 3.5 of \cite{Lnamj}, we conclude that
\begin{equation*}
{{\tau\mapsto}}\widehat{h_{j,n}}(\tau){{(=\tau(h_{j,n}))}}~~\in~ {{R_{\td{\phi}_1, \td{\phi}_2}(K_1(M_{\td{\phi}_1, \td{\phi}_2}))|_{T(A)}}}\subset \rho_A(K_0({ A}))\,\,\,\,\,\, {{{\rm on}\,\,\, T(A).}}
\end{equation*}
{{Note also that $\widehat{h_{j,n}}(\tau_{\C})=0.$}}

{{Now}} define $\af_n': K_1(C_n)\to \rho_{A}(K_0(A))$ by
\beq\label{71-15}
{{\hspace{-0.2in}\af_n'(z_{j,n})(\tau)=\widehat{h_{j,n}}(\tau)=\tau(h_{j,n}),\,\,\, 1\le j\le f(n)
\andeqn
\af_n'(z_{j,n}')=0,\,\,\, 1\le j\le t(n),}}
\eneq
{{for all $n.$}}
 Since $\af_n'(K_1(C_n))$ is free,
there is a \hm\, $\af_n^{(1)}: K_1(C_n)\to K_0(A)$
such that
\beq\label{71-16}
&&(\rho_A\circ \af_n^{(1)}(z_{j,n}))(\tau)=\tau(h_{j,n}),\,\,\, j=1,2,...,f(n),\,\,{{\tau\in T(A)}} \,\,\andeqn\\
&&{{\af_{n}^{(1)}(z_{j,n}')=0,\,\,\, j=1,2,...,t(n).}}
\eneq
Define
$\af_n^{(0)}: K_0(C_n)\to K_1(A)$ by $\af_n^{(0)}=0.$ By the UCT,
there is $\kappa_n\in KL({ S}C_n, A)$ such that $\kappa_n|_{K_i(C_n)}=\af_n^{(i)},$ $i=0,1,$
where $SC_n$ is the suspension of $C_n$ (here, we also {identify} $K_i(C_n)$ with $K_{i+1}(SC_n)$).

By the UCT again, there is $\af_n\in KL(C_n\otimes C(\T), A)$ such that
$\af_n\circ \boldsymbol{\bt}|_{\underline{K}(C_n)}=\kappa_n.$ In particular,
$\af_n\circ \boldsymbol{\bt}|_{K_1(C_n)}=\af_n^{(1)}.$   By
\ref{TTbote} and (\ref{09-07-2021}),
there exists a unitary
$U_n\in U_0({\tilde A})$ such that
\beq\label{71-17}
\|[\phi_{ 2}(c),\, U_n]\|<\ep_n''\rforal c\in {\cal F}_n\andeqn\\\label{71-17+}
\rho_A({\text{bott}_1}(\phi_2,\, U_n)(\iota_{n*1}(z_{j,n})))=-\rho_A\circ \af_n(z_{j,n}),
\eneq
$j=1,2,...,f(n).$ We also have
\beq\label{71-18}
\rho_A({\text{bott}_1}(\phi_2,\, U_n)(\iota_{n,*1}(z_{j,n}')))=0,\,\,\, j=1,2,...,t(n).
\eneq
By the Exel trace formula {{(see, for example, \cite{HL})}}, {{(\ref{71-16}) and (\ref{71-17+}),}}
 we have
\beq\label{71-19}
\hspace{-0.2in}\tau(h_{j,n})=-\rho_A({\rm bott}_1(\phi_2, U_n)(\iota_{n*1}(z_{j,n}))(\tau)
    =-\tau({1\over{2\pi i}}\log(U_n\phi_2^\sim(w_{j,n})U_n^*\phi_2^\sim(w_{j,n}^*)))
\eneq
for all $\tau\in T(A),$ $j=1,2,...,f(n).$
Define $u_n=v_nU_n,$ $n=1,2,....$ By 6.1 of \cite{Lnind},  {{\eqref{71-13} and  (\ref{71-19}),}}
we compute that
\beq\label{71-21}
&&\hspace{-0.9in}\tau({1\over{2\pi i}}\log(\phi_2^\sim(w_{j,n}){\rm Ad} u_n(\phi_{ 1}^\sim(w_{j,n}^*)))))\\
&&\hspace{-0.6in}=\tau({1\over{2\pi i}}\log(U_n\phi_2^\sim(w_{j,n})U_n^*v_n^*\phi_1^\sim(w_{j,n}^*)v_n)))\\
&&\hspace{-0.6in}=\tau({1\over{2\pi i}}\log(U_n\phi_2^\sim(w_{j,n})U_n^*\phi_2^\sim(w_{j,n}^*)\phi_2^\sim(w_{j,n})v_n^*
\phi_1(w_{j{{,n}}}^*)v_n)))\\
&&\hspace{-0.6in}=\tau({1\over{2\pi i}}\log(U_n\phi_2^\sim(w_{j,n})U_n^*\phi_2^\sim(w_{j,n}^*))))
+\tau({1\over{2\pi i}}\log(\phi_2^\sim(w_{j,n})v_n^*\phi_1^\sim(w_{j,n}^*)v_n)))\\
&&\hspace{-0.6in}={{\rho_A({\rm bott}_1(\phi_2, U_n)(z_{j,n}))(\tau)}}+\tau(h_{j,n})=0
\eneq
for all $\tau\in T(A),$ $j=1,2,...,f(n)$ and $n=1,2,....$
By (\ref{71-13+}) and (\ref{71-18}),
\begin{equation}\label{71-22}
\tau({1\over{2\pi i}}\log(\phi_2^\sim(w_{j,n}'){\rm Ad} u_n(\phi_{1}^\sim((w_{j,n}')^*))))=0,
\end{equation}
$j=1,2,...,t(n)$ and $n=1,2,....$
Let
\beq\label{71-23}
b_{j,n}&=&{1\over{2\pi i}}\log(u_n\phi_2^\sim(w_{j,n})u_n^*\phi_1^\sim(w_{j,n}^*)),\\
b_{j,n}'&=&{1\over{2\pi i}}\log(\phi_2^\sim(w_{j,n})u_n^*u_{n+1}\phi_2^\sim(w_{j,n}^*)u_{n+1}^*u_n)\andeqn\\
b_{j,n+1}''&=&{1\over{2\pi i}}\log(u_{n+1}\phi_2^\sim(w_{j,n})u_{n+1}^*\phi_1^\sim(w_{j,n}^*)).
\eneq
$j=1,2,...,f(n)$ and $n=1,2,....$ We have, by (\ref{71-21}),
\beq\label{71-24}
\tau(b_{j,n})&=&\tau({1\over{2\pi i}}\log(u_n\phi_2^\sim(w_{j,n})u_n^*\phi_1^\sim(w_{j,n}^*)))\\
&=&\tau({1\over{2\pi i}}\log(\phi_2^\sim(w_{j,n})u_n^*\phi_1^\sim(w_{j,n}^*)u_n))=0
\eneq
for all $\tau\in T(A),$ $j=1,2,..., f(n)$ and $n=1,2,....$ Note that
$\tau(b_{j, n+1})=0$ for all $\tau\in T(A),$ $j=1,2,...,f(n+1).$
It follows from  27.2 of \cite{GLN2}
and \eqref{71-9} that
\begin{equation*}
\tau(b_{j,n+1}'')=0\tforal \tau\in T(A),\,\,\,j=1,2,...,f(n),\,\,\,n=1,2,....
\end{equation*}
Note  also
that
\begin{equation*}
u_ne^{2\pi i b_{j,n}'}u_n^*=e^{2\pi i b_{j,n}}\cdot e^{-2\pi i b_{j,n+1}''},\,\,\,j=1,2,...,f(n).
\end{equation*}
Thus, by 6.1 of \cite{Lnind}, we compute that
\beq\label{71-26}
\tau(b_{j, n}')=\tau(b_{j,n})-\tau(b_{j,n+1}'')=0\tforal \tau\in T(A).
\eneq
By the Exel formula {{(see \cite{HL})}} and (\ref{71-26}),
\beq\label{71-27}
&&\hspace{-0.8in}\rho_A({\rm bott}_1(\phi_2, u_n^*u_{n+1}))(w_{j,n}^*)(\tau)
=\tau({1\over{2\pi i}}\log(u_n^*u_{n+1}\phi_2^\sim(w_{j,n})u_{n+1}^*u_n\phi_2^\sim(w_{j,n}^*)))\\
&=&{{\tau({1\over{2\pi i}}\log(\phi_2^\sim(w_{j,n})u_{n}^*u_{n+1}\phi_2^\sim(w_{j,n}^*)u_{n+1}^*u_n))=0}}
\eneq
for all $\tau\in T(A)$ and $j=1,2,...,f(n).$
Thus
\beq\label{71-28}
\rho_A({\rm bott}_1(\phi_2,u_n^*u_{n+1})(w_{j,n})(\tau)=0\tforal \tau\in T(A),
\eneq
$j=1,2,...,f(n)$ and $n=1,2,....$ We also have
\beq\label{71-29}
\rho_A({\rm bott}_1(\phi_2, u_n^*u_{n+1})(w_{j,n}')(\tau)=0\tforal \tau\in T(A),
\eneq
$j=1,2,...,f(n)$ and $n=1,2,....$ By  27.2 of \cite{GLN2},
we have that
\beq\label{71-30}
\rho_A({\rm bott}_1(\phi_2, u_n^*u_{n+1})(z)=0\tforal z\in {\cal P}_n,
\eneq
$n=1,2,....$
 \end{proof}

 \begin{thm}\label{Tm72}
 Let $C_0$ be a  simple \CA\,  in   ${\cal M}_1$
 with continuous scale,  {{$A_1$}} be a  separable simple \CA\, in {${\cal D}$} with continuous scale, and {{$U_1$}} and $U_2$ two UHF-algebras of infinite type. Let $C={{C_0}}\otimes U_1$ and $A=A_1\otimes U_2.$
 Suppose that $\phi_1,\, \phi_2: C\to A$ are two monomorphisms which {{map}} strictly positive elements
 to strictly positive elements.
Then they are {{strongly}} asymptotically unitarily equivalent  {{(see \ref{Dstrongaue})}} if and only if
\beq\label{Tm72-1}
[\phi_1]=[\phi_2]\,\,\,{\rm in}\,\,\, KK(C,A),\\
\phi_1^{\dag}=\phi_2^{\dag},\,\,\,(\phi_1)_T=(\phi_2)_T \andeqn {\overline{R_{\phi_1,\phi_2}}}=0.
\eneq
\end{thm}

\begin{proof}
We will prove the ``if " part only. The ``only if" part follows
from 4.3 of \cite{Lininv}.
%
%
Let $C=\overline{\cup_{n=1}^\infty C_n}$ with $\lim_{n\to\infty}{\rm dist} (x, C_n)=0$ for any $x\in C,$
where $C_n$ is in ${\cal M}_1$ such that $K_i(C_n)$ is finitely generated
as \ref{Lfinitelimit} 
 and
\ref{Lappinductivelimit}.
By passing to a subsequence, applying \ref{Lappinductivelimit},  we may assume
that $\underline{K}(C)=\lim_{n\to\infty}(\underline{K}(C_n), j_n),$ where
$j_n\in {{{\rm Hom}_{\Lambda}(\underline{K}(C_n), \underline{K}(C_{n+1}))}},$
$j_{n, \infty}=[\imath_n],$ $[\imath_n](\underline{K}(C_n))\subset [\imath_{n+1}](\underline{K}(C_{n+1}))$ and   $\imath_n: C_n\to C$ is the embedding, $n=1,2,....$
Let ${\cal F}_n\subset C$ be
a sequence of subsets of $C$ such that
${\cup_{n=1}^{\infty}{\cal F}_n}$ is dense in $C.$
Put
$$
M_{\phi_1, \phi_2}=\{(f, {  c})\in C([0,1], A){ \oplus C}: f(0)=\phi_1(c)\andeqn f(1)=\phi_2(c)\}.
$$
Since $C$ satisfies the UCT, the assumption that $[\phi_1]=[\phi_2]$ in
$KK(C,A)$ implies that the following exact sequence splits:
\beq\label{72-1}
0\to \underline{K}(SA)\to \underline{K}(M_{\phi_1, \phi_2}) \rightleftharpoons_{\theta}^{\pi_e} \underline{K}(C)\to 0
\eneq
for some $\theta\in Hom(\underline{K}(C), \underline{K}(A)),$ where
$\pi_e: M_{\phi_1,\phi_2}\to C$ is { the projection  to $C$ (see
 \ref{Mappingtorus}}).
Furthermore, since $\tau\circ \phi_1=\tau\circ \phi_2$ for all $\tau\in T(A)$ and ${\overline{R_{\phi_1,\phi_2}}}=0,$ we may also assume that
\beq\label{72-2}
R_{\phi_1,\phi_2}(\theta(x))=0\rforal x\in K_1(C).
\eneq
Recall  {{that}}  $\lim_{n\to\infty}(\underline{K}(C_n),[\imath_n])=\underline{K}(C).$
Since $K_i(C_n)$ is finitely generated, there exists $K(n)\ge 1$ such that
\beq\label{72-3+}
Hom_{\Lambda}(F_{K(n)}\underline{K}(C_n),\, F_{K(n)}\underline{K}(A))=Hom_{\Lambda}(\underline{K}(C_n), \underline{K}(A))
\eneq
(see also \cite{DL} for the notation $F_m$ there).

Let $\dt_n'>0$ (in place of $\dt$), $\sigma_n'>0$ (in place of $\sigma$), ${\cal G}_n'\subset C$ (in place of ${\cal G}$), \linebreak ${\{p_{1,n}', p_{2,n}',...,p'_{I(n),n)}, q_{1,n}', q_{2,n}',...,q_{I(n),n}'\}}$ (in place of
$\{p_1, p_2,...,p_k, q_1, q_2,...,q_k\}),$ ${\cal P}_n'\subset  {\underline{K}(C)}$
(in place of ${\cal P}$) corresponding to $1/2^{n+2}$ (in place of $\ep$) and ${\cal F}_n$ (in place of ${\cal F}$)
{{be}} as required by {{14.8 of \cite{GLII}.}}
Note that, by the choice above  as in
{{14.8 of}} \cite{GLII}, 
 $G_{u,n}',$ the subgroup generated by
${\{[p'_{i,n}]-[q'_{i,n}]: 1\le i\le I(n)\}}$ is  {{assumed to be free.}}

Without loss of generality, we may assume that
${\cal G}_n'\subset \imath_n({\cal G}_n)$
and ${\cal P}_n'\subset [\imath_n]({\cal P}_n)$ for some finite subset ${\cal G}_n\subset C_n,$ and for some finite subset ${\cal P}_n\subset \underline{K}(C_n),$ we may assume that $p_{i,n}'=\imath_n^\sim (p_{i,n})$ and $q_{i,n}'=\imath_n^\sim(q_{i,n})$
for some projections $p_{i,n}, q_{i,n}\in M_{R(n)}(\td C_n),$ $i=1,2,...,I(n).$
Write  $p_{i,n}={\bar p}_{i,n}+x(p_{i,n})$ and $q_{i,n}={\bar q}_{i,n}+x(q_{i,n}),$
where ${\bar p}_{i,n}$ and ${\bar q}_{i,n}$ are scalar projections in $M_{R(n)},$ $x(p_{i,n})$ and
$x(q_{i,n})\in  M_{R(n)}(C_n)_{s.a.}.$ Also, we may assume
that $x(p_{i,n}), x(q_{i,n})\subset \{(a_{k,j})_{R(n)\times R(n)}: a_{k,j}\in {\cal G}_n\}.$

\Wlog\, we may {{also}} assume that the subgroup $G_{n,u}$  generated by $\{[p_{i,n}]-[q_{i,n}]: 1\le i\le I(n)\}$
is free.

 We may assume that ${\cal P}_n$ contains a set of generators
of $F_{K(n)}\underline{K}(C_n),$ ${\cal F}_n\subset {\cal G}_n'$ and
$\dt_n'<1/2^{n+3}.$
We may also assume that
${\rm Bott}(h', u')|_{{\cal P}_n}$ is well defined whenever $\|[h'(a),\, u']\|<\dt_n'$ for all $a\in {\cal G}_n'$ and for any \hm\, $h'$ from $C_n$ and a unitary $u'$ in the target algebra. Note that ${\rm Bott}(h',\,u')|_{{\cal P}_n}$ defines ${\rm Bott}(h'\, u').$

We further assume that
\beq\label{72-4}
{\rm Bott}(h,\, u)|_{{\cal P}_n}={\rm Bott}(h',u)|_{{\cal P}_n}
\eneq
provided that $h\approx_{\dt_n'}h'$ on ${\cal G}_n'.$
We may also assume
that $\dt_n'$ is smaller than {$\dt/16$} for the $\dt$ defined in 2.15 of \cite{Lininv} for $C_n$ (in place of $A$) and ${\cal P}_n$ (in place of ${\cal P}$).
{Let $k(n)\ge n$ (in place of $n$), $\eta_n'>0$ (in place of $\eta$) and ${\cal Q}_{k(n)}\subset K_1(C_{k(n)})$
(in place of ${\cal P}^{(0)}$) be required by 
\ref{TTboteplus}
for
$\dt_{k(n)}'/4$ (in place of $\ep$),
$\imath_{n}({\cal G}_{k(n)})$ (in place of ${\cal F}$),  ${\cal P}_{k(n)}$ (in place of ${\cal P}$) and
$[p_{j,n}]-[q_{j,n}]$ (in place ${\bar s}_j$) ($j=1,2,...,k(n)$),
{and $\sigma_{k(n)}'/16$ (in place of $\sigma$)}.
 We may assume that ${\cal Q}_{k(n)}$ forms a  free
 generator set  for the free part of
 $K_1(C_{k(n)})$ (see \ref{TTboteplus}).
Since  ${\cal P}_{k(n)}$
generates $F_{K(n)}\underline{K}(C_{k(n+1)}),$ we may assume that ${\cal Q}_{k(n)}\subset {\cal P}_{k(n)}.$}

For $C_n,$ since $K_i(C_n)$ ($i=0,1$) is finitely generated, by
(\ref{72-3+}), we may further assume that
$j_{k(n), \infty}$ is injective on $j_{n, k(n)}(\underline{K}(C_n)),$ $n=1,2,....$ By passing to a subsequence, to simplify notation, we may also assume that
$k(n)=n+1.$ Let $\dt_n=\min\{\eta_n, \sigma_n',\dt_n'/2\}/R(n)^2.$ By \ref{inv71}, there are unitaries $v_n\in {{U_0({\tilde A})}}$ such that
\beq\label{72-5}
&&{\rm Ad}\, v_n\circ \phi_1\approx_{\dt_{n+1}/4} \phi_2\,\,\,{\rm on}\,\,\,
\imath_{n+1}({\cal G}_{n+1}),\\
&&\rho_A({\rm bott}_1(\phi_2, v_n^*v_{n+1}))(x)=0
\,\,\,\rforal x\in
[\imath_{n+1}](K_1(C_{n+1}))\andeqn\\
&&\|[\phi_2(c),\, v_n^*v_{n+1}]\|<\dt_{n+1}/2\,\,\rforal a\in
\imath_{n+1}({\cal G}_{n+1})
\eneq
(note that $K_1(C_{n+1})$ is finitely generated).
Note that, by (\ref{72-4}),
we may also assume that
\beq\label{72-7}
{\rm Bott}(\phi_1,\, v_{n+1}v_n^*)|_{[\imath_{n}]({\cal P}_n)}&=&{\rm Bott}(v_n^*\phi_1v_n, \, v_n^*v_{n+1})|_{[\imath_{n}]({\cal P}_n)}\\
&=&{\rm Bott}(\phi_2,\, v_n^*v_{n+1}])|_{[\imath_{n}]({\cal P}_n)}.
\eneq
In particular,
\beq\label{72-8}
{\rm bott}_1(v_n^*\phi_1v_n,\, v_n^*v_{n+1})(x)={\rm bott}_1(\phi_2,\, v_n^*v_{n+1})(x)
\eneq
for all $x\in {{(\imath_{n+1})_{*1}}}(K_1(C_{n+1})).$
%
%
%
By applying {{\ref{D642cbms},}}
without loss of generality,
we may assume that the triple $\phi_1, \phi_2$ and $v_n$ defines
an element $\gamma_n\in
KL_{loc}^{C}([\imath_{n+1}](\underline{K}(C_{n+1})),\underline{K}(M_{\phi_1, \phi_2}))$
and
{{$[\pi_e]\circ\gamma_n\circ {{[\imath_{n+1}]}}=
[\imath_{n+1}]$ (see \eqref{2020-826-1}).}}
Moreover, by 10.4 and 10.5 of \cite{Lnamj} (see also the end of
6.4.2 of \cite{Lncbms}), we may assume, without loss of {generality}, that
\beq\label{72-9}
|\tau(\log(\phi_2^\sim\circ \imath_{n}^\sim(z_j^*)v_n\phi_1^\sim\circ \imath_{n}^\sim(z_j)v_n))|<{{\dt_{n+1}/2}},
\eneq
$j=1,2,..., {{r(n+1)}}$, where $\{z_1,z_2,...,z_{r(n)}\} \subset U(\td C_{n+1})$
which forms a set of generators of $K_1(C_{n+1})$ (recall  that $C_{n+1}$ is a simple \CA\, of stable rank one).
{We may assume that
$[z_j]\in {\cal Q}_n\subset {\cal P}_n,$ $j=1,2,...,r(n).$}

Let $H_n=[\imath_{n+1}](\underline{K}(C_{n+1})).$  Since $\bigcup_{n=1}[\imath_n](\underline{K}(C_n))=\underline{K}(C)$
and $[\pi_e]\circ \gamma_n\circ [\imath_{n+1}]=[\imath_{n+1}],$ we conclude that
\beq\label{72-10}
\underline{K}(M_{\phi_1, \phi_2})=\underline{K}(SA)+\cup_{n=1}^{\infty}\gamma_n(H_n).
\eneq
{{Recall also $H_n=[\imath_{n+1}](\underline{K}({{C_{n+1}}}))\subset {{[\imath_{n+2}]}}(\underline{K}(C_{n+2})).$
So $\gamma_n-\gamma_{n+1}$ is defined on $H_n.$}}
Thus, by passing to a subsequence, we may further assume that
\beq\label{72-11}
\gamma_{n+1}(H_n)\subset \underline{K}(SA)+\gamma_{n+2}(H_{n+1}),\,\,\,n=1,2,....
\eneq

By identifying $H_n$ with $\gamma_{n+1}(H_n),$ we may write $j_n': \underline{K}(SA)\oplus H_n\to \underline{K}(SA)\oplus H_{n+1}$ {for the inclusion in  (\ref{72-11}).} By
(\ref{72-10}), the inductive limit is $\underline{K}(M_{\phi_1,\phi_2}).$
From the definition of $\gamma_n,$ we note that $\gamma_n\circ [\iota_{n+1}]-\gamma_{n+1}\circ [\imath_{n+2}]$ maps $\underline{K}(C_{n+1})$ into $\underline{K}(SA).$
By 10.6 of \cite{Lnamj} {{(see also 27.4 of \cite{GLII}, in particular, (e.27.54) there),}}
$$
\Gamma({\rm Bott}(\phi_1,\, v_nv_{n+1}^*))|_{H_n}=(\gamma_{n+1}-\gamma_{n+2})|_{H_n}
$$
(see 10.4, 10.5 and 10.6  of \cite{Lnamj}, {{and also 27.4 of \cite{GLN}}} for the definition of $\Gamma({\rm Bott}(\cdot,\cdot))$) gives a \hm\, $\xi_n: H_n\to \underline{K}(SA).$ {{Let}} $\zeta_n=\gamma_{n+1}|_{H_n}$  and $j^n:H_n\to H_{n+1}$ {{be the}}
embedding. Then
\beq\label{72-13}
j_n'(x,y)=(x+\xi_n(y), j^n(y))
\eneq
for all $(x,y)\in \underline{K}(SA)\oplus H_n.$ Thus we obtain the following diagram:
\beq\label{72-14}
\begin{array}{ccccccc}
0 \to  & \underline{K}(SA)  &\to & \underline{K}(SA)\oplus H_n &\to
& H_n &\to 0\\\nonumber
 &\| & &\hspace{0.4in}\| \hspace{0.15in}\swarrow_{\xi_n} \hspace{0.05in}\downarrow_{j^{n}}
 &&
 \hspace{0.2in}\downarrow_{j^{n}}
 &\\\label{72-15}
 0 \to  & \underline{K}(SA)  &\to & \underline{K}(SA)\oplus H_{n+1} &\to & H_{n+1} &\to 0\\\nonumber
 &\| & &\hspace{0.4in}\| \hspace{0.1in}\swarrow_{\xi_{n+1}}\downarrow_{j^{n+1}}
 &&
\hspace{0.2in} \downarrow_{j^{n+1}}
&\\
 0 \to  & \underline{K}(SA)  &\to & \underline{K}(SA)\oplus H_{n+2} &\to & H_{n+2} &\to {{0.}}\\
 \end{array}
\eneq
By the assumption that ${{\overline{ R_{\phi_1, \phi_2}}}}=0,$ map $\theta$ also gives the following decomposition:
\beq\label{72-16}
{\rm ker}R_{\phi_1,\phi_2}={\rm ker}\rho_A\oplus K_1(C).
\eneq

Define $\theta_n=\theta\circ [\imath_{n+1}]$ and $\kappa_n=\zeta_n\circ [\imath_{n+1}]-\theta_n.$
Note that
\beq\label{72-17}
\theta_n=\theta_{n+1}\circ j_{n+1}
\eneq
(recall $j_{n+1}\in {\rm Hom}_{\Lambda}(\underline{K}(C_{n+1}), \underline{K}(C_{n+2}))$.) We also have that
\beq\label{72-18}
\zeta_n-\zeta_{n+1}\circ j^{n}
=\xi_n.
\eneq
Since $[\pi_e]\circ (\zeta_n\circ [\imath_{n+1}]-\theta_n)|_{\underline{K}(C_{n+1})}=0,$ $\kappa_n$ maps
$\underline{K}(C_{n+1})$ into $\underline{K}(SA).$ It follows that
\beq\label{72-19}
\kappa_n-\kappa_{n+1}\circ  j_{n+1}
&=&
\zeta_n\circ [\imath_{n+1}]-\theta_n-\zeta_{n+1}\circ
[ \imath_{n+2}]\circ j_{n+1}
+\theta_{n+1}\circ j_{n+1}
\\{\label{72-19+}}
&=&\zeta_n\circ [\imath_{n+1}]-\zeta_{n+1}\circ[\imath_{n+2}]\circ j_{n+1}=\xi_n\circ [\imath_{n+1}].
\eneq
It follows from \ref{VuV} that there {{is an integer}}
$N_1\ge 1,$  {{an $
\imath_{n+1}({\cal
G}_{n+1})$-${\dt_{n+1}\over{4}}$}}-multiplicative \morp\, $L_n: \imath_{n+1}
(C_{n+1})\to  M_{1+N_1}(M_{\phi_1,\phi_2}),$ a  \hm\,
$h_0: \imath_{n+1}(C_{n+1})\to {{{\cal W}\to M_{N_1}(A),}}$
and a continuous path of unitaries $\{V_n(t): t\in [0,3/4]\}$ of
$M_{1+N_1}(A)$ such that $[L_n]|_{{\cal P}_{n+1}'}$ is well
defined, $V_n(0)=1_{M_{1+N_1}(A)},$
\begin{equation*}
[L_n\circ \imath_{n+1}]
|_{{\cal P}_n}=(\theta\circ
[\imath_{n+1}])|_{{\cal P}_n},
\end{equation*}
and
\begin{equation*}
\pi_t\circ L_n\circ
\imath_{n+1}\approx_{\dt_{n+1}/4} {{\rm Ad}}\, V_n(t)\circ
((\phi_1\circ \imath_{n+1})\oplus (\phi_1\circ h_0\circ
\imath_{n+1}))
\end{equation*}
on $\imath_{n+1}({\cal G}_{n+1})$ for all $t\in (0,3/4],$
\begin{equation*}
\pi_t\circ L_n\circ
\imath_{n+1}\approx_{\dt_{n+1}/4} {{\rm Ad}}\, V_n(3/4)\circ
((\phi_1\circ \imath_{n+1})\oplus (\phi_1\circ h_0\circ
\imath_{n+1}))
\end{equation*}
on $\imath_{n+1}({\cal G}_{n+1})$ for all $t\in (3/4,1),$ and
\begin{equation*}
\pi_1\circ L_n\circ
\imath_{n+1}\approx_{\dt_{n+1}/4}\phi_2\circ
\imath_{n+1}\oplus \phi_2\circ h_0\circ \imath_{n+1}
\end{equation*}
on $\imath_{n+1}({\cal G}_{n+1}),$ where $\pi_t:
M_{\phi_1,\phi_2}\to A$ is the {{point evaluation}} at $t\in (0,1).$

Note that $R_{\phi_1,\phi_2}(\theta(x))=0$ for all $x\in
[\imath_{n+1}](K_1(C_{n+1})).$  As computed in 10.4 of
\cite{Lnamj},
\beq\label{72-24}
\tau(\log((\phi_2^\sim(x)\oplus h_{0,2}^\sim(x)^*{{V^*_n(3/4)}}(\phi_1^\sim (x)\oplus
h_{0,1}^\sim(x))V_n(3/4)))=0
\eneq
for $x\in \imath_{n+1}^\sim(\{z_1,z_2,...,z_{r(n)}\})$
and for all $\tau\in T(A),$  {{where $h_{0,i}:=\phi_i\circ h_0.$}}
{{Since {{$h_0$}} factors through ${\cal W}$ and {{$(\phi_1)_T=(\phi_2)_T,$}}
there is a unitary $S_n\in U(M_{N_1}(A)^\sim)$
such that $S_n^*h_{0,1}(c)S_n\approx_{\dt_{n+2}/16} h_{0,2}(c)$
for all $c\in \imath_n({\cal G}_n)\cup \imath_{n+1}({\cal G}_{n+1}).$}}
{{We may even assume that, for $1\le i\le r(n),$
$$\tau(\log(h_{0,2}(z_j)S_n^*h_{0,1}S_n))<\dt_{n+1}/4\, \rforal \tau\in T(A).$$}}
Define $W_n'={{{\rm diag}}}(v_{n+1}, S_n)\in M_{1+N_1}(A).$ Then
$\text{Bott}((\phi_1\oplus h_{0,1})\circ \imath_{n+1},\,
W_n'({{V^*_n(3/4)}})$ defines a \hm\, ${\tilde \kappa}_n\in
{\rm Hom}_{\Lambda}(\underline{K}(C_{n+1}),\underline{K}(SA)).$
Then, by
(\ref{72-9}),
for $\tau\in T(A),$
\beq\label{72-25}
|\tau(\log((\phi_2\oplus h_{0,2})^\sim)\circ \imath_{n+1}^\sim(z_j)^*(W_n')^*(\phi_1\oplus h_{0,1})^\sim\circ \imath_{n+1}^\sim(z_j)W_n'))|<\dt_{n+1},
\eneq
$j=1,2,...,r(n).$
  {{Since $[h_{0,i}]=0,$}} one computes {{(see the late part of \ref{D642cbms})}} that
\beq\label{72-25+}
\hspace{-0.2in}\Gamma(\text{Bott}(\phi_1\circ {\imath}_{n+1}\oplus h_{0,1}\circ \imath_{n+1},\, W_n'{{V^*_n(3/4)}}))|_{{\cal P}_n}=
{{(\gamma_n-[L])[\imath_n]|_{{\cal P}_n}=
(\gamma_{n}-\theta)[\imath_n]|_{{\cal P}_n}.}}
\eneq
Put ${\tilde V}_n=V_n(3/4).$
Let 
\beq\label{72-26}
&&\hspace{-0.7in}b_{j,n}={1\over{2\pi i}}\log({\tilde V}_n^*(\phi_1\oplus h_{0,1})^\sim \circ \imath_{n+1}^\sim(z_j){\tilde V}_n (\phi_2\oplus h_{0,2})^\sim\circ \imath_{n+1}^\sim(z_j)^*),\\
&&\hspace{-0.7in}b_{j,n}'={1\over{2\pi i}}\log((\phi_1\oplus h_{0,1})^\sim\circ \imath_{n+1}^\sim(z_j){\tilde V}_n(W_n')^*(\phi_1\oplus h_{0,1})^\sim\circ \imath_{n+1}^\sim(z_j)^*W_n'{\tilde V}_n^*)\andeqn\\
&&\hspace{-0.7in}b_{j,n}''={1\over{2\pi i}}\log((\phi_2\oplus h_{0,2})^\sim\imath_{n+1}^\sim (z_j)(W_n')^* (\phi_1\oplus h_{0,1})^\sim\circ \imath_{n+1}^\sim(z_j)^*W_n'),
\eneq
$j=1,2,...,r(n).$
By (\ref{72-24}) and (\ref{72-25}),
\beq\label{72-27}
\tau(b_{j,n})=0 \andeqn |\tau(b_{j,n}'')|<\dt_{n+1}
\eneq
for all $\tau\in T(A).$
Note that
\beq\label{72-28}
{\tilde V}_n^*e^{2\pi i b_{j,n}'}{\tilde V}_n=
e^{2\pi i b_{j,n}} e^{2\pi i b_{j,n}''}{{.}}
\eneq
Then, by 6.1 of \cite{Lnind} and  by (\ref{72-27})
\beq\label{72-29}
\tau(b_{j,n}')&=&\tau(b_{j,n})-\tau(b_{j,n}'')
=\tau(b_{j,n}'')
\eneq
for all $\tau\in T(A).$
It follows from this and (\ref{72-7}) that
\beq\label{72-30}
|\rho_A({{{\td \kappa}_n}}(z_j))(\tau)| 
<\dt_{n+1},\,\,\,j=1,2,{{...,}}
\eneq
for all $\tau\in T(A).$
It follows from \ref{TTboteplus}
that there is a unitary $w_n'\in {{U_0({\tilde A})}}$
such that
\beq\label{72-31}
&&\|[\phi_1(a), w_n']\|<\dt_{n+1}'/4\rforal a\in \imath_{n+1,\infty}({\cal G}_{n+1})\andeqn\\
&&\text{Bott}(\phi_1\circ \imath_{n+1},\, w_n')=-{\tilde
\kappa}_n\circ[\imath_{n+1}].
\eneq
By (\ref{72-4}),
\beq\label{72-32}
\text{Bott}(\phi_2\circ \imath_{n+1},\,v_n^*w_n'v_n)|_{{\cal
P}_n}=-{\tilde \kappa}_n\circ [\imath_{n+1}]|_{{\cal P}_n}.
\eneq
 It follows from 10.6 of \cite{Lnamj}  and \eqref{72-25+} that
\beq\label{72-33}
\Gamma(\text{Bott}(\phi_1\circ \imath_{n+1}, w_n'))&=&-\kappa_n\circ [\imath_{n+1}]
\andeqn\\
 \Gamma(\text{Bott}(\phi_1\circ \imath_{n+2},
w_{n+1}'))&=&-\kappa_{n+1}\circ [\imath_{n+2}].
\eneq
We also have
\beq\label{72-34}
\Gamma(\text{Bott}(\phi_1\circ \imath_{n+1},
v_nv_{n+1}^*))|_{\underline{K}(C_{n+1})}=(\zeta_n-\zeta_{n+1}\circ j^n)\circ[\imath_{n+1}]
=\xi_n\circ [\imath_{n+1}].
\eneq
But, by (\ref{72-19}) {and (\ref{72-19+})},
\beq\label{72-35}
(-\kappa_n +\xi_n\circ[\imath_{n+1}]+\kappa_{n+1}\circ j_{n+1})=0.
\eneq
By 10.6 of \cite{Lnamj} ({{see also 27.4 of \cite{GLN2}}}), $\Gamma({\rm Bott}(.,.))=0$ if and only
if ${\rm Bott}(.,.)=0.$ Thus, by (\ref{72-32}), (\ref{72-33}) and (\ref{72-34}),
\beq\label{72-36}
{\hspace{-0.5in}} -\text{Bott}(\phi_1\circ \imath_{n+1},\,w_n')
+\text{Bott}(\phi_1\circ \imath_{n+1}, \,v_nv_{n+1}^*)
+\text{Bott}(\phi_1\circ\imath_{n+1}, {w_{n+1}'}) =0.
\eneq
Put $w_n=v_n^*(w_n')v_n$ and  $u_n=v_nw_n^*,$ $n=1,2,....$
{{Since $v_n, w_n'\in U_0({\tilde A}),$ $u_n\in U_0({\tilde A}),$ $n=1,2,....$}}
 Then, by (\ref{72-5}) and
(\ref{72-31}),
\beq\label{72-37}
{{\rm Ad}}\, u_n\circ \phi_1\approx_{\dt_n'/2} \phi_2\rforal a\in
\imath_{n+1}({\cal G}_{n+1}).
\eneq
From (\ref{72-7}), (\ref{72-4}) and (\ref{72-36}), we compute
that
\beq\label{72-38}
&&\hspace{-1.3in}\text{Bott}(\phi_2\circ \imath_{n+1},u_n^*u_{n+1})
= \text{Bott}(\phi_2\circ \imath_{n+1}, w_nv_n^*v_{n+1}w_{n+1}^*)\\
&&\hspace{0.3in} {{=\text{Bott}(\phi_2\circ \imath_{n+1}, w_n)+\text{Bott}(\phi_2\circ \imath_{n+1}, v_n^*v_{n+1})}}\\
&&\hspace{0.6in}
+\text{Bott}(\phi_2\circ \imath_{n+1}, w_{n+1}^*)\\
&&\hspace{0.3in}=\text{Bott}(\phi_1\circ \imath_{n+1},w_n')+\text{Bott}(\phi_1\circ \imath_{n+1}, {{v_{n+1}v_{n}^*}})\\
&&\hspace{0.6in}+\text{Bott}(\phi_1\circ \imath_{n+1}, (w_{n+1}')^*)\\
&&\hspace{0.3in}=-[-\text{Bott}(\phi_1\circ \imath_{n+1}, w_n')+\text{Bott}(\phi_1\circ \imath_{n+1}, v_nv_{n+1}^*)\\\label{Add-403-1}
&&\hspace{2.1in}+\text{Bott}(\phi_1\circ \imath_{n+1}, {w_{n+1}'})]=0.
\eneq
Let $x_{i,n}=[p_{i,n}]-[q_{i,n}], \,\, 1\le i\le I(n).$ Note that
$G_{u,n}$ is {{assumed to be}} {{a free abelian}} group generated by  $\{x_{i,n}: 1\le i\le I(n)\}.$ \Wlog, we may assume that these are free generators.

{{Let $e_{i,n}=\phi_2^\sim\circ \imath_{n+1}^\sim(p_{i,n}),$
$e'_{i,n}=\phi_2^\sim\circ \imath_{n+1}^\sim(q_{i,n}),$
$i=1,2,..., I(n).$  Define $s_1=1$  and ${\td u}_1=u_1s_1^*=u_1.$}}
{{By \eqref{Add-403-1},
$$
\Lambda_1(x_{i,1})=\overline{\langle ((1-e_{i,1})+e_{i,1}{{{\tilde u}_1^*}}u_2)((1-e'_{i,1})+e'_{i,1}u_2^*{{{\tilde u}_1}})\rangle}\in U_0(A)/CU(A).
$$}}
Define a \hm\, $\Lambda_1: G_{u,1}\to U_0(A)/CU(A)$   {{ by $\Lambda_1(x_{i,1})$ as above ($i=1,2,...,I(1)$).}}
%
Since $\Lambda$ factors through $G_{u,n}',$ applying Theorem \ref{TTboteplus}
(with $\af=0$)
to $\phi_2\circ\iota_{2, \infty}$, one obtains  a unitary $s_{2}\in {{CU({\tilde A})}}$ such that
\beq
||[\phi_2\circ\iota_{2, \infty}(f), {{s_{2}}}]||<\delta'_{2}/4\rforal f\in\mathcal G_{2},\\
\mathrm{Bott}(\phi_2\circ\iota_{2, \infty}, s_2)|_{\mathcal P_2}=0, \andeqn
\eneq
\beq\label{n-72-39+}
\hspace{-0.2in}{\rm dist}(\overline{\langle ((1-e_{i,2})+e_{i,2}{{s_{2}^*}})((1-e'_{i,2})+
e'_{i,2}s_{2})\rangle},
\Lambda_{1}(-x_{i,2}))<\sigma_2'/16.
\eneq
Define $\td u_2=u_2s_2^*.$
%
%
In what follows, we will construct unitaries $s_2, ..., {s_n,...}$ in {{$CU({\tilde A})$}} such that
\beq\label{n-72-37}
||[\phi_2\circ\iota_{j+1, \infty}(f), s_{j+1}]||<\delta'_{n+1}/4\rforal f\in\mathcal G_{j+1},\\
\label{n-72-38}
\mathrm{Bott}(\phi_2\circ\iota_{j+1, \infty}, s_{j+1})|_{\mathcal P_j}=0, \andeqn\\
\label{n-72-39}
\hspace{-0.2in}{\rm dist}(\overline{\langle ((1-e_{i,j})+{{e_{i,j}}}s_j)((1-e'_{i,j})+e'_{i,j}s_j^*)\rangle},
\Lambda_j(-x_{i,j}))<\sigma_j'/16,
\eneq
where $\Lambda_j: G_{u,j}\to U_0(A)/CU(A)$  is a \hm\, {{given}}  by
\beq\label{2019-sept-10-2}\Lambda_n(x_{i,n})=(\overline{\langle ((1-e_{i,j})+{{e_{i,j}}}{{\td u_j^*u_{j+1}}})((1-{{e'_{i,j}}})+e'_{i,j}{{u_{j+1}^*\td u_j}})\rangle}\eneq
{{(see \eqref{Add-403-1}) and $\widetilde{u}_j=u_js_j^*,$  $j=1,2,....$}}

%
%
{{Assume}} that $s_2, s_3, ..., s_n$ are already constructed. Let us construct $s_{n+1}$. Note that by \eqref{72-38}--\eqref{Add-403-1}, the $K_1$ class of the unitary $u_n^*u_{n+1}$ is trivial. In particular, the $K_1$ class of $s_nu_n^*u_{n+1}$ is trivial.    {{By \eqref{Add-403-1} and \eqref{n-72-38}, $\Lambda_n(x_{i,n})\in U_0(A)/CU(A).$}}
 {Since $\Lambda$ factors through $G_{u,n}',$ applying Theorem \ref{TTboteplus} (with $\af=0$)
 to $\phi_2\circ\imath_{{n+1}}$, one obtains  a unitary {{$s_{n+1}\in CU({\tilde A})$}} such that
 (see also \ref{Djcubt})
\beq
&&||[\phi_2\circ\iota_{{n+1}}(f), s_n]||<\delta'_{n+1}/4\rforal f\in\mathcal G_{n+2},\\
&&\mathrm{Bott}(\phi_2\circ\iota_{n+1}, s_{n+1})|_{\mathcal P_n}=0\andeqn\\
\label{n-72-39+}
&&\hspace{-0.92in}{\rm dist}(\overline{\langle ((1-e_{i,n})+e_{i,n}s_{n+1})((1-e'_{i,n})+e'_{i,n}s_{n+1}^*)\rangle},
\Lambda_{{{n+1}}}(-x_{i,n}))<\sigma_n'/16,
\eneq
$i=1,2,...,I(n+1).$
Then  $s_1, s_2, ..., s_{n+1}$ satisfy \eqref{n-72-37}, \eqref{n-72-38} and \eqref{n-72-39}.

Put $\widetilde{u_n}=u_ns_n^*\in {{U_0({\tilde A})}}.$ Then by \eqref{72-37} and \eqref{n-72-37}, one has
\beq\label{n-72-40}
{{\rm Ad}}\, \widetilde{u_n}\circ \phi_1\approx_{\dt_n'} \phi_2\rforal a\in
\imath_{n+1}({\cal G}_{n+1}).
\eneq
By \eqref{72-38} -- \eqref{Add-403-1}  and \eqref{n-72-38}, one has
\begin{equation}\label{n-72-41}
\mathrm{Bott}(\phi_2\circ\imath_{n+1}, (\widetilde{u_{n}})^* \widetilde{u_{n+1}})|_{\mathcal P_n}=0.
\end{equation}
Note that
\beq\label{72-n2014103}
\overline{\langle (1-e_{i,n})+e_{i,n}\widetilde{u_n}^* \widetilde{u_{n+1}}\rangle \langle(1-e_{i,n}')+e_{i,n}'{\widetilde{u_{n+1}}^*}{\widetilde{u_{n}}}\rangle}
=\overline{c_1c_2c_4c_3}=
\overline{c_1c_3c_2c_4},
\eneq
where
{{\beq
c_1=\langle(1-e_{i,n})+e_{i,n}{\widetilde{u_n}}^*u_{n+1}\rangle,
\,\, c_2=\langle (1-e_{i,n})+e_{i,n}s_{n+1}^*\rangle.\\
c_3=\langle (1-e'_{i,n})+{e_{i,n}'}u_{n+1}^*\widetilde{u_n}\rangle,\,\,
c_4=\langle (1-e'_{i,n})+e'_{i,n}s_{n+1}\rangle.
\eneq}}
Therefore,  by \eqref{n-72-39+} {{and by \eqref{2019-sept-10-2},}}
one has
{{\beq\label{n-72-42}
&&\hspace{-0.2in}\mathrm{dist}(\overline{\langle ((1-e_{i,n})+e_{i,n}{\widetilde{u_n}}^*{\widetilde{u_{n+1}}})((1-e_{i,n}')+e_{i,n}'{\widetilde{u_{n+1}}^*{\widetilde{u_n}})}\rangle}), {\bar 1})\\
&&<{\rm dist}(\Lambda_n(x_{i,n})\Lambda_n(-x_{i,n}), {\bar 1})+\sigma_n'/16
={{\sigma_{n}'/16,}}
\eneq}}
$i=1,2,...,I(n).$
%
%
%
Therefore, by 
{{14.8 of \cite{GLII},}}
there exists a piece-wise smooth and
continuous path of unitaries $\{z_n(t): t\in [0,1]\}$ of $A$ such
that
\beq\label{72-39}
&&z_n(0)=1,\,\,\, z_n(1)=(\widetilde{u_n})^*\widetilde{u_{n+1}}\andeqn\\\label{72-40}
&&\|[\phi_2(a),\, z_n(t)]\|<1/2^{n+2}\rforal a\in {\cal F}_n\andeqn
t\in [0,1].
\eneq
Define
$$
u(t+n-1)=\widetilde{u_n}z_{n+1}(t)\,\,\,t\in (0,1].
$$
Note that $u(n)=\widetilde{u_{n+1}}$ for all integer $n$ and $\{u(t):t\in [0,
\infty)\}$ is a continuous path of unitaries in ${{U_0({\tilde A})}}.$ One estimates
that, by (\ref{72-37}) and (\ref{72-40}),
\beq\label{72-41}
{{\rm Ad}}\, u(t+n-1)\circ \phi_1 \approx_{\dt_n'} {{\rm Ad}}\,z_{n+1}(t)\circ \phi_2
 \approx_{1/2^{n+2}} \phi_2
\,\,\,\,\,\,\,\text{on}\,\,\,\,{\cal
F}_n
\eneq
 for all $t\in (0,1).$
It then follows that
\beq\label{72-42}
\lim_{t\to\infty}u^*(t)\phi_1(a) u(t)=\phi_2(a)\rforal a\in C.
\eneq
}
\end{proof}

\section{KK-Lifting and rotation maps}

Let {{us}} begin this section with  Theorem \ref{Tcbms748} which is the same as
Theorem 3.17 of \cite{LN}. But  we will use the form of
Theorem 7.4.8 of \cite{Lncbms}.
To state it, we will refer {{to the}}  property (B2) {{in Definition \ref{DB2}.}}

Let $A$ be a \CA\, and $B
\subset A$ be a \SCA, and {{$\af: B\to A$}} a \hm. We write $\af\in \overline{\rm Inn}(B,A),$
if there exists a sequence of unitaries $u_n\in {\tilde A}$ such that
$\af(b)=\lim_{n\to\infty}u_n^*bu_n$ for all $b\in {{B}}$ (converges in norm).

{{In general, let $\iota: B\to A$ be the embedding and $\af:B\to A$ be a monomorphism.
Denote $M_\af:=M_{\iota, \af},$ the mapping torus.  Suppose that $\af_{*i}=\iota_{*i},$ $i=0,1.$
Let $E_i:=K_i(M_\af),$  $i=0,1.$   Then we have group extensions:
\beq
0\to K_i(SA)\to E_i\to K_i(B)\to 0\,\,\, (i=0,1).
\eneq
Write  $\eta_0(M_\af):=E_0\in {{{\rm ext}(K_0(B),K_1(A))}}$ and $\eta_1(M_\af):=E_1\in {{{\rm ext}(K_1(B),K_0(A))}}.$

Let   $G_1$ and $G_2$ be abelian groups. Denote by {{${\rm Pext}(G_2, G_1)$}}
the (equivalence classes) of extensions  $0\to G_1\to E\to G_2\to 0$ such that
every finitely generated subgroup of $G_2$ splits.
}}

\begin{df}\label{DmathttP}
Let $0<\dt^{\mathtt p}<1/4$ be the constant described in  Section 3 of \cite{LN} before 3.1.
Let $A$ be a unital \CA, $l\in  \N,$ $U(t)\in C([0,1], M_l(A))$ be a continuous path of unitaries, and $z\in U(M_l(A))$
be a unitary such that $U(0)=1,$ $\|[U(1),\, z]\|<\dt^{\mathtt p}.$ So ${\rm bott}_1(U(1),z)$ is well defined.
Define a loop of unitaries
${\mathtt p}(U, z)(t)$  just as in the first paragraph of Section 3  of \cite{LN}.

Denote by $K_1^{{\mathtt P} -}(SA)$ the subset  consisting of $[{\mathtt p}(U,z)(t)]$ as described above (see
3.4 of \cite{LN}).  Define a map ${\mathtt P}: [{\mathtt p}(U,z)]\to {\rm bott}_1(U(1),z).$
Then, by Lemma 3.3 and 3.4 of \cite{LN} (see also Remark 3.3.4, Lemma 3.3.5 and 3.36 of \cite{Lncbms}),
${\mathtt P}$ is an injective \hm.  In fact, by Theorem 3.3.7 and 3.3.10 of \cite{Lncbms} that
$K_1^{{\mathtt P} -}(SA)=K_1(SA)\cong K_0(A)$ and ${\mathtt P}$ is an isomorphism.

\end{df}

\begin{thm}[Theorem 3.17 of \cite{LN} and Theorem 6.4.8 of \cite{Lncbms}]\label{Tcbms748}
Let $A$ be an
infinite dimensional simple \CA\,  and {{$B$}} a separable \SCA.
Suppose that $A$ has  {{the}}
property (B2) associated with $B$ and  certain $\DT_B$ here {{(see Definition \ref{DB2})}}
and
$K_0(A)$ is tracially approximately divisible (see 3.15 of \cite{LN}).
For any $E_0\in {{\mathrm{Pext}(K_0(B),K_1(A))}}$ and
$E_1\in {{\mathrm{Pext}(K_1(B),K_0(A))}}$, there exists
$\alpha\in\overline{\mathrm{Inn}}(B,A)$ such that
$\eta_0(M_\alpha)=E_0$ and $\eta_1(M_\alpha)=E_1,$
or equivalently,
$$
[\af]-[\imath]=(E_0, E_1)\in{\mathrm{Pext}}(K_1(A), K_0(B)){\oplus} {\mathrm{Pext}}(K_0(A), K_1(B)),
$$
where $\imath: B\to A$ is the {inclusion}.
\end{thm}

\begin{proof}
Write
\beq\nonumber
&&0\to K_1(A)\to E_0\,{\stackrel{\pi^{(0)}}{\longrightarrow} }\, K_0(B)\to 0\andeqn\\
&&0\to K_1(SA)\to E_1\,{\stackrel{\pi^{(1)}}{\longrightarrow} }\, K_1(B)\to 0
\eneq
for the two given extensions of abelian groups.
Write  (with $p_1=1_{\td B}$)
\beq\nonumber
K_0(\td B)_+=\{[p_1], [p_2],  ... , [p_n], ... \}\andeqn
 K_1(B)=\{[z_1], [z_2],  ... , [z_n],  ... \}.
\eneq
Put ${\cal P}^{(0)}=\{\{[p_1]-k_1[1_{\td B}], [p_2]-k_2[1_{\td B}],...,[p_n]-k_n[1_{\td B}],...\},$
where $k_i=[\pi_\C^B(p_i)]$ and  $\pi_\C^B: \widetilde B\to \C$
is the quotient map. Let $y_i:=[p_i]-k_i[1_{\td B}],$ $i\in \N.$

Denote by $G_n^{(0)}$ the subgroup of $K_0(B)$ generated by $\{y_1, y_2,...,y_n\},$
$\td G_n^{(0)}$  {{the subgroup}} of $K_0(\td B)$ generated by $\{[p_1], [p_2],...,[p_n]\}$ and
by $G_n^{(1)}$ the subgroup of $K_1(B)$ generated by $\{[z_1], [z_2],...,[z_n]\}.$
Denote also by $\iota_{n, n+1}^i: G_n^i\to G_{n+1}^i$ the embedding ($i=0,1$).
We also write  $\td \iota_{n,n+1}^0: \td G_n^{(0)}\to \td G_{n+1}^{(0)}$ for the extension
(and $\td \iota_{n, n+1}^1=\iota_{n, n+1}$).
Let $\{x_1, x_2,...,x_n,...\}\subset B$ be a dense sequence in the unit ball of $B.$
Let $\{{\cal F}_i\}$ be an increasing family of finite subsets in the unit ball of $B$
such that
$
\{x_1, x_2,...,x_n\}{\subset} {\cal F}_n
$
and, assume, for
for each $i$, $p_i\in\mathrm{M}_{r_i}({\cal F}_i+\C1_{\td B})$ is a projection
and $z_i\in\mathrm{M}_{r_i}({\cal F}_i{{+\C1_{\td B}}})$ is a unitary.
We may assume that $r_i\le r_{i+1},$ $i\in\mathbb N.$
Denote by $\td {\cal F}_n=\{1_{\td B}+x: x\in {\cal F}_n\}\cup {\cal F}_n.$
In what follows, if $v\in A,$ by $v^{(m)},$ we mean
$v^{(m)}=v\otimes 1_{M_m}.$

We claim that there are unitaries $\{u_n\}$  with $[u_n]=0$ in $K_1(A)$ and diagrams, for $i=0,1,$
(where $A^0=A$ and $A^1=SA$),

\vskip 3mm
\begin{xy}
(0,0)*{~}="0";
(20,45)*{0}="1uu";
(40,45)*{K_1(A^i)}="2uu";
(78,45)*{E_i}="4uu";
(110,45)*{K_i(B)}="6uu";
(130, 45)*{0}="7uu";
(20,15)*{0}="1";
(40,15)*{K_1(A^i)}="2";
(70,15)*{K_1(A^i)}="3";
(78,15)*{\oplus}="4";
(86,15)*{G_{n+1}^{(i)}}="5";
(110,15)*{G_{n+1}^{(i)}}="6";
(130, 15)*{0}="7";
(20,0)*{0}="1l";
(40,0)*{K_1(A^i)}="2l";
(70,0)*{K_1(A^i)}="3l";
(78,0)*{\oplus}="4l";
(86,0)*{G_n^{(i)}}="5l";
(110,0)*{G_n^{(i)}}="6l";
(130, 0)*{0}="7l";
(40,-15)*{\vdots}="2ll";
(70,-15)*{\vdots}="3ll";
(86,-15)*{\vdots}="5ll";
(110,-15)*{\vdots}="6ll";
(40,30)*{\vdots}="2u";
(70,30)*{\vdots}="3u";
(78,30)*{~}="4u";
(86,30)*{\vdots}="5u";
(110,30)*{\vdots}="6u";
{\ar "1";"2"};
{\ar "2";"3"};
{\ar^{\pi_n^{(i)}}@<.5ex> "5";"6"};
{\ar^{\theta_{n+1}^i}@<.5ex> "6";"5"};
{\ar "6";"7"};
{\ar "1uu";"2uu"};
{\ar "2uu";"4uu"};
{\ar^{\pi^{(i)}} "4uu";"6uu"};
{\ar "6uu";"7uu"};
{\ar "1l";"2l"};
{\ar "2l";"3l"};
{\ar^{\pi_n^{(i)}}@<.5ex> "5l";"6l"};
{\ar^{\theta_n^i}@<.5ex> "6l"; "5l"};
{\ar "6l";"7l"};
{\ar@{=} "2l";"2"};
{\ar@{=} "3l";"3"};
{\ar_{\iota_{n, n+1}^i} "5l";"5"};
{\ar_{\iota_{n, n+1}^i} "6l";"6"};
{\ar_{\gamma^{i}_n} "5l";"3"};
{\ar@{=} "2u";"2uu"};
{\ar "4u";"4uu"};
{\ar "6u";"6uu"};
{\ar@{=} "2ll";"2l"};
{\ar@{=} "3ll";"3l"};
{\ar "5ll";"5l"};
{\ar "6ll";"6l"};
{\ar "5ll"; "3l"};
{\ar@{=} "2";"2u"};
{\ar@{=} "3";"3u"};
{\ar "5";"5u"};
{\ar "6";"6u"};
{\ar "5"; "3u"};
\end{xy}
\vskip 3mm
%
%
\noindent such that
$$
\|[u_{n+1},\, a]\|\leq \dt_{n+1}/r_n^2
$$
for any $a\in\mathrm{M}_{r_n}(w_n^*{\cal F}_nw_n)$, where $w_n=u_1\cdots u_n$ and $u_1=1$
and where $\dt_n$ is as chosen in  Lemma 3.13 of \cite{LN} (see Lemma 6.45 of \cite{Lncbms})
corresponding to
$\{\td {\cal F}_n\},$  $\td G_n^{(0)}$ and $G_n^{(1)}$ as well as $\td A$ and $\td B.$  Note that we assume that
$\dt_n\le {\dt^{\mathtt p}\over{2^{n+4}}}.$  Moreover,
%
%
$$
\mathrm{bott}_1({{(w^*_n)^{(r_i)}}}z_iw^{(r_i)}_n, u_{n+1}^{(r_i)})=
{\mathtt P}\circ
\gamma^1_n([z_i])\quad\textrm{and}\quad\mathrm{bott}_0({{(w^*_n)^{(r_i)}}}p_iw_n^{(r_i)}, u_{n+1}^{(r_i)})=\gamma^0_n([p_i]),
$$
where ${\mathtt P}$ is defined in \ref{DmathttP}.
Note that, if $q\in M_{r_i}(\C)$ is a scalar projection of rank $R_q,$
there exists a scalar unitary $Z\in M_{r_i}(\C)$ such that
$Z^*qZ=\diag(1,1,..,1,0,...,0),$ where $1$ repeats $R_q$ times.
Since $u_n\in U(\td A),$   then $Zu_n^{(r_i)}Z^*=u_n^{(r_i)}.$ Therefore
$$[{{(w^*_{n-1})^{(r_i)}}} q w_{n-1}^{(r_i)}u_n^{(r_i)}{{(w^*_{n-1})^{(r_i)}}} q w_{n-1}^{(r_i)}]=[qu_nq]=[Z^*qZuZ^*qZ]=[u^{(R_q)}\oplus 1_{r_i-R_q}].$$
It follows that, if $[u_n]=0$ in $K_1(A),$  then
\beq\label{pi=xi}
{\rm bott}_0({\rm Ad}(w_{n-1}), u_n)(x_i)={\rm bott}_0({{(w^*_{n-1})^{(r_i)}}}p_iw_{n-1}^{(r_i)}, u_n).
\eneq

As $A$ has property  (B2) associated with $B$ and   $\DT_B$ as in Definition \ref{DB2}, let $G_i^n$ (in place of $G_i,$ $i=0, 1$)
 and
 ${\cal Q}_n\subset G_1^n$ (in place of ${\cal Q}$)  given in Definition \ref{DB2} with respect to $\dt_{n+1}/(2r_{n+1}^2)$(in place of $\ep$), ${\cal F}_n$ (in place of ${\cal F}$), ${\cal P}_n^{(0)}:=\{y_1, y_2,...,y_n\}, {\cal P}_n^{(1)}:=\{[z_1], [z_2],...,[z_n]\}$ (in place of ${\cal P}_0, {\cal P}_1$), and $\imath$ (in place of $h$).
  Each partial splitting map $\theta_n^{i}$ ($i=0, 1$)  can be extended to a partial splitting map $\tilde{\theta}_n^i$ ($i=0, 1$) defined on the subgroup generated by $G_n^{(0)}\cup \mathcal G_0^n$ or by $G_n^{(1)}\cup \mathcal G_1^n\cup {{{\cal Q}_n}}$, where ${\cal  G}_i^n$  is the set of generators of {{$G_i^n$}}. Denote {the} subgroups generated by {{$G_n^{(0)}\cup \mathcal G_0^n$}} 
and  {{$G_n^{(1)}\cup \mathcal G_1^n\cup {{{\cal Q}_n}}$}} 
by
$\bar{G}_n^{(0)}$ and $\bar{G}_n^{(1)}$ respectively.

Assume that we have constructed the unitaries $\{u_1=1, u_2,  ... , u_{n}\}$ with $[u_j]=0$ in $K_1(A)$
($1\le j\le n$) and the diagrams
\vskip 3mm
\begin{xy}
(0,0)*{~}="0";
(20,30)*{0}="1uu";
(40,30)*{K_1(A^i)}="2uu";
(78,30)*{E_i}="4uu";
(110,30)*{K_i(B)}="6uu";
(130, 30)*{0}="7uu";
(20,15)*{0}="1";
(40,15)*{K_1(A^i)}="2";
(70,15)*{K_1(A^i)}="3";
(78,15)*{\oplus}="4";
(86,15)*{G_n^{(i)}}="5";
(110,15)*{G_n^{(i)}}="6";
(130, 15)*{0}="7";
(20,0)*{0}="1l";
(40,0)*{K_1(A^i)}="2l";
(70,0)*{K_1(A^i)}="3l";
(78,0)*{\oplus}="4l";
(86,0)*{G_{n-1}^{(i)}}="5l";
(110,0)*{G_{n-1}^{(i)}}="6l";
(130, 0)*{0}="7l";
(40,-15)*{\vdots}="2ll";
(70,-15)*{\vdots}="3ll";
(86,-15)*{\vdots}="5ll";
(110,-15)*{\vdots}="6ll";
{\ar "1";"2"};
{\ar "2";"3"};
{\ar^{\pi_n^{(i)}}@<.5ex> "5";"6"};
{\ar^{\theta_{n}^i}@<.5ex> "6";"5"};
{\ar "6";"7"};
{\ar "1uu";"2uu"};
{\ar "2uu";"4uu"};
{\ar^{\pi^{(i)}} "4uu";"6uu"};
{\ar "6uu";"7uu"};
{\ar "1l";"2l"};
{\ar "2l";"3l"};
{\ar^{\pi_{n-1}^{(i)}}@<.5ex> "5l";"6l"};
{\ar^{\theta_{n-1}^i}@<.5ex> "6l";"5l"};
{\ar "6l";"7l"};
{\ar@{=} "2l";"2"};
{\ar@{=} "3l";"3"};
{\ar_{\iota_{n-1, n}} "5l";"5"};
{\ar_{\iota_{n-1, n}} "6l";"6"};
{\ar_{\gamma^{i}_{n-1}} "5l";"3"};
{\ar@{=} "2";"2uu"};
{\ar "4";"4uu"};
{\ar "6";"6uu"};
{\ar@{=} "2ll";"2l"};
{\ar@{=} "3ll";"3l"};
{\ar "5ll";"5l"};
{\ar "6ll";"6l"};
{\ar "5ll"; "3l"};
\end{xy}
\vskip 3mm
%
%
\noindent satisfying the above claim.

We note that $[p_i]=[(w_n^{(r_n)})^*p_iw_n^{(r_n)}]$ and $[z_i]=[(w_n^{(r_n)})^*z_iw_n^{(r_n)}],$ $i=1,2, ... ${.}

Choose (since $A$ has  property (B2) associated with $B$ and $\Delta_B$)
$$
\sigma_n=\Delta_B({\dt_{n+1}\over{2 r_{n+1}^2}}, {\cal F}_{n}, {\cal P}_{n}^{(0)}, {\cal P}_{n}^{(1)}, {\rm Ad}\, w_n\circ \imath).
$$


Since $E_0$ and $E_1$ are pure extensions, there are partial splitting maps $\bar{\theta}_{n+1}^0:\bar{G}_{n+1}^{(0)}\to E_0$ and $\bar{\theta}_{n+1}^1: \bar{G}_{n+1}^{(1)}\to E_1$. Since
$K_0(A)$ is tracially approximately divisible, by
Lemma 3.16 of \cite{LN},
the partial splitting map $\bar{\theta}^1_{n+1}$
can be chosen {{such that,}} for any $g\in\mathcal G_1^{n}\cup\{[z_1],...,[z_n]\}\cup Q$,
$$
|\rho_A((\bt^{(0)})^{-1}\circ{\gamma}^1_{n}(g))(\tau)|<
\sigma_n
\quad \tforal \tau\in T(A),
$$
where  $\gamma_n^i={\bar \theta_{n+1}^i}|_{G_n^{(i)}}-\bar \theta_n^i$ ($i=0,1$) and
$\bt^{(0)}:K_0(A)\to K_1(SA)$ is  the Bott isomorphism (specified, for example, as in {{2.1.21 of}} \cite{Lncbms}).
Note that
${\mathtt{P}}$ is defined on $K_1(SA).$
By Lemma 3.5 of \cite{LN} (see also 3.7 of \cite{HL}).
$$
|\tau({\mathtt{P}}(\gamma^1_{n})(g))|<\sigma_n\quad \tforal \tau\in T(A).
$$
Note also that the map  $\gamma_n^i$ is   defined on ${\bar G}_n^{(i)},$ $i=0,1.$

Put $b_0= \gamma_n^{0}$ and $b_1={\mathtt P}\circ \gamma_n^1.$
Using the property (B2),
one obtains
a unitary $u_{n+1}\in A$  with $[u_{n+1}]=0$ in $K_1(A)$ such that
$$
\|[u_{n+1}, a]\|\leq\frac{\dt_{n+1}}{r_{n+1}^2}
$$
for {{all}} $a\in\mathrm{M}_{r_n}(w^*_n\mathcal F_nw_n)$ and  ($1\le i\le n$)
$$
\mathrm{bott}_1({\rm Ad}\, w_n, u_{n+1})([z_i])=
{\mathtt P}\circ
\gamma^1_{n}([z_i])\quad\textrm{and}\quad\mathrm{bott}_0({\rm Ad}\, w_n, u_{n+1})(y_i)
=\gamma^0_{n}(y_i).
$$
Denote by $\theta_{n+1}^i$  {{the restriction}} of $\tilde{\theta}_{n+1}^i$ to {{$G_{n+1}^{(i)},\, i=0,1.$}}
Repeating this procedure, one obtains a sequence of unitaries $\{u_n\}$ and diagrams satisfying the claim.

We extend $\td \gamma_n^{(0)}: \td G_n\to K_1(S\td A)$ by $\td \gamma_n^{(0)}([1_{\td B}])=0.$
By Lemma 3.13 of \cite{LN},
the inner automorphisms
$\{\mathrm{Ad}(u_1\cdots u_n)\}$ converge on $\td B$ to a monomorphism $\td \alpha$
(and on $B$ to $\af$),
and the extension $\eta_0(M_{\td \alpha})$ and $\eta_1(M_{\td \alpha})$ are
determined by the inductive limits of
\vskip 3mm
\xy
(0,0)*{~}="0";
(20,15)*{0}="1";
(40,15)*{K_1(A)}="2";
(70,15)*{K_1(A)}="3";
(78,15)*{\oplus}="4";
(86,15)*{\td G_{n+1}^{(0)}}="5";
(110,15)*{\td G_{n+1}^{(0)}}="6";
(130, 15)*{0}="7";
(20,0)*{0}="1l";
(40,0)*{K_1(A)}="2l";
(70,0)*{K_1(A)}="3l";
(78,0)*{\oplus}="4l";
(86,0)*{\td G_n^{(0)}}="5l";
(110,0)*{\td G_n^{(0)}}="6l";
(130, 0)*{0}="7l";
{\ar "1";"2"};
{\ar "2";"3"};
{\ar "5";"6"};
{\ar "6";"7"};
{\ar "1l";"2l"};
{\ar "2l";"3l"};
{\ar "5l";"6l"};
{\ar "6l";"7l"};
{\ar@{=} "2l";"2"};
{\ar@{=} "3l";"3"};
{\ar_{\td \iota_{n, n+1}^0} "5l";"5"};
{\ar_{\td \iota_{n, n+1}^0} "6l";"6"};
{\ar_{\tilde{\gamma}^{0}_n} "5l";"3"};
\endxy
\vskip 3mm
\noindent and
\vskip 3mm
\xy
(0,0)*{~}="0";
(20,15)*{0}="1";
(40,15)*{K_1(SA)}="2";
(70,15)*{K_1(SA)}="3";
(78,15)*{\oplus}="4";
(86,15)*{G_{n+1}^{(1)}}="5";
(110,15)*{G_{n+1}^{(1)}}="6";
(130, 15)*{0}="7";
(20,0)*{0}="1l";
(40,0)*{K_1(SA)}="2l";
(70,0)*{K_1(SA)}="3l";
(78,0)*{\oplus}="4l";
(86,0)*{G_{n}^{(1)}}="5l";
(110,0)*{G_{n}^{(1)}}="6l";
(130, 0)*{0}="7l";
{\ar "1";"2"};
{\ar "2";"3"};
{\ar "5";"6"};
{\ar "6";"7"};
{\ar "1l";"2l"};
{\ar "2l";"3l"};
{\ar "5l";"6l"};
{\ar "6l";"7l"};
{\ar@{=} "2l";"2"};
{\ar@{=} "3l";"3"};
{\ar_{\iota_{n, n+1}^1} "5l";"5"};
{\ar_{\iota_{n, n+1}^1} "6l";"6"};
{\ar_{\tilde{\gamma}^{1}_n} "5l";"3"};
\endxy
\vskip 3mm
%
\noindent respectively, where (for $1\le i\le n$)
\beq\nonumber
\tilde{\gamma}_n^{0}([p_i])=\tilde{\gamma}_n^{0}([p_i])=[{{((w^*_n)^{(r_i)}p_i
w_n^{(r_i)})}}u_{n+1}({{(w^*_n)^{(r_i)}}}p_iw_n^{(r_i)})+(1-{{(w^*_n)^{(r_i)}}}p_iw_n^{(r_i)})]\andeqn\\\nonumber
\tilde{\gamma}_n^{1}([z_i])=\tilde{\gamma}_n^{1}([z_i])=[{\mathtt{p}}(R^*(u_{n+1}, t), w_n^*z_iw_n)]:=
[{\mathtt{p}}(R^*(u_{n+1}^{(r_i)}, t), \diag({{(w^*_n)^{(r_i)}}}z_iw_n^{(r_i)}), 1_{r_i})]
\eneq
(see \eqref{DRU(t)} and \eqref{DmathttP}).
Thus $\tilde{\gamma}_n^{0}([p_i])=\mathrm{bott}_0({{(w^*_n)^{(r_i)}}}p_iw_n^{(r_i)}, u_{n+1})$ ($1\le i\le n$).
It follows that, for $1\le i\le n,$  $\tilde{\gamma}_n^{0}(y_i) =\gamma_n^{(0)}(y_i)$
(see \eqref{pi=xi}).   Moreover,  for $1\le i\le n,$
$$
{\mathtt P}\circ
\tilde{\gamma}_n^1([z_i])=\mathrm{bott}_1({{(w^*_n)^{(r_i)}}}z_iw_n^{(r_i)}, R(u_{n+1}^{(r_i)}, 1))=\mathrm{bott}_1({{(w^*_n)^{(r_i)}}}z_iw_n^{(r_i)}, u_{n+1}^{(r_i)})={\mathtt P}\circ\gamma_n^1([z_i]),
$$
Since ${\mathtt P}$ is injective, we have that $\tilde{\gamma}_n^1={\gamma}_n^1$. Hence, one has that $\eta_0(M_\alpha)=E_0$ and $\eta_1(M_\alpha)=E_1$, as desired.

\end{proof}

%
%

\begin{cor}\label{L86}
Let $B$ and $A_1$ be separable amenable simple \CA s with continuous {{scales,}}
let $C=B\otimes U_1\in {\cal D},$  {{$A=A_1\otimes U_2\in {\cal D},$}} where $U_1$ and $U_2$ are
UHF-algebras of infinite type. Suppose that $B$ satisfies the UCT {{and}}  $\kappa\in
{{KK(C,A)}},$
$\gamma: T(A)\to T(C)$ is a continuous affine map, and $\af: U({{\td C}})/CU({{\td C}})\to
U({{\td A}})/CU({{\td A}})$ is a continuous \hm\, for which $\gamma,\, \af,$ and $\kappa$ are compatible {{(see Definition \ref{Dcompatible})}}.
Then, there exists a monomorphism $h: C\to A$ such that
\begin{enumerate}
\item $[h]=\kappa$ in $KK(C,A)$,
\item $h_T=\gamma$ and $h^{\dag}=\af.$
\end{enumerate}
\end{cor}

\begin{proof}
The proof follows the same lines as that of Theorem 8.6 of \cite{Lininv}, following
the proof of Theorem 3.17 of \cite{LN}.
First note that, by  Theorem \ref{Misothm} and {{(3) of Remark 4.32 (see Theorem 4.31 also) of \cite{GLrange}}},
 $C$ is isomorphic to a \CA\, in ${\cal M}_1\cap {\cal D}^d$ {{which satisfies the condition of Theorem \ref{Text1} for algebra $A$ there}} {{(see \ref{DdertF} {{and Theorem \ref{L215}}}).}}
Denote by $\overline{\kappa}\in KL(C, A)$
 the image of $\kappa$. It follows from Theorem \ref{Text1}
 that there is a  monomorphism $\phi: C\to A$ such that
$$[\phi]=\overline{\kappa}, \quad \phi^\dag=\alpha,\quad\textrm{and}\quad (\phi)_T=\gamma.$$
Note that it follows from the UCT that (as an element  of $KK(C,A)$)
$$\kappa-[\phi]\in\textrm{Pext}(K_*(C), K_{*+1}(A)).$$
By Lemmas \ref{TTbote},
$A$ has
{{the}} property (B2) associated with $C$
in the sense  of  Definition \ref{DB2}.
Note that
$A$ is approximately divisible.
It follows from   Theorem \ref{Tcbms748}
that  there is a monomorphism $\psi_0: A\to A$ which is approximately inner and such that
$$[\psi_0\circ\phi]-[\phi]=\kappa-[\phi]\quad\textrm{in $KK(C, A)$}.$$ Then the map
$$h:=\psi_0\circ\phi$$
satisfies the requirements of the corollary.
\end{proof}


\begin{lem}\label{L92}
Let {{$A$}}  be a separable \CA\, such that $T(A)$ is a compact.
{Suppose that $B$ is a
separable \CA\,}
and suppose that $\phi,\,\psi: B \to A$ are two  monomorphisms such that
\beq\label{L92-1}
[\phi]=[\psi]\,\,\,{\textrm in}\,\,\,KK(B,A),\,\,\,
\phi_T=\psi_T\tand \phi^{\dag}=\psi^{\dag}.
\eneq
Then
\beq\label{L92-3}
R_{\phi, \psi}{{(K_1(B))}}
\subset {{\overline{\rho_A(K_0( A))}}} 
\eneq

\end{lem}

\begin{proof}
Let $z\in K_1(B)$ be represented by the unitary $u\in {U(M_m({\tilde B}))}$ with $u=1_{M_m({\tilde B})}+x$ for $x\in M_m(B)$ and for some integer $m$.
Then, by {{(\ref{L92-1}),}} 
$$
{(\phi\otimes {\rm id}_{M_m})(u)(\psi\otimes {\rm id}_{M_m})(u)^*\in CU(M_m({\tilde A})).}
$$
Suppose that $\{u(t): t\in [0,1]\}$ is a {{continuous and}} piecewise smooth path in ${M_m(U({\tilde A}))}$ such that
$u(0)={(\phi\otimes {\rm id}_{M_m})}(u)$ and $u(1)=(\psi\otimes{\rm id}_{M_m})(u).$ Without lose of generality, we can assume that $\pi_\C(u(t))=1_{M_m(\C)}$ for all $t\in [0,1]$, where $\pi_\C:M_m({\tilde B}) \to M_m(\C)$ is the quotient map. Put $w(t)={(\psi\otimes {\rm id}_{M_m})(u)^*u(t)}.$  Then $w(0)={(\psi\otimes {\rm id}_{M_m})(u)^*(\phi\otimes {\rm id}_{M_{m}})(u)}\in CU(M_m({\tilde A}))$ and $w(1)=1_{M_m(\tilde A)}.$ Thus,
\beq\label{L92-4}
R_{\td \phi, \td \psi}(z)(\tau)&=&{1\over{2\pi i}}\int_0^1 \tau({du(t)\over{dt}}u^*(t))dt
= {1\over{2\pi i}}\int_0^1 \tau(\psi(u)^*{du(t)\over{dt}}u^*(t)\psi(u))dt\\
&=& {1\over{2\pi i}}\int_0^1 \tau({dw(t)\over{dt}}w^*(t))dt
\eneq
for all $\tau\in T({\tilde A}).$ By 3.1 and 3.2 of \cite{Thomsen},
$R_{\td \phi, \td \psi}(z)\in \overline{\rho_{\td A}(K_0(\td A))}.$ Furthermore ${{R_{\td \phi, \td \psi}(z)(\tau_{\C})=0}}$ as $\pi(u(t))=1_{M_m(\C)}$.
{{Let $0<\ep<1/2.$}}
{{Choose $x\in K_0(\td A)$ such that $\|\rho_{\td A}(x)-R_{\td\phi, \td\psi}(z)\|<\ep.$
Then $|x(\tau_\C)|<\ep<1/2.$  Since $\tau_\C(K_0(\td A))\subset \Z,$
$(\pi_\C)_{*0}(x)=0.$ It follows that $x\in  K_0(A).$
Note that $R_{\phi, \psi}=R_{\td \phi, \td\psi}|_{T(A)}.$}}
Hence $R_{\phi, \psi}(z)\in \overline{\rho_{A}(K_0( A))}$.
It follows that
\beq\nonumber
R_{\phi, \psi}\in {\rm Hom}(K_1(B), \overline{\rho_{A} (K_0( A))}).
\eneq
\end{proof}

\begin{lem}\label{L714}
Let $A$ be a  \CA\, with  $T(A)\not=\emptyset$ and let $H$ be a finitely generated
abelian group. Let $\psi\in {\rm Hom}(H, \overline{\rho_A(K_0(A))}).$ Fix $\{g_1, g_2,...,g_n\}\subset H.$
Then, for any $\ep>0,$ there exists a \hm\, $h: H\to K_0(A)$ such that
\beq
|\psi(g_i)-\rho_A(h(g_i))|<\ep,\,\,\, 1\le i\le n.
\eneq
\end{lem}

\begin{proof}
Let $H_0:=H/{\rm Tor}(H)$ and $q: H\to H_0$ be the quotient map.
Note that ${\rm Tor}(H)\subset {\rm ker}\psi.$ There is a \hm\, $\psi_0: H_0\to \overline{\rho_A(K_0(A))}$
such that $\psi=\psi_0\circ q.$
\Wlog, we may assume that $\{q(g_1), q(g_2),...,q(g_k)\}$ is a set of free generators of $H_0$
(for some $0\le k\le n$).
Since ${\rm im}\psi\subset \overline{\rho_A(K_0(A))},$ there are $a_1, a_2,...,a_k\in K_0(A)$ such that
\beq
|\psi(g_i)-\rho_A(g_i)|<\ep,\,\,\,i=1,2,...,k.
\eneq
Define $h_0: H_0\to K_0(A)$ by $h(q(g_i))=a_i$ ($1\le i\le k$) and $h:=h_0\circ q.$ Lemma follows.
\end{proof}

%



\begin{thm}\hspace{-0.05in}{\rm (Lemma 4.2 of \cite{LN})}\label{rotation-maps}
Let $A$ be an infinite dimensional simple \CA\, {{with $T(A)\not=\emptyset,$}}
$B\subseteq A$ be a  {{\SCA}}\,  and  $\iota$ the inclusion map. Suppose that $A$ has
{{the}} property $\mathrm{(B2)}$ associated with $B$ and certain $\Delta_B.$
For any {{$\psi\in \mathrm{Hom}( K_1(B), \overline{\rho_A(K_0(A))}),$}}
there exists $\alpha\in\overline{\mathrm{Inn}}(B, A)$ such that there are {{\hm s}} $\theta_i:{K}_i(B)\to{K}_i({{M_{\iota, \alpha}}})$ with $\pi_e\circ\theta_i=\mathrm{id}_{{K}_i(B)}$, $i=0, 1$,
where $\pi_e: M_{\iota, \af}\to B$ {{is}} the quotient map,
and  the rotation map $R_{\iota, \alpha}:  K_1(M_{\iota, \alpha})\to\mathrm{Aff}(\mathrm{T}(A))$ is given by
$$R_{\iota, \alpha}(c)=\rho_A(c-\theta_1([\pi_e]_1(c)))+\psi([\pi_e]_1(c)) \tforal c\in K_1(M_{\iota, \alpha}(A)).$$ In other words, $$[\alpha]=[\iota]$$ in ${KK}(B, A)$, and the rotation map $R_{\iota, \alpha}:  K_1(M_{\iota,\alpha})\to\mathrm{Aff}{{(T(A))}}$ is given by
$$R_{\iota, \alpha}(a, b)=\rho_A(a)+\psi(b)$$ for some identification of $ K_1(M_{\iota,\af})$ with $ K_0(A)\oplus K_1(B)$.
\end{thm}

\begin{proof}
The proof is exactly the same as that of  Theorem 4.2 of \cite{LN} (see also Theorem {{7.4.1}} of \cite{Lncbms}).
The first part of the proof  {{is}} similar to the proof of Theorem \ref{Tcbms748}.
Let $\{{\cal F}_i\}$
be an increasing family of finite subsets {of $B$ with dense union}. 
Assume that for each $i$, there is a unitary $z_i$ and a
projection $p_i$ in $\mathrm{M}_{r_i}({\cal F}_i+\C\cdot 1_{\td B})$ for some
natural number $r_i$ such that
\beq\nonumber
K_0(\td B)_+\setminus \{0\}=\{[p_1], [p_2],...,[p_n],\cdots\}
\andeqn K_1(B)=\{[z_1], [z_2],  ... , [z_n],  ... \}
\eneq
(It might be helpful to note that, if $B$ has stable rank one, we may assume that $z_i\in B$).
Put ${\cal P}^{(0)}=\{\{[p_1]-k_1[1_{\td B}], [p_2]-k_2[1_{\td B}],...,[p_n]-k_n[1_{\td B}],...\},$
where $k_i=[\pi_\C^B(p_i)]$ and  $\pi_\C^B: \widetilde B\to \C$
is the quotient map. Put $x_i=[p_i]-k_i[1_{\td B}]$ ($i\in \N$).
Note that ${\cal P}^{(0)}$ generates $K_0(B).$
Denote by $G_n$ the subgroup generated by $\{x_1, x_2,...,x_n\}$
and
by $H_n$ the subgroup generated by $\{[z_1], [z_2],  ... , [z_n]\}.$
Without loss of generality, we may assume that $r_i\le r_{i+1},$
$i\in\mathbb N.$  In what follows, if $v\in A,$ by $v^{(m)},$ we mean
$v^{(m)}=v\otimes 1_{M_m}.$ Let $\dt_n$ be a sequece of positive number such that $\sum_{n=1}^{\infty} \dt_n< \infty$ and $\dt_{n+1} <\dt_n$.

We assert that there are \hm {s} $h_i: {\tilde{H}}_i\to K_0(A)$  ($i=1,2,...$),
where ${\tilde H}_i$ is the subgroup generated by $ H_n\cup \cup_{i=1}^n{\cal Q}^i$
(${\cal Q}^n$ is specified be{l}ow){,}
and unitaries $\{u_i\,{:}\, i\in\mathbb N\}$ in $A$
with $[u_n]=0$ in $K_1(A)$ such that for any $n$, {writing} $w_n=u_1\cdots u_{n-1}$ (assume $u_{-1}=u_0=1$), one has
\begin{enumerate}
\item {{for}} any $x\in\{[z_1], [z_2], ... , [z_n]\}\cup\cup_{i=1}^n\mathcal Q^i_1$,
$$
|\rho_A\circ h_n(x)-\psi(x)|<\frac{\sigma_n}{2^{n+1}},
$$
where
$$
\sigma_n= \Delta_B(\dt_n/2r_n^2, {\mathcal F_n}, {\mathcal P_n^{(0)}}, {\mathcal P_n^{(1)}}, \mathrm{Ad}(w_{n-1})\circ\iota)
$$
for $\mathcal P_n^{(0)}=\{x_1, x_2,...,x_n\}$
and $\mathcal P_n^{(1)}=\{[z_1], [z_2], ... , [z_n]\}$,
and $\mathcal Q^n$ is the finite subset ${\cal Q}$ {in} Definition \ref{DB2}
with respect to ${\dt_n\over{2r_n^2}}$(in place of $\ep$) , ${\mathcal F_n}$ (in place of ${\cal F}$, ${\mathcal P_n^{(0)}}, {\mathcal P_n^{(1)}},$ (in place of ${\cal P}_0, {\cal P}_1$), and  $\iota$ (in place of $h$). 

\item {{for}} any $a\in w_{n-1}^*{\mathcal F_n}w_{n-1}$,
$$
\|[u_{n}, a]\|<{\dt_n\over{2r_n^2}},
$$
where $w_n=u_1 ...  u_n$, and
$$
\mathrm{bott}_1({\rm Ad}(w^{(r_{n})}_{n-1})(z_i), u_{n})
={\mathtt P}\circ\bt^{(0)}\circ \phi_{n}([z_i])\quad\andeqn\quad
\mathrm{bott}_0(
(w^{(r_{n})}_{n-1})^*p_iw^{(r_{n})}_{n-1}, u_{n})=0,
$$
 where $\phi_n:=h_{n+1}|_{{\tilde H_n}}-h_n,$ for any $i=1,  ... , n$.
\end{enumerate}
%
%
%
Note that $\bt^{(0)}: K_0(A)\to K_1(SA)$ is  the Bott map  {{as mentioned  in the proof of \ref{Tcbms748}.}}

If $n=1$, since $\psi\in \text{Hom}(K_1(B), \overline{\rho_A( K_0(A))}),$   by Lemma \ref{L714},
there is a \hm\, $h_1:
{\tilde H}_1\to K_0(A)$ such that, for any $z\in \{[z_1]\}\cup\mathcal Q^1,$
$$
|\psi(z)-\rho_A(h_1(z))|<\frac{\sigma_1}{2^2},
$$
and a \hm\, $h_2: {\tilde H}_2
\to K_0(A)$ such that
$$
|\psi(z)-\rho_A(h_2(z))|<\frac{\sigma_2}{2^3}
$$
for {{all}}
$z\in \{[z_1], [z_2]\}\cup\mathcal Q^2$. For $\phi_1=h_2|_{{\tilde H}_1}-h_1,$ we have that
$|\tau(\phi_1(z)|<\dt_1/2$ for all $z\in \{[z_1]\}\cup\mathcal Q^1$ and
for any $\tau\in\mathrm{T}(A).$
By Lemma 3.5 of \cite{LN} (see also 3.7 of \cite{HL}),
$$|\tau({\mathtt P}\circ\bt^{(0)}(\phi_1(z)))|=|\tau(\phi_1(z))|<\sigma_1/2$$
for {{all}} $\tau\in\mathrm{T}(A)$ and {{all}} $z\in\mathcal \{[z_1]\}\cup Q^1$. 
Since $A$ has Property (B2) associated with $B$
and $\Delta_B$, there is a unitary $u_1\in  {{U(\td A)}}$  with $[u_1]=0$ in $K_1(A)$ such
that
$$\|[u_1, a]\|<{\dt_1\over{2r_1^2}}
\rforal a\in {\cal F}_1,$$ and
\beq\label{u=0}
\mathrm{bott}_1(z_1, u_1)={\mathtt P}\circ\bt^{(0)}\circ \phi_1([z_1])\quad\textrm{and}\quad \mathrm{bott}_0(p_1, u_1)=0
\,\,\,({\rm recall}\,\, \eqref{pi=xi}).
\eneq
(${\mathtt P}: K_1(SA)\to K_0(A)$ is defined in 3.3.6 of \cite{Lncbms} which is denoted by
$\Lambda$ in 3.4 of \cite{LN}.  By Theorem 3.3.10  of \cite{Lncbms}, ${\mathtt P}$ is an isomorphism).

Assume that we have constructed the \hm s $h_i: {\tilde H}_i\to K_0(A),$ $i=1,2,...,n,$
and unitaries $\{u_i\,{:}\, i=1,  ... , n-1\}$ satisfying the assertion above. By  Lemma \ref{L714},
for the subgroup ${\tilde H}_{n+1}$  of $K_1(B)$ generated by $H_{n+1}\cup\cup_{i=1}^{n+1} {\cal Q}_1^{i},$
there is a function $h_{n+1}: {\tilde H}_{n+1}
\to K_0(A), $ such that for any $z\in\{[z_1], [z_2],  ... , [z_{n+1}]\}\cup\cup_{i=1}^{n+1}\mathcal Q^{n+1}$,
\beq\label{714-n1}
|\rho_A\circ h_{n+1}(z)-\psi(z)|<\frac{\sigma_{n+1}}{2^{n+2}}
\eneq
where
$$
\sigma_{n+1}=\Delta_B( {\dt_{n+1}\over{2r_{n+1}^2}},
{\mathcal F_{n+1}}, {\mathcal P_{n+1}^{(0)}}, {\mathcal P_{n+1}^{(1)}}, \mathrm{Ad}(w_{n})\circ\iota)
$$
for $\mathcal P_{n+1}^{(0)}=\{x_1, x_2, ... , x_{n+1}\}$ and $\mathcal P_{n+1}^{(1)}=\{[z_1], [z_2], ... , [z_{n+1}]\}$. 

Recall that $\phi_n=h_{n+1}|_{{\tilde H}_n}-h_n$. Then, by Lemma 3.5 of \cite{LN} (see also 3.7 of \cite{HL}),
and \eqref{714-n1},
 for any $\tau\in T(A)$, \beq\label{rotationmap-2}
|\tau{(}({\mathtt P}\circ \bt^{(0)}\circ\phi_{n})(x){)}|=|\tau(\phi_n(x))|<\sigma_n/2^n
\eneq
for any $x\in\mathcal  {\td H}_n.$
Since $A$ {{has}} Property (B2) associated with $B$ and $\Delta_B$,  there is a unitary $u_{n}\in U(\td A)$
(with $[u_n]=0$ in $K_1(A)$) such that
$$
\|[u_{n}, w_{n-1}^*aw_{n-1}]\|<{\dt_{n+1}\over{2r_{n+1}^2}}
\rforal a\in {\cal F}_n,
$$
and
\beq\label{rotationmap-3}
&&\mathrm{bott}_1(\mathrm{Ad}(w_{n-1})\circ\iota, u_{n})|_{{\cal P}_{n}^{(1)}}=
{\mathtt P}\circ\bt^{(0)}\circ \phi_n|_{{\cal P}_{n}^{(1)}}\quad\textrm{and}\\\label{bott0-e}
&&\quad \mathrm{bott}_0(\mathrm{Ad}(w_{n-1})\circ\iota, u_{n})|_{{\cal P}_{n}^{(0)}}=0.
\eneq
 This proves the assertion.

 Note that  \eqref{bott0-e}  and \eqref{pi=xi}  imply that (for $1\le i\le n$)
 \beq\label{bott0-ei}
\hspace{-0.25in} {{[\la {{(w^*_{n-1})^{(r_i)}}} p_i w_{n-1}^{(r_i)}u_n^{(r_i)}{{(w^*_{n-1})^{(r_i)}}} p_i w_{n-1}^{(r_i)}+(1- {{(w^*_{n-1})^{(r_i)}}} p_i w_{n-1}^{(r_i)}\ra ]=[1]=0\in
 K_1(S{\td A}).}}
 \eneq

Next, we consider $\td A$ and $\td B.$

By Lemma 3.13 of \cite{LN},
$\textrm{Ad}(w_n)$ converges to a monomorphism $\td \alpha: \td B\to \td A$
and $\af: B\to A.$  Moreover, the extension $\eta_0(M_{\td \alpha})$ is trivial, and $\eta_1(M_{\td \alpha})$
 is determined by the inductive limit of
\vskip 3mm
\xy
(0,0)*{~}="0";
(20,15)*{0}="1";
(40,15)*{K_1(S\td A)}="2";
(70,15)*{K_1(S\td A)}="3";
(78,15)*{\oplus}="4";
(86,15)*{H_{n+1}}="5";
(110,15)*{H_{n+1}}="6";
(130, 15)*{0}="7";
(20,0)*{0}="1l";
(40,0)*{K_1(S\td A)}="2l";
(70,0)*{K_1(S\td A)}="3l";
(78,0)*{\oplus}="4l";
(86,0)*{H_{n}}="5l";
(110,0)*{H_{n}}="6l";
(130, 0)*{0}="7l";
{\ar "1";"2"};
{\ar "2";"3"};
{\ar "5";"6"};
{\ar "6";"7"};
{\ar "1l";"2l"};
{\ar "2l";"3l"};
{\ar "5l";"6l"};
{\ar "6l";"7l"};
{\ar@{=} "2l";"2"};
{\ar@{=} "3l";"3"};
{\ar_{\iota_{n, n+1}} "5l";"5"};
{\ar_{\iota_{n, n+1}} "6l";"6"};
{\ar_{{\gamma}_n} "5l";"3"};
\endxy
\vskip 3mm
%
\noindent where
$\gamma_n([z_i])=[{\mathtt p}(R^*(u_{n}^{r_i}, t), \diag({{(w_{n-1}^{(r_i)})^*}}z_iw_{n-1}^{(r_i)},1_{r_i})]$ (see the notation explanation
at the end of page 178 of \cite{Lncbms}) and $R(u,t)$ is as defined \eqref{DRU(t)}.
The assertion that $\eta_0(M_{\td \af})$ is trivial follows from \eqref{bott0-ei}.
Note that $\pi_\C^A(u_n), \pi_\C^A(w_n)\in \C.$ Thus
$((\pi_\C^A)_{*0}\circ \gamma_n)|_{H_n}=0.$
Therefore  $\eta_0(M_\af)$ is   trivial and $\eta_1(M_\af)$ is determined  by
\vskip 3mm
\xy
(0,0)*{~}="0";
(20,15)*{0}="1";
(40,15)*{K_1(SA)}="2";
(70,15)*{K_1(SA)}="3";
(78,15)*{\oplus}="4";
(86,15)*{H_{n+1}}="5";
(110,15)*{H_{n+1}}="6";
(130, 15)*{0}="7";
(20,0)*{0}="1l";
(40,0)*{K_1(SA)}="2l";
(70,0)*{K_1(SA)}="3l";
(78,0)*{\oplus}="4l";
(86,0)*{H_{n}}="5l";
(110,0)*{H_{n}}="6l";
(130, 0)*{0}="7l";
{\ar "1";"2"};
{\ar "2";"3"};
{\ar "5";"6"};
{\ar "6";"7"};
{\ar "1l";"2l"};
{\ar "2l";"3l"};
{\ar "5l";"6l"};
{\ar "6l";"7l"};
{\ar@{=} "2l";"2"};
{\ar@{=} "3l";"3"};
{\ar_{\iota_{n, n+1}} "5l";"5"};
{\ar_{\iota_{n, n+1}} "6l";"6"};
{\ar_{{\gamma}_n} "5l";"3"};
\endxy

However, since
$$
{\mathtt P}\circ \gamma_n([z_i])=\mathrm{bott}_1(\textrm{Ad}(w_{n-1})\circ\iota, u_{n})([z_i])={\mathtt P}\circ\bt^{(0)}\phi_n([z_i])
$$
and ${\mathtt P}$ is an isomorphism, we have that
\beq\label{rotationmap-1}
\gamma_n=\bt^{(0)}\circ \phi_n,\,\,\, n=1,2,....
\eneq
We assert that $\eta_1(M_\alpha)$ is also trivial.

Write
$$
K_1(M_\af)=\lim(K_1(SA)\oplus H_n, j_n),
$$
$j_n: K_1(SA){\oplus}  H_n\to  K_1(SA){\oplus}  H_{n+1}$ be defined
by
$$j_n(x, y)=(x+ \gamma_n(y), \imath_{n,n+1}(y))$$
as determined by the above
diagram.

For any $n$, define a map $\theta'_n:H_{n}\to K_1(M_{\alpha})$,
for each $g\in H_n,$ {by}
 $$
 \theta'_n(g)=(\bt^{(0)}\circ h_{n}(g), g{)}.
 $$
 We then have, for $g\in H_n,$  since $\phi_n=h_{n+1}|_{H_n}-h_n,$ by \eqref{rotationmap-1},
\begin{eqnarray*}
&&\theta'_{n+1}\circ\iota_{n, n+1}(g)-\theta'_n(g)\\
&=&(\bt^{(0)}\circ h_{n+1}\circ\iota_{n, n+1}(g), \iota_{n, n+1}(g))-(\bt^{(0)}\circ h_n(g)+\bt^{(0)}\circ \phi_n(g)), \iota_{n, n+1}(g))\\
&=&0,
\end{eqnarray*}
and hence $(\theta'_n)$ define{s} a homomorphism $\theta_1: K_1(B) \to K_1(M_\alpha)$. Moreover, since $\pi\circ\theta_1=\mathrm{id}_{K_1(B)}$, the extension $\eta_1(M_\alpha)$ splits. Therefore, $[\alpha]=[\iota]$ in ${KK}(B, A)$.

It remains to calculate the rotation map $R_{\iota, \af}.$  But that calculation  follows exactly
word by word as in the proof of Theorem 4.2 of \cite{LN} starting the last paragraph of page 1761 of the
proof of Theorem 4.2 of \cite{LN} (see also  the last paragraph of p.207  to p.208 of the proof of Theorem 7.41 of \cite{Lncbms}).

\end{proof}

\begin{cor}\label{T94}
Let $C_1, C_2$ be amenable separable simple \CA s, $A=C_1\otimes U_1\in {\cal D}$ $B=C_2\otimes U_2\in {\cal D},$ where $U_1$ and $U_2$ are  UHF-algebras of infinite type, {{and
$B$ satisfies the UCT.}}  Suppose that
$A$ has continuous scale {{and}} $B$ is a \SCA\, of $A,$ and denote by $\imath$ the embedding. For any $\lambda\in
{{{\rm Hom}(K_{{1}}(B), \overline{\rho_A(K_0(A))})}},$
{{there exists $\phi\in {\overline{{\rm Inn}}}(B,A)$}}
such that there are \hm s $\theta_i: K_i(B)\to K_i(M_{\imath, \phi})$ with $(\pi_e)_{*i}\circ \theta_i={\rm id}_{K_i({{B}})},$ $i=0,1,$ and
the rotation map $R_{\imath, \phi}: K_1({{M_{\imath, \phi}}})\to {\rm Aff}(T(A))$ given by
\beq\label{T94-1}
R_{\imath, \phi}(x)=\rho_A(x-\theta_1(\pi_e)_{*1}(x))+\lambda\circ (\pi_e)_{*1}(x))
\eneq
for all $x\in K_1(M_{\imath, \phi}).$ In other words,
\beq\label{T94-2}
[\phi]=[\imath]\,\,\,{\rm in}\,\,\, KK(B,A)
\eneq
and the rotation map $R_{\imath, \phi}: K_1(M_{\imath, \phi})\to
{\rm Aff}(T(A))$ is given by
\beq\label{T94-3}
R_{\imath, \phi}(a,b)=\rho_A(a)+\lambda(b)
\eneq
for some identification of $K_1(M_{\imath, \phi})$ with $K_0(A)\oplus K_1(B).$
\end{cor}

\begin{proof}
By  \ref{Misothm}, $B\in {\cal M}_1\cap {\cal D}^d.$
For each $\ep>0$ and {{finite subsets}}  ${\cal F},$ ${\cal P}_0\subset K_0(B),$
${\cal P}_1\subset K_1(B),$
and ${\cal P}={\cal P}_0\cup {\cal P}_1,$ choose
${{\Delta_B}}(\ep, {\cal F}, {\cal P}_0, {\cal P}_1):=\eta$ {{to}} be given by
Theorem \ref{TTbote}.
Then,   by Theorem \ref{TTbote}, $B$ has {{the}} property ($B_2$) associated with $B$ and
${{\Delta_B}}.$ Therefore, this corollary follows from \ref{rotation-maps}.

\end{proof}

\begin{thm}\label{T96}
Let $C$ and $A$ be two  separable  amenable \CA s such that
$T(C)$ and $T(A)$ are compact.  Suppose that $\phi_1, \phi_2, \phi_3: {C\to A}$
are three  monomorphisms for which
\beq\label{Mul-1}
[\phi_1]=[\phi_2]=[\phi_3]\,\,\,{\textrm in}\,\,\, KK({C,A}))\tand
(\phi_1)_T=(\phi_2)_T=(\phi_3)_T.
\eneq
Then
\beq\label{Mul-2}
\overline{R}_{\phi_1, \phi_2}+\overline{R}_{\phi_2,
\phi_3}=\overline{R}_{\phi_1,\phi_3}.
\eneq
\end{thm}

\begin{proof}
The proof is exactly the same as that of Theorem 9.6 of \cite{Lininv}.
\end{proof}

\begin{lem}\label{L113}
Let $A_1$ and $B$ be two separable simple \CA s  with continuous {{scales}}
which satisfy the UCT,
${{A}}=A_1\otimes U_1$
where $U_1$
is a
UHF-algebras of infinite type.
Let $\phi: A\to B$ be an isomorphism and let $\bt:
B\otimes M_{\mathfrak{p}}\to B\otimes M_{\mathfrak{p}}$ be an
automorphism such that $\bt_{*1}={\rm id}_{K_1(B\otimes
M_{\mathfrak{p}})}$ for some supernatural number $\mathfrak{p}$ of infinite type.
Then
$$\psi^{\dag}(U(\td A)/CU(\td A))=
(\phi_0)^{\dag}(U(\td A)/CU(\td A))=U(\td B)/CU(\td B),
$$
where $\phi_0=\imath\circ \phi,$ $\psi=\bt\circ \imath\circ \phi,$
and where $\imath: B\to B\otimes M_{\mathfrak{p}}$ is defined by
$\imath(b)=b\otimes 1$ for all $b\in B.$ Moreover, there is an
isomorphism $\mu: U(\td B)/CU(\td B)\to U(\td B)/CU(\td B)$ with
$\mu(U_0(\td B)/CU(\td B))\subset U_0(\td B)/CU(\td B)$ such that
$$
\imath^{\dag}\circ \mu\circ \phi^{\dag}=\psi^{\dag}\andeqn
q_1\circ \mu=q_1,
$$
where $q_1: U(\td B)/CU(\td B)\to K_1(B)$ is the quotient map.
\end{lem}

\begin{proof}
We first note that $T({\tilde A})=T(\widetilde{A\otimes M_{\p}})=T({\widetilde{B\otimes M_\p}}).$
Applying the K{\"u}nneth  formula, we compute
that
$$
{{\overline{\rho_A^\sim (K_0({{\td A}}))}=\overline{\{r \cdot x: r\in \R \andeqn x\in \rho_A(j_*(K_0( A_1)))\}}+
{{\Z\cdot 1_{T(A)}}},}}
$$
where $j: A_1\to A_1\otimes U_1=A$ is the embedding $x\mapsto x\otimes 1_{U_1}.$
It follows that $\iota: A\to A\otimes M_\p$ defined by
$\iota(a)=a\otimes 1_{M_\p}$ (for all $a\in A$) induces an isomorphism
$U_0({\tilde A})/CU({\tilde A})\cong U_0({\widetilde{A\otimes M_\p}})/CU({\widetilde{A\otimes M_\p}}).$
We also note that, since $U_0({\tilde A})/CU({\tilde A})\cong \Aff(T({\tilde A}))/\overline{\rho_{\td A} (K_0(\td A))},$
it is a divisible group.
{{With  these facts,}} the proof of this lemma is exactly the same as that of Lemma 11.3 of \cite{Lininv}.
\end{proof}

\begin{lem}\label{L114}
Let $A_1$ and $B_1$ be   separable simple amenable \CA s  which {{satisfy}} the UCT,
Suppose that  $A=A_1\otimes U_1\in {\cal D}$, {{and}}
$B=B_1\otimes U_2\in {\cal D},$ where $U_1$ and $U_2$ are
UHF-algebras {{of infinite type}}, and suppose that
$\phi_1, \phi_2: A\to B$ are two isomorphisms such that
$[\phi_1]=[\phi_2]$ in $KK(A,B).$ Then there exists an
automorphism $\bt: B\to B$ such that $[\bt]=[{\rm id}_B]$ in
$KK(B,B)$ and $\bt\circ \phi_2$ is strongly asymptotically unitarily
equivalent  (see \ref{Dstrongaue}) to $\phi_1.$
\end{lem}

\begin{proof}
It follows from  Corollary
\ref{L86}
 that there is an automorphism $\bt_{1}:
B\to B$ satisfying the following {{conditions}}:
\beq
[\bt_1]=[{\rm id}_B]\,\,\,{\rm in}\,\,\,KK(B,B),\\
\bt_1^{\dag}=\phi_1^{\dag}\circ (\phi_2^{-1})^{\dag}\andeqn
(\bt_1)_T=(\phi_1)_T\circ (\phi_2)_T^{-1}.
\eneq
Let {{$\lambda:=R_{\phi_1, \bt_1\circ \phi_2}\circ (\phi_2)_{*1}^{-1}$.}} {{By Lemma \ref{L92}, $\ld \in {\rm Hom}(K_1(B), \overline{\rho_B (K_0( B))}).$}}
By \ref{T94},
there is {{an}} automorphism $\bt_2\in \mathrm{Aut}(B)$ such that
\beq
[\bt_2]=[{\rm id}_B]\,\,\, {\rm in}\,\,\, KK(B,B),\\
\bt_2^{\dag}={\rm id}_B^{\dag},\,\,\,(\bt_2)_T=({\rm
id}_B)_T, \andeqn\\\label{10-6-1-n2021}
 \overline{R}_{{\rm id}_B,\bt_2}=-\overline{R}_{\phi_1,\bt_1\circ \phi_2}\circ
(\phi_2)_{*1}^{-1}.
\eneq
Put $\bt=\bt_2\circ \bt_1.$
 It follows that
\beq
[\bt\circ \phi_2]=[\phi_1]\,\,\,{\rm in}\,\,\, KK(A,B),
(\bt\circ \phi_2)^{\dag}=\phi_1^{\dag}, \andeqn (\bt\circ
\phi_2)_T=(\phi_1)_T.
\eneq
Moreover, by \ref{T96}, {{the fact ${\bt_1\circ \phi_2}_{*1}={\phi_2}_{*1}$ and \eqref{10-6-1-n2021}}}
\beq
\overline{R}_{\phi_1,\bt\circ \phi_2}&=&{{\overline{R}_{\bt_1\circ \phi_2, \bt\circ \phi_2}+\overline{R}_{\phi_1, \bt_1\circ \phi_2}}}=
\overline{R}_{{\rm
id}_B,\bt_2}\circ
(\phi_2)_{*1}+\overline{R}_{\phi_1,\bt_1\circ\phi_2}\\
&=&(-\overline{R}_{\phi_1,\bt_1\circ \phi_2}\circ
(\phi_2)_{*1}^{-1})\circ (\phi_2)_{*1}
+\overline{R}_{\phi_1,\bt_1\circ\phi_2}=0.
\eneq
It follows from \ref{Tm72} that $\bt\circ \phi_2$ and $\phi_1$ are  strongly
asymptotically unitarily equivalent.

\end{proof}

\section{Isomorphism theorem}

\begin{lem}\label{L115}
Let $A_1$ and  $B_1$ be separable simple amenable \CA s
which satisfy the UCT and have
 continuous scale, and let
$A=A\otimes U_1, B=B_1\otimes U_2\in {\cal D}$
for  UHF-algebras $U_1$ and $U_2$ of infinite type.  Suppose that  there is an isomorphism  $\phi: A\to B$
and
$\bt\in \mathrm{Aut}(B\otimes M_{\mathfrak{p}})$ such that
$$[\bt]=[{\rm id}_{B\otimes
M_{\mathfrak{p}}}]\,\,\,{\textrm in}\,\,\,KK(B\otimes
M_{\mathfrak{p}},B\otimes M_{\mathfrak{p}}) \tand \bt_T=({\rm
id}_{B\otimes M_{\mathfrak{p}}})_T
$$
for some supernatural number $\mathfrak{p}$ of infinite type.

Then there exists an automorphism $\af\in Aut(B)$ with
$[\af]=[{\rm id}_{B}]$ in $KK(B,B)$ such that $ \imath\circ
\af\circ \phi $ and $\bt\circ \imath\circ \phi$ are asymptotically
unitarily equivalent, where $\imath: B\to B\otimes
M_{\mathfrak{p}}$ is defined by $\imath(b)=b\otimes 1$ for all
$b\in B.$
\end{lem}

\begin{proof}
It follows from Lemma \ref{L113} that there is an isomorphism $\mu:
U(\td B)/CU(\td B)\to U(\td B)/CU(\td B)$ such that
$$
\imath^{\dag}\circ \mu\circ \phi^{\dag}=(\bt\circ \imath\circ
\phi)^{\dag}.
$$
Note that $\imath_T: T(B\otimes M_{\mathfrak{p}})\to T(B)$ is an
affine homeomorphism.

It follows from \ref{T94}
that there is an automorphism $\af: B\to
B$ such that
\beq\label{l1-1}
&&[\af]=[{\rm id}_B]\,\,\,{\rm in}\,\,\,KK(B,
B),\\
&&\af^{\dag}=\mu,\,\,\, \af_T=(\bt\circ \imath\circ \phi)_T\circ
((\imath\circ
\phi)_T)^{-1}=({\rm id}_{B\otimes M_{\mathfrak{p}}})_T\andeqn\\\label{10-6-2-n2021}
&&\overline{R}_{{\rm
id}_B,\af}(x)(\tau)=-\overline{R}_{\bt\circ\imath\circ \phi,\,
\imath\circ \phi}(\phi_{*1}^{-1}(x))(\imath_T(\tau))\tforal x\in
K_1(A)
\eneq
and for all $\tau\in T(B).$

Put $\psi=\imath\circ \af\circ \phi.$ Then we  compute, applying \ref{T96},
\beq
&&[\psi]= [\imath\circ \phi]=[\bt\circ\imath\circ \phi]\,\,\,{\rm
in}\,\,\, KK(A, B\otimes M_{\mathfrak{p}})\\
&&\psi^{\dag}=\imath^{\dag}\circ\mu\circ \phi^{\dag}=(\bt\circ
\imath\circ \phi)^{\dag}, \, \mathrm{and}\\
&&\psi_T=(\imath\circ \af\circ \phi)_T=(\imath\circ
\phi)_T=(\bt\circ \imath\circ \phi)_T.
\eneq
Moreover, for any $x\in K_1(A)$ and $\tau\in T(B\otimes
M_{\mathfrak{p}}),$  {{by \eqref{T96} and \eqref{10-6-2-n2021},}}
\beq
\overline{R}_{\bt\circ \imath\circ \phi, \psi}(x)(\tau)&=&
\overline{R}_{\bt\circ \imath\circ\phi,
\imath\circ\phi}(x)(\tau)+\overline{R}_{\imath, \imath\circ \af}\circ
\phi_{*1}(x)(\tau)\\
&=&\overline{R}_{\bt\circ \imath\circ\phi,
\imath\circ\phi}(x)(\tau)+\overline{R}_{{{{\rm id}_B, \af}}}\circ
\phi_{*1}(x)(\imath_T^{-1}(\tau))\\
&=&\overline{R}_{\bt\circ \imath\circ\phi,
\imath\circ\phi}(x)(\tau)-\overline{R}_{\bt\circ \imath\circ\phi,
\imath\circ\phi}(\phi_{*1}^{-1})(\phi_{*1}(x))(\tau)=0.
\eneq
It follows from Theorem \ref{Tm72}  that $\imath\circ \af\circ \phi$ and
$\bt\circ \imath\circ \phi$ are asymptotically unitarily
equivalent.
\end{proof}

The following is a restatement of an important  result of W. Winter (\cite{Wlocal}).
\begin{thm}[Proposition 4.5  of \cite{Wlocal}]\label{TWinter}
{{Let  $\p$ and $\q$ be relatively prime supernatural numbers.
Suppose that $A$ and $B$ are non-unital separable
${\cal Z}$-stable
\CA s  and  $\phi: A\otimes {{\cal Z}_{\p,\q}}\to B\otimes {{\cal Z}_{\p,\q}}$
is a unitarily suspended $C([0,1])$-isomorphism.
Then, there is an isomorphism
${\bar \phi}: A\to B\otimes {\cal Z}.$
Moreover,
${\bar \phi}$ is approximately unitarily equivalent to the \hm\,
$(\id_B\otimes {\bar \sigma_{\p, \q}})\circ \phi\circ (\id_A\otimes 1_{{{\cal Z}_{\p,\q}}}),$
where ${{\id_A\otimes 1_{{{\cal Z}_{\p,\q}}}(a}})=a\otimes 1_{{{\cal Z}_{\p,\q}}}$ 
and
${\bar \sigma_{\p,\q}}: {{\cal Z}_{\p,\q}}\to {\cal Z}$ is the standard embedding (see 3.4
of \cite{Wlocal}).}}
\end{thm}

\begin{proof}

The unital version of this is  stated as Proposition 4.5 of \cite{Wlocal}
which follows from Proposition 4.4 and 4.3 {{of \cite{Wlocal}.}}
One notes that, both $A\otimes {{\cal Z}_{\p,\q}}$ and $B\otimes {{\cal Z}_{\p,\q}}$
are $C([0,1])$-algebras as $C([0,1])$ {{can be}} embedded unitally into
the centers of  both $M(A\otimes {\cal Z}_{\p, \q})$ and $M(B\otimes {\cal Z}_{\p, \q}).$
The unitarily suspended $C([0,1])$-isomorphism of Definition 4.2 of \cite{Wlocal}
remains the same in the non-unital case except that the continuous path of unitaries {{$\{u(t): t\in [0,1)\}$
in Definition 4.2 of \cite{Wlocal}  has the property that $u(t)$ is in}}
$(B\otimes M_\p\otimes M_\q)^\sim$
instead in $B\otimes M_\p\otimes M_\q$   for each $t\in [0,1)$
(we still have $u(0)=1$).
 Moreover, by multiplying  a continuous path of scalar  unitaries $\{v(t): t\in [0,1)\}$
 (with $v(t)\in \C$),
we may always assume that  $\pi(u(t))=1,$ for all $t\in [0,1),$ where $\pi: (B\otimes M_\p\otimes M_\q)^\sim
\to \C$ is the quotient map.
%
%

We claim {{that}} Lemma 4.3 of \cite{Wlocal} holds for the case that $A$ and $B$ are not unital,
but the unitary {{$u$ is in}} $(B\otimes {\cal Z}_{\p,\q})^\sim.$
To see this,
extends $\phi$ to an isomorphism, still denoted by  $\phi,$
from
$M(A\otimes {{\cal Z}_{\p,\q}})$ onto $M(B\otimes {{\cal Z}_{\p,\q}})$ (see 3.12.10 of \cite{Pbook}).
Note that ${{\cal Z}_{\p,\q}}$ is unital. Maps with the form
$\id_{A\otimes {\cal Z}}\otimes \rho_m$ and $\id_{B\otimes {{\cal Z}}}\otimes \rho_m$
can also {{be extended}} so they {{are defined}} on  {{$(A\otimes {\cal Z}\otimes {\cal Z}_{\p,\q})^\sim.$}}


In what follows we will retain the notation used in the proof of Lemma 4.3 of \cite{Wlocal}.
Let us make some modification.
For any finite subset ${\cal F}_0'\subset (A\otimes {\cal Z})_+,$
with an arbitrarily small error, \wilog, we may
choose $e_0, e_1\in (A\otimes {\cal Z})_+^{\bf 1}$ such that
$e_0e_1=e_1e_0=e_1$ and
$e_1a=ae_1=a$ for all $a\in {\cal F}_0'.$  Choose
${\cal F}'=\{e_1\}\cup {\cal F}_0'.$
In the definition  of  ${\cal F}$ after equation  (10) of the proof of Lemma 4.3 of \cite{Wlocal}, we  may
assume  that ${\cal F}'$ is as we described.

{{The only other additional modification is for the choice of the unitary $u$ in (35) of the proof of Lemma 4.3 of \cite{Wlocal}.}}
Let $z\in U(Z_{P_{2k_m+1}, Q_{2k_m+1}}\otimes Z_{P_{2k_m+1}, Q_{2k_m+1}})$ be as in  (11) of the proof of  Lemma 4.3 of \cite{Wlocal}.   Recall that this $z$ is a small perturbation of a unitary in ${\cal Z}\otimes {\cal Z}={\cal Z}$ (pull back from
the inductive limit to finite stage as described at the end of the proof  Lemma 3.11 of \cite{Wlocal}).
Since $U({\cal Z}\otimes {\cal Z})=U_0({\cal Z}\otimes {\cal Z}),$ we may assume that $z\in U_0(Z_{P_{2k_m+1}, Q_{2k_m+1}}\otimes Z_{P_{2k_m+1}, Q_{2k_m+1}}).$ {{(In fact, one can directly prove that $U(Z_{P_{2k_m+1}, Q_{2k_m+1}}\otimes Z_{P_{2k_m+1}, Q_{2k_m+1}})=U_0(Z_{P_{2k_m+1}, Q_{2k_m+1}}\otimes Z_{P_{2k_m+1}, Q_{2k_m+1}})$.)}}
We may write $z=\exp(i H_1)\exp (i H_2)\cdots \exp(i  H_l)$
for some integer $l$ (we can make $l=2,$ {{see Corollary 3.11 of \cite{LinZexp})
and  for some $H_j\in (Z_{P_{2k_m+1}, Q_{2k_m+1}}\otimes Z_{P_{2k_m+1}, Q_{2k_m+1}})_{s.a.}.$}}
Define $z_0:=\exp(i(e_0\otimes H_1))\exp (i(e_0\otimes  H_2))\cdots \exp(i(e_0\otimes  H_l)).$
Note that $$z_0\in ((A\otimes {\cal Z})\otimes Z_{P_{2k_m+1}, Q_{2k_m+1}}\otimes Z_{P_{2k_m+1}, Q_{2k_m+1}})^\sim.$$
Put
 \beq
 z'&:=&({\bf 1}_{(A\otimes {\cal Z})^\sim}\otimes \gamma_{2k_m+1})\otimes \theta^{[3]})(z)\in  {\bf 1}_{(A\otimes {\cal Z})^\sim}\otimes {\cal Z}\otimes {\cal Z}_{\p, \q}\andeqn\\
     z_0'&:=&(\id_{(A\otimes {\cal Z})^\sim} \otimes \gamma_{2k_m+1}\otimes \theta^{[3]})(z_0)\in (A\otimes {\cal Z}\otimes
     {\cal Z}\otimes {\cal Z}_{\p,\q})^\sim
     \\
    z_0''&:=&  (\phi^{[1,2,4]}\otimes \id_{\cal Z}^{[3]})(z_0').
   \eneq
%
We have,  by the definition of $z_0,$ for $a\otimes c\in {\cal F}$ (see also (12) in the proof of \cite{Wlocal}),
\beq\label{43u'}
&&\hspace{-0.4in}{\rm Ad}\, z_0'' ((\phi
^{[1,2,4]})(a\otimes {\bf 1}_{{\cal Z}_{\p,\q}})
\otimes (\theta\circ {\bar \varrho}_{k_m}\circ \kappa_{2k_m})^{[3]}(c))\\
&&{{=(\phi^{[1,2,4]}
\otimes \id_{\cal Z}^{[3]})\circ {\rm Ad}\, z_0'((a\otimes \gamma_{2k_m+1}({\bf 1}_{C_{2k_m+1}})\otimes (\theta\circ {\bar \varrho_{k_m}}\circ \kappa_{2k_m})^{[3]}(c)))}}\\\nonumber
&&{{=(\phi^{[1,2,4]}\otimes \id_{\cal Z}^{[3]})\circ
(\id_{A\otimes {\cal Z}}\otimes (\theta\otimes \gamma_{2k_m+1})^{[3,4]})({\rm Ad}\, z_0((a\otimes {\bf 1}_{C_{2k_m+1}}) \otimes {\bar \rho}_{k_m}\circ \kappa_{2k_m}^{[3]}(c)))}}\\\nonumber
&&=(\phi^{[1,2,4]}\otimes \id_{\cal Z}^{[3]})\circ
(\id_{A\otimes {\cal Z}}\otimes (\theta\otimes \gamma_{2k_m+1})^{[3,4]})({\rm Ad}\, z((a\otimes {\bf 1}_{C_{2k_m+1}}) \otimes {\bar \rho}_{k_m}\circ \kappa_{2k_m}^{[3]}(c)))\\\nonumber
&&\approx_\eta(\phi^{[1,2,4]}\otimes \id_{\cal Z}^{[3]})\circ
(\id_{A\otimes {\cal Z}}\otimes (\theta\otimes \gamma_{2k_m+1})^{[3,4]})((a\otimes {\bf 1}_{C_{2k_m+1}}) \otimes {\bar \rho}_{k_m}\circ \kappa_{2k_m}(c)))\\\label{NN-12}
&&= (\phi^{[1,2,4]}\otimes \id_{\cal Z}^{[3]}) (a\otimes  {\bf 1}_{\cal Z}\otimes  \gamma_{2k_m+1}\circ {\bar \varrho}_{k_m}\circ \kappa_{2k_m}(c)).
\eneq
We then, for large $n,$  let (new)
\beq
{\tilde u}:=(\id_{B\otimes {\cal Z}}\otimes \varrho_{k_n})^{[1,2,4]}\circ(\phi^{[1,2,4]}\otimes  \id_{\cal Z}^{[3]})(z_0')\in (B\otimes {\cal Z}\otimes {\cal Z}\otimes
Z_{\p,\q})^\sim.
\eneq
As in the proof of \cite{Wlocal}, for large $n,$ we also obtain a unitary $u'\in (B\otimes {\cal Z}\otimes {\cal Z}\otimes {\cal Z}_{\p,\q})^\sim$ such
that
\beq\label{new17}
&&\hspace{-3in}(17')\hspace{1.5in} \|u'-{\tilde u}\|<\eta.
\eneq
Note that, for all $a\otimes c\in {\cal F}$ (for large enough $n$ as in  (19) of the proof of Lemma 4.3 of \cite{Wlocal}),
\beq\nonumber
&&\hspace{-0.4in}(19')\hspace{0.1in}{\rm Ad}\, \td u \circ (\id_{B\otimes  {\cal Z}}\otimes {\varrho}_{k_n})^{[1,2,4]}\circ \phi\circ (\id_{A\otimes {\cal Z}}\otimes {\bf 1}_{{\cal Z}_{\p,\q}})\otimes (\theta\circ {\bar \varrho}_{k_m}\circ \kappa_{2k_m})^{[3]}(a\otimes c)\\
\nonumber
&&\approx_\eta (\id_{B\otimes {\cal Z}}\otimes \varrho_{k_n})^{[1,2,4]}\circ \phi)\otimes \id_{\cal Z}^{[3]}\circ
{\rm Ad}\, z' {{((a\otimes \gamma_{2k_m+1}({\bf 1}_{C_{2k_m+1}})\otimes (\theta\circ {\bar \varrho_{k_m}}\circ \kappa_{2k_m})^{[3]}(c)))}}\\\nonumber
&&=(\phi^{[1,2,4]}\otimes \id_{\cal Z}^{[3]})\circ
(\id_{A\otimes {\cal Z}}\otimes (\theta\otimes \gamma_{2k_m+1})^{[3,4]})({\rm Ad}\, z((a\otimes {\bf 1}_{C_{2k_m+1}}) \otimes {\bar \rho}_{k_m}\circ \kappa_{2k_m}^{[3]}(c)))
\eneq
(This serves as new (19) in the proof of Lemma 4.3 of \cite{Wlocal}).

Denote $u''':=(\phi^{[1,2,4]}\otimes \id_{\cal Z}^{[3]})\circ (\id_{A\otimes {\cal Z}}\otimes (\theta\otimes \gamma_{2k_m+1})^{[3,4]})(z_0)\in (B\otimes {\cal Z}\otimes {\cal Z}\otimes {\cal Z}_{\p,\q})^\sim$
(as $u'''$ in the proof of Lemma 3.4 of \cite{Wlocal}).
By the definition of $z_0$ and recalling the fact that $\theta({\bf 1}_{C_{k_m+1}})={\bf 1}_{\cal Z},$
 for all $a\otimes c\in {\cal F}$  (see [12] of the proof of Lemma 3.4 of \cite{Wlocal}),
\beq\nonumber
&&\hspace{-0.6in}{\rm Ad}\, u'''\circ \phi^{[1,2,4]}\circ (\id_{A\otimes {\cal Z}}\otimes \rho_{k_m})(a\otimes c)\\\nonumber
&&={\rm Ad}\, u'''\circ (\phi^{[1,2,4]}\otimes \id_{\cal Z}^{[3]})(a\otimes \theta ({\bf 1}_{C_{2k_m+1}}) \otimes \gamma_{2k_m+1}\circ {\bar \rho}_{k_m}\circ \kappa_{2k_m}(c))\\\nonumber
&&\hspace{-0.3in}=(\phi^{[1,2,4]}\otimes \id_{\cal Z}^{[3]})\circ
(\id_{A\otimes {\cal Z}}\otimes (\theta\otimes \gamma_{2k_m+1})^{[3,4]})({\rm Ad}\, z_0((a\otimes {\bf 1}_{C_{2k_m+1}}) \otimes {\bar \rho}_{k_m}\circ \kappa_{2k_m}(c)))\\\nonumber
&&\hspace{-0.3in}=(\phi^{[1,2,4]}\otimes \id_{\cal Z}^{[3]})(
(\id_{A\otimes {\cal Z}}\otimes (\theta\otimes \gamma_{2k_m+1})^{[3,4]})({\rm Ad}\, z((a\otimes {\bf 1}_{C_{2k_m+1}}) \otimes {\bar \rho}_{k_m}\circ \kappa_{2k_m}(c)))\\
&&\hspace{-0.3in}=(\phi^{[1,2,4]}\otimes \id_{\cal Z}^{[3]})(
(\id_{A\otimes {\cal Z}}\otimes (\theta\otimes \gamma_{2k_m+1})^{[3,4]})(a \otimes {\bar \rho}_{k_m}\circ \kappa_{2k_m}(c)\otimes {\bf 1}_{C_{2k_m+1}})\\
&&\approx_{\eta}(\phi^{[1,2,4]}\otimes \id_{\cal Z}^{[3]})(a\otimes \theta\circ {\bar \rho}_{k_m}\circ \kappa_{2k_m}(c)\otimes \gamma_{2k_m+1}({\bf 1}_{C_{2k_m+1}}))\\\label{New 12}
&&=\phi^{[1,2,4]}(a\otimes {\bf 1}_{{\cal Z}_{\p,\q}})\otimes ((\theta\circ {\bar \rho}_{k_m}\circ \kappa_{2k_m})^{[3]}(c)).
\eneq

Exactly the same perturbation as   in the proof of Lemma 4.3 of \cite{Wlocal}, the unitary $u$ in equation (14) of that proof
can be chosen  to be in $C([0,1], B\otimes M_{P_{2k_n}}\otimes M_{Q_{2k_n}})^\sim$ (see the first paragraph of this proof).
Therefore  the unitary $w'$ (just above (28) of that proof) can be chosen to be in $(B\otimes {\cal Z}\otimes C_{2k_n}\otimes C_{2k_n})^\sim.$
The corresponding $w''$ is in  $(B\otimes {\cal Z}\otimes C_{2k_n})^\sim.$
Consequently, we obtain a unitary $w'''\in (B\otimes {\cal Z}\otimes {\cal Z}_{\p,\q})^\sim$
such that (30) in the proof of Lemma 4.3 of \cite{Wlocal} holds.
We then retain the notation $u'':=\id_{B\otimes {\cal Z}\otimes {\cal Z}_{\p,\q}}^{[1,2,4]}(w''').$
Define, as in  (35) of the proof of Lemma 4.3 of \cite{Wlocal}, $U:=v_B\otimes \id_{{\cal Z}_{\p,\q}}(u'u''u''')^*\in (B\otimes {\cal Z}_{\p,\q})^\sim$
(as $u$ in (35) there).

We now proceed the same estimates as the ones between (35) and (36) of the proof of Lemma 4.3
of \cite{Wlocal}, using $U$ instead of $u.$  For the second inequality in these estimates which uses  original (12) for
the first time, we will use the estimate \eqref{New 12}
instead.
We also replace original (17) by  (17') above, and  original (19) by (19') above, and, for the  second use of original (12)
in that estimates,  we use
\eqref{NN-12}.
As in the proof of Lemma 4.3 of \cite{Wlocal}, we conclude that, for all $f\in {\cal F},$
\beq
\phi\circ (\id_A\otimes \rho_m)(f)\approx_{\ep} {\rm Ad}\, U\circ (\id_B\otimes \varrho_n)\circ \phi\circ (\id_A\otimes \varrho_m)(f).
\eneq
This proves the claim.

Then the conclusion {of} Proposition 4.4 of \cite{Wlocal} holds under the assumption
that both $A$ and $B$ are non-unital
with  $w_i^A\in (A\otimes {{\cal Z}_{\p,\q})^\sim}$
and $w_i^B\in (B\otimes {{\cal Z}_{\p,\q}})^\sim.$
The same proof works for the non-unital case.
In the proof, applying the non-unital version of Lemma 4.3 as just shown,
the unitary $x_i$ is in $(B\otimes {\cal Z}_{\p,\q})^\sim.$
{{So  $(\phi^\sim)^{-1}(x_i)\in   (A\otimes Z_{\p,\q})^\sim.$
Therefore, the same perturbation  as in the proof of Proposition 4.4 of \cite{Wlocal}, one can choose
${\bar x}_i\in (A\otimes \gamma_{2{\bar m_i''}}(Z_{P_{2{\bar m_i}''},Q_{2{\bar m_i}''}}))^\sim,$ and, exactly
the same way,}}  ${\bar y}_i\in (B\otimes \gamma_{2{\bar n_i''}}(Z_{P_{2{\bar n_i''}},Q_{2{\bar n_i''}}}))^\sim$ as desired.
%
With {{the}} modified version of 4.3 and 4.4 of \cite{Wlocal} described above, Proposition
4.5 of \cite{Wlocal} holds as stated in the current lemma.

\end{proof}

  \begin{thm}\label{CMT1}
Let $A$ and $B$ be  separable  simple 
\CA s  which have continuous {{scale}} and
 satisfy the UCT.
Suppose that there is an isomorphism
$$
\Gamma: {\rm Ell}(A)={{((K_0(A), T(A), \rho_A), K_1(A))\to {\rm Ell}(B)=((K_0(B), T(B), \rho_B), K_1(B))}}
$$
{{({{see}} the {{end of \ref{DElliott+} and \ref{DElliott}}}).}}
Suppose also that, for some pair of relatively prime supernatural
numbers $\mathfrak{p}$ and $\mathfrak{q}$ of infinite type such
that $M_{\mathfrak{p}}\otimes M_{\mathfrak{q}}\cong Q,$ we have
$A\otimes M_{\mathfrak{p}}\in {\cal D},$
$B\otimes
M_{\mathfrak{p}}\in {\cal D},$
$A\otimes M_{\mathfrak{q}}\in {\cal D},$
and
$B\otimes M_{\mathfrak{q}}\in {\cal D}.$
Then
$$
A\otimes {\cal Z}\cong B\otimes {\cal Z}.
$$
{{Moreover, the isomorphism induces $\Gamma.$}}
\end{thm}

\begin{proof}
The proof is almost identical to that of 11.7 of \cite{Lininv}, with a few necessary modifications.
Note that $\Gamma$ induces an isomorphism
$$
\Gamma_{\mathfrak{p}}: {\rm Ell}(A\otimes M_{\mathfrak{p}})\to
{\rm Ell}(B\otimes M_{\mathfrak{p}}).
$$
 Since $A\otimes M_{\mathfrak{p}},
B\otimes M_{\mathfrak{p}}\in {\cal D},$ we have, by Theorem \ref{Misothm},
$A\otimes M_{\mathfrak{p}},
B\otimes M_{\mathfrak{p}}\in {\cal M}_1\cap {\cal D}^d.$
By Theorem \ref{Misothm},
there is an isomorphism $\phi_{\mathfrak{p}}:
A\otimes M_{\mathfrak{p}}\to B\otimes M_{\mathfrak{p}}.$ Moreover
$\phi_{\mathfrak{p}}$ carries $\Gamma_{\mathfrak{p}}.$  In the same way, $\Gamma$ induces an isomorphism
$$
\Gamma_{\mathfrak{q}}:{\rm Ell}(A\otimes M_{\mathfrak{q}})\to
{\rm Ell}(B\otimes M_{\mathfrak{q}})
$$
and there is an isomorphism $\psi_{\mathfrak{q}}: A\otimes
M_{\mathfrak{q}}\to B\otimes M_{\mathfrak{q}}$ which induces
$\Gamma_{\mathfrak{q}}.$
{{Denote by}} $\phi_\p^\sim: {\widetilde{A\otimes M_\p}}\to {\widetilde{B\otimes M_\p}}$
and $\psi_\q^\sim: {\widetilde{A\otimes M_\q}}\to {\widetilde{B\otimes M_\q}}$
the unital extensions of $\phi_\p$ and $\psi_\q,$ respectively.
 Put $\phi=\phi_{\mathfrak{p}}\otimes {\rm
id}_{M_{\mathfrak{q}}}: A\otimes Q\to B\otimes Q$ and
$\psi=\psi_{\mathfrak{q}}\otimes {\rm
id}_{{M_{\mathfrak{p}}}}: A\otimes Q\to B\otimes Q.$
Also $\phi^\sim:=\phi^\sim_{\mathfrak{p}}\otimes {\rm
id}_{M_{\mathfrak{q}}}$ and $\psi^\sim:=\psi^\sim_\q\otimes {\rm id}_{M_\p}.$
Note that
$$
(\phi)_{*i}=(\psi)_{*i}\,\,{\rm (} i=0,1 {\rm )} \andeqn
\phi_T=\psi_T
$$
(all four of these maps are induced by $\Gamma$). Note that $\phi_T$ and $\psi_T$
are affine homeomorphisms. Since $K_{*i}(B\otimes Q)$ is
divisible, we in fact have $[\phi]=[\psi]$ (in $KK(A\otimes Q,
B\otimes Q)$). It follows from Lemma \ref{L114} that there is an
automorphism $\bt: B\otimes Q\to B\otimes Q$ such that
$$
[\bt]=[{\rm id}_{B\otimes Q}]\,\,\, \mathrm{in}\  KK(B\otimes Q, B\otimes Q)
$$
{and} such that $\phi$ and $\bt\circ \psi$ are strongly asymptotically unitarily
equivalent.
We will {{write}} $\bt^\sim: {\widetilde{B\otimes Q}}\to {\widetilde{B\otimes Q}}$ for the unital extension.
Note {{that}} also in this case,
$$
\bt_T=({\rm id}_{B\otimes Q})_T.
$$
Let $\imath: B\otimes M_{\mathfrak{q}}\to B\otimes Q$ be defined by
$\imath(b)=b\otimes 1$ for $b\in B\otimes M_\p.$ We consider the pair
$\bt\circ \imath\circ {\psi}_{\mathfrak{q}}$ and $\imath \circ
{\psi}_{\mathfrak{q}}.$ Applying Lemma \ref{L115}, we obtain an
automorphism $\af: B\otimes M_{\mathfrak{q}}\to B\otimes
M_{\mathfrak{q}}$ such that $\imath\circ \af\circ
\psi_{\mathfrak{q}}$ and $\bt\circ \imath\circ
\psi_{\mathfrak{q}}$ are  strongly asymptotically unitarily equivalent (in
$B\otimes Q$).
Moreover,
$$
[\af]=[{\rm id}_{B\otimes M_{{\mathfrak{q}}}}]\,\,\,{\rm in}\,\,\,
KK(B\otimes M_{\mathfrak{q}},B\otimes M_{\mathfrak{q}}).
$$

We will show that $\bt\circ \psi$ and $(\af\circ{\psi}_{\mathfrak{q}})\otimes {\rm id}_{M_{\mathfrak{p}}}$
 are
strongly asymptotically unitarily equivalent. Define
$\bt_1=(\bt\circ \imath\circ\psi_{\mathfrak{q}})
\otimes {\rm
id}_{M_{\mathfrak{p}}}: B\otimes M_\q\otimes M_{\mathfrak{p}}\to
B\otimes Q\otimes M_{\mathfrak{p}}.$ Let $j: M_\p\to M_\p\otimes
M_{\mathfrak{p}}$ be defined by $j(b)=b\otimes 1$ for $b\in M_\p.$
Let $s: M_\p\otimes M_\p\to M_\p$ be {{an isomorphism}} such that
$s\circ j$ is strongly asymptotically unitarily equivalent to ${\rm id}_{M_\p}.$
Let ${\bar s}: B\otimes M_\q\otimes M_\p\otimes M_\p\to B\otimes M_\q\otimes M_\p$ be
given by ${\bar s}(b\otimes a)=b\otimes s(a)$ for all $b\in B\otimes M_\q$ and $a\in M_\p\otimes M_\p,$
and ${\bar j}: B\otimes M_\q\otimes M_\p\to B\otimes M_\q\otimes M_\p\otimes M_\p$
defined by {{${\bar j}(b\otimes c)=b\otimes  j(c)$}}
for all $b\in B\otimes M_\q$ and $c\in M_\p.$

It follows that $
(\af\circ\psi_{\mathfrak{q}})
\otimes {\rm id}_{M_{\mathfrak{p}}}$ and
$
{\bar s}\circ\bt_1$
are strongly asymptotically unitarily
equivalent (note that $\imath\circ \af\circ
\psi_{\mathfrak{q}}$ and $\bt\circ \imath\circ
\psi_{\mathfrak{q}}$  are
strongly asymptotically unitarily equivalent).

Let $\psi^\sigma: A\otimes M_\q\otimes M_\p\to B\otimes M_\q\otimes M_\p\otimes M_\p$ {{be}}
defined by $\psi^\sigma(a\otimes b)=\psi_\q(a)\otimes j_r(b)$
for all $a\in A\otimes M_\q$ and $b\in M_\p,$ where $j_r: M_\p\to  M_\p\otimes M_\p$ {{is}} defined by
$j_r(b)=1_{M_\p}\otimes b$  for all $b.$
Define $\bt^\sigma:=\bt\otimes \id_{M_\p}: B\otimes M_\q\otimes M_\p\otimes M_\p
\to B\otimes M_\q\otimes M_\p\otimes M_\p.$
Note that
\beq
\bt^\sigma(\psi_q(a)\otimes b\otimes 1_{M_\p})=\bt(\psi_q(a)\otimes b)\otimes 1_{M_\p}\rforal a\in A\otimes M_\q\andeqn b\in M_\p.
\eneq
By Theorem \ref{Tm72},
there is a continuous path of unitaries $\{v(t): t\in [0,1)\}\subset (B\otimes M_\q\otimes M_\p\otimes M_\p)^\sim$
with $v(0)=1$ such that
\beq
\hspace{-0.2in}\lim_{t\to 1}v(t)^* (\psi_q(a)\otimes b\otimes 1_{M_\p})v(t)=\psi_q(a)\otimes 1_{M_\p}\otimes b\rforal a\in A\otimes M_\q\andeqn b\in M_\p.
\eneq
Therefore
\beq
v_1(t)^*\bt^\sigma(\psi_q(a)\otimes b\otimes 1_{M_\p})v_1(t)=\bt(\psi_q(a)\otimes 1_{M_\p})\otimes b=\bt(\iota \circ \psi_q(a))\otimes b
\eneq
for all $ a\in A\otimes M_\q\andeqn
b\in M_\p,$ where $v_1(t)=\bt^{\sigma\sim}(v(t)).$

It follows that
\beq
{{\lim_{t\to 1}}}~ {\bar s}(v_1(t))^*({\bar s}({\bar j}(\bt\circ \psi(a\otimes b)))){\bar s}(v_1(t))=
{\bar s}(\bt\circ \iota\circ \psi_\q(a)\otimes b)
\eneq
for all $a\in A\otimes M_\q$ and $b\in M_\p.$
Since ${\bar s}\circ {\bar j}$ is strongly asymptotically unitarily equivalent to ${\rm id}_{B\otimes Q},$
$\bt\circ \psi$ and ${\bar s}\circ \bt_1$ are strongly asymptotically unitarily equivalent.
Hence $\bt\circ \psi$ and
$
(\af\circ \psi_{\mathfrak{q}})\otimes {\rm id}_{M_{\mathfrak{p}}}$
are strongly asymptotically unitarily equivalent. Finally, we
conclude that $
(\af\circ \psi_{\mathfrak{q}})\otimes {\rm
id}_{M_{\mathfrak{p}}}$ and $\phi$ are strongly asymptotically
unitarily equivalent. Note that $\af\circ \psi_{\mathfrak{q}}$ is
an isomorphism which induces $\Gamma_{\mathfrak{q}}.$

Let $\{u(t): t\in [0,1)\}$ be a continuous path of unitaries in
$B\otimes Q$ with $u(0)=1_{\widetilde{B\otimes Q}}$ such that
$$
{{\lim_{t\to 1}}}~{\rm Ad}\, u(t)\circ \phi(a)=\af\circ
\psi_{\mathfrak{q}}\otimes {\rm id}_{{M_{\mathfrak{p}}}}(a)\tforal
a\in A\otimes Q.
$$
One then obtains a unitary suspended  $C([0,1])$-isomorphism which lifts
$\Gamma$ along $Z_{p,q}$ (see \cite{Wlocal}). It follows from Theorem \ref{TWinter}
(4.5  of \cite{Wlocal}) that $A\otimes {\cal Z}$ and $B\otimes {\cal Z}$
are isomorphic, {{and, by the exactly the same computation  in the proof of
7.1 of \cite{Wlocal}, the isomorphism induces $\Gamma.$}}
\end{proof}

{{\begin{rem}\label{RTThom}
Suppose that $\Gamma: {\rm Ell}(A)\to {\rm Ell}(B)$ is  a \hm\, (see \ref{DElliott}).
In the proof above, we obtain \hm s  $\phi_{\mathfrak{d}}:
A\otimes M_{\mathfrak{d}}\to B\otimes M_{\mathfrak{d}}$  which carries $\Gamma_{\mathfrak{d}},$
$\mathfrak{d}=\mathfrak{p}, \mathfrak{q},$ by applying \ref{Text1} (as well as \ref{Misothm}){{---note that, in Theorem \ref{Text1}, one does not assume the invariants of $A\otimes M_{\mathfrak{d}}$ and $B\otimes M_{\mathfrak{d}}$ to be isomorphic}}.
A modification of the proof would give a \hm\,  $\Phi: A\otimes {\cal Z}\to B\otimes {\cal Z}$ which carries $\Gamma.$
\end{rem}}}


\begin{thm}\label{TisomTconscale}
Let $A$ and $B$ be   separable amenable simple ${\cal Z}$-stable  \CA s
 satisfying the UCT.  Suppose that $A$ and $B$
  have continuous {{scale}}
 and $A\otimes Q, B\otimes Q\in {\cal D}.$
Suppose that there is an isomorphism
\beq
\Lambda: ((K_0(A), T(A), \rho_A), K_1(A))
 \cong ((K_0(B), T(B), \rho_B), K_1(B)).
\eneq
{{Then there exists an isomorphism $\phi: A\cong B$ which induces $\Lambda.$}}
\end{thm}
\begin{proof}
We only need to prove the ``if" part of the statement.
It follows from Theorem {{6.13 of \cite{GLrange}}} 
that $A\otimes M_\p, B\otimes M_\p\in {\cal D}$ for any supernatural number $\p.$
Thus, by Theorem \ref{CMT1},  $A\cong B.$

\end{proof}

\begin{cor}\label{TisomTconscale-c}
Let $A$ and $B$ be   separable amenable simple   \CA s
in ${\cal D}$
 satisfying the UCT.  Suppose that $A$ and $B$
  have continuous {{scales.}}
Then $A\cong B$ if and only if
\beq
 ((K_0(A), T(A), \rho_A), K_1(A))
 \cong ((K_0(B), T(B), \rho_B), K_1(B)).
\eneq
\end{cor}

\begin{proof}
It follows from Theorem 16.10 of \cite{GLII} that $A$ and $B$ are ${\cal Z}$-stable.
Since $A, B\in {\cal D},$ $A\otimes Q, {{B\otimes Q}}\in {\cal D}.$
\end{proof}

\begin{cor}\label{TisomTconscale-2}
Let $A$ and $B$ be stably projectionless  separable simple \CA s which {{have}} 
{{finite nuclear dimension}}
 and satisfying the UCT.  Suppose that $A$ and $B$ have continuous scale.
Then $A\cong B$ if and only if
\beq
 ((K_0(A), T(A), \rho_A), K_1(A))
 \cong ((K_0(B), T(B), \rho_B), K_1(B)).
\eneq
\end{cor}

\begin{proof}
We only need to  prove the ``if" part of the statement.
By {{Theorem 6.14 of \cite{GLrange},}} 
 $A\otimes Q, B\otimes Q\in {\cal D}.$
Thus,  the corollary follows from Theorem \ref{TisomTconscale}.
%
\end{proof}

\begin{thm}\label{TisomorphismT}
Let $A$ and $B$ be separable stably projectionless simple amenable \CA s  which have generalized tracial rank at most one (see \ref{DD0})
 and satisfying the UCT.
 Suppose that there is an isomorphism $\Lambda:$
 \beq\label{TisomT-1}
 &&((K_0(A),  \{0\},  {\tilde T}(A), {{\widehat{\la e_A\ra},}}\, \rho_A), K_1(A))\\
 &&\hspace{0.3in}\cong ((K_0(B), \{0\}, {\tilde T}(B), {{\widehat{\la e_B\ra}}}, \rho_B), K_1(A)).
 \eneq
 Then there is an isomorphism $\phi: A\cong B$ which induces $\Lambda.$
 Moreover, for any  simple ordered group paring $(G_0, \{0\}, T, s, \rho)$ with
 $\rho(G_0)\cap \Aff_+(T)=\{0\}$ and any
  countable abelian group $G_1,$
  there is a stably projectionless simple amenable \CA\, $A$ with $gTR(A)\le 1$
  such that
  $$
  ((K_0(A),  \{0\}, {\tilde T}(A), {{\widehat{\la e_A\ra},}} \rho_A), K_1(A))=(G_0, \{0\}, T, s, \rho, G_1)
  $$
  which is stably isomorphic to a simple \CA\, constructed in {{Theorem 4.31 of \cite{GLrange}.}} 

\end{thm}
{\em Note that, in this case $\Sigma(K_0(A))=\Sigma(K_0(B))=\{0\}$ is automatic.}

\begin{proof}
First we note  that  {{the}} ``Moreover" part follows from {{Theorems 5.3 and  4.31 of \cite{GLrange}.}}

So we will prove the isomorphism part of the theorem.
We first note that  by {{A.10 of \cite{GLII}}} (see also 11.7 of \cite{eglnp}),
for any $a\in {\rm Ped}(A)_+\setminus \{0\},$  $\overline{aAa}$ is ${\cal Z}$-stable.
Also {{$\overline{bBb}$ is ${\cal Z}$-stable  for any $b\in {\rm Ped}(B)_+\setminus \{0\}.$}}

By 9.1 of \cite{eglnp},
${\tilde T}(A)\not=\{0\}$ and ${\tilde T}(B)\not=\{0\}.$
Let
\hspace{-0.12in}\beq
\Gamma:((K_0(A), \{0\}, {\tilde T}(A), {\widehat{\la e_A\ra}}, \rho_A), K_1(A))\to  ((K_0(B), \{0\}, {\tilde T}(B), {\widehat{\la e_B\ra}},  \rho_B), K_1(B))
\eneq
be an isomorphism.
Let $\Gamma_T: {\tilde T}(A)\to {\tilde T}(B)$ be the cone  homeomorphism
such that
\beq
\Sigma_B(\Gamma_T(\tau))=\Sigma_A(\tau)\rforal \tau\in {\tilde T}(A).
\eneq
{{Let $a\in {\rm Ped}(A)_+\setminus \{0\}.$ Then, there exists
$e_a\in  (\overline{aAa})_+$
{{with $\|e_a\|=1$}}  such that
$A_0:=\overline{e_aAe_a}$  has continuous scale (see
5.2 of
\cite{eglnp}).}}
{{Note that ${{e_a}}\in {\rm Ped}(A).$ In particular, $A_0$ is ${\cal Z}$-stable.}}
By  Proposition 11.11  (and Theorem 9.4) of \cite{eglnp},
there {{exists}} $e_b\in B_+$ such that $\|e_b\|=1$ and
$d_{\Gamma(t)}(e_b)=d_t(e_a)$
for all $t\in \td T(A).$  It follows that $\widehat{e_b}$ is continuous and
(by Proposition 5.4 of \cite{eglnp}) $B_0:=\overline{e_bBe_b}$ has continuous scale.
{{Since $\overline{bBb} \in {\cal D}$ is ${\cal Z}$-stable, $e_b\in {\rm Ped}(B)_+$ and
$B_0$ is also ${\cal Z}$-stable.}}
Then $T(A_0)$  and $T(B_0)$ are  metrizable Choquet {{simplexes}}.
Moreover $T(A_0)$  and $T(B_0)$ can be identified with
\beq
T_A=\{\tau\in {\tilde T}(A): d_\tau(e_a)=1\}\andeqn  T_B=\{s\in {\tilde T}(B): d_s(e_b)=1\},
\eneq
respectively.
%
%
%
%
%
It follows that $\Gamma$ induces the following isomorphism
\beq
((K_0(A_0), T(A_0), \rho_{A_0}), K_1(A_0))\cong ((K_0(B_0), T(B_0), \rho_{B_0}), K_1(B_0)).
\eneq
Note {{that}} now both $A_0$ and $B_0$  have continuous {{scales.}}
It follows from    {{Theorem \ref{TisomTconscale}}}
that there is an isomorphism  $\phi_0: A_0\to  B_0$
which induces $\Gamma$ on $((K_0(A_0), K_1(A_0), T(A_0), r_{A_0}), K_1(A_0).$
({\bf Remark}:  Note that we only need $A_0$ and $B_0$ {{to be}} ${\cal Z}$-stable and
$A_0\otimes Q, B_0\otimes Q\in {\cal D}$
{{to apply}} Theorem \ref{TisomTconscale}).
Then $\phi_0$ gives an isomorphism from $A_0\otimes {\cal K}$ onto
$B_0\otimes {\cal K}.$
By \cite{Br}, we may identify $A$ with a hereditary \SCA\, of $A_0\otimes {\cal K}.$
With this {{identification}}, choose
$b\in (B_0\otimes {\cal K})_+$ such that $\phi(a)=b.$
Then
\beq
d_t(b)=\lim_{n\to\infty} t\circ \phi(a^{1/n})\tforal   t\in {\tilde T}(B).
\eneq
Note $\Sigma_B(t)=d_t(b).$ Since $B$ is simple and has stable rank one, this implies  that
$B\cong \overline{b(B_0\otimes {\cal K})b}.$  The theorem follows.
\end{proof}

{{Now we present a unified form of  isomorphism theorem
for the class of {{finite}} simple separable \CA s {{of}} 
{{finite nuclear dimension}}
which satisfy the UCT.
The following isomorphism theorem {{combined}} with {{the Elliott}} range  theorems (Theorem 5.3 and Theorem 5.2 of \cite{GLrange}) 
gives a complete classification of {class of \CA s mentioned above.}}

{{Note also  {{that}} {{$\widehat{\la e_A\ra}$}} is sometime written as the scale function $\Sigma_A\in {\rm LAff}_+(\td T(A))$
since  it is independent of the choice of the strictly positive element $e_A.$
In fact $\Sigma_A(t)=\sup\{\tau(a): a\in A_+ \,\,{\rm with}\,\, \|a\|\le 1\}.$
}}

\begin{thm}\label{TisomorphismC}
Let $A$ and $B$ be finite separable simple  \CA s
{{with}} 
finite nuclear dimension which
satisfy the UCT.
Then $A\cong B$ if and only if
 \beq\label{TisomT-1.1}
 &&((K_0(A),  \Sigma(K_0(A)),  {\tilde T}(A), {{\widehat{\la e_A\ra}}}, \rho_A), K_1(A))\\
 &&\hspace{0.3in}\cong ((K_0(B), \Sigma(K_0(B)), {\tilde T}(B), {{\widehat{\la e_B\ra}}}, \rho_B), {{K_1(B)}}).
 \eneq
\end{thm}

\begin{proof}
Since isomorphic \CA s have the same Elliott invariant,
we will prove the ``if" part only.

Denote  {{by}} $\Gamma$ the isomorphism between the Elliott invariant.
If $\Sigma(K_0(A))$ has a unit, i.e.,
there is $u\in \Sigma(K_0(A))$ such that
 $\rho_A(u){{(\tau)}}=\sup\{\tau(a): a\in {\rm Ped}(A)_+,\,\,\|a\|\le 1\}$ {{for any $\tau\in {\tilde T}(A)$,}}
so does $\Sigma(K_0(B)).$
{{In this case, let $p\in A$ be a projection such that
$[p]=u.$
Consider $A_1={{(1_{\tilde A}-p)A(1_{\tilde A}-p)}}.$ 
{{If $A_1\not=\{0\},$}}
let $a_0\in {A_1}_+$ {{be}} with $\|a_0\|\le 1$ and $a_1=a_0+p.$
Then $a_1\in A_+$ with $\|a_1\|=1.$ But ${\rm d}_\tau(a_1)>\tau(p)$
for all $\tau\in \td T(A).$  This {{contradicts}} the assumption that $\tau(p)=\Sigma_A(\tau)$
for all $\tau\in \td T(A).$
It follows that $A_1=\{0\},$ whence $p$ is the unit of $A.$}}
{{The same argument applies to $B.$}}
In other words, both $A$ and $B$ are unital.
The unital case has been established
which can be quoted from \cite{EGLN}.

If $K_0(A)_+\not=\{0\},$ i.e.,
$\rho_A(K_0(A))\cap \Aff_+(\td T(A))\not=\{0\},$ then
$\rho_A(K_0({{A}}))\cap \Aff_+(\td T(A))\not=\{0\}.$
Pick a projection $p\in M_n(A)$
for some $n.$ Put $A_1=pM_n(A)p.$ Let $y=\Gamma([p])$ and $q\in M_m(B)$ be a projection
for some $m$ such that $[q]=y.$  Put $B_1=qM_m(B)q.$ Then both $A_1$ and $B_1$ are
unital.  Using $\Gamma,$ by the classification theorem stated in \cite{EGLN},
there is an isomorphism $\phi: A_1\to B_1$ which is consistent with $\Gamma.$
Let $\psi: A_1\otimes {\cal K}\cong B_1\otimes {\cal K}$ {{be}} induced by $\phi.$
{{Choose}} a positive element $e_a\in A_1\otimes {\cal K}$ such
that {{$\widehat{\la e_a\ra }=\widehat{\la e_A\ra}$}} in ${\rm LAff}_+(\td T(A)).$
Then, since $A$ has (almost) stable rank one, $A\cong \overline{e_a(A_1\otimes {\cal K})e_a}$
(see Theorem 1.2 \cite{Rlz}). On the other hand, {{$\widehat{\la \psi(e_a)\ra }=\widehat{\la e_B\ra}.$}} Applying \cite{Rlz} again,
{{we have}}  $\psi(A)\cong B.$
Note that $\psi$ carries $\Gamma.$

Now we consider the case   that $K(A)_+=\{0\}.$ Then, $K_0(B)_+=\{0\}.$
So both
$A$ and $B$ are stably projectionless.   It follows from \cite{T-0-Z} that
both $A$ and $B$ are ${\cal Z}$-stable.  The proof of \ref{TisomorphismT} produces the
hereditary \SCA s  $A_0$ of $A$ and $B_0$ of $B$ which have continuous scale as described in
the proof \ref{TisomorphismT},
respectively. As in the proof of \ref{TisomorphismT}, it suffices to show that $A_0$ and $B_0$ are isomorphic.
Note that $A_0$ and $B_0$ {{have}} 
finite nuclear dimension {{(see
Corollary 3.6 of \cite{BP} and 2.8 of \cite{WZ-ndim}).}}
By  Theorem {{6.2 of \cite{GLrange},}} 
$A_0\otimes Q$ and $B_0\otimes Q$ are in ${\cal D}.$ Using the remark within
the proof of \ref{TisomorphismT}, $A_0\cong B_0.$
\end{proof}

\begin{rem}\label{Rfinal} 
 {{In the case}} that $\rho_A(K_0(A))\cap \Aff_+(\td T(A))=\{0\},$ when $A$ is ${\cal Z}$-stable
and has rationally generalized tracial rank at most one,
by  Theorem \ref{TisomTconscale} and  \ref{TisomorphismT}, if, {{in addition,}}  $A$ is amenable and satisfies the UCT,
then $A$ actually has generalized tracial rank at most one. Note that, in  the unital case,
there are simple separable amenable \CA s in the UCT class which have rationally generalized tracial rank
one but do not have generalized tracial rank one.


\end{rem}









 \providecommand{\href}[2]{#2}

 \end{document}